\documentstyle{amsart}

\newtheorem{Theorem}{Theorem}[section]

\newtheorem{Lemma}[Theorem]{Lemma}

\newtheorem{Corollary}[Theorem]{Corollary}

\newtheorem{Remark}[Theorem]{Remark}

\newtheorem{Definition}[Theorem]{Definition}

\begin{document}

\title{Toroidalization of birational morphisms of 3-folds}

\author{Steven Dale Cutkosky}
\thanks{Research   partially supported by NSF}

\maketitle
Suppose that $f:X\rightarrow Y$ is a morphism of algebraic varieties, over a field $\bold k$ of characteristic zero.
The structure of such morphisms is quite rich. The simplest class of such morphisms is the toroidal morphisms.
If $X$ and $Y$ are nonsingular,
$f:X\rightarrow Y$ is toroidal if there are simple normal crossing divisors $D_X$ on $X$ and $D_Y$ on $Y$ such that
$f^*(D_Y)=D_X$, and $f$ is  locally given by monomials in appropriate etale local parameters on $X$.
The precise definition of this concept is in \cite{AK} (see also \cite{KKMS}).  We state the Definition of toroidal  in \ref{Def274}. The problem of toroidalization is to determine, given a dominant morphism $f:X\rightarrow Y$,
 if there exists a commutative diagram 
 \begin{equation}\label{eq393}
 \begin{array}{rcl}
 X_1&\stackrel{f_1}{\rightarrow}&Y_1\\
 \Phi\downarrow&&\downarrow\Psi\\
 X&\stackrel{f}{\rightarrow}&Y
 \end{array}
 \end{equation}
 such that $\Phi$ and $\Psi$ are products of blow ups of nonsingular subvarieties, $X$ and $Y$ are nonsingular,
 and there exist simple normal crossing divisors $D_Y$ on $Y$ and $D_X=f^*(D_Y)$ on $X$ such that $f_1$ is toroidal
 (with respect to $D_X$ and $D_Y$).     This  is stated in Problem 6.2.1. of \cite{AKMW}.
 Toroidalization, and related concepts, have been considered earlier in different contexts,
 mostly for morphisms of surfaces. Torodialization is the strongest structure theorem which could be true for
 general morphisms.  The concept of torodialization fails completely in positive characteristic. A simple example is
 shown in \cite{C2}.
 
 In the case when $Y$ is a curve, toroidalization follows from embedded resolution of singularities (\cite{H}).
 When $X$ and $Y$ are surfaces, there are several proofs in print (\cite{AkK}, Corollary 6.2.3 \cite{AKMW}, \cite{Mat}).  They all make use of special properties of the
 birational geometry of surfaces.  An outline of proofs of the above cases can be found in the introduction to
 \cite{C2}.

 In \cite{C2}, the toroidalization problem is solved in the case when $X$ is a 3-fold and $Y$ is a surface.
 In this paper, we prove toroidalization for birational morphisms of 3-folds. 
 
 \begin{Theorem}\label{TheoremA} Suppose that $f:X\rightarrow Y$ is a birational morphism of 3-folds which are proper over an
algebraically closed field $\bold k$ of characteristic 0. Then there exists a commutative diagram of morphisms
$$
\begin{array}{rll}
X_1&\stackrel{f_1}{\rightarrow}&Y_1\\
\Phi\downarrow&&\downarrow\Psi\\
X&\stackrel{f}{\rightarrow}&Y
\end{array}
$$
where $\Phi,\Psi$ are products of blow ups of nonsingular curves and points, and there exists a simple normal crossings divisor $D_{Y_1}$ on $Y_1$
such that $D_{X_1}=f_1^{-1}(D_{Y_1})$ is a simple normal crossings  divisor and $f_1$ is toroidal with respect to $D_{X_1}$ and $D_{Y_1}$.
\end{Theorem}

If we relax some of the restrictions in the definition of toroidalization, there are other constructions
producing a toroidal morphism $f_1$, which 
are valid for arbitrary dimensions of $X$ and $Y$.
In \cite{AK} it is shown that a diagram (\ref{eq393}) can be constructed where $\Phi$ is weakened to being a
modification (an arbitrary birational morphism).  In \cite{C1} and \cite{C4}, it is shown that a diagram (\ref{eq393})
can be constructed where $\Phi$ and $\Psi$ are locally products of blow ups, but the morphisms $\Phi$, $\Psi$ and $f_1$
may not be separated.  This construction is obtained by patching local solutions valid for any given valuation.

It has been shown in \cite{AKMW} and \cite{W1} that weak factorization of birational morphisms holds in characteristic zero,
and arbitrary dimension.  That is, birational morphisms of complete varieties can be factored by an alternating
sequence of blow ups and blow downs of non singular subvarieties.  Weak factorization of birational (toric) morphisms
of toric varieties, (and of birational toroidal morphisms) has been proven by Danilov \cite{D1} and Ewald \cite{E} (for 3-folds), and by Wlodarczyk \cite{W}, Morelli \cite{Mo} and Abramovich, Matsuki and Rashid \cite{AMR} in general dimensions.

Our Theorem \ref{TheoremA}, when combined with  weak factorization for toroidal morphisms (\cite{AMR}),
gives a new proof of weak factorization of birational morphisms of 3-folds. 
We point out that our proof uses an  analysis of the structure as power series of local germs of a mapping, as opposed to the entirely different proof of weak factorization, using geometric invariant theory, of \cite{AKMW} and \cite{W1}.

\begin{Corollary} Suppose that $f:X\rightarrow Y$ is a birational morphism of 3-folds which are proper over an algebraically closed field $\bold k$ of characteristic
zero. Then there exists a commutative diagram of morphisms factoring $f$,
$$
\begin{array}{llllllllllllll}
&&X_1&&&&X_3&&&&&X_n&&\\
&\swarrow&&\searrow&&\swarrow&&\searrow&&\cdots&\swarrow&&\searrow&\\
X&&&&X_2&&&&X_4&&&&&Y
\end{array}
$$
where each arrow is a product of blow ups of points and nonsingular curves.
\end{Corollary}
 
 The problem of strong factorization, as proposed by Abhyankar \cite{Ab2} and Hironaka \cite{H}, is to factor a birational morphism
$f:X\rightarrow Y$ by constructing a diagram
$$
\begin{array}{lllll}
&&Z&&\\
&\swarrow&&\searrow&\\
X&&\stackrel{f}{\rightarrow}&&Y
\end{array}
$$
where $Z\rightarrow X$ and $Z\rightarrow Y$ factor as products of blow ups of nonsingular subvarieties.
Oda \cite{O} has proposed the analogous problem for (toric) morphisms of toric varieties.

A birational morphism $f:S\rightarrow Y$ of (nonsingular) surfaces can be directly factored by blowing up points
(Zariski \cite{Z3} and Abhyankar \cite{Ab2}), but there are examples showing that a direct factorization is not possible in general
for 3-folds (Shannon \cite {Sh} and Sally\cite{S}).

We also obtain as an immediate corollary the following new result,
which reduces the problem of strong factorization of 3-folds to the case of toroidal morphisms.
\begin{Corollary}
Suppose that the Oda conjecture on strong factorization of birational toroidal morphisms of 3-folds is true.
Then the Abhyankar, Hironaka strong factorization conjecture of birational morphisms of complete (characteristic zero) 3-folds is true.
\end{Corollary}

Abhyankar's local factorization conjecture \cite{Ab2}, which is ``strong factorization'' along a valuation,  follows from 
local monomialization (Theorem A \cite{C1}), to reduce to a locally toroidal morphism, and local factorization for
toroidal morphisms along a valuation Christensen  \cite{Ch} (for 3-folds), and Karu \cite{K} in general dimensions.

\section{An outline of the proof}
We will say a few words about the structure of the proof of Theorem \ref{TheoremA}.
Most of the proof of Theorem \ref{TheoremA} is valid for generically finite morphisms of 3-folds. We use the fact that a birational morphism
$f:X\rightarrow Y$ of complete 3-folds is an isomorphism in codimension 1 in $Y$, and has a very simple structure in codimension 2 in $Y$ (is a product of blowups of nonsingular curves above the base \cite{Ab1} or \cite{D}). The structure in codimension 2 is however not that much more difficult in codimension 2 for arbitrary dominant morphisms of 3-folds (in characteristic zero).
We hope to develop this theme in a later paper.

If $X$ is a nonsingular variety and $D_X$ is a SNC divisor on $X$, then $D_X$ defines a toroidal structure on $X$. If $V$ is a nonsingular subvariety of $X$, which is supported on $D_X$ and makes SNCs with $D_X$, then $V$ is a possible center. Let $\Phi:X_1\rightarrow X$ be the blow up of $V$. Then $X_1$ is nonsingular with toroidal structure $D_{X_1}=\Phi^{-1}(D_X)$.

Suppose that $f:X\rightarrow Y$ is a birational morphism of nonsingular complete 3-folds of characteristic zero, and
suppose that $D_Y$ is a simple normal crossings divisor on $Y$ such that $D_X=f^*(D_Y)$ is also a simple  normal 
crossings divisor, defining toroidal structures on $X$ and $Y$. Further suppose that the locus  of points in $X$ where $f$ is not smooth is contained in $D_X$.  We will refer to points where three components of $D_Y$ (or $D_X$) intersect as 3-point, points
where 2-components intersect as 2-points, and the remaining points of $D_Y$ (and $D_X$) as 1-points.

We develop a series of algorithms to manipulate local germs of our mappings, expressed in terms of series and polynomials.
We are required to blow up in both the domain and target. 

 The main result of resolution of singularities \cite{H} tells us that
we can construct a diagram 
$$
\begin{array}{rll}
X_1&&\\
\downarrow&\searrow\overline f&\\
X&\stackrel{f}{\rightarrow}&Y
\end{array}
$$
where for all points $p\in X$ and $q=f(p)\in Y$ there are  regular parameters $x,y,z$ at $p$ and $u,v,w$ at $q$,
which contain local equations of components of $D_X$ and $D_Y$ passing through these points, such that we have an expression
\begin{equation}\label{eq394}
\begin{array}{ll}
u&=x^ay^bz^c\gamma_1\\
v&=x^dy^ez^f\gamma_2\\
w&=x^gy^hz^i\gamma_3
\end{array}
\end{equation}
where $\gamma_1,\gamma_2,\gamma_3$ are units at $p$.  The algorithms of resolution give us no information about the
structure of the units $\gamma_1,\gamma_2,\gamma_3$. In general, we will have that the monomials
$x^ay^bz^c,x^dy^ez^f,x^gy^hz^i$ are algebraically dependent. We can then hope to blow up nonsingular subvarieties of $Y$,
leading to new regular parameters at a point $q'$ above $q$ with regular parameters $u',v',w'$ which are obtained by dividing  $u,v,w$ by each other. we will eventually obtain an expression such as $w'=\gamma_3-\gamma_3(p)$ by this procedure,
which need not be any better than the expressions we started with for our original map $f:X\rightarrow Y$.

In attempting to obtain a sufficiently deep understanding of map germs to prove toroidalization, we must understand 
information about germs such as (\ref{eq394}), up to the level of a series expansion of the units $\gamma_i$.
In considering this problem, we are led to generalized notions of multiplicity, which can actually increase after
making blow ups above $X$. In summary, the problem of toroidalization  requires new methods, which are not contained
in any  proofs of resolution of singularities.

In Definition \ref{Def1} we define a prepared morphism $f:X\rightarrow Y$ of 3-folds. This is related to the notion of a prepared
morphism from a 3-fold to a surface (Definition 6.5 \cite{C2}).  The basic idea of a prepared morphism of 3-folds is that locally 
an appropriate projection of $Y$ onto a surface $S$ is a toroidal morphism.  Some special care is required in
handling 3-points of $Y$.
If $f$ is prepared, by a simple local calculation involving Jacobian determinants, we have very simple expressions of the form Definition \ref{torf} and Lemma \ref{Lemmatorf} at all points of $X$. A typical case, when $p$ is a 1-point and $q$ is a 2-point, is
$$
u=x^a, v=x^b(\alpha+y), w=g(x,y)+x^cz
$$
where $0\ne\alpha\in \bold k$ and $g(x,y)$ is a series.

In Section 4 it is shown that we can construct from our given birational morphism $f:X\rightarrow Y$ a 
commutative diagram

$$
\begin{array}{rll}
X_1&\stackrel{f_1}{\rightarrow}&Y_1\\
\Phi\downarrow&&\downarrow\Psi\\
X&\stackrel{f}{\rightarrow}&Y
\end{array}
$$
where $\Phi,\Psi$ are products of blow ups of nonsingular curves and points, and there exists a SNC divisor $D_{Y_1}$ on $Y_1$
such that $D_{X_1}=f_1^{-1}(D_{Y_1})$ is a SNC divisor and $f_1$ is prepared with respect to $D_{X_1}$ and $D_{Y_1}$. We also make $D_{X_1}$ cuspidal for $f_1$ (Definition \ref{Def247}). That is, $f_1$ is toroidal in a neighborhood of all components of $D_{X_1}$ which do not contain a 3-point, and in a neighborhood of all 2-curves of $D_{X_1}$ which do not contain a 3-point.

It may appear that the local forms of a prepared morphism are very simple, and we can easily modify them to obtain
a toroidal form.  However, a little exploration with these local forms will reveal that the notion of being prepared is stable under blow ups of 2-curves (on $Y$), but is in general not stable under blow ups of points and curves which are not 2-curves (on $Y$).
It also will quickly become apparent that it is absolutely necessary to blow up subvarieties of $Y$ other than 2-curves to
toroidalize.  This leads to the notion of super parameters (Definition \ref{Def357}), which is necessary for all the blow ups which we consider to preserve the notion of being prepared, and for a global invariant $\tau$ to behave well under blow ups.

To prove Theorem \ref{TheoremA}, we may assume that $f$ is prepared, and $D_X$ is cuspidal for $f$. These conditions are preserved throughout the proof.

We define the $\tau$-invariant of a 3-point $p\in X$ (Definition \ref{Def221}). Since $f$ is prepared, $f(p)=q$ is a 2-point or a 3-point.
There are regular parameters $u,v,w$ in ${\cal O}_{Y,q}$ and $x,y,z$ in $\hat{\cal O}_{X,p}$ such that $xyz=0$ is a local equation of $D_X$,
$uv=0$ or $uvw=0$ is a local equation of $D_Y$ and there is an expression
\begin{equation}\label{N3}
\begin{array}{ll}
u&=x^ay^bz^c\\
v&=x^dy^ez^f\\
w&=\sum_{i\ge 0}\alpha_iM_i+N
\end{array}
\end{equation}
with $\alpha_i\in {\bf k}$, $M_i=x^{a_i}y^{b_i}z^{c_i}$, $N=x^gy^hz^i$,
$$
\text{rank}\left(\begin{matrix} a&b&c\\ d&e&f\end{matrix}\right)=2,
\text{det}\left(\begin{matrix} a&b&c\\ d&e&f\\ a_i&b_i&c_i\end{matrix}\right)=0\text{ for all }i,
$$
$$
\text{det}\left(\begin{matrix} a&b&c\\ d&e&f\\ g&h&i\end{matrix}\right)\ne 0.
$$
If $q$ is a 3-point, then
$$
w=\text{unit series }N
$$
 if and only if $f$ is toroidal at $p$. In this case define $\tau_f(p)=-\infty$.

Otherwise, define
$$
H_p={\bold Z}(a,b,c)+{\bold Z}(d,e,f)+\sum_i {\bold Z}(a_i,b_i,c_i),
$$
$$
A_p=\left\{\begin{array}{ll} {\bold Z}(a,b,c)+{\bold Z}(d,e,f)+{\bold Z}(a_0,b_0,c_0)&\text{ if $q$ is a 3-point (we have $w=\text{unit series }M_0$)}\\
{\bold Z}(a,b,c)+{\bold Z}(d,e,f)&\text{ if $q$ is a 2-point}.\end{array}\right.
$$
Now define
$$
\tau_f(p)=|H_p/A_p|.
$$
We define
$$
\tau_f(X)=\text{max}\{\tau_f(p)\mid p\in X\text{ is a 3-point}\}.
$$

\vskip .2truein

 We show 
in Theorem \ref{Theorem391}, 
 that $f$ is toroidal if
$\tau_f(X)=-\infty$.

We have that $\tau_f(X)\ge 0$ or $\tau_f(X)=-\infty$. The proof of Theorem \ref{TheoremA} is by descending induction on $\tau_f(X)$.
In our proof of Theorem \ref{TheoremA} we may thus assume that $\tau=\tau_f(X)\ne -\infty$ (so that $\tau\ge 0)$.
\vskip .2truein
\noindent{\bf Step 1.} (Theorem \ref{Theorem268})  There exist sequences of blow ups of 2-curves
$$
\begin{array}{lll}
X_1&\stackrel{f_1}{\rightarrow}&Y_1\\
\downarrow&&\downarrow\\
X&\stackrel{f}{\rightarrow}&Y
\end{array}
$$
such that $f_1$ is prepared, $D_{X_1}$ is cuspidal for $f_1$, $\tau_f(X_1)=\tau$, and $\tau_{f_1}(p)=\tau$ implies that $f_1(p)$ is a 2-point.
Theorem \ref{Theorem268} is a consequence of Theorem \ref{Theorem137} and the concept of 3-point relation (Definition \ref{Def160}).
\vskip .2truein
\noindent{\bf Step 2.} (Theorem \ref{Theorem269})
In this step we construct a commutative diagram of morphisms
$$
\begin{array}{lll}
X_1&\stackrel{f_1}{\rightarrow}&Y_1\\
\Phi\downarrow&&\downarrow\Psi\\
X&\stackrel{f}{\rightarrow}&Y
\end{array}
$$
such that
\begin{enumerate}
\item[1.]
$\Phi$ and $\Psi$ are products of blow ups of possible centers.
\item[2.]
$\tau_{f_1}(X_1)=\tau$, and if $p \in X$ is a 3-point such that  $\tau_{f_1}(p)=\tau$ then $f_1(p)$ is a 2-point.
\item[3.]
 $D_{X_1}$ is cuspidal for $f_1$.
\item[4.]
$f_1$ is $\tau$-very-well prepared.
\end{enumerate}
\vskip .2truein
Step 2 is the most difficult step technically. It is the content of Sections 7 and 8.

 The  definition of $\tau$-very-well prepared is given in Definition \ref{Def130}. It  uses the concept of 2-point relation (Definition \ref{Def156}), and requires the preliminary definitions
of $\tau$-quasi-well prepared (Definition \ref{Def128}) and $\tau$-well prepared (Definition \ref{Def65}).

By virtue of the result of this step, we can  assume that $f$ is $\tau$-very-well prepared. We will also assume that $\tau>0$. The case when
$\tau=0$ is actually a little easier, but the definition is a bit different.

We now summarize some of the properties of a $\tau$-very-well-prepared morphism (Definition \ref{Def130}).

There exists a finite, distinguished set of nonsingular algebraic surfaces $\Omega(\overline R_i)$ in $Y$, with a SNC divisor $F_i$ on $\Omega(\overline R_i)$
such that the intersection graph of $F_i$ is a tree.

Suppose that $p\in X$ is a 3-point with $\tau_f(p)=\tau$ (so that $q=f(p)$ is 2-point).
Then the following  conditions hold.
\begin{enumerate}
\item[1.]
The expression (\ref{N3}) has the form
\begin{equation}\label{N4}
w=\gamma M_0
\end{equation}
where $\gamma$ is a unit series, $M_0^e=u^{\overline a}v^{\overline b}$, with $\overline a,\overline b,e\in{\bold Z}$, $e>1$, and $\text{gcd}(a,b,e)=1$.
Observe that we cannot have both $\overline a<0$ and $\overline b<0$, since $M_0,u,v$ are all monomials in $x,y,z$.
\item[2.] Suppose that $V$ is the curve in $Y$ with local equations $u=w=0$ (or $v=w=0$) at $q$. Then
$V$ is a $*$-permissible center (Definition \ref{Def219}). That is, $V$ is a possible center for $D_Y$ and there exists a commutative diagram of morphisms
\begin{equation}\label{eqN5}
\begin{array}{lll}
X_1&\stackrel{f_1}{\rightarrow}&Y_1\\
\Phi_1\downarrow&&\downarrow\Psi_1\\
X&\stackrel{f}{\rightarrow}&Y
\end{array}
\end{equation}
where $\Psi_1$ is the blow up of $V$ (possibly followed by blow ups of some special 2-points), such that $f_1$ and $\overline f=\Psi_1\circ f_1:X_1\rightarrow Y$ are prepared, $\tau_{f_1}(X_1)\le\tau$
and $\Phi_1$ is toroidal at 3-points $p_1\in\Phi_1^{-1}(p)$. Further, $f_1$ is $\tau$-very-well prepared.

\item[3.] There exists a surface $\Omega(\overline R_i)$ such that
\begin{enumerate}
\item[a.] $f(p)=q\in \Omega(\overline R_i)$.
\item[b.]  The $w$ of (\ref{N4}) gives a local equation $w=0$ of $\Omega(\overline R_i)$ at $q$.
\item[c.] $uv=0$ is a local equation of $F_i$ (on the surface $\Omega(\overline R_i)$) at $q$.
\end{enumerate}
\end{enumerate}
\vskip .2truein

The necessity of several different surfaces $\Omega(\overline R_i)$ arises because of the possibility that there may be several 3-points $p_j$
with $\tau_f(p_j)=\tau$  which map to $q$, and require different $w$ in their expressions (\ref{N4}). We require that the surfaces $\Omega(\overline R_i)$
intersect in a controlled way.

The first step in the construction of a $\tau$-very well prepared morphism is the construction of a morphism such that for all 3-points $p$ with
$\tau_f(p)=\tau$, an expression (\ref{N4}) holds
for some possibly formal $w$.

\vskip .2truein
\noindent{\bf Step 3.} (Theorem \ref{Theorem270}) We construct a commutative diagram
$$
\begin{array}{lll}
X_n&\stackrel{f_n}{\rightarrow}&Y_n\\
\downarrow&&\downarrow\\
X&\stackrel{f}{\rightarrow}&Y
\end{array}
$$
such that $\tau_{f_n}(X_n)<\tau$.  By induction on $\tau$, we then obtain the proof of Theorem \ref{TheoremA}.

We fix an index $i$ of the surfaces $\Omega(\overline R_i)$.
A curve $E$ on $Y$ is good if it is a component of $F_i$, and if $j$ is such that $E\cap \Omega(\overline R_j)\ne\emptyset$,
then $E$ is a component of $F_j$.

In our construction we begin with $i=1$, and blow up a good curve $V$ on $Y$, by a morphism (\ref{eqN5}). Part of the definition of $\tau$-very-well prepared implies
the existence of a good curve. Suppose that $p\in X$ is a 3-point with $\tau_f(X)=\tau$
and $q=f(p)\in V$.
Suppose that $p_1\in\Phi_1^{-1}(p)$ is a 3-point. Set $q_1=f_1(p_1)$. If $V$ has local equations $u=w=0$ at $q$, then $q_1$ has regular parameters $u_1,v,w_1$ with
\begin{equation}\label{N6}
u=u_1w_1, w=w_1
\end{equation}
or
\begin{equation}\label{N7}
u=u_1,w=u_1(w_1+\alpha)
\end{equation}
and $\alpha\in {\bf k}$.

If (\ref{N6}) holds then $q_1$ is a 3-point. Since $e>1$, we have
$$
\tau_{f_1}(p_1)=|H_p/A_p+M_0{\bold Z}|<|H_p/A_p|=\tau.
$$
If (\ref{N7}) holds, then $e>1$ implies $\alpha=0$. Thus $f_1$ has the form (\ref{N3}), (\ref{N4}) at $p_1$ with $(\overline a, \overline b,e)$ changed to
$(\overline a-e,\overline b,e)$. If $V$ has local equations $v=w=0$, then $f_1$ has the form (\ref{N3}), (\ref{N4}) at $p_1$ with
$(\overline a, \overline b,e)$ changed to
$(\overline a,\overline b-e,e)$.

We have SNC divisors $\Phi_1^{-1}(F_i)$ on the surfaces $\Phi_1^{-1}(\Omega(\overline R_i))$).
If there are no 3-points $p_1$ in $X_1$ satisfying 1, 2 and 3 of Step 2 for $\Phi_1^{-1}(\Omega(\overline R_1))$, then we increase $i$ to 2.

Otherwise, there  exists a good curve on $Y_1$ for the SNC divisor $\Phi^{-1}(F_1)$ on the surface $\Phi_1^{-1}(\Omega(\overline R_1))$).
We continue to iterate, blowing up good curves. If we always have a 3-point satisfying 1, 2 and 3 for the preimage of $\Omega(\overline R_1)$,
then we eventually obtain a form (\ref{N4}) with both $\overline a<0$ and $\overline b<0$ which is impossible.

We then continue this algorithm for the preimages of all of the surfaces $\Omega(\overline R_i)$. The algorithm terminates in the construction of a
morphism with a drop in $\tau$ as desired.

 The final proof of Theorem \ref{TheoremA} is given after Theorem \ref{Theorem391}.

\section{Notation}
Throughout this paper, $\bold k$ will be an algebraically closed field of characteristic zero. A curve, surface or 3-fold is
a quasi-projective variety over $\bold k$ of respective dimension 1, 2 or 3.
If $X$ is a variety, and $p\in X$ is a nonsingular point, then regular parameters at $p$ are regular parameters in ${\cal O}_{X,p}$.
Formal regular parameters at $p$ are regular parameters in $\hat{\cal O}_{X,p}$.
 If $X$ is a variety and $V\subset X$ is a subvariety, then
${\cal I}_V\subset {\cal O}_X$ will denote the ideal sheaf of $V$.
If $V$ and $W$ are subvarieties of a variety $X$, we denote the scheme theoretic intersection 
$Y=\text{spec}({\cal O}_X/{\cal I}_V+{\cal I}_W)$ by $Y=V\cdot W$.

 Suppose that $a,b,c,d\in{\bf Q}$. Then we will write $(a,b)\le (c,d)$ if
$a\le b$ and $c\le d$.

A toroidal structure on a nonsingular variety $X$ is a simple normal crossing divisor (SNC divisor) $D_X$ on $X$.

We will say that a nonsingular curve $C$ which is a subvariety of a nonsingular 3-fold $X$ with toroidal structure
$D_X$ makes simple normal crossings (SNCs) with $D_X$ if for all $p\in C$, there exist regular parameters
$x,y,z$ at $p$ such that $x=y=0$ are local equations of $C$, and $xyz=0$ contains the support of $D_X$ at $p$.

Suppose that $X$ is a nonsingular 3-fold with toroidal structure $D_X$. If $p\in D_X$ is on the intersection of three components of $D_X$ then $p$ is called a 3-point. If $p\in D_X$
is on the intersection of two components of $D_X$ (and is not a 3-point) then $p$ is called a 2-point. If $p\in D_X$
is not a 2-point or a 3-point, then $p$ is called a 1-point. If $C$ is an irreducible component of the intersection of two
components of $D_X$, then $C$ is called a 2-curve. $\Sigma(X)$ will denote the closed locus of 2-curves on $X$.

By a general point $q$ of a variety $V$, we will mean a point $q$ which 
satisfies conditions which hold on some nontrivial open subset of $V$.
The exact open condition which we require will generally be clear from context.
By a general section of a coherent sheaf ${\cal F}$ on a projective variety $X$, we mean the section corresponding 
to a general point of the $\bold k$-linear space $\Gamma(X,{\cal F})$.

If $X$ is a variety, ${\bold k}(X)$ will denote the function field of $X$. A 0-dimensional valuation $\nu$ of ${\bold k}(X)$ is
a valuation of ${\bold k}(X)$ such that $\bold k$ is contained in the valuation ring $V_{\nu}$ of $\nu$ and the residue field of
$V_{\nu}$ is $\bold k$.  If $X$ is a projective variety which is birationally equivalent to $X$, then there exists a
unique (closed) point $p_1\in X_1$ such that $V_{\nu}$ dominates ${\cal O}_{X_1,p_1}$. $p_1$ is called the center of $\nu$
on $X_1$. If $p\in X$ is a (closed) point, then there exists a 0-dimensional valuation $\nu$ of ${\bold k}(X)$ such that $V_{\nu}$ dominates ${\cal O}_{X,p}$
(Theorem 37, Section 16, Chapter VI \cite{ZS}). For $a_1,\ldots, a_n\in {\bold k}(X)$, $\nu(a_1),\ldots, \nu(a_n)$ are rationally dependent if there exist
$\alpha_1,\ldots, \alpha_n\in{\bf Z}$ which are not all zero, such that $\alpha_1\nu(a_1)+\cdots\alpha_n\nu(a_n)=0$ (in the value group of $\nu$).
Otherwise, $\nu(a_1),\ldots, \nu(a_n)$ are rationally independent.

If $f:X\rightarrow Y$ is a morphism of varieties, and $D$ is  a Cartier divisor on $Y$, then $f^{-1}(D)$ will denote
the reduced divisor $f^*(D)_{red}$.

\section{Prepared, monomial and toroidal morphisms}

Throughout this section we assume that
 $f:X\rightarrow Y$ is a dominant morphism of nonsingular 3-folds,  $D_Y$
is SNC divisor on $Y$ such that $D_X=f^{-1}(D_Y)$ is a SNC divisor, and the singular locus of $f$ is contained in $D_X$.
$D_X$ and $D_Y$ define toroidal structures on $X$ and $Y$.

A possible center on a nonsingular 3-fold $X$ with toroidal structure defined by a SNC divisor $D_X$, is a point
on $D_X$ or a nonsingular curve in $D_X$ which makes SNCs with $D_X$. A possible center on a nonsingular surface $S$
with toroidal structure defined by a SNC divisor $D_S$  is a point on $D_S$.

Observe that if $\Phi:X_1\rightarrow X$ is the blow up of a possible center, then $D_{X_1}=\Phi^{-1}(D_X)$ is a SNC
divisor on $X_1$. Thus $D_{X_1}$ defines a toroidal structure on $X_1$. All blow ups $\Phi:X_1\rightarrow X$
considered in this paper will be of possible centers, and we will impose the toroidal structure on $X_1$ defined by
$D_{X_1}=\Phi^{-1}(D_X)$.

Suppose that $x,y,z$ are indeterminants, and
$M_1,M_2,\ldots,M_r$ are Laurent monomials in  $x,y,z$, so that there are expressions
$$
M_i=x^{a_{i1}}y^{a_{i2}}z^{a_{i3}}
$$
with all $a_{ij}$ in $\bf Z$.  We define
$$
\text{rank}_{(x,y,z)}(M_1,\ldots,M_r)=n
$$
if the matrix $(a_{ij})$ has rank $n$. If there is no danger of confusion, we will denote
$$
\text{rank}(M_1,\ldots,M_r)=\text{rank}_{(x,y,z)}(M_1,\ldots,M_r).
$$

Suppose that $q\in Y$. We say that $u,v,w$ are (formal) permissible
parameters at $q$ if $u,v,w$ are regular parameters in  $\hat {\cal O}_{Y,q}$ such that
\begin{enumerate}
\item[1.] If $q$ is a 1-point, then $u\in{\cal O}_{Y,q}$ and $u=0$ is a local equation of
$D_Y$ at $q$.
\item[2.] If $q$ is a 2-point then $u,v\in {\cal O}_{Y,q}$ and $uv=0$ is a local equation of $D_Y$ at $q$.
\item[3.] If $q$ is a 3-point then $u,v,w\in{\cal O}_{Y,q}$ and $uvw=0$ is a local equation of
$D_Y$ at $q$.
\end{enumerate}
$u,v,w$ are algebraic permissible parameters if we further have that $u,v,w\in{\cal O}_{Y,q}$.

\begin{Definition}\label{torf} Suppose that $u,v,w$ are (possibly formal) permissible parameters
at $q\in Y$. Then
$u,v$ are {\bf toroidal forms} at $p\in f^{-1}(q)$ if there exist regular parameters $x,y,z$
in $\hat{\cal O}_{X,p}$ such that
$x=0$, $xy=0$ or $xyz=0$ are local equations of $D_X$ and
\begin{enumerate}
\item[1.]  $q$ is a 2-point or a 3-point, $p$ is a 1-point and 
\begin{equation}\label{eqTF1}
u=x^a,
v=x^b(\alpha+y)
\end{equation}
where $0\ne \alpha\in \bold k$.
\item[2.] $q$ is 2-point or a 3-point, $p$ is a 2-point and 
\begin{equation}\label{eqTF21}
u=x^ay^b,
v=x^cy^d
\end{equation}
with $\text{rank}(u,v)=2$.
\item[3.] $q$ is a 2-point or a 3-point, $p$ is a 2-point and 
\begin{equation}\label{eqTF22}
u=(x^ay^b)^k,
v=(x^ay^b)^t(\alpha+z)
\end{equation}
where $0\ne\alpha\in \bold k$, $a,b,k,t>0$ and $\text{gcd}(a,b)=1$.
\item[4.] $q$ is a 2-point or a 3-point, $p$ is a 3-point and 
\begin{equation}\label{eqTF3}
u=x^ay^bz^c,
v=x^dy^ez^f
\end{equation}
where $\text{rank}(u,v)=2$.
\item[5.] $q$ is a 1-point, $p$ is a 1-point and 
\begin{equation}\label{eqTF01}
u=x^a,
v=y
\end{equation}
\item[6.] $q$ is a 1-point, $p$ is a 2-point and 
\begin{equation}\label{eqTF02}
u=(x^ay^b)^k,
v=z
\end{equation}
with $a,b,k>0$ and $\text{gcd}(a,b)=1$
\end{enumerate}
\end{Definition}

Regular parameters $x,y,z$ as in Definition \ref{torf} will be called permissible
parameters for $u,v,w$ at $p$.

\begin{Lemma}\label{Lemmatorf}
Suppose that $q\in Y$, $p\in f^{-1}(q)$ and $u,v,w$ are permissible
parameters at $q$ such that $u,v$ are toroidal forms at $p$. Then there exist
permissible parameters $x,y,z$ for $u,v,w$ at $p$ such that
an expression of
Definition \ref{torf} holds for $u$ and $v$, and one of the following respective forms for $w$ holds at $p$.
\begin{enumerate}
\item[1.]  $q$ is a 2-point or a 3-point, $p$ is a 1-point, $u,v$ satisfy (\ref{eqTF1})  and 
\begin{equation}\label{eqTF1w}
w=g(x,y)+ x^cz
\end{equation}
where $g$ is a series.
\item[2.] $q$ is 2-point or a 3-point, $p$ is a 2-point, $u,v$ satisfy (\ref{eqTF21})  
and 
\begin{equation}\label{eqTF21w}
w=g(x,y)+x^ey^fz
\end{equation}
where $g$ is a series.
\item[3.] $q$ is a 2-point or a 3-point, $p$ is a 2-point, $u,v$ satisfy (\ref{eqTF22}) and 
\begin{equation}\label{eqTF22w}
w=g(x^ay^b,z)+x^cy^d
\end{equation}
where $g$ is a series and $\text{rank}(u,x^cy^d)=2$.
\item[4.] $q$ is a 2-point or a 3-point, $p$ is a 3-point, $u,v$ satisfy (\ref{eqTF3}) and 
\begin{equation}\label{eqTF3w}
w=g(x,y,z)+N
\end{equation}
where $g$ is a series in monomials $M$ in $x,y,z$ such that $\text{rank}(u,v,M)=2$,
and $N$ is a monomial in $x,y,z$ such that $\text{rank}(u,v,N)=3$.
\item[5.] $q$ is a 1-point, $u,v$ satisfy (\ref{eqTF01})  and 
\begin{equation}\label{eqTF01w}
w=g(x,y)+ x^cz
\end{equation}
where $g$ is a series.
\item[6.] $q$ is a 1-point, $p$ is a 2-point, $u,v$ satisfy (\ref{eqTF02})  and 
\begin{equation}\label{eqTF02w}
w=g(x^ay^b,z)+x^cy^d
\end{equation}
where $g$ is a series and $\text{rank}(u,x^cy^d)=2$.
\end{enumerate}
\end{Lemma}

\begin{pf} Choose permissible parameters $x,y,z$ for $u,v,w$ at $p$. The Lemma follows from an explicit calculation of the jacobian determinant
$$
J=\text{Det}\left(\begin{array}{ccc}
\frac{\partial u}{\partial x} \frac{\partial u}{\partial y} \frac{\partial u}{\partial z}\\
\frac{\partial v}{\partial x} \frac{\partial v}{\partial y} \frac{\partial v}{\partial z}\\
\frac{\partial w}{\partial x} \frac{\partial w}{\partial y} \frac{\partial w}{\partial z}
\end{array}\right),
$$
and a change of variables of $x,y,z$. Observe that $J=0$ is supported on $D_X$, since the singular locus of $f$
is contained in $D_X$.

We indicate the proof if (\ref{eqTF3}) holds at $p$. There exists a unit series $\gamma$ in $x,y,z$ and $l,m,n\in{\bf N}$
such that $J=\gamma x^ly^mz^n$. Write $w$ as a series
$$
w=\sum c_{ijk}x^iy^jz^k
$$
with $c_{ijk}\in \bold k$. We compute
$$
J=\sum c_{ijk}\text{Det}\left(\begin{array}{lll}
a&b&c\\
d&e&f\\
i&j&k
\end{array}\right)
x^{a+d+i-1}y^{b+e+j-1}z^{c+f+k-1}=\gamma x^ly^mz^n
$$
from which we obtain forms (\ref{eqTF3}) and (\ref{eqTF3w}), after making a change of variables in $x,y,z$,
multiplying $x,y,z$ by appropriate unit series.
\end{pf}

\begin{Definition}\label{Def198}
Let notation be as in  Lemma \ref{Lemmatorf}.
If $p\in X$ is a 3-point, we will say that permissible parameters $u,v,w$ at $q=f(p)$ have a monomial form at $p$ if
there exist permissible parameters $x,y,z$ for $u,v,w$ at $p$ such that there is an expression
$$
\begin{array}{ll}
u&=x^ay^bz^c\\
v&=x^dy^ez^f\\
w&=x^gy^hz^i
\end{array}
$$
(with $\text{rank}(u,v,w)=3$).
\end{Definition}

\begin{Definition}\label{Def1}
A birational morphism $f:X\rightarrow Y$ of nonsingular 3-folds with toroidal structures determined by SNC divisors $D_Y$, $D_X=f^{-1}(D_Y)$ such that the singular locus of $f$ is contained in $D_X$ is {\bf prepared} if:
\begin{enumerate}
\item[1.] If $q\in Y$ is a 3-point,  $u,v,w$ are  permissible parameters at $q$
and $p\in f^{-1}(q)$, then $u,v$ and $w$ are each a unit (in $\hat{\cal O}_{X,p}$) times a monomial in local equations of the toroidal
structure $D_X$ at $p$ . Furthermore, there exists a permutation of $u,v,w$ such that
$u,v$ are toroidal forms  at $p$.
\item[2.] If $q\in Y$ is a 2-point, $u,v,w$ are permissible parameters at $q$ and $p\in f^{-1}(q)$, then either
\begin{enumerate}
\item $u,v$ are toroidal forms  at $p$ or
\item $p$ is a 1-point and there exist regular parameters $x,y,z\in\hat{\cal O}_{X,p}$ such that there is
 an expression 
$$
\begin{array}{ll}
u&=x^a\\
v&=x^c(\gamma(x,y)+x^dz)\\
w&=y
\end{array}
$$
where $\gamma$ is a unit series and $x=0$ is a local equation of $D_X$, or
\item $p$ is a 2-point and there exist regular parameters $x,y,z$ in $\hat{\cal O}_{X,p}$ such that there is
 an expression 
$$
\begin{array}{ll}
u&=(x^ay^b)^k\\
v&=(x^ay^b)^l(\gamma(x^ay^b,z)+x^cy^d)\\
w&=z
\end{array}
$$
where $a,b>0$, $\text{gcd}(a,b)=1$, $ad-bc\ne 0$, $\gamma$ is a unit series and $xy=0$ is a local equation of $D_X$.
\end{enumerate}
\item[3.] If $q\in Y$ is a 1-point,  and $p\in f^{-1}(q)$,
 then there exist  permissible parameters $u,v,w$ at $q$ such that 
  $u,v$ is
a toroidal form  at $p$.
\end{enumerate}
\end{Definition}

We call $x,y,z$ in 2 (b) or 2 (c) of Definition \ref{Def1} permissible parameters for $u,v,w$ at $p$.

\begin{Lemma}\label{Lemma141} Suppose that $X$, $Y$ are projective, $f:X\rightarrow Y$ is birational, prepared and
$q\in Y$ is a 1-point such that $f$ is not an isomorphism over $q$. Then the fundamental locus of $f$ contains a single curve $C$ passing through $q$ and $C$ is nonsingular at $q$.

Furthermore, there exist algebraic permissible parameters $u,v,w$ at $q$ such that a form (\ref{eqTF01}) or (\ref{eqTF02}) of 
Definition \ref{torf} (for this fixed choice of $u,v,w$) holds  at $p$ for all $p\in f^{-1}(q)$.
\end{Lemma}

\begin{pf} 
Since  for all $p\in f^{-1}(q)$, there exist permissible parameters at $q$ such that a  form (\ref{eqTF01}) or (\ref{eqTF02}) of 
Definition \ref{torf}  holds  at $p$, we have that $\text{dim }f^{-1}(q)=1$. Thus if $E'$ is an exceptional
component of $f$ such that $q\in f(E')$, then $f(E')$ is a curve. 

Let $D$ be the component of $D_Y$ containing $q$. Let $F$ be the strict transform of $D$ on $X$ and let $p\in f^{-1}(q)\cap F$.
Then $p$ must be a 2-point, and since $f$ is birational, there exist permissible parameters $\overline u,\overline v,\overline w$ at $q$ such that there
is an expression in $\hat{\cal O}_{X,p}$
$$
\begin{array}{ll}
\overline u&=x^ay\\
\overline v&=z\\
\overline w&=\phi(x^ay,z)+x^dy^e
\end{array}
$$
where $y=0$ is a (formal) local equation of $F$, and $x=0$ is a (formal) local equation of the other component $E$ of $D_X$
containing $p$. Computing the Jacobian determinant of $f$ at $p$, we see that $e=0$. We consider the morphism $f^*:{\cal O}_{D,q}\rightarrow {\cal O}_{F,p}$.  $\hat f^*:\hat{\cal O}_{D,q}\rightarrow
\hat{\cal O}_{F,p}$ is the $\bold k$-algebra homomorphism $\hat f^*:{\bold k}[[\overline v,\overline w]]\rightarrow {\bold k}[[x,z]]$ given by
$\overline v=z$, $\overline w=\phi(0,z)+x^d$. Thus $\hat{\cal O}_{F,p}$ is finite over $\hat{\cal O}_{D,q}$. It follows that
$f^*:{\cal O}_{D,q}\rightarrow {\cal O}_{F,p}$ is quasi-finite, and thus ${\cal O}_{F,p}\cong{\cal O}_{D,q}$ by Zariski's Main Theorem.
In particular, $\{p\}=f^{-1}(q)\cap F$. We thus have that the only component of the fundamental locus of $f$ through $q$ is the
algebraic curve $C=f(E)$, which has  analytic local equations $\overline u=\overline w-\phi(\overline u,\overline v)=0$ at $q$. Thus $C$ is nonsingular at $q$.

If $E'$ is a component of the exceptional locus of $f$ such that $q\in f(E')$, we must have that $f(E')=C$. Now let $u,v,w$ be
permissible parameters at $q$ such that $u=w=0$ are local equations of $C$ at $q$. We see that $u,v,w$ must have a form (\ref{eqTF01})
or (\ref{eqTF02}) of Definition \ref{torf} for all $p\in f^{-1}(q)$.

\end{pf}

\begin{Lemma}\label{Remark368} Let notation be as in Definition \ref{Def1}, and suppose that $X$ and $Y$
are projective.
\begin{enumerate}
\item[1.] If $q\in Y$ is a 2-point (or a 3-point) and for all $p\in f^{-1}(q)$ there exist permissible parameters $u_p,v_p,w_p$ at $q$ which satisfy one of 2 (a) -- 2 (c)   of Definition
\ref{Def1} (or one of (\ref{eqTF1}) -- (\ref{eqTF3}) of Definition \ref{torf}) at $p$, then any
permissible parameters $u,v,w$  at $q$ satisfy  2  (or 1) of Definition \ref{Def1} for all $p\in f^{-1}(q)$.
\item[2.] If $f:X\rightarrow Y$ is prepared and $q\in Y$ is a 1-point,  
 then there exist permissible parameters $u,v,w$ at $q$ such that 
 for all $p\in f^{-1}(q)$, $u,v$ are
 toroidal forms in local equations at $p$ of the toroidal structure.
\end{enumerate}
\end{Lemma}
\begin{pf}  1 follows from the definitions, and a local calculation.
The most difficult case to verify is when $q\in Y$ is a 2-point, $u,v,w$ are permissible parameters at $q$,  $p\in f^{-1}(p)$ is a 2-point, and $u_p,v_p,w_p$ have a form 2 (c) at $p$,
$$
u_p=(\overline x^a\overline y^b)^k,
v_p=(\overline x^a\overline y^b)^l(\gamma(\overline x^a\overline y^b,\overline z)+\overline x^c\overline y^d),
w_p=\overline z.
$$
First observe that if $u=v_p,v=u_p, w=w_p$, we can find regular parameters $\overline x,\overline y,\overline z$ in $\hat{\cal O}_{X,p}$ such that $u,v,w$ have a form 2 (c) with respect to $\overline x, \overline y, \overline z$.

By the formal implicit function theorem, we have reduced to the case where there exist unit series $\gamma_1,\gamma_2,\gamma_3$ in the variables $u,v,w$ and a series $\lambda$ in $u,v$ (with no constant term)
such that 
$$
u=\gamma_1u_p, v=\gamma_2v_p, w=\gamma_3(w_p+\lambda).
$$
There exist unit series $\alpha_1,\alpha_2\in\hat{\cal O}_{X,p}$ such that if $x_1=\alpha_1\overline x$, $y_1=\alpha_2\overline y$, $z_1=w$, then $x_1,y_1,z_1$ are regular parameters in $\hat{\cal O}_{X,p}$ such that 
\begin{equation}\label{eq407}
u=\gamma_1u_p=(x_1^ay_1^b)^k, 
v=(x_1^ay_1^b)^l(\gamma_2\gamma_1^{-\frac{l}{k}}\gamma(\gamma_1^{-\frac{1}{k}}x_1^ay_1^b,\overline z)+x_1^cy_1^d),
w=z_1.
\end{equation}
We have
$$
\overline x^a\overline y^b=\gamma_1^{-\frac{1}{k}}x_1^ay_1^b
$$
and
$$
\overline x^c\overline y^d=\gamma_2^{-1}\gamma_1^{\frac{l}{k}}x_1^cy_1^d.
$$
There exist series $g,h,h_1$ such that
$$
\begin{array}{ll}
\overline z&=w_p=g(u,w)+vh(u,v,w)\\
&=g(u,w)+(x_1^ay_1^b)^lh_1(x_1^ay_1^b,x_1^cy_1^d,z_1,\overline x^a\overline y^b,\overline x^c\overline y^d,\overline z)
\end{array}
$$
and there exists a series $h_2$ such that 
$$
v=(x_1^ay_1^b)^lh_2(x_1^ay_1^b,x_1^cy_1^d,z_1,\overline x^a\overline y^b,\overline x^c\overline y^d,\overline z).
$$
By iteration, we see that there exists a series  $g_2$ such that 
\begin{equation}\label{eq395}
\overline z\equiv g_2(x_1^ay_1^b,z_1)\text{ mod }(x_1^{c+1}y_1^{d+1}),
\end{equation}
and there exists a series $h_3$ such that 
\begin{equation}\label{eq406}
v\equiv  h_3(x_1^ay_1^b,z_1)\text{ mod }(x_1^{c+1}y_1^{d+1}).
\end{equation}
Substituting equations (\ref{eq395}) and (\ref{eq406}) into
$\gamma_2\gamma_1^{-\frac{l}{k}}\gamma(\gamma_1^{-\frac{1}{k}}x_1^ay_1^b,\overline z)$ in (\ref{eq407}), we see that there exist unit series
$\overline\gamma$ and $\overline\tau$ such that 
$$
v=(x_1^ay_1^b)^l(\overline\gamma(x_1^ay_1^b,z_1)+x_1^cy_1^d\overline\tau).
$$ 
Finally, we can find  unit series $\overline \alpha_1,\overline \alpha_2$ in $x_1,y_1,z_1$ such that if we set
$$
x=x_1\overline \alpha_1, y=y_1\overline \alpha_2, z=z_1,
$$
then $u,v,w$ have an expression of the form 2 (c) in terms of $x,y,z$.

 2 follows from Lemma \ref{Lemma141}.
\end{pf}

\begin{Definition}\label{Def274} (\cite{KKMS}, \cite{AK})
A normal variety $\overline X$ with a SNC divisor $D_{\overline X}$ on $\overline X$ is called toroidal if for every point $p\in \overline X$ there exists an affine
toric variety $X_{\sigma}$, a point $p'\in X_{\sigma}$ and an isomorphism of $\bold k$-algebras
$$
\hat{\cal O}_{\overline X,p}\cong \hat{\cal O}_{X_{\sigma},p'}
$$
such that the ideal of $D_{\overline X}$ corresponds to the ideal of $X_{\sigma}-T$ (where $T$ is the  torus in $X_{\sigma}$). Such a pair
$(X_{\sigma},p')$ is called a local model at $p\in \overline X$. $D_{\overline X}$ is called a toroidal structure on $\overline X$.

A dominant morphism $\Phi:\overline X\rightarrow \overline Y$ of toroidal varieties with SNC divisors $D_{\overline Y}$
on $\overline Y$ and  $D_{\overline X}=\Phi^{-1}(D_{\overline Y})$ on $\overline X$, 
is called toroidal at $p\in\overline X$, and we will say that $p$ is a toroidal point of $\Phi$ if with $q=\Phi(p)$, there exist local models
$(X_{\sigma},p')$ at $p$, $(Y_{\tau},q')$ at $q$ and a toric morphism $\Psi:X_{\sigma}\rightarrow Y_{\tau}$ such that the following
diagram commutes:
$$
\begin{array}{rll}
\hat{\cal O}_{\overline X,p}&\leftarrow &\hat{\cal O}_{X_{\sigma},p'}\\
\hat\Phi^*\uparrow&&\hat\Psi^*\uparrow\\
\hat{\cal O}_{\overline Y,q}&\leftarrow&\hat{\cal O}_{Y_{\tau},q'}.
\end{array}
$$
$\Phi:\overline X\rightarrow \overline Y$ is called toroidal (with respect to $D_{\overline Y}$ and $D_{\overline X}$) if $\Phi$ is toroidal at all $p\in \overline X$.
\end{Definition}
The following is the list of toroidal forms for a dominant morphism $f:X\rightarrow Y$
of nonsingular 3-folds with toroidal structure $D_Y$ and $D_X=f^{-1}(D_X)$. Suppose that $p\in D_X$, $q=f(p)\in D_Y$,
and $f$ is toroidal at $p$.
Then there exist permissible parameters $u,v,w$ at $q$ and permissible parameters $x,y,z$ for $u,v,w$ at $p$ such that
one of the following forms hold:
\begin{enumerate}
\item[1.] $p$ is a 3-point and $q$ is a  3-point,
$$
\begin{array}{ll}
u&=x^ay^bz^c\\
v&=x^dy^ez^f\\
w&=x^gy^hz^i,
\end{array}
$$
where $a,b,d,e,f,g,h,i\in{\bf N}$ and
$$
\text{Det}\left(\begin{array}{lll}
a&b&c\\
d&e&f\\
g&h&i
\end{array}\right)\ne 0.
$$
\item[2.] $p$ is a 2-point and $q$ is a 3-point,
$$
\begin{array}{ll}
u&=x^ay^b\\
v&=x^dy^e\\
w&=x^gy^h(z+\alpha)
\end{array}
$$
with $0\ne \alpha\in \bold k$ and $a,b,d,e,g,h\in{\bf N}$ satisfy $ae-bd\ne 0$.
\item[3.] $p$ is a 1-point and $q$ is a 3-point,
$$
\begin{array}{ll}
u&=x^a\\
v&=x^d(y+\alpha)\\
w&=x^g(z+\beta)
\end{array}
$$
with $0\ne \alpha,\beta\in \bold k$, $a,d,g>0$.
\item[4.] $p$ is a 2-point and $q$ is a 2-point,
$$
\begin{array}{ll}
u&=x^ay^b\\
v&=x^dy^e\\
w&=z
\end{array}
$$
with $ae-bd\ne 0$
\item[5.] $p$ is a 1-point and $q$ is a 2-point,
$$
\begin{array}{ll}
u&=x^a\\
v&=x^d(y+\alpha)\\
w&=z
\end{array}
$$
with $0\ne \alpha\in \bold k$, $a,d>0$.
\item[6.] $p$ is a 1-point and $q$ is a 1-point,
$$
\begin{array}{ll}
u&=x^a\\
v&=y\\
w&=z
\end{array}
$$
with $a>0$.
\end{enumerate}

\begin{Definition}\label{Def247} Suppose that $f:X\rightarrow Y$ is a
prepared morphism. Then $D_X$ is  cuspidal for $f$ if:
\begin{enumerate}
\item[1.] If $E$ is a component of $D_X$ which does not contain a 3-point then $f$ is toroidal
in a Zariski open neighborhood of $E$.
\item[2.] If $C$ is a 2-curve of $X$ which does not contain a 3-point then $f$ is toroidal
in a Zariski open neighborhood of $C$.
\end{enumerate}
\end{Definition}

\begin{Definition}\label{Def221}
 Suppose that $f:X\rightarrow Y$ is prepared, and $p\in X$ is a 3-point.
Suppose that $u,v,w$ are permissible parameters at $q=f(p)$.   Then   there is an expression (after possibly permuting
$u,v,w$ if $q$ is a 3-point) 
\begin{equation}\label{eq16}
\begin{array}{ll}
u&=x^ay^bz^c\\
v&=x^dy^ez^f\\
w&=\sum_{i\ge 0} \alpha_iM_i +N
\end{array}
\end{equation}
where $x,y,z$ are permissible parameters at $p$ for $u,v,w$, $\text{rank}(u,v)=2$,
 the sum in $w$ is over (possibly infinitely many) monomials $M_i$ in $x,y,z$ such that
$\text{rank}(u,v,M_i)=2$,
$\text{deg}(M_i)\le\text{deg}(M_j)$ if $i<j$, $N$ is a monomial in $x,y,z$ such that 
$\text{rank}(u,v,N)=3$
and $N\not\,\mid M_i$ for any $M_i$ in the series $\sum\alpha_i M_i$.

If $q$ is a 3-point and $u,v,w$ is not a monomial form
(at $p$), we necessarily have (since $f$ is prepared) that 
\begin{equation}\label{eq224}
\sum \alpha_iM_i=M_0\gamma
\end{equation}
 where
 $\gamma$
is a unit series in the monomials $\frac{M_i}{M_0}$ (in $x,y,z$)  such that $$
\text{rank}(u,v,M_0)=\text{rank}(u,v,\frac{M_i}{M_0})=2
$$
for all $i$, and $M_0\mid N$.

If $q$ is a 3-point and $u,v,w$ have a monomial form at $p$, so that $w=N$, define $\tau(p)=\tau_f(p)=-\infty$. Otherwise,
define a group $H_p=H_{f,p}$ as follows.
The Laurent monomials  in $x,y,z$ form a group under multiplication.
We define $H_p=H_{f,p}$ to be the subgroup generated by $u,v$ and the terms $M_i$
appearing in the expansion (\ref{eq16}). We will write the group $H_p$
additively as:
$$
H_p=H_{f,p}={\bf Z}u+{\bf Z}v+\sum{\bf Z}M_i.
$$
Define a subgroup $A_p$ of $H_p$ by:
$$
A_p=A_{f,p}=\left\{\begin{array}{ll} {\bf Z}u+{\bf Z}v+{\bf Z}M_0&\text{ if $q$ is 3-point}\\
{\bf Z}u+{\bf Z}v&\text{if $q$ is a 2-point}.
\end{array}\right.
$$
Define
$$
L_p=L_{f,p}=H_p/A_p,
$$
$$
\tau(p)=\tau_f(p)=\mid L_p\mid.
$$
\end{Definition}

Observe that $\tau(p)<\infty$ in Definition \ref{Def221}, since $H_p$ is a finitely generated
group, and $H_p/A_p$ is a torsion group.

We define 
$$
\tau(X)=\tau_f(X)=\text{max}\{\tau_f(p)\mid p\in X\text{ is a 3-point}\}.
$$
\begin{Lemma}\label{Lemma0}
$\tau_f(p)$ is independent of choice of permissible  parameters $u,v,w$ at $q=f(p)$ and permissible parameters $x,y,z$ at $p$ for $u,v,w$. 
\end{Lemma}

\begin{pf} 
Suppose that $q\in Y$ is a 3-point.

The condition $\tau(p)=-\infty$ is independent of permuting
$u,v,w$, multiplying $u,v,w$ by units in $\hat{\cal O}_{Y,q}$ and multiplying $x,y,z$ by units in $\hat{\cal O}_{X,p}$ so that the conditions
of (\ref{eq16}) hold. Thus $\tau(p)=-\infty$  is independent of choice of
permissible parameters at (the 3-points) $q$ and $p$.

Suppose that $u,v,w$ are permissible parameters at $q$ and $x,y,z$ are permissible parameters
at $p$ for $u,v,w$ satisfying (\ref{eq16}). Let $\tau$ be the computation of $\tau(p)$
for these variables.

Suppose that $\tilde u,\tilde v,\tilde w$ is another set of permissible parameters at $q$,
and $\tilde x,\tilde y,\tilde z$ are permissible parameters for $\tilde u,\tilde v,\tilde w$ at
$p$,  satisfying (\ref{eq16}). Let $\tau_1$ be the computation for $\tau(p)$ with respect
to these variables. We must show that $\tau=\tau_1$. We may assume that $\tau\ge 0$ and 
$\tau_1\ge 0$.

Since $p$ and $q$ are 3-points, $\tilde u,\tilde v,\tilde w$ can be obtained from $u,v,w$ by
permuting the variables $u,v,w$ and then multiplying $u,v,w$ by unit series (in $u,v,w$). $\tilde x,\tilde y,\tilde z$ can be obtained from $x,y,z$ by permuting $x,y,z$ and then multiplying $x,y,z$ by unit series (in $x,y,z$).

We then reduce to proving the following:
\begin{enumerate}
\item[1.] Suppose that $\tilde u,\tilde v,\tilde w$ is a permutation of $u,v,w$ such that
$\tilde u,\tilde v$ are toroidal forms at $p$. Then there exist permissible parameters
$\tilde x,\tilde y,\tilde z$ at $p$ for $\tilde u,\tilde v,\tilde w$ such that a form
(\ref{eq16}) holds, and $\tau_1=\tau$.
\item[2.] Suppose that $\tilde u,\tilde v,\tilde w$ are obtained from $u,v,w$ by multiplying $u,v,w$
by unit series. Then there exist permissible parameters $\tilde x,\tilde y,\tilde z$ at $p$
for $\tilde u,\tilde v,\tilde w$ such that a form (\ref{eq16}) holds, and $\tau_1=\tau$.
\item[3.] Suppose that $u=\tilde u$, $v=\tilde v$ and $w=\tilde w$ and $x,y,z$, $\tilde x,\tilde y,\tilde z$ are two sets of permissible parameters for $u,v,w$. Then $\tau_1=\tau$.
\end{enumerate}

We now verify 1. The case when $\tilde u=v$, $\tilde v=u$, $\tilde w=w$ is immediate. We
will verify the case when $\tilde u=w$, $\tilde v=v$ and $\tilde w=u$.  Since the
symmetric group $S_3$ is generated by the permutations $(12)$ and $(13)$, the remaining 
cases of 1 will follow. 

Since $\tau\ge 0$, and $\tilde u,\tilde v,\tilde w$ have a form (\ref{eq16}) at $p$, we have $\text{rank}(v,M_0)=2$.
Since $w,v$ are toroidal forms at $p$,   there exist permissible parameters $\tilde x,\tilde y,
\tilde z$ at $p$ such that $w,v,u$ (in this order) have an expression of the form of
(\ref{eq16}) in terms of $\tilde x,\tilde y,\tilde z$. We will show that there is an isomorphism of the corresponding group
$\tilde H$ (computed for these variables) and the group $H$ computed for $u,v,w$ and
$x,y,z$ which takes the corresponding group $\tilde A$ to $A$.

With the notation of (\ref{eq224}), we have
$$
w=M_0(\gamma+\overline N_0)
$$
where $\overline N_0=\frac{N}{M_0}$ is a monomial in $x,y,z$.

Set $\overline M_0=M_0$, $\overline M_i=\frac{M_i}{M_0}$ for 
all $M_i$ appearing in the series $\sum_{i\ge 1}\alpha_i M_i$.
Thus $\gamma=\alpha_0+\sum_{i\ge 1}\alpha_i\overline M_i$.
There exist $a_i,b_i,c_i\in{\bf N}$ such that
$$
\overline M_i=x^{a_i}y^{b_i}z^{c_i}
$$
for $0\le i$.
$$
\begin{array}{ll}
H&={\bf Z}u+{\bf Z}v+\sum_{i\ge 0}{\bf Z}M_i\\
&={\bf Z}u+{\bf Z}v+\sum_{i\ge 0}{\bf Z}\overline M_i.
\end{array}
$$
Define a finite sequence
$$
1=\mu(1)<\mu(2)<\cdots<\mu(\overline r)
$$
(for appropriate $\overline r$) so that
$$
\sum_{i=1}^{\mu(\overline r)}{\bf Z}\overline M_i=\sum_{i=1}^{\infty}{\bf Z}\overline M_i,
$$
$$
\sum_{i=1}^{\mu(j)}{\bf Z}\overline M_i=\sum_{i=1}^n{\bf Z}\overline M_i
$$
if $\mu(j)<n<\mu(j+1)$ and
$$
\sum_{i=1}^{\mu(j)}{\bf Z}\overline M_i\ne\sum_{i=1}^{\mu(j+1)}{\bf Z}\overline M_i.
$$
Set 
$$
G_j=\sum_{i=1}^{\mu(j)}{\bf Z}\overline M_i.
$$

There exist $k_i,l_i \in {\bf Z}$ and $e_i\in {\bf N}$ with $\text{gcd}(e_i,k_i,l_i)=1$ such that
$\overline M_i^{e_i}=u^{k_i}v^{l_i}$ for $i\in\{0,\mu(1),\ldots,\mu(\overline r)\}$, and there exist $g,h,i\in {\bf N}$
such that $\overline N_0=x^gy^hz^i$. Thus $H$ is (isomorphic to) the subgroup of ${\bf Z}^3$
generated by 
$$
\begin{array}{l}
\delta_1=(a,b,c),\delta_2=(d,e,f),\epsilon_0=(a_0,b_0,c_0),\\
\epsilon_{\mu(1)}=(a_{\mu(1)},b_{\mu(1)},c_{\mu(1)}),\ldots,
\epsilon_{\mu(\overline r)}=(a_{\mu(\overline r)},b_{\mu(\overline r)},c_{\mu(\overline r)}).
\end{array}
$$ 
 $A$ is the subgroup with generators $\delta_1,\delta_2$ and $\epsilon_0$.

Since $\text{rank}(v,M_0)=2$, we can make a change of variables
$$
x=\overline x\lambda_1,
y=\overline y\lambda_2,
z=\overline z\lambda_3
$$
where $\lambda_i=(\gamma+\overline N_0)^{\beta_i}$ for some $\beta_i\in{\bf Q}$,
 so that
$$
\begin{array}{ll}
\lambda_1^{a_0}\lambda_2^{b_0}\lambda_3^{c_0}&=(\gamma+\overline N_0)^{-1}\\
\lambda_1^{d}\lambda_2^e\lambda_3^f&=1\\
\lambda_1^a\lambda_2^b\lambda_3^c&=(\gamma+\overline N_0)^t
\end{array}
$$
for some $t\in{\bf Q}$.
Since $u,v,w$ are algebraically independent in $\hat{\cal O}_{X,p}$ (by
Zariski's subspace theorem, Theorem 10.14 \cite{Ab}), we have that $t\ne 0$, and 
$$
\begin{array}{ll}
w&=\overline x^{a_0}\overline y^{b_0}\overline z^{c_0}\\
v&=\overline x^d\overline y^e\overline z^f\\
u&=\overline x^a\overline y^b\overline z^c(\gamma^t+\overline x^g\overline y^h\overline z^i\gamma_2)
\end{array}
$$
where $\gamma_2$ is a unit series in $\overline x,\overline y,\overline z$.

Set $\tilde M_i=\overline x^{a_i}\overline y^{b_i}\overline z^{c_i}$ for $i\ge 1$.
$$
\gamma(x,y,z)^t\equiv \sum_{j=0}^{\infty}q_j\text{ mod }(\overline x,\overline y,\overline z)\overline N_0
$$
where each $q_j$ is a series in monomials of degree $j$ in $\{\tilde M_i\mid i\ge 1\}$, and
$$
q_0=\alpha_0^t,q_1=\sum_{i=1}^{\infty}\sigma_i\tilde M_i
$$
with 
$$
 \sigma_i=t\alpha_0^{t-1+a_i\beta_1+b_i\beta_2+c_i\beta_3}\alpha_i
$$
for $i\ge 1$. Let
$$
\omega=\alpha_0^t+\sum_{i=1}^{\infty}\sigma_i\tilde M_i.
$$
There exists a unit series $\gamma_2$ such that 
$$
u=\overline x^a\overline y^b\overline z^c(\omega+\overline x^g\overline y^h\overline z^i\gamma_2).
$$
We now see that the coefficient of $\tilde M_{\mu(j)}$ for $1\le j\le\overline r$ in the
expansion of $\omega$ as a series in $\overline x,\overline y,\overline z$ is $\sigma_{\mu(j)}$.
If not, there would exist a relation $\tilde M_{\mu(j)}=\tilde M_{i_1}\cdots\tilde M_{i_n}$ for some $n>1$. Thus $i_1,\ldots,i_n<\mu(j)$ and $G_{\mu(j)}=G_{\mu(j-1)}$, a contradiction.

 Since 
 $$
 \text{Det}\left(\begin{array}{lll}
 a&b&c\\
 d&e&f\\
 g&h&i
 \end{array}
 \right)\ne 0,
 $$
  there exists
a change of variables
$$
\overline x=\tilde x\phi_1,
\overline y=\tilde y\phi_2,
\overline z=\tilde z\phi_3
$$
where $\phi_1,\phi_2,\phi_3$ are unit series in $\overline x,\overline y,\overline z$ such that
$$
\begin{array}{ll}
\phi_1^{\frac{a}{ e_0}}\phi_2^{\frac{b}{ e_0}}\phi_3^{\frac{ c}{ e_0}}&=1\\
\phi_1^{\frac{d}{ e_0}}\phi_2^{\frac{e}{ e_0}}\phi_3^{\frac{f}{e_0}}&=1\\
\phi_1^g\phi_2^h\phi_3^i&=\gamma_2^{-1}.
\end{array}
$$
Since $e_i(a_i,b_i,c_i)=k_i(a,b,c)+l_i(d,e,f)$ for $i\ge 0$,
we have an expression
$$
\begin{array}{ll}
w&=\tilde x^{a_0}\tilde y^{b_0}\tilde z^{c_0}\\
v&=\tilde x^d\tilde y^e\tilde z^f\\
u&=\tilde x^a\tilde y^b\tilde z^c(\omega+\tilde x^g\tilde y^h\tilde z^i)
\end{array}
$$
and $\tilde M_i=\tilde x^{a_i}\tilde y^{b_i}\tilde z^{c_i}\eta_i$
for some $e_i$-th root of unity $\eta_i$ in $\bold k$
for all $1\le i$.

 $\tilde H$ is thus (isomorphic to) the subgroup of ${\bf Z}^3$ with generators 
$$
\overline \delta_1,\overline \delta_2,
\overline \epsilon_0,\overline \epsilon_{\mu(1)},\ldots,\overline \epsilon_{\mu(\overline r)},
$$
defined by 
$$
\begin{array}{l}
\overline\delta_1=(a_0,b_0,c_0),\overline \delta_2=(d,e,f),\overline \epsilon_0=(a,b,c),\\
\overline \epsilon_{\mu(1)}=(a_{\mu(1)},b_{\mu(1)},c_{\mu(1)}),\ldots,
\overline \epsilon_{\mu(\overline r)}=(a_{\mu(\overline r)},b_{\mu(\overline r)},c_{\mu(\overline r)}).
\end{array}
$$ 

$\tilde A$ is the subgroup of $\tilde H$ generated by $\overline\delta_1,\overline\delta_2$ and
$\overline\epsilon_0$.
Thus, we  have an isomorphism of $H$ with $\tilde H$ which takes $A$ to $\tilde A$ by mapping
$\delta_1$ to $\overline \epsilon_0$, $\delta_2$ to $\overline\delta_2$,
$\epsilon_0$ to $\overline \delta_1$ and $\epsilon_i$ to $\overline\epsilon_i$ for
$i\in\{\mu(1),\ldots,\mu(\overline r)\}$.

We have thus completed the verification of 1. The verification of 2 and 3 follows from simpler calculations. We thus obtain
the conclusions of the lemma when $q$ is a 3-point.

Now assume that $q\in Y$ is a 2-point. $\tau(p)$ is independent of interchanging $u$ and $v$, multiplying $u$ and $v$ by unit series in $\hat{\cal O}_{Y,q}$, and permuting $x,y,z$ and multiplying $x,y,z$ by unit series in $\hat{\cal O}_{X,p}$, so that the conditions of (\ref{eq16}) hold. If we replace
$w$ by $w'\in\hat{\cal O}_{Y,q}$ so that $u,v,w'$ are permissible
parameters at $q$, then by the formal implicit function theorem, there exists a unit series $\alpha(u,v,w)\in\hat{\cal O}_{Y,q}$ and a series $\beta(u,v)\in {\bold k}[[u,v]]$ such that $w=\alpha^{-1}(w'-\beta(u,v))$.

There exists a series $\phi$ in $x,y,z$ such that
$$
\alpha(u,v,w)=\alpha(u,v,\sum\alpha_iM_i+N)=\alpha(u,v,\sum\alpha_iM_i)+N\phi.
$$
Then
$$
w'=\beta(u,v)+\alpha(u,v,\sum\alpha_iM_i)(\sum\alpha_iM_i)+
N[\alpha(u,v,\sum\alpha_iM_i)+(\sum\alpha_iM_i)\phi+N\phi].
$$
Now as in the calculation we make in the verification that $\tau_1=\tau$ in the case when
$q$ is a 3-point, we see that there exist permissible parameters $\tilde x,\tilde y,\tilde z$ for $u,v,w'$ at $p$ such that $\tau_1=\tau$ (where $\tau$ is computed for $u,v,w$ and $x,y,z$ and $\tau'$ is computed for $u,v,w'$ and $\tilde x,\tilde y,\tilde z$).
Thus $\tau(p)$ is
independent of choice of permissible parameters at $q$ and $p$ when $q$ is a 2-point.

\end{pf}

\begin{Lemma}\label{Lemma353} Suppose that $X$ is a nonsingular 3-fold with
SNC divisor $D_X$, defining a toroidal structure on $X$. Suppose that ${\cal I}$ is an
ideal sheaf on $X$ which is locally generated by monomials in local equations of components
of $D_X$. Then there exists a sequence of blow ups of 2-curves $\Phi_1:X_1\rightarrow X$ such
that ${\cal I}{\cal O}_{X_1}$ is an invertible ideal sheaf. If ${\cal I}$ is locally generated by two equations, then $\Phi_1$ is an isomorphism away
from the support of ${\cal I}$.
\end{Lemma}

This lemma is an extension of Lemma 18.18 \cite{C2}, and is generalized to all dimensions in \cite{G}.  

\begin{pf} 
$X$ has a cover by affine open sets $U_1,\ldots, U_n$ such that there exist $g_{i,1},\ldots, g_{i,l}\in\Gamma(U_i,{\cal O}_X)$
such that $g_{i,j}=0$ are local equations in $U_i$ of irreducible components of $D_X$, and there exist $f_{i,1},\ldots, f_{i,m(i)}\in \Gamma(U_i,{\cal O}_X)$ such that the $f_{i,j}$ are monomials in the $g_{i,k}$ and $\Gamma(U_i,{\cal I})=(f_{i_,1},\ldots,f_{i,m(i)})$.

Let $D_{ij}$ be an effective divisor supported on the components of $D_X$ such that there is equality of divisors $D_{ij}\cap U_i=(f_{i,j})\cap U_i$.
Let ${\cal I}_i\subset{\cal O}_X$ be the ideal sheaf which is locally generated by local equations of $D_{i,1},\ldots, D_{i,m(i)}$.
By construction, ${\cal I}_i\mid U_i={\cal I}\mid U_i$ for all $i$.

We will show that for an ideal sheaf of the form ${\cal I}_1$, there exists a sequence of blow ups of 2-curves, $\pi:X_1\rightarrow X$
such that ${\cal I}_1{\cal O}_{X_1}$ is invertible. Since ${\cal I}_i{\cal O}_{X_1}$ are locally generated by local equations of
$\pi^*(D_{i,1}),\ldots,\pi^*(D_{i,m(i)})$, there exists $\pi_2:X_2\rightarrow X$ which is a sequence of blow ups of 2-curves such that
${\cal I}_i{\cal O}_{X_2}$ is invertible for all $i$. Since ${\cal I}_i{\cal O}_{X_2}\mid \pi_2^{-1}(U_i)={\cal I}{\cal O}_{X_2}\mid\pi_2^{-1}(U_i)$ for all $i$, ${\cal I}{\cal O}_{X_2}$ is invertible.

We may now suppose that
there exists $n>0$ and effective divisors $D_1,\ldots, D_n$ on $X$ whose supports are
unions of components of $D_X$, such that ${\cal I}$ is locally generated by local equations of
$D_1,\ldots, D_n$.

First suppose that $n=2$. Suppose that $p\in X$ is a general point of a 2-curve $C$. Let 
$x=0$, $y=0$ be local equations of the components of $D_X$ containing $p$. $x=y=0$ are local equations of $C$ at $p$. Then there exist $a,b,c,d\in{\bf N}$ such that $D_1$ is
defined near $p$ by the divisor of $x^ay^b$, and $D_2$ is defined near $p$ by the divisor of
$x^cy^d$. Define
$$
\omega(C)=\left\{
\begin{array}{ll}
\text{max}\{(|a-c|,|b-d|),(|b-d|,|a-c|)\}&\text{if $a-c$, $b-d$ are nonzero}\\
&\text{and have opposite
signs},\\
-\infty&\text{otherwise}
\end{array}
\right.
$$
Here the maximum  is computed in the lexicographic order. We see that the stalk ${\cal I}_p$
is invertible if and only if $\omega(C)=-\infty$.

Further, if $\omega(C)=-\infty$ for all 2-curves $C$ of $X$, then ${\cal I}$ is invertible, as
follows since the divisors  $D_1$ and $D_2$ are given locally at a 3-point $p$ by the divisors of monomials
$x^ay^bz^c$ and $x^dy^ez^f$ where $xyz=0$ is a local equation of $D_X$ at $p$. $a-d$ and $b-e$ have the same signs,  $a-d$, $c-f$ have the same signs, and $b-e$, $c-f$ have the same signs, so $x^ay^bz^c|x^dy^ez^f$ or
$x^dy^ez^f\mid x^ay^bz^c$.

Now define
$$
\overline\omega(X)=\text{max}\{\omega(C)\mid C\text{ is a 2-curve of $X$}\}.
$$
We have seen that ${\cal I}$ is invertible if and only if $\overline\omega(X)=-\infty$.
Suppose that $\overline\omega(X)\ne -\infty$ and $C$ is a 2-curve of $X$ such that $\omega(C)=\overline\omega(X)$. Let $\pi:X_1\rightarrow X$ be the blow up of $C$. Let
$D_{X_1}=\pi^{-1}(D_X)=\pi^*(D_X)_{red}$, $D_1'=\pi^*(D_1)$, $D_2'=\pi^*(D_2)$.

We can define the function $\omega$ for 2-curves on $X_1$, relative to $D_1'$ and $D_2'$, and define
$\overline\omega(X_1)$.

By a local calculation (as shown in the proof of Lemma 18.18 \cite{C2}) we see that $\omega(C_1)<\overline\omega(X)$ if $C_1$ is a 2-curve which is contained in the exceptional
divisor of $\pi$.

Suppose that $C_1,\ldots, C_r$ are the 2-curves $C$ on $X$ such that $\omega(C)=\overline \omega(X)$. We obtain a reduction $\overline\omega(X_1)<\overline\omega(X)$ after blowing up 
(the strict transforms of) these $r$ curves. By induction on $\overline\omega(X)$, we must obtain
that ${\cal I}{\cal O}_{X_2}$ is invertible after an appropriate sequence of blow ups of 2-curves $X_2\rightarrow X$.

Now suppose that ${\cal I}$ is locally generated by local equations of $D_1,\ldots, D_n$ (with $n>2$). Let ${\cal I}_1\subset {\cal I}$
be the ideal sheaf which is locally generated by local equations of $D_1$ and $D_2$. We have seen that there exists a sequence of
blow ups of 2-curves $\pi_1:X_1\rightarrow X$ such that ${\cal I}_1{\cal O}_{X_1}$ is invertible. Thus there exists a divisor $\overline D$
on $X_1$ whose support is a union of components of $D_{X_1}$ such that ${\cal I}_1{\cal O}_{X_1}$ is locally generated by a local equation
of $\overline D$.

Let $\overline D_i=\pi_1^*(D_i)$ for $3\le i\le n$. Then ${\cal I}{\cal O}_{X_1}$ is locally generated by local equations of the $n-1$
divisors $\overline D,\overline D_3,\ldots, \overline D_n$. By induction, there exists a sequence of blow ups of 2-curves $X_2\rightarrow X$ such that ${\cal I}{\cal O}_{X_2}$ is invertible.

\end{pf}

\begin{Lemma}\label{Lemma423} Suppose that $f:X\rightarrow Y$ is a prepared morphism, and $\overline C$ is a 2-curve in $Y$. then there
exists a sequence of blowups of 2-curves $\Phi:X_1\rightarrow X$ such that ${\cal I}_{\overline C}{\cal O}_{X_1}$ is invertible and $\Phi$
is an isomorphism over $f^{-1}(Y-\overline C)$.
\end{Lemma}
\begin{pf} The Lemma is a consequence of  Lemma \ref{Lemma353}.
\end{pf}

\begin{Lemma}\label{Lemma419} Suppose that $X$ is a nonsingular 3-fold with
SNC divisor $D_X$, defining a toroidal structure on $X$. Suppose that ${\cal I}$ is an
ideal sheaf on $X$ which is locally generated by monomials in local equations of components
of $D_X$. Then there exists a sequence of blow ups of 2-curves and 3-points $\Phi_1:X_1\rightarrow X$ such
that ${\cal I}{\cal O}_{X_1}$ is an invertible ideal sheaf and $\Phi_1$ is an isomorphism away from the support of ${\cal I}$. 
\end{Lemma}

The proof of this lemma follows from the proof of principalization of ideals, as in \cite{BEV} or \cite{BrM} (cf. the proof of Theorem 6.3 \cite{C3})
in the case when the ideal to be principalized is locally generated by monomials in the toroidal structure.

\section{Preparation}
The following theorem is  Theorem 19.11 \cite{C2}, with the additional conclusions that
all 2-curves of $X_2$ contain a 3-point, and all components of $D_{X_2}$ contain a 3-point.

\begin{Theorem}\label{Theorem183}  Suppose that $\Phi:X\rightarrow S$ is a dominant
morphism from a nonsingular 3-fold $X$ to a nonsingular surface $S$ and $D_S$
is a SNC divisor on $S$ such that $D_X=\Phi^{-1}(D_S)$ is a SNC divisor which
contains the singular locus of $\Phi$. Further suppose that every component of $D_X$
contains a 3-point and every 2-curve of $X$ contains a 3-point.
\begin{enumerate}
\item[1.]
Then there exists a sequence of blow ups of possible centers
 $\alpha_1:X_1\rightarrow X$ such that
\begin{enumerate}
\item The fundamental locus of $\alpha_1$ is contained in the union of irreducible
components $E$ of $D_X$ such that $E$ contains a point $p$ such that $\Phi$ is not prepared
at $p$ (Definition 6.5 \cite{C2}).
\item $\Phi_1=\Phi\circ\alpha_1:X_1\rightarrow S$ is prepared
(Definition 6.5 \cite{C2})
\item Each 2-curve of $X_1$ contains a 3-point and each component of $D_{X_1}=\alpha_1^{-1}(D_X)$
contains a 3-point.
\end{enumerate}
\item[2.] Further, there exist sequences of blow ups of possible centers, $\alpha_2:X_2
\rightarrow X_1$ and $\beta:S_1\rightarrow S$ such that:
\begin{enumerate}
\item There is a commutative diagram
$$
\begin{array}{rll}
X_2&\stackrel{\Phi_2}{\rightarrow}&S_1\\
\alpha_2\downarrow&&\downarrow \beta\\
X_1&\stackrel{\Phi_1}{\rightarrow}&S
\end{array}
$$
such that $\Phi_2$ is toroidal,
\item $\beta$ is an isomorphism away from $\beta^{-1}(\Phi_1(Z))$
where $Z$ is the locus where $\Phi_1$ is not toroidal.
\item $\alpha_2$ is an isomorphism away from $\alpha_2^{-1}(\Phi_1^{-1}(\Phi_1(Z)))$
\item Each 2-curve of $X_2$ contains a 3-point and each component of $D_{X_2}$
contains a 3-point.
\end{enumerate}
\end{enumerate}
\end{Theorem}

\begin{pf}  For the proof we need only make some small modifications in the proof of Theorem
19.11 \cite{C2}.

We first prove 1 of the theorem. By Lemma 6.2 \cite{C2}, there exists a commutative diagram
$$
\begin{array}{rll}
X_0&&\\
\alpha_0\downarrow&\searrow\Phi_0&\\
X&\stackrel{\Phi}{\rightarrow}&S
\end{array}
$$
such that $\Phi_0:X_0\rightarrow S$ is a weakly
prepared morphism (Definition 6.1 \cite{C2}), and the fundamental locus of $\alpha_0$ is contained in the locus where
$\Phi$ is not weakly prepared, (which is contained in the locus where $\Phi$ is not
prepared). It is not necessary to blow up points on $S$ since $D_S$, $D_X$ are  SNC divisors. Let $D_{X_0}=\alpha_0^{-1}(D_X)$. By further blowing up of points in the exceptional locus of $\alpha_0$, we may assume
that all components of $D_{X_0}$ contain a 3-point, and all 2-curves of $X_0$ contain a 3-point.

By Theorem 17.2 \cite{C2} there exists a commutative diagram
$$
\begin{array}{rll}
X_1&&\\
\overline \alpha\downarrow&\searrow\Phi_1&\\
X_0&\stackrel{\Phi_0}{\rightarrow}&S
\end{array}
$$
such that $\Phi_1$ is prepared. The algorithm consists of a sequence 
$$
X_1=Y_n\stackrel{\overline\alpha_n}{\rightarrow}Y_{n-1}\rightarrow
\cdots\rightarrow Y_1\stackrel{\overline\alpha_1}{\rightarrow}Y_0=X_0
$$
of blow ups of
points and nonsingular curves  which are possible centers; that is make SNCs with the preimage $D_{Y_i}$ of $D_{X_0}$ and are
contained in a  component $E$ of $D_{Y_i}$ such that $E$ contains a point which is
not prepared for 
$$
\overline\Phi_i=\Phi_0\circ\overline\alpha_1\circ\cdots\circ\overline\alpha_i:Y_i\rightarrow S.
$$
 If $p\in Y_i$ is prepared for $\overline \Phi_i$ then all points
of $\overline\alpha_{i+1}^{-1}(p)$ are prepared for $\overline\Phi_{i+1}$. Thus conditions (a) and (b) of 1 hold.

The strict transforms on $X_1$ of all components of $D_X$ must contain a 3-point
and if $C$ is a 2-curve of $X$ which is not contained in the fundamental locus of
$\alpha_1=\alpha_0\circ\overline\alpha$, then its strict transform on $X_2$ must also contain a
3-point. Thus any components of $D_{X_2}$ which do not contain a 3-point, and 
2-curves of $X_2$ which do not contain a 3-point, are contained in the exceptional locus of $\alpha_1$.

To complete the proof of 1 of the theorem, we need only show that if $p\in X_1$ is a
point in the exceptional locus of $\alpha_1$, and $\alpha^*:\overline X\rightarrow X_1$
is the blow up of $p$, then $\Phi_1\circ\alpha^*:\overline X\rightarrow S$ is prepared.
 This can directly be seen by substituting  local equations for
the blow up of a point into (17)--(20) of Definition 6.5 \cite{C2}.
We can thus construct $\alpha_1:X_1\rightarrow X$ such that $\Phi_1:X_1\rightarrow S$
is prepared and all conditions of 1 hold.

We now verify 2 of the theorem. We first observe that $\Phi_1$ (from the conclusions
of part 1 of this theorem) is strongly prepared (Definition 18.1 \cite{C2}). We now
examine the monomialization algorithm of Chapter 18 \cite{C2} and the toroidalization
algorithm of Chapter 19 \cite{C2}, applied to $\Phi_1:X_1\rightarrow S$.

The monomialization algorithm of Theorem 18.19 and Theorem 18.21 \cite{C2} consists in constructing a commutative diagram 
\begin{equation}\label{eq*1}
\begin{array}{rll}
\tilde X_2&\stackrel{\tilde \pi_2}{\rightarrow}&X_1\\
\tilde \Phi_2\downarrow&&\downarrow\Phi_1\\
\tilde S_1&\stackrel{\tilde \Psi_1}{\rightarrow}&S
\end{array}
\end{equation}
such that $\tilde\Phi_2$ is monomial (all points of $\tilde X_2$ are good for $\tilde \Phi_2$
as defined in Definition 18.5 \cite{C2}). 

(\ref{eq*1}) has a factorization 
\begin{equation}\label{eq*3}
\begin{array}{rrrrrrrrrrr}
\tilde X_2&=&Z_l&\stackrel{\epsilon_l}{\rightarrow}&\cdots&\rightarrow&
Z_1&\stackrel{\epsilon_1}{\rightarrow}&Z_0&=&X_1\\
\tilde\Phi_2\downarrow&&\Omega_l\downarrow&&&&\Omega_1
\downarrow&&\Omega_0\downarrow&&\Phi_1\downarrow\\
\tilde S_1&=&T_l&\stackrel{\delta_l}{\rightarrow}&\cdots&\rightarrow&
T_1&\stackrel{\delta_1}{\rightarrow}&T_0&=&S
\end{array}
\end{equation}
where each $\Omega_i$ is strongly prepared,
each $\delta_{i+1}$ is the   blow up of a point
$q_i$ such that $\Omega_i^{-1}(q_i)$ contains a point $p_i$ at which $\Omega_i$ is not monomial and
$\epsilon_{i+1}$ is a sequence of blow ups of curves which are exceptional to such $q_i$. This step is accomplished by
performing the algorithms of Lemmas 18.16, 18.17 and 18.18 of \cite{C2}.
We will make a minor modification in the algorithm of Theorem 18.19, which will ensure that
all 2-curves of $\tilde X_2$ contain a 3-point, and all components of $D_{\tilde X_2}$ contain a 3-point.

Each map $Z_{i+1}\rightarrow Z_i$ has a factorization 
\begin{equation}\label{eq*2}
Z_{i+1}=\overline Z_m\stackrel{\lambda_m}{\rightarrow}\overline Z_{m-1}
\rightarrow\cdots\stackrel{\lambda_1}{\rightarrow}\overline Z_0=Z_i
\end{equation}
where each $\lambda_{j+1}$ is a blow up of a 2-curve or of a curve $C_j$ which  contains a 1-point, makes SNCs with the preimage $D_{\overline Z_j}$ of $D_X$ on $\overline Z_j$, and is contained in a component of $D_{\overline Z_j}$. To construct (\ref{eq*2}) we successively apply Lemmas 18.16, 18.17, 18.18 of \cite{C2}.

The algorithms of Lemma 18.16 and Lemmas 18.18 \cite{C2} consist of blow ups of 2-curves and the
condition that all 2-curves contain a 3-point and all components of $D_{\overline Z_j}$ contain a 3-point
is preserved by this condition. 

The algorithm of Lemma 18.17 \cite{C2} consists of a sequence of blow ups of curves 
$\lambda_{j+1}:\overline Z_{j+1}\rightarrow \overline Z_j$ of
$C_j\subset D_{\overline Z_j}$
which are not 2-curves, and are contained in the locus where $m_{q_i}{\cal O}_{\overline Z_j}$
is not invertible.  Let $p\in C_j$ be a general point, so that $p$
is a 1-point. There exist permissible parameters $(u,v)$ at $q_i$ and regular parameters
$x,y,z$ in $\hat{\cal O}_{\overline Z_j,p}$ such that a form (185) of Lemma 18.12 \cite{C2} holds, 
\begin{equation}\label{eq181}
u=x^k,
v=x^cy
\end{equation}
with $c<k$. $x=0$ is a local equation of $D_{\overline Z_j}$ at $p$, and $x=y=0$ are local equations
of $C_j$. $u=v=0$ are local equations of $q$ in $T_i$.
 The exceptional divisor of $\lambda_{j+1}$ contains a 2-curve
which is a section over $C_j$. At the two point $p_1\in\lambda_{j+1}^{-1}(p)$, we have regular
parameters $x_1,y_1,z_1$ in $\hat{\cal O}_{\overline Z_{j+1},p}$ such that 
\begin{equation}\label{eq179}
u=x_1^ky_1^k, v=x_1^cy_1^{c+1}.
\end{equation}
If $C_j$ contains a 2-point then all components of $D_{\overline Z_{j+1}}$
contain a 3-point, and all 2-curves of $\overline Z_{j+1}$ contain a 3-point.

Suppose that $C_j$ does not contain a 2-point. Then $u,v$ have an expression of the form (\ref{eq181}) at all points $\overline p\in C_j$.

If $C_j$ does not contain a 2-point,  we modify the algorithm of Lemma 18.17 \cite{C2}, inserting an extra step here,
by performing the blow up $\lambda_{j+1}':\overline Z_{j+1}'\rightarrow \overline Z_{j+1}$ of the point $p_1$.
In this case, all points $p\in C_j$ are general points, and we we may choose $p=\lambda_{j+1}(p_1)$ to be any point of $C_j$ which is convenient. We will make use of this observation in the proof of Theorem \ref{Theorem190}.
 Points $p_2$
above $p_1$ have regular parameters $(x_2,y_2,z_2)$ such that 
\begin{equation}\label{eq177}
x_1=x_2, y_1=x_2(y_2+\alpha), z_1=x_2(z_2+\beta)
\end{equation}
with $\alpha,\beta\in \bold k$, 
\begin{equation}\label{eq178}
x_1=x_2y_2, y_1=y_2, z_1=y_2(z_2+\beta)
\end{equation}
with $\beta\in \bold k$, or 
\begin{equation}\label{eq228}
x_1=x_2z_2, y_1=y_2z_2, z_1=z_2
\end{equation}
Substituting (\ref{eq177}) into (\ref{eq179}) we have 
\begin{equation}\label{eq180}
u=x_2^{2k}(y_2+\alpha)^k, v=x_2^{2c+1}(y_2+\alpha)^{c+1}.
\end{equation}
Thus if $\alpha\ne 0$ we have a good point of the form (183) of \cite{C2} and
$m_{q_i}{\cal O}_{\overline Z_{j+1}',p_2}$ is invertible.  If $\alpha= 0$, then
$$
u=x_2^{2k}y_2^k, v=x_2^{2c+1}y_2^{c+1}
$$
which is a good point of the form (179) of \cite{C2} and is a form (187) of \cite{C2} if $m_{q_i}{\cal O}_{\overline Z_{j+1}',p_2}$ is not invertible.

Under substitution of (\ref{eq178}) into (\ref{eq179}), we see that
$$
u=x_2^ky_2^{2k}, v=x_2^{c}y_2^{2c+1}
$$
which is a good point of the form (179) of \cite{C2} and is a form (187) of \cite{C2} if $m_{q_i}{\cal O}_{\overline Z_{j+1}',p_2}$ is not invertible.

Under substitution of (\ref{eq228}) into (\ref{eq179}), we obtain
$$
u=x_2^ky_2^kz_2^{2k}, v=x_2^cy_2^{c+1}z_2^{2c+1}
$$
which is a good point of the form (193) of \cite{C2} if $m_{q_i}{\cal O}_{\overline Z_{j+1}',p_2}$ is not invertible.

Observe that the locus of points in $(\lambda_{j+1}')^{-1}(p_1)$ where $m_{q_i}{\cal O}_{\overline Z_{j+1}'}$ is not invertible is a union of 2-curves.

 We now
continue the algorithm as in the the proof of Lemma 18.17 \cite{C2}. As the invariant
$\Omega(C_j)=k-c$ of Lemma 18.17 \cite{C2} which is decreased in the algorithm of Lemma 18.17
\cite{C2} is computed at generic points of curves $C_j$ (which contain a 1-point) and for
which $m_{q_i}{\cal O}_{\overline Z_j,p}$ is not invertible, these invariants are not affected by inserting these
new blow ups of points $\lambda_{j+1}'$ into (\ref{eq*2}). Thus the conclusions of Theorem 18.19 \cite{C2} will
hold, for the modified $\tilde X_2\rightarrow \tilde S_1$, but we may further assume that 
each 2-curve of $\tilde X_2$ contains a 3-point and each component of $D_{\tilde X_2}$ contains
a 3-point.

Theorem 19.9 \cite{C2} and Theorem 19.10 \cite{C2} imply there exists a commutative diagram 
\begin{equation}\label{eq*4}
\begin{array}{rll}
\tilde X_3&\stackrel{\tilde\pi_3}{\rightarrow}&\tilde X_2\\
\tilde \Phi_3\downarrow&&\downarrow\tilde \Phi_2\\
\tilde S_2&\stackrel{\tilde \psi_2}{\rightarrow}&\tilde S_1
\end{array}
\end{equation}
such that $\tilde X_3\rightarrow \tilde S_2$ is toroidal, and 2 (b), 2 (c) of the conclusions of the theorem hold. We will indicate  how we can modify
the proof slightly to ensure that 2 (d) of the conclusions of the theorem holds for $\tilde X_3\rightarrow \tilde S_2$.

The algorithm of Theorem 19.9 \cite{C2} consists of a sequence of blow ups of curves
above $\tilde X_2$ and finitely many blow ups of points over $\tilde S_1$.

In the algorithm, we first construct a diagram  
\begin{equation}\label{eq*5}
\begin{array}{rll}
\tilde X_3'&\stackrel{\pi_3'}{\rightarrow}&\tilde X_2\\
\tilde\Phi_3'\downarrow&&\downarrow\tilde\Phi_2\\
\tilde S_2'&\stackrel{\Psi_2'}{\rightarrow}&\tilde S_1
\end{array}
\end{equation}
which has a factorization by a diagram of the form (\ref{eq*3}), so that
a global invariant $I(\tilde\Phi_3')\le 0$ (this invariant is defined on
page 227 of \cite{C2}).
The factorizations (\ref{eq*2}) of the morphisms of (\ref{eq*3}) consist of a sequence of blow ups of
curves, using first Lemma 18.17 \cite{C2} to blow up curves which contain a
1-point and are possible centers (make SNCs with $D_{\overline Z_j}$ and are contained in a component of $D_{\overline Z_j}$) and then Lemma 18.18 \cite{C2}
to blow up 2-curves.

If a 2-curve is blown up, then
the condition that all 2-curves contains a 3-point is preserved.

Suppose that a curve $C_j$ is blown up which contains a 1-point by $\lambda_{j+1}:\overline Z_{j+1}\rightarrow\overline Z _j$ (in (\ref{eq*2})). This is analyzed 
in Lemma 19.8 \cite{C2}. Let  $p\in C_j$ be a general point. Then a form (185) of \cite{C2} (as in (\ref{eq181}) of our analysis of monomialization)
holds at $p$, and if $p_1\in \lambda_{j+1}^{-1}(p)$ is the 2-point, then a form (\ref{eq179})
holds at $p_1$. 

Assuming that $C_j$ does not contain a 2-point, We now modify the algorithm of Theorem 19.9 \cite{C2} by blowing up
the point $p_1$. Let $\lambda_{j+1}':\overline Z_{j+1}'\rightarrow \overline Z_{j+1}$ be this map. Let $\overline E$ be the
exceptional divisor of $\lambda_{j+1}'$. We see (from (\ref{eq180})) that a form
$$
u=\overline x_2^{2k}, v=\overline x_2^{2c+1}(\overline \alpha+\overline y_2)
$$
with $\overline \alpha\ne 0$ holds at a general point $p_2$ of $\overline E$.
Let 
$$
\tilde \Omega_j= \Omega_i\circ\lambda_1\circ\ldots\circ\lambda_{j+1}
\circ\lambda_{j+1}'.
$$
We have 

\begin{equation}\label{eq182}
I(\tilde\Omega_j,\overline E)=(2c+1)-2k=2(c-k)+1<0
\end{equation}
since $c<k$. We may thus continue the algorithm of Theorem 19.9 \cite{C2}.
We modify (\ref{eq*5}) by adding in these blow ups, $\lambda_{j+1}'$, to achieve the reduction $I(\tilde\Phi_3',E)\le 0$ for all components $E$ of $D_{\tilde X_3'}$ which contain
a 1-point mapping to a 1-point, and so that all components of $D_{\tilde X_3'}$ contain a 3-point, and all 2-curves of $\tilde X_3'$ contain a 3-point.

The algorithm of Theorem 19.10 \cite{C2} consists of a sequence of 
blow ups of curves over $\tilde X_3'$ and points over $\tilde S_2'$. We
construct a commutative diagram 
\begin{equation}\label{eq*6}
\begin{array}{rll}
\tilde X_3&\stackrel{\tilde \pi_3'}{\rightarrow}&\tilde X_3'\\
\tilde\Phi_3\downarrow&&\downarrow\tilde\Phi_3'\\
\tilde S_2&\stackrel{\tilde\Psi_2'}{\rightarrow}&\tilde S_2'
\end{array}
\end{equation}
such that $\tilde \pi_3'$ is toroidal, which has a factorization by a diagram of the form (\ref{eq*3}).
The factorization (\ref{eq*2}) of the morphisms in (\ref{eq*3}) consists of a sequence of blow ups of curves
$C_j$ which are possible centers (make SNCs with $D_{\overline Z_j}$ and are contained in $D_{\overline Z_j}$),
using Lemma 18.17 \cite{C2}.

Suppose that a curve $C_j$ is blown up by $\lambda_{j+1}:\overline Z_{j+1}\rightarrow \overline Z_j$ in equation (\ref{eq*2}).
If $C_j$ contains a 2-point then (assuming that all components of $D_{\overline Z_j}$ contain a 3-point and all 2-curves of $\overline Z_j$ contain a 3-point) all components  of $D_{\overline Z_{j+1}}$ contain a 3-point, and all 2-curves of $\overline Z_{j+1}$ contain a 3-point. 

Suppose that $C_j$ does not contain a 2-point. Then $u,v$ have an expression of the form (\ref{eq181}) at all $p\in C_j$.
Let $p\in C_j$. If
$p_1\in\lambda_{j+1}^{-1}(p)$ is the 2-point, then (\ref{eq179}) holds at
$p_1$. We now modify the algorithm of Theorem 19.10 \cite{C2} by  blowing up  the
2-point $p_1$. Let $\lambda_{j+1}':\overline Z_{j+1}'\rightarrow\overline Z_{j+1}$ be this map. Let $\overline E$ be the exceptional divisor of $\lambda_{j+1}'$. Let
$$
\overline\Omega_j=\Omega_i\circ\lambda_1\circ\cdots\circ \lambda_{j+1}\circ\lambda_{j+1}'.
$$
We have $I(\overline\Omega_j,\overline E)<0$ (as shown in (\ref{eq182})).
Now by Lemma 19.6 \cite{C2} we can continue the algorithm of Theorem 19.10
\cite{C2} to achieve the conclusions of Theorem 19.10 \cite{C2}, with the conclusions
2 (a) - 2 (d) of the conclusions of this theorem.
\end{pf}

\begin{Lemma}\label{Theorem185} Suppose that $f:X\rightarrow Y$ is a birational
projective morphism of nonsingular 3-folds with toroidal structure, defined by SNC divisors $D_Y$ and $D_X=f^{-1}(D_Y)$, there exists
$q\in D_Y$ such that
there exist uniformizing parameters $u,v,w$ on $Y$ (an etale morphism $Y\rightarrow \text{spec}({\bold k}[u,v,w])$) such that
$u=v=w=0$ are equations of $q$ in $Y$,
$$
D_Y=\{uv=0\},
$$
 and
the fundamental locus of $f$ is $C_1\cup C_2$ where $u=w=0$ are equations of $C_1$ in $Y$, $v=w=0$ are equations of $C_2$ in $Y$.
Suppose that $u=0$, $v=0$ are integral surfaces in $Y$, and  $C_1$, $C_2$ are irreducible.
Let $\pi:Y\rightarrow S=\text{spec}({\bold k}[u,v])$ be the projection. 
Let $\overline q=\pi(q)$ and $\gamma=\pi^{-1}(\overline q)\subset Y$. Let $D_S=\{uv=0\}$, a SNC divisor on $S$.
Assume that
$g=\pi\circ f:X\rightarrow S$ is toroidal  away from $f^{-1}(q)$, and prepared away from $f^{-1}(\gamma)$.

Then there exists a sequence of blow ups 
\begin{equation}\label{eq190}
X_n\stackrel{\Psi_n}{\rightarrow} X_{n-1}\rightarrow\cdots\stackrel{\Psi_1}{\rightarrow} X
\end{equation}
where each $\Psi_i:X_i\rightarrow X_{i-1}$ is the blow up of a possible center (a point or a nonsingular curve
contained in
$D_{X_{i-1}}=\Phi_{i-1}^{-1}(D_{X_{i-2}})$ which makes SNCs with $D_{X_{i-1}}$) which is supported over $f^{-1}(q)$
 such that if $F$ is a component of $D_X$
which dominates a component of $D_Y$ (or dominates $C_1$ or $C_2$), and $F_n$
is the strict transform of $F$ on $X_n$, then $X_n\rightarrow S$ is toroidal
 in a
neighborhood of $F_n$, and $X_n\rightarrow S$  is prepared on $F_n$ away from the strict transform of
$\gamma$. Further, $X_n\rightarrow S$ is prepared away from the preimage of $\gamma$, and is toroidal away from the
preimage of $q$.

Further assume that every irreducible component of $D_X$ contains a 3-point and every 2-curve of $D_X$ contains a 3-point. Then every
irreducible component of $D_{X_n}=(\Psi_1\circ\cdots\circ\Psi_n)^{-1}(D_X)$ contains a 3-point and every 2-curve of $D_{X_n}$ contains a 3-point.
\end{Lemma}

\begin{pf} After possibly blowing up points and curves over $X$ which
are supported over $f^{-1}(q)$,
 we may assume
that $f^{-1}(q)$ is a divisor, and if $F$ is a component of $D_X$ which dominates a component $E$ of $D_Y$, and $L$ is
an exceptional component of $f$ which intersects $F$, then $f(L)\subset E$. Let $E_1$ be the component of $D_Y$
with local equation $u=0$, $E_2$ be the component of $D_Y$ with local equation $v=0$. We may further assume that if
$F$ is a component of $D_X$ which dominates $C_1$ (respectively $C_2$) and $L$ is an exceptional component of $f$
which intersects $F$, then $f(L)\subset E_1$ (respectively $f(L)\subset E_2$).
Finally, if every irreducible component of $D_X$ contains a 3-point and every 2-curve of $D_X$ contains a 3-point, we may assume that this condition is
preserved.
 Let $G=f^{-1}(q)$.

Suppose that $F$ is a component of $D_X$ which dominates a component $E$ of $D_Y$. Without
loss of generality, $E$ has the equation $u=0$. By assumption, $g$ is toroidal  at points of $F-f^{-1}(q)$. Suppose that $p\in F\cap f^{-1}(q)$. Then $p$ must be a 2-point or a 3-point.
Recall that $uv=0$ is an equation of the SNC divisor $D_X$ on $X$.

If $p\in F\cap f^{-1}(q)$ is a 2-point then there exist regular parameters $x,y,z$
at $p$ such that $xy=0$ is a local equation of $D_X$, $x=0$ is a local equation of $F$
and there is an expression
$$
u=xy^g\lambda_1, v=y^c\lambda_2
$$
where $\lambda_1,\lambda_2$ are units, $g>0$ and $c>0$. 

Thus $g$ is toroidal and prepared in a neighborhood of $p$.

If $p\in F\cap f^{-1}(q)$ is a 3-point and $p$ is not on the strict transform of the component $E'$ of $D_Y$ with local equation $v=0$, then there exist regular parameters $x,y,z$
at $p$ such that $xyz=0$ is a local equation of $D_X$, $x=0$ is a local equation of $F$,
and there is an expression
$$
u=xy^bz^c\lambda_1,
v=y^dz^e\lambda_2
$$
where $\gamma_1,\gamma_2$ are units, $b,c>0$ and $d+e>0$. Thus $g$ is toroidal and prepared in a neighborhood of $p$.

If $p\in F\cap f^{-1}(q)$ is a point on  the strict transform on $X$ of the component $E'$ of $D_Y$ with local equation $v=0$,
then $p$ is a 3-point and there is an expression 
\begin{equation}\label{eq288}
\begin{array}{ll}
u&=xz^a\lambda_1\\
v&=yz^b\lambda_2
\end{array}
\end{equation}
at $p$ where $a,b>0$, $\lambda_1,\lambda_2$ are units,
$x=y=0$ are local equations of the strict transform $\gamma_1$ of $\gamma=\pi^{-1}(q)$,  and $g$ is toroidal in a neighborhood
of $p$.

Suppose that $F$ is a component of $D_X$ which dominates $C_1$. By assumption, $g$ is toroidal at points of $F-f^{-1}(q)$,
and on points of the strict transform of a component of $D_Y$.

Suppose that $p\in F\cap f^{-1}(q)$ is not on the strict transform of a component of $D_Y$. by construction, $p$ must
be a 2 or 3-point. If $p$ is a 2-point, then there exist regular parameters $x,y,z$ at $p$ and unit series $\lambda_1,\lambda_2$ such that 
$$
u=x^ay^b\lambda_1, v=y^c\lambda_2,
$$
$x=0$ is a local equation of $F$, $y=0$ is a local equation of a component of $D_X$ which maps to $E_1$, so that
$a,b,c>0$. Thus $g$ is toroidal and prepared at $p$. 

If $p$ is a 3-point, then there exist regular parameters $x,y,z$ at $p$ such that 
$$
u=x^ay^bz^c\lambda_1, v=y^dz^e\lambda_2
$$
where $x=0$ is a local equation of $F$, $y=0$, $z=0$ are local equations of exceptional components of $D_X$ which
map into $E_1$, $d+e>0$ and $\lambda_1,\lambda_2$ are unit series. Thus $a,b,c>0$ and $f$ is toroidal and prepared at $p$.

The same analysis applies if $F$ is a component of $D_X$ which dominates $C_2$.

\end{pf} 

\begin{Lemma}\label{Lemma280}
Suppose that $f:X\rightarrow Y$ is a birational morphism of nonsingular projective 3-folds. Then there exists
 a commutative diagram
$$
\begin{array}{rll}
X_1&\stackrel{f_1}{\rightarrow}&Y_1\\
\downarrow&&\downarrow\\
X&\stackrel{f}{\rightarrow}&Y
\end{array}
$$
where the vertical arrows are products of blow ups of nonsingular subvarieties
such that the (reduced) fundamental locus $\Gamma$ of  $f_1$ is a union of nonsingular curves
and points such that two curves of $\Gamma$ intersect in at most one point, and this
intersection is transversal (the two curves have distinct tangent directions). Further, the intersection of any three curves of $\Gamma$
is empty.
\end{Lemma}

\begin{pf}
Let $S$ be a reduced (but  not necessarily irreducible) surface in $Y$ containing the fundamental
locus of $f$. By the standard theorems of resolution of singularities (\cite{H}, Section 6.8 \cite{C3}), there exists a 
commutative diagram
$$
\begin{array}{rll}
X_1&\stackrel{f_1}{\rightarrow}&Y_1\\
\Phi_1\downarrow&&\downarrow\Psi_1\\
X&\stackrel{f}{\rightarrow}&Y
\end{array}
$$
where $\Phi_1$ and $\Psi_1$ are products of blow ups of points and nonsingular curves, 
such that $\Psi_1^{-1}(S)$ is a  divisor whose irreducible components are nonsingular, which necessarily contains the  fundamental
curve $\Gamma_1$ of $f_1$ (the reduced 1-dimensional scheme consisting of the 1-dimensional components of the fundamental locus of $f_1$). Let $S_1,\ldots, S_n$ be the irreducible components of $\Psi_1^{-1}(S)$.

Let $H$ be a hyperplane section of $Y_1$. Let $q\in \Gamma_1$ be a singular point. Let
$m_q\subset {\cal O}_{Y_1,q}$ be the ideal sheaf of the point $q$, $R={\cal O}_{\Gamma_1,q}$,
$\overline m_q=m_qR$ be the maximal ideal of $R$.

Suppose that $\vec \alpha, \vec \beta\in m_q/m_q^2$ are linearly independent over $\bold k$.
There exist regular parameters $u,v,w$ in ${\cal O}_{Y_1,q}$ such that $\vec\alpha=[u],\vec\beta=[v]\in m_{q}/m_{q}^2$. Let
${\cal I}_{\vec\alpha}\subset{\cal O}_{Y_1}$ be the ideal sheaf defined by
$$
{\cal I}_{\vec\alpha,p}=\left\{\begin{array}{ll}
{\cal O}_{Y_1,p}&\text{if }p\ne q\\
(u)+m_q^2&\text{if }p=q.
\end{array}\right.
$$
Let $Y_{\vec\alpha}$  be the blow up of ${\cal I}_{\vec\alpha}$, with projection $\pi_{\vec\alpha}:Y_{\vec\alpha}\rightarrow Y_1$.
Let ${\cal L}_{\vec\alpha}={\cal I}_{\vec\alpha}{\cal O}_{Y_{\vec \alpha}}$. Then 
${\cal M}_m^{\vec\alpha}=\pi_{\vec\alpha}^*{\cal O}_{Y_1}(mH)\otimes{\cal L}_{\vec\alpha}$ is very ample for $m \gg 0$. ${\cal I}_{\vec\alpha,p}$ is a complete ideal, so
$$
{\cal N}_m^{\vec\alpha}={\cal O}_{Y_1}(mH)\otimes{\cal I}_{\vec\alpha}\cong (\pi_{\vec\alpha})_*({\cal M}_m^{\vec\alpha})
$$
and 
$$
\Gamma(Y_{\vec\alpha},{\cal M}_m^{\vec\alpha})=\Gamma(Y_1,{\cal N}_m^{\vec\alpha}).
$$
Since ${\cal M}_m^{\vec\alpha}$ is generated by global sections, the divisor of a general section of $\Gamma(Y_{\vec\alpha},{\cal M}_m^{\vec\alpha})$ is irreducible and nonsingular away from $\pi_{\vec\alpha}^{-1}(q)$
by Bertini's theorem, (cf. Theorems 7.18 and 7.19 \cite{I}). 
Thus the divisor of a general section of $\Gamma(Y_1,{\cal N}_m^{\vec\alpha})$ is irreducible and is nonsingular away from $q$.
Consider the exact sequence
$$
0\rightarrow m_q^2\rightarrow {\cal I}_{\vec\alpha}\rightarrow {\cal I}_{\vec \alpha}/m_q^2\rightarrow 0.
$$
Tensoring with ${\cal O}_{Y_1}(mH)$ for $m\gg 0$, we see that there is a surjection 
$$
\Gamma(Y_1,{\cal N}_m^{\vec\alpha})\rightarrow {\cal I}_{\vec\alpha}/m_q^2\cong \bold k.
$$
Thus a general section $\sigma\in\Gamma(Y_1,{\cal N}_m^{\vec\alpha})$ is such that
$$
\sigma\equiv \lambda_1 u\mod m_q^2
$$
for some $0\ne\lambda_1\in \bold k$. In particular, the divisor $D_{\sigma}$ of a general section $\sigma$ is irreducible, nonsingular
 and $u=0$ is a local equation of its tangent space at $q$.

In an analogous way, we can define an ideal sheaf ${\cal I}_{\vec\beta}\subset {\cal O}_{Y_1}$ by
$$
{\cal I}_{\vec\beta,p}=\left\{\begin{array}{ll}
{\cal O}_{Y_1,p}&\text{if }p\ne q\\
(v)+m_q^2&\text{if }p=q,
\end{array}
\right.
$$
and show that for $m\gg0$, if ${\cal N}_m^{\vec\beta}={\cal O}_{Y_1}(mH)\otimes{\cal I}_{\vec\beta}$, the divisor $D_{\tau}$ of a
general section $\tau$ of $\Gamma(Y_1,{\cal N}_m^{\vec\beta})$ is irreducible, nonsingular, and $v=0$ is a local equation
of its tangent space at $q$.

Since the base locus of the linear system of $\Gamma(Y_1,{\cal N}_m^{\vec\alpha})$ is the point $q$,
 the divisor $D_{\sigma}$ of a general section $\sigma$ of $\Gamma(Y_1,{\cal N}_m^{\vec\alpha})$ intersects $\Gamma_1$ at $q$
 and at finitely  many other points. Since the base locus of the linear system $\Gamma(Y_1,{\cal N}_m^{\vec\beta})$ is $q$, if $D_{\tau}$ is the divisor of a general section $\tau$ of $\Gamma(Y_1,{\cal N}_{m'}^{\vec \beta})$, (for appropriate $m'\ge m$) then $\Gamma_1$ intersects the scheme $D_{\sigma}\cdot D_{\tau}$ at the point $q$ only, and by Bertini's theorem, $\gamma=D_{\sigma}\cdot D_{\tau}$ is an irreducible curve 
which is nonsingular away from $q$. Since $\sigma\equiv \lambda_1 u \mod m_q^2$ and $\tau\equiv \lambda_2v\mod m_q^2$ for some $0\ne\lambda_1,\lambda_2\in \bold k$, we have that $\gamma$ is nonsingular at $q$, and $u=v=0$ are local equations of the tangent space to $\gamma$ at $q$. By Bertini's theorem applied to the surfaces $S_1,\ldots, S_n$, we see that the divisor $D_{\sigma}$ of a general section  $\sigma$ of $\Gamma(Y_1,{\cal N}_m^{\vec\alpha})$ intersects the $S_i$ transversally away from $q$ (the tangent spaces of $D_{\sigma}$ and $S_i$ intersect in a line), and a further application
of Bertini's theorem shows that the divisor $D_{\tau}$ of a general section $\tau$ of $\Gamma(Y_1,{\cal N}_{m'}^{\vec\beta})$ is such that $\gamma=D_{\sigma}\cdot D_{\tau}$ intersects each $S_i$ transversally away from $q$.

In summary, for $m'\gg m\gg 0$ there exist $\sigma_1\in
\Gamma(Y_1,{\cal O}_{Y_1}(mH))$, and $\sigma_2\in
\Gamma(Y_1,{\cal O}_{Y_1}(m'H))$, with respective divisors $H_1$ and $H_2$
 such that $\gamma=H_1\cdot H_2$
is a nonsingular curve which intersects $\Gamma_1$ in the point $q$ only, $\gamma$ intersects
each $S_i$ transversally at all points other than $q$, and  for some nonzero $\lambda_1,\lambda_2\in \bold k$, 
$\sigma_1$ has image
$\lambda_1\vec\alpha$  and $\sigma_2$ has image $\lambda_2\vec \beta$  in $m_q/m_q^2$ under the natural maps 
$\Gamma(Y_1,{\cal O}_{Y_1}(mH)\otimes m_q)\rightarrow m_q/m_q^2$,
and $\Gamma(Y_1,{\cal O}_{Y_1}(m'H)\otimes m_q)\rightarrow m_q/m_q^2$.
Let $f_i=0$ be a local equation of $S_i$ at $q$ for all $i$ such that $q\in S_i$.
We have a surjection $m_q/m_q^2\rightarrow \overline m_q/\overline m_q^2$ and
$\text{dim}_k\overline m_q/\overline m_q^2>1$ since $q$ is a singular point of $\Gamma_1$.
Choose $\vec\alpha,\vec\beta\in m_q/m_q^2$ so that the $\bold k$-span of $\vec\alpha$ and $\vec\beta$ does not contain
the class of $f_i$ for any $S_i$ containing $q$, and so that the images of $\vec\alpha$
and $\vec\beta$ in $\overline m_q/\overline m_q^2$ are linearly independent over $\bold k$.
Now choose $H_1, H_2$ and $\gamma=H_1\cdot H_2$ as above.

Let $\Psi_2:Y_2\rightarrow Y_1$ be the blow up of $\gamma$. Let $\overline\Gamma_1$ be
the strict transform of $\Gamma_1$ on $Y_2$, $\overline S_i$ for $1\le i\le n$ be the
strict transforms of $S_i$ on $Y_2$. We can assume that  $u=0$ is a local equation of $H_1$ at $q$,
$v=0$ is a local equation of $H_2$ at $q$. By assumption, $u,v,f_i$ is a regular system
of parameters  in ${\cal O}_{Y_1,q}$ for all $i$ such that $q\in S_i$, so $\gamma$ intersects $S_i$
transversally at $q$. Thus for all $i$ such that $q\in S_i$, $\overline S_i$ is the
blow up of (the ideal sheaf of) $q$ and a finite number of other points on $S_i$,
which are disjoint from $\Gamma_1$. In particular, each $\overline S_i$ is nonsingular.

Let $C_1,\ldots, C_m$ be the irreducible components of $\Gamma_1$ containing $q$ and let
$K_i$ be the function field of $C_i$ for $1\le i\le m$. Let $\overline K=
K_1\oplus \cdots \oplus K_m$. Since $\Gamma_1$ is reduced, we have
 natural inclusions
$$
R\rightarrow \overline A= A_1\oplus\cdots\oplus A_m\rightarrow \overline K
$$
where $A_i$ is the integral closure of ${\cal O}_{C_i,q}$ in $K_i$.
$\overline A$ is finite over $R$, since $\overline A$ is the normalization of $R$. Let $s(q)$ 
be the length of $\overline A$ as an $R$-module.

Suppose that $q_1\in\overline\Gamma_1 \cap \Psi_2^{-1}(q)$. Let $R_1={\cal O}_{\overline\Gamma_1,q_1}$. $R_1$ is a local ring of the blow up of $(u,v)R$ in $R$. $(u,v)R$
is not a principal ideal in $R$, since $\vec \alpha$, $\vec\beta$ are linearly independent in $\overline m_q/\overline m_q^2$,  and by construction,
$\overline m_qR_1$ is principal. Thus we have inclusions
$$
R\rightarrow R_1\rightarrow \overline K
$$
with $R\ne R_1$. Since $R_1$ is finite over $R$, we have an inclusion $R_1\subset \overline A$,
and
$\overline A$ has length $s_1<s(q)$ as an $R_1$ module.

The strict transform $\gamma'$ of $\gamma$ on $X_1$ is necessarily a nonsingular curve.
Let $\Phi_2':X_2'\rightarrow X_1$ be the blow up of $\gamma'$. 
${\cal I}_{\gamma}{\cal O}_{X_2'}$
is invertible, except possibly over $q$. There exists a sequence of blow ups $X_2\rightarrow X_2'$ supported over $q$ such that ${\cal I}_{\gamma}{\cal O}_{X_2}$ is invertible. Thus we have a natural morphism
$f_2:X_2\rightarrow Y_2$. The fundamental curve of $f_2$ is contained in the union of the strict transform $\overline\Gamma_1$ of $\Gamma_1$ and the nonsingular curve $l_2=\Psi_2^{-1}(q)$. $l_2$ lies
on all surfaces $\overline S_i$ such that $S_i$ contains $q$.

By induction on
$$
\text{max}\{s(q)\mid q\in \Gamma_1\text{ is a singular point}\}
$$
we can iterate the blowups of such curves $\gamma$ until we construct a commutative diagram
$$
\begin{array}{rll}
X_3&\stackrel{f_3}{\rightarrow}&Y_3\\
\Phi_3\downarrow&&\downarrow\Psi_3\\
X_1&\stackrel{f_1}{\rightarrow}&Y_1
\end{array}
$$
such that
\begin{enumerate}
\item[1.] If $\Gamma'$ is the strict transform of $\Gamma_1$ on $Y_3$ then $\Gamma'$
is a disjoint union of nonsingular curves.
\item[2.] The strict transform $S_i'$ of $S_i$ on $Y_3$
  is nonsingular for all $i$.
\item[3.] The fundamental curve $\Gamma_3$ of $f_3$ is the union of $\Gamma'$ and a union
of curves which are contained in the exceptional loci of  the morphisms
$\lambda_i=(\Psi_3\mid S_i'):S_i'\rightarrow S_i$ for $1\le i\le n$.
\end{enumerate}

For $q\in\Gamma_1$, and $S_i$ such that $q\in S_i$, $\Psi_3^{-1}(q)=\lambda_i^{-1}(q)$ is a SNC divisor on $S_i'$.
Thus $\Gamma_3$ can only fail to satisfy the conclusions of the theorem at a finite
number of points $q'$ such that $q'$ is contained in a (unique) component $C_i'$ of $\Gamma'$
which is contained in  some $S_i'$, and there exists a neighborhood $U$ of $q'$ in $Y_3$
such that $\Gamma_3\cap U$ is a union of components of 
$$
(\lambda_i^{-1}(q) \cup C_i')\cap U,
$$
which is a divisor on the surface $U\cap S_i'$.
Now we can choose a nonsingular curve $\tilde\gamma$ on $Y_3$ which intersects $\Gamma_3$ at
 $q'$ only, and  intersects the surfaces $S_1',\ldots, S_n'$ transversally.
 Let $\Psi_4:Y_4\rightarrow Y_3$ be the
blow up of $\tilde\gamma$, $X_4'\rightarrow X_3$ be the blow up of the strict transform
of $\tilde\gamma$ on $X_3$, and let $X_4\rightarrow X_4'$ be a prinicipalization of
${\cal I}_{\tilde\gamma}{\cal O}_{X_4'}$ (obtained by blowing up points and nonsingular
curves)  which is an isomorphism away from points above $q'$.
Thus we have a commutative diagram
$$
\begin{array}{rll}
X_4&\stackrel{f_4}{\rightarrow}&Y_4\\
\Phi_4\downarrow&&\downarrow\Psi_4\\
X_3&\stackrel{f_3}{\rightarrow}&Y_3
\end{array}
$$
such that the fundamental curve of $f_4$ is contained in the union of the strict transform
of $\Gamma_3$ on $Y_4$ and the curve $l_4=\Psi_4^{-1}(q')$.

Let $S_i''$ be the strict transform of $S_i'$ on $Y_4$ for $1\le i\le n$.
For $i$ such that $q'\in S_i'$, $S_i''\rightarrow S_i'$ is the blow up of $q'$ on $S_i'$,
with exceptional divisor $l_4$. Thus by embedded resolution of plane curve singularities (cf. 
Section 3.4, Exercise 3.13 \cite{C3}),
after a finite number of blow ups of such curves $\tilde \gamma$ we obtain a diagram
$$
\begin{array}{rll}
X_5&\stackrel{f_5}{\rightarrow}&Y_5\\
\downarrow&&\downarrow\\
X_4&\stackrel{f_4}{\rightarrow}&Y_4
\end{array}
$$
such that the fundamental locus of $f_5$ satisfies the conclusions of this lemma.
\end{pf}

\begin{Remark}\label{Remark390} In the conclusions of Lemma \ref{Lemma280}, we can assume that the fundamental locus of $f_1$ has no isolated points. To see this, we make the
following construction. Let $A$ be the isolated points in the fundamental locus of $f_1$.
Let $\gamma$ be a general curve on $Y_1$ through $A$ (an intersection of two general
hypersurface sections through $A$). Then $\gamma\cap\Gamma=A$. Let $\Phi_2:X_2\rightarrow X_1$ be the blow up of the strict transform of $\gamma$ on $X_1$.  The fundamental locus of the resulting map
$X_2\rightarrow Y_1$ is $\gamma\cup\Gamma_1$, which satisfies the conclusions of Lemma
\ref{Lemma280}, and has no isolated points.
\end{Remark}

\begin{Theorem}\label{Theorem190} Suppose that $f:X\rightarrow Y$ is a
birational morphism of nonsingular projective 3-folds and the fundamental locus
$\Gamma$ of $f$ is a union of  nonsingular curves such
that two curves of $\Gamma$ intersect in most one point, and this intersection is
transversal. Further assume that the intersection of any three curves of $\Gamma$ is empty. Then there
exists a commutative diagram
$$
\begin{array}{rll}
X_1&\stackrel{f_1}{\rightarrow}&Y_1\\
\downarrow&&\downarrow\\
X&\stackrel{f}{\rightarrow}&Y
\end{array}
$$
such that the vertical arrows are products of blow ups of points and nonsingular curves, $X_1$, $Y_1$ have toroidal structures $D_{Y_1}$, $D_{X_1}=f_1^{-1}(D_{X_1})$,
$f_1$ is prepared (Definition \ref{Def1}),  every 2-curve of $X_1$ contains a 3-point
and every component of $D_{X_1}$ contains a 3-point.
\end{Theorem}

\begin{pf}
Let $\{C_1,\ldots,C_n\}$ be the irreducible components of $\Gamma$. 

Let $H$ be a hyperplane
section of $Y$. 
Whenever $j\ne i$ and $C_i\cap C_j$ is nonempty, let $C_i\cap C_j=\{q_{ij}\}$.
For $m\gg 0$, and $1\le i\le n$, let $H_i$ be divisors of general sections $\sigma_i$ of $\Gamma(Y,{\cal I}_{C_i}\otimes{\cal O}_{Y_i}(mH))$.
By the arguments using Bertini's theorem of Lemma \ref{Lemma280}, we conclude the following:
\begin{enumerate}
\item[1.] Each $H_i$ is a nonsingular irreducible surface and
$\overline D_Y=H_1+\cdots+H_n$
is a SNC divisor on $Y$.
 \item[2.] For $i\ne j$,
if $C_i\cap C_j\ne\emptyset$, then $H_i$ intersects $C_j$ transversally at $q_{ij}$ plus a sum
of general points of $C_{j}$.
If $C_i\cap C_j=\emptyset$, then $H_i$ intersects $C_j$
transversally at a sum of general points of $C_j$.
\item[3.]  $H_i\cap H_j\cap H_k$ is disjoint from
$\Gamma$ for $i,j,k$ distinct.   
\end{enumerate}

Let $\overline D_X=f^{-1}(\overline D_Y)$. Away from $f^{-1}(\Gamma)$, $\overline D_X$ is a SNC divisor and 
$f$ is  prepared (Definition \ref{Def1}). Suppose that $\eta\in C_j$ is a general point.
Then $f^{-1}(H_j)$ is a SNC divisor over $\eta$ and  $f$ has fiber dimension 1 over $\eta$.
Further, $f$ is a product of blow ups of sections over $C_j$ above $\eta$
which make SNCs with the preimage of $H_j$ (by \cite{Ab1} or \cite{D}, since $f$ is birational and $\eta\in C_j$ is a general point). If $u=0$ is a local equation of $H_j$
at $\eta$ and $v=0$ is a local equation of a nonsingular surface transversal to $C_j$ at $\eta$ then $u,v$ are
toroidal forms (Definition \ref{torf}) at all points of $f^{-1}(\eta)$. In fact we have a form 
\begin{equation}\label{eq283}
u=x^a,
v=y
\end{equation}
 or 
\begin{equation}\label{eq284}
u=x^ay^b,
v=z
\end{equation}
 at all points $p\in f^{-1}(\eta)$.

If $\eta\in H_i$ for some $i\ne j$ (and $\eta\ne q_{ij}$), then we can take $v=0$ to be a local equation of $H_i$ at $\eta$. 

Thus $f^{-1}(\overline D_Y)$ is a SNC divisor except possibly over  $A=\{q_{ij}\}$ and over a finite
number of 1-points $B=\{q_k\}$ of $\overline D_Y$ (contained in $\Gamma$).
After possibly extending $B$ by adding a finite number of points
which are 1-points of $\overline D_Y$, we have that $f$
is prepared away from the points of $A\cup B$. 

Index $B$ as $B=\{q_{n+1},\ldots,q_r\}$. For $q_i\in B$, let $H_i$ be the divisor of a general
section $\sigma_i\in \Gamma(Y,{\cal I}_{q_i}\otimes{\cal O}_Y(mH))$.
For $i\ge n+1$, $H_i$ intersects $\Gamma$ at $q_i$ plus a sum of general points of
$C_1,\ldots, C_n$ (we can make our initial choice of $m$ so that this property holds), 
$$
D_Y=H_1+\cdots+H_r=\overline D_Y+H_{n+1}+\cdots+H_r
$$
is a SNC divisor on $Y$, and $D_X=f^{-1}(D_Y)$ is a SNC divisor on $X$, except possibly over
points of $A\cup B$. Thus after blowing up points and nonsingular curves supported above
$f^{-1}(A\cup B)$, we may assume that $D_X=f^{-1}(D_Y)$ is a SNC divisor, and every irreducible component of $D_X$ contains a 3-point, every 2-curve
of $D_X$ contains a 3-point.

Observe that the points where $f$ is not prepared are intersection points $H_i\cdot H_j\cdot\Gamma$ for $i\ne j$. We may assume that $r\ge 3$. For $i\ne j$ let $\Gamma_{ij}=H_i\cdot H_j$. 
$\Gamma_{ij}$ are nonsingular irreducible curves (by Bertini's theorem).

We now apply for $i=1$ and $j=2$ a general construction that we will iterate for
all $i<j$.    We will assume that $n\ge 2$ and $C_1\cap C_2\ne\emptyset$. The case when $C_1\cap C_2=\emptyset$ (or $n=1$) is simpler. 
 Let $\gamma=\Gamma_{12}$,
$q=q_{12}$, $D_{12}=\Gamma_{12}\cap (C_1\cup C_2)=\Gamma_{12}\cap \Gamma$.

$H_1,H_2,H_3$ are the divisors of sections  $\sigma_1,\sigma_2,\sigma_3\in \Gamma(Y,{\cal O}_Y(mH))$ which define a rational map $\pi:Y\rightarrow {\bold P}^2$ by
$Q\mapsto (\sigma_1:\sigma_2:\sigma_3)(Q)$ for closed points $Q\in Y$. Let $U=Y-( \cup_{j>2}H_j)$, an affine
neighborhood of $q$ in $Y$
on which $\pi$ is a morphism. Let
$$
f_1=\frac{\sigma_1}{\sigma_3}, f_2=\frac{\sigma_2}{\sigma_3}\in\Gamma(U,{\cal O}_Y).
$$
\begin{equation}\label{eq286}
\pi:U\rightarrow {\bold A}^2
\end{equation}
  is defined by the inclusion of $\bold k$-algebras
${\bold k}[u,v]\rightarrow \Gamma(U,{\cal O}_Y)$ given by
$$
u=f_1, v=f_2.
$$
$\overline q=\pi(q)$ has  equations $u=v=0$ in ${\bf A}^2$.
$f_1=f_2=0$ are  equations of $\overline \gamma=\gamma\cap U$ in $U$.  $\overline\gamma$ contains $D_{12}$. $\overline\gamma$ is a nonsingular curve. Hence
$f_1,f_2$ are part of a regular system of parameters at all points of $\overline \gamma$.
Thus $\pi:U\rightarrow {\bold A}^2$ is smooth in a neighborhood of $\overline\gamma$.
Since $C_1$ and $C_2$ intersect $\overline\gamma$ transversally, there exists
an open neighborhood $\overline U$ of $\overline\gamma$ in $U$
such that 
\begin{equation}\label{eq389}
\pi:\overline U\rightarrow{\bold A}^2
\end{equation}
 is smooth
and $\pi\mid C_1\cap \overline U$, $\pi\mid C_2\cap \overline U$ are unramified. ${\bold A}^2$ has toroidal
structure $uv=0$, and $\overline U$ has toroidal structure
$D_Y\cap\overline U$ which is defined by $f_1f_2=0$.
Let $\overline X=f^{-1}(\overline U)$,
$\overline f=f\mid\overline X$, $g=\pi\circ \overline f:\overline X\rightarrow {\bold A}^2$.
$u,v-v(g(\overline p))$ must have a form (\ref{eq283}) or (\ref{eq284}) at points $\overline p$ of $\overline X$
above $C_1-q$, and $v,u-u(g((\overline p))$ have a form (\ref{eq283}) or (\ref{eq284}) at points $\overline p$ of $\overline X$ above  $C_2-q$, so $g$ is toroidal away from $\beta=f^{-1}(q)$
and is prepared away from $f^{-1}(\overline\gamma)$.

Let $R={\cal O}_{Y,q}$ with maximal ideal $m$, $I_{C_1}={\cal I}_{C_1,q}$,
$I_{C_2}={\cal I}_{C_2,q}$. $I_{C_1}$ has generators $u,w_1$, $I_{C_2}$ has generators
$v,w_2$ where $u,v,w_1$ and $u,v,w_2$ are  bases of $m/m^2$ since $H_2$ intersects $C_1$ transversally at $q$ and $H_1$ intersects $C_2$ transversally at $q$. By the formal
implicit function theorem in $\hat R$, $w_2=\phi(w_1-\psi(u,v))$ where $\phi$ is a unit series, and $\psi$ is a series. 
$$
\hat I_{C_2}=(v,w_1-\psi(u,v))=(v,w_1-\psi(u,0)).
$$
Set $\overline w=w_1-\psi(u,0)$. $\hat I_{C_1}=(u,\overline w)$ and $\hat I_{C_2}=(v,\overline w)$. Thus
$$
\widehat{I_{C_1}\cap I_{C_2}}=\hat I_{C_1}\cap \hat I_{C_2}=(uv,\overline w).
$$
There exists $w\in I_{C_1}\cap I_{C_2}$ such that $w\equiv \overline w\text{ mod }m^2\hat R$.
We have $I_{C_1}\cap I_{C_2}=(uv,w)$.
Thus there exists an affine neighborhood $\overline U_1$
of $q$ in $\overline U$ and $w\in \Gamma(\overline U_1,{\cal O}_Y)$ such that $u,v,w$
are uniformizing parameters in $\overline U_1$,  $u=w=0$ are  equations
of $C_1$, and $v=w=0$ are local equations of $C_2$ in $\overline U_1$.

After possibly blowing up points supported above $q$, we may suppose that every irreducible component of $D_{\overline X}$ contains a 3-point and
every 2-curve of $D_{\overline X}$ contains a 3-point.

 By Lemma \ref{Theorem185} (applied to $f^{-1}(\overline U_1)\rightarrow \overline U_1$,
and extending trivially to $\overline X\rightarrow \overline U$),

there exists a commutative diagram
$$
\begin{array}{rll}
\overline X_0&&\\
\overline\Phi_0\downarrow&\searrow&g_0\\
\overline X&\stackrel{g}{\rightarrow}&{\bold A}^2.
\end{array}
$$
such that  $\overline\Phi_0$ is an isomorphism away from $\overline f^{-1}(q)$,
$g_0$ is a toroidal in a neighborhood of the strict transform of $D_Y$ on $\overline X_0$, and is toroidal in a neighborhood of all components of $D_{\overline X_0}$
which dominate $C_1$ or $C_2$, and is prepared in a neighborhood of all components of $D_{\overline X_0}$ which dominate a component of $D_{\overline U}$ or
dominate $C_1$ or $C_2$, away from the strict transform of $\overline\gamma$ on $\overline X_0$.
Further, every irreducible component of $D_{\overline X_0}=\overline\Phi_0^{-1}(D_{\overline X})$ contains a 3-point and every 2-curve of $D_{\overline X_0}$
contains a 3-point.

We now apply the algorithm of the proof of Lemma \ref{Theorem185} to the other points of $D_{12}$ (the conclusions of the
lemma hold if $C_1=\emptyset$ or $C_2=\emptyset$). We construct $\overline\Phi_0':\overline X_0'\rightarrow\overline X_0$
such that
\begin{enumerate}
\item[1.] $\overline\Phi_0\circ\overline\Phi_0'$ is an isomorphism away from $(\overline f\circ\overline\Phi_0\circ\overline\Phi_0')^{-1}(D_{12})$.
\item[2.] $g_0'=g_0\circ\overline\Phi_0'$ is toroidal in a neighborhood of the strict transform of $D_Y$ on $\overline X_0'$,
and in a neighborhood of all components of $D_{\overline X_0'}$ which dominate $C_1$ or $C_2$ and 
is toroidal away from $(\overline f\circ\overline\Phi_0\circ\overline\Phi_0')^{-1}(D_{12})$.
\item[3.] $g_0'$ is prepared in a neighborhood of all components of $D_{\overline X_0'}$ which dominate a component of $D_{\overline U}$ or dominate $C_1$ or $C_2$,
away from the strict transform $\gamma_1$ of $\overline \gamma$ on $\overline X_0'$.
\item[4.] $D_{\overline X_0'}=(\overline \Phi_0')^{-1}(D_{\overline X_0})=(\overline f\circ\overline\Phi_0\circ\overline\Phi_0')^{-1}(D_{Y})$  and $(f\circ\overline\Phi_0\circ\overline\Phi_0')^{-1}(D_{12})$ are SNC divisors.
\item[5.] Every irreducible component of $D_{\overline X_0'}$ contains a 3-point and every 2-curve of $D_{\overline X_0'}$ contains a 3-point.
\end{enumerate}

Let $\overline\Phi_1:\overline X_1\rightarrow\overline X_0'$ be the blow up of the
strict transform $\gamma_1$ of $\overline\gamma$ on $\overline X_0'$.
Let $E_1=\overline\Phi_1^{-1}(\gamma_1)$. Since $\gamma_1$ is a 2-curve of $D_{\overline X_0'}$, we have that every irreducible component of
$D_{\overline X_1}=\overline\Phi_1^{-1}(D_{\overline X_0'})$ contains a 3-point and every 2-curve of $D_{\overline X_0}$ contains a 3-point.

If $p\in\gamma_1-(\overline f\circ\overline\Phi_0\circ\overline\Phi_0')^{-1}(D_{12})$ then $u,v$ are part of a regular
system of parameters at $p$, and $u=v=0$ are local equations of $\gamma_1$
on $\overline X_0'$
at $p$. If $p\in \gamma_1\cap
(\overline f\circ\overline\Phi_0\circ\overline\Phi_0')^{-1}(D_{12})$, then $u,v$ have a toroidal form of the type of equation (\ref{eq288}) at $p$, 
\begin{equation}\label{eq387}
u=xz^a, v=yz^b
\end{equation}
where $a,b>0$, $x=y=0$ are (formal) local equations of $\gamma_1$, $xyz=0$ is a local equation of $D_{\overline X_0'}$.

Suppose that $p'\in E_1\cap (\overline f\circ\overline\Phi_0\circ\overline\Phi_0'\circ\overline\Phi_1)^{-1}(D_{12})$.
Then $p=\overline\Phi_1(p')\in\gamma_1\cap (\overline f\circ\overline\Phi_0\circ\overline\Phi_0')^{-1}(D_{12})$, and $u,v$ have 
a toroidal form (\ref{eq387}) at $p$. Thus $p'$ has (formal) regular parameters $x_1,y_1,z_1$ such that
$$
x=x_1, y=x_1(y_1+\alpha), z=z_1
$$
for some $\alpha\in \bold k$, or
$$
x=x_1y_1, y=y_1, z=z_1.
$$

Substituting into (\ref{eq387}), we see that $g_1=g_0'\circ\overline\Phi_1:\overline X_1\rightarrow {\bf A}^2$ is toroidal and prepared in a neighborhood of the strict transform
of $D_Y$ and in a neighborhood of $E_1$.
Furthermore, $g_1$ is prepared and toroidal in a neighborhood of all components $D_{\overline X_1}$ which dominate $C_1$ or $C_2$ and $g_1$ is prepared and toroidal away from $(\overline f\circ\overline\Phi_0\circ\overline\Phi_0'\circ\overline\Phi_1)^{-1}(D_{12})$.

By 1 of Theorem \ref{Theorem183} there exists a morphism $\overline\Phi_2:\overline X_2
\rightarrow \overline X_1$ such that $\overline \Phi_2$ is an isomorphism away from
$(\overline f\circ\overline\Phi_0\circ\overline\Phi_0'\circ\overline\Phi_1\circ\overline\Phi_2)^{-1}(D_{12})$,
$g_2=g_1\circ\overline\Phi_2:\overline X_2\rightarrow {\bf A}^2$ 
 is prepared,
all 2-curves of $\overline X_2$ contain a 3-point, all components of $D_{\overline X_2}$
contain a 3-point. 

Now by 2 of Theorem \ref{Theorem183}, there exists a commutative diagram
$$
\begin{array}{rll}
\overline X_3&\stackrel{g_3}{\rightarrow}&S_1\\
\overline\Phi_3\downarrow&&\downarrow\overline\Psi_1\\
\overline X_2&\stackrel{g_2}{\rightarrow}&{\bf A}^2
\end{array}
$$
such that $\overline \Phi_3$ is an isomorphism away from
$(g_2\circ\overline\Phi_3)^{-1}(\overline q)$,
$\overline\Psi_1$ is an isomorphism away from $\overline q$ and $g_3$ is toroidal.
We further have that all 2-curves of $\overline X_3$ contain
a 3-point and all components of $D_{\overline X_3}$ contain a 3-point.

In the sequences (\ref{eq*2}) of the proof of Theorem \ref{Theorem183}, we only insert point blow ups $\lambda_{j+1}':\overline Z_{j+1}'\rightarrow
\overline Z_{j+1}$ at a point above the curve $C_j\subset\overline Z_j$
if $C_j$ contains no 2-points. In this case, we are free to blow up the 2-point $p_1$ above any point $p\in C_j$ that we wish. We can thus assume
that $p_1$ lies above 
$\overline\beta_3=(\overline f\circ\overline\Phi_0\circ\overline\Phi_0'\circ\overline\Phi_1\circ
\overline\Phi_2)^{-1}(D_{12})$.

$\overline\Psi_1:S_1\rightarrow {\bf A}^2$ is the blow up of an ideal sheaf ${\cal I}\subset {\cal O}_{{\bf A}^2}$.  We have (by the algorithm of 2 of Theorem \ref{Theorem183})
that $\overline\Phi_3$ is the blow up of ${\cal I}{\cal O}_{\overline X_2}$
away from $\overline\beta_3$.

We can assume that $\overline\Psi_1$ is nontrivial so that $m_{\overline q}{\cal O}_{\overline X_3}$
is an invertible ideal sheaf (where $m_{\overline q}=(u,v)\subset {\bold k}[u,v]$).
 
Let $\overline Y_3=\overline U\times_{{\bf A}^2}S_1$.  $\overline Y_3$ is obtained from $\overline U$
by blowing up sections over $\overline\gamma=\pi^{-1}(\overline q)$, which make SNCs with  the preimage of
$D_{\overline U}$.
Let $\overline\Psi_3:\overline Y_3\rightarrow \overline U$ be the first
projection. 
$\overline Y_3$ is nonsingular and
$D_{\overline Y_3}=\overline\Psi_3^{-1}(D_{\overline U})$ is a SNC divisor.
Let
$$
\tilde\Phi=\overline \Phi_0\circ\cdots \circ\overline\Phi_3:\overline X_3\rightarrow \overline X.
$$
There is a natural birational morphism $\overline f_3:\overline X_3\rightarrow\overline Y_3$ where
$\overline f_3=(\overline f\circ\tilde\Phi)\times g_3$.
By our construction,
$\overline f_3$ is an isomorphism away from $(\overline f\circ\tilde\Phi)^{-1}(C_1\cup C_2)$.
Recall that $\overline\Phi_1:\overline X_1\rightarrow \overline X_0'$ is the blow up of $m_{\overline q}{\cal O}_{\overline X_0'}$ away from $(\overline f\circ\overline\Phi_0\circ\overline\Phi_0'\circ\overline\Phi_1)^{-1}(C_1\cup C_2)$.

Now consider the commutative diagrams
$$
\begin{array}{rll}
X&\stackrel{f}{\rightarrow}&Y\\
\uparrow&&\uparrow\\
\overline X&\stackrel{\overline f}{\rightarrow}&\overline U
\end{array}
$$
where the vertical arrows are the natural inclusions,
and 
\begin{equation}\label{eq272}
\begin{array}{rll}
\overline X_3&\stackrel{\tilde \Phi}{\rightarrow}
&\overline X\\
\overline f_3\downarrow&&\downarrow \overline f\\
\overline Y_3&\stackrel{\overline\Psi_3}{\rightarrow}&\overline U.
\end{array}
\end{equation}

 $\overline\Psi_3$ is constructed by blowing up sections
over $\overline\gamma$ (which make SNCs with the preimage of $D_{\overline U}$). Over $\overline X-\overline f^{-1}(D_{12})$, $\tilde\Phi$ is constructed by
blowing up the isomorphic sections over $\overline \gamma$
and $\overline f_3\mid \overline X_3-\tilde\Phi^{-1}(\overline f^{-1}(C_1\cup C_2))$ is an isomorphism.

Recall that $f$ is an
isomorphism over  $H_j\cap \gamma$ for $j>2$, and $(\cup_{j>2}H_j)\cap\gamma=\gamma-\overline \gamma$.

 There exists a factorization of the above
diagram (\ref{eq272}) over $\overline U-D_{12}$ by a commutative diagram 
\begin{equation}\label{eq271}
\begin{array}{lllllll}
\hat X_m&\stackrel{\Omega_m}{\rightarrow}&\cdots&\stackrel{\Omega_2}{\rightarrow}
&\hat X_1&\stackrel{\Omega_1}{\rightarrow}&\hat X\\
\downarrow&&&&\downarrow&&\downarrow\\
\hat Y_m&\stackrel{\Lambda_m}{\rightarrow}&\cdots&\stackrel{\Lambda_2}{\rightarrow}
&\hat Y_1&\stackrel{\Lambda_1}{\rightarrow}&\hat Y\\
\end{array}
\end{equation}

where $\hat Y=\overline U-D_{12}$,
 $\hat X=\overline X-\overline f^{-1}(D_{12})$, $\hat Y_m=\overline Y_3-\overline\Psi_3^{-1}(D_{12})$,
$\hat X_m=\overline X_3-\tilde\Phi^{-1}
(\overline f^{-1}(D_{12}))$.
The vertical arrows are isomorphisms away from the preimage of $C_1\cup C_2$, each map $\hat Y_{i+1}\rightarrow\hat Y_i$ is the
blow up of a section $\tilde\gamma_i$ over $\tilde\gamma=\gamma\cap \hat Y$ which
makes SNCs with the toroidal structure ($\tilde\gamma_i$ is a possible blow up).

We can thus extend (\ref{eq271}) to a commutative diagram of projective morphisms 
\begin{equation}\label{eq273}
\begin{array}{rll}
\tilde X_3&\rightarrow& X-f^{-1}(D_{12})\\
\downarrow&&\downarrow\\
\tilde Y_3&\rightarrow&Y-D_{12}
\end{array}
\end{equation}
such that the vertical arrows are isomorphisms away from $\Gamma$, and the horizontal arrows are isomorphisms
over $Y-\gamma$,
by performing the blow up 
$\tilde Y_1\rightarrow Y-D_{12}$ of $\gamma\cap (Y-D_{12})$, then blowing up points
over  $H_j\cap \gamma$ for $j>2$ to ensure that the closure of $\tilde\gamma_1$  in $\tilde Y_1$ makes SNCs with the toroidal structure, and
continue to construct $\tilde Y_3\rightarrow Y-D_{12}$ which extends $\hat Y_m\rightarrow \hat Y$. Now we can perform
the corresponding point blow ups over $X$ ($f$ is an isomorphism over 
 $ H_j\cap\gamma$ for $j>2$) and blow ups of the closures of $\tilde\gamma_i$
to get  $\tilde X_3\rightarrow X- f^{-1}(D_{12})$ which extends $\hat X_m\rightarrow \hat X$.
Patching (\ref{eq272}) and (\ref{eq273}) we  have a commutative diagram 
\begin{equation}\label{eq285}
\begin{array}{rll}
X_3&\stackrel{f_3}{\rightarrow}&Y_3\\
\Phi_3\downarrow&&\downarrow \Psi_3\\
X&\stackrel{f}{\rightarrow}&Y
\end{array}
\end{equation}
of projective morphisms 
where $\Phi_3$ and $\Psi_3$ are products of blow ups of possible centers (points and
nonsingular curves which make SNCs with the preimages of $D_X$ and $D_Y$)
such that $f_3$ is prepared in a neighborhood of $\Psi_3^{-1}(\gamma)$, $\Psi_3$
is an isomorphism away from $\gamma$, $\Phi_3$ is an isomorphism away from $f^{-1}(\gamma)$,
and $f_3$ is an isomorphism away from $\Psi_3^{-1}(\Gamma)$.
We further have that every 2-curve of $X_3$ contains a 3-point and every
component of $D_{X_3}$ contains a 3-point.

We now iterate this construction for the other pairs $H_i,H_j$ with $1\le i<j\le r$, first continuing  the procedure for $H_1$
and $H_3$. We will indicate this next step.

Let $H_i'=\Psi_3^*(H_i)$ for $1\le i\le r$.
$H_1', H_3',H_2'$ are divisors of sections $\sigma_1',\sigma_3',\sigma_2'
\in H^0(Y_3,{\cal O}_{Y_3}(m\Psi_3^*(H))=H^0(Y,{\cal O}_Y(mH))$.
Let $U'=Y_3-(\cup_{j\not\in\{1,3\}}H_j')$. Since $\Psi_3$ is an
isomorphism away from $\gamma=\Gamma_{12}\subset H_2$,
$U'\cong Y-(\cup_{j\not\in \{1,3\}}H_j)$. Thus $U'$ is affine and we have
a morphism $\pi':U'\rightarrow{\bf A}^2$ defined by the inclusion of $\bold k$-algebras
${\bold k}[u,v]\rightarrow\Gamma(U',{\cal O}_Y)$ given by
$$
u=f_1=\frac{\sigma_1'}{\sigma_2'}, v=f_2=\frac{\sigma_3'}{\sigma_2'}.
$$
Let $\gamma'$ be the strict transform of $\Gamma_{13}$ in $Y_3$ and let
$\overline\gamma'=\gamma'\cap U'$. Since $\Gamma_{13}\cap\Gamma\cap H_i=\emptyset$ for $i\not\in\{1,3\}$, 
 and the fundamental locus of $f_3$ is contained in $\Psi_3^{-1}(\Gamma)$,
$f_3$ is an isomorphism over $\gamma'-\overline\gamma'$.

As in the construction of (\ref{eq272}), there exists an open neighborhood $\overline U'$
of $\overline\gamma'$ in $U'$ such that $\pi':\overline U'\rightarrow {\bf A}^2$
is smooth, and enjoys the properties of the morphism $\pi:\overline U\rightarrow {\bf A}^2$
constructed in (\ref{eq389}). Thus there exists a commutative diagram
$$
\begin{array}{rll}
\overline X_3'&\stackrel{\Phi'}{\rightarrow}&\overline X'\\
f_3'\downarrow&&\downarrow \overline f' \\
\overline Y_3'&\stackrel{\Psi_3'}{\rightarrow}&\overline U'
\end{array}
$$
having the properties of (\ref{eq272}) (where $\overline X'=f_3^{-1}(\overline U')$, $\overline f'=f_3\mid \overline X$). We can thus (as in the construction of (\ref{eq285}))
construct a commutative diagram
$$
\begin{array}{rll}
X_4&\stackrel{f_4}{\rightarrow}&Y_4\\
\Phi_4\downarrow&&\downarrow\Psi_4\\
X_3&\stackrel{f_3}{\rightarrow}&Y_3\\
\Phi_3\downarrow&&\downarrow\Psi_3\\
X&\stackrel{f}{\rightarrow}&Y
\end{array}
$$
where $\Phi_4$ and $\Psi_4$ are products of blow ups of points and
nonsingular curves which are possible centers (make SNCs with the preimages of $D_X$ and $D_Y$),
such that $f_4$ is prepared in a neighborhood of $(\Psi_3\circ\Psi_4)^{-1}(\Gamma_{12}\cup\Gamma_{13})$, $\Psi_3\circ\Psi_4$ is an isomorphism away from $\Gamma_{12}\cup\Gamma_{13}$,
$\Phi_3\circ\Phi_4$ is an isomorphism away from $f^{-1}(\Gamma_{12}\cup\Gamma_{13})$
and $f_4$ is an isomorphism away from $(\Psi_3\circ\Psi_4)^{-1}(\Gamma)$. We further have that every 2-curve of $X_4$ contains
a 3-point and every component of $D_{X_4}$ contains a 3-point.

We can continue in this way to construct, by induction, a commutative diagram
$$
\begin{array}{rll}
X_5&\stackrel{f_5}{\rightarrow}&Y_5\\
\Phi_5\downarrow&&\downarrow\Psi_5\\
X&\stackrel{f}{\rightarrow}&Y
\end{array}
$$
such that $f_5$ is prepared in a neighborhood of $\Psi_5^{-1}(\cup_{i\ne j}\Gamma_{ij})$,
$\Psi_5$ is an isomorphism away from $\cup_{i\ne j}\Gamma_{ij}$, $\Phi_5$ is an isomorphism
away from $f^{-1}(\cup_{i,\ne j}\Gamma_{ij})$ and $f_5$ is an isomorphism away from $\Psi_5^{-1}(\Gamma)$. Thus  $f_5$ is prepared.

\end{pf}

\section{First Properties of Prepared Morphisms} 

In this section, we suppose that $f:X\rightarrow Y$ is a birational morphism of nonsingular projective 3-folds, with
toroidal structures $D_Y$ and $D_X=f^{-1}(D_Y)$ such that $D_Y$ contains the fundamental locus of $f$.

\begin{Lemma}\label{Lemma143} Suppose that $f:X\rightarrow Y$ is prepared,
$C\subset Y$ is a curve such that $C$ is contained in the fundamental locus of $f$ and $C$
contains a 1-point. Let $\Psi_1:Y_1\rightarrow Y$ be the blow up of $C$, and let $q\in C$
be a 1-point. Then there exists an affine neighborhood $\overline Y$ of $q$ such that
$C\cap \overline Y$ is nonsingular, and if $\overline Y_1=\Psi_1^{-1}(\overline Y)$ and $\overline X=f^{-1}(\overline Y)$, then there exists
a factorization $f_1:\overline X\rightarrow \overline Y_1$ such that $\Psi_1\circ f_1=f\mid \overline X$. Further,
$f_1$ is prepared.
\end{Lemma}

\begin{pf} Since $f$ is birational, $f^{-1}(q)$ has dimension 1 and $C$ is the only
component of the fundamental locus of $f$ through $q$
(by Lemma \ref{Lemma141}), there exists a neighborhood $\overline
Y$ of $q$ and a factorization  
$$
f_1=\Psi_1^{-1}\circ f: \overline X=f^{-1}(\overline Y)\rightarrow \overline Y_1=\Psi_1^{-1}(\overline Y)
$$
of the blow up $\Psi_1$ of $C$
by a morphism $f_1$ (by \cite{D}). We will show that $f_1$ is prepared.  As shown in the proof of Lemma \ref{Lemma141},
there exist permissible parameters $(u,v,w)$ at $q$ such that 
for all $p\in f^{-1}(q)$, $u,v$ are toroidal forms at $p$ and $u=w=0$ are local
equations of $C$ at $q$. Thus if $p\in f^{-1}(q)$ we have permissible parameters $x,y,z$ for $u,v,w$
at $p$ such that if $p$ is a 1-point, 
\begin{equation}\label{eq147}
\begin{array}{ll}
u&=x^a\\
v&=y\\
w&=\lambda(x,y)+x^cz,
\end{array}
\end{equation}
and if $p$ is a 2-point, 
\begin{equation}\label{eq148}
\begin{array}{ll}
u&=(x^ay^b)^k\\
v&=z\\
w&=\lambda(x^ay^b,z)+x^cy^d
\end{array}
\end{equation}
with $ad-bc\ne 0$. Further, $(u,w)\hat{\cal O}_{X,p}$ is invertible. In case (\ref{eq147}), if $u\mid w$, we have that $f_1(p)$ is a 1-point, and there exists $\alpha\in \bold k$ such that if 
$$
u_1=u, v_1=v, w_1=\frac{w}{u}-\alpha,
$$
then $u_1,v_1,w_1$ are regular parameters at $f_1(p)$ and there exists a
series $\overline\lambda$ such that 
\begin{equation}\label{eq300}
\begin{array}{ll}
u_1&=x^a\\
v_1&=y\\
w_1&=\overline\lambda(x,y)+x^{c-a}z
\end{array}
\end{equation}
of the form (\ref{eqTF01}) of Definition \ref{torf} at $p$. If $w\mid u$ (and $u\not\,\mid w$) in Case (\ref{eq147}), we have that
$w=x^n\gamma$, where $\gamma$ is a unit series and $n<a$. Thus after computing a Taylor series expansion of an appropriate root of
$\gamma$, we see that (by a similar argument to Case 1.1 of Lemma 57 of \cite{C2}) that there exist
regular parameters $\overline x,y,\overline z$ at $p$ and a unit series $\overline\lambda$ such that
$$
\begin{array}{ll}
u&=\overline x^a(\overline\lambda(\overline x,y)+\overline x^{ c-n}\overline z)\\
v&=y\\
w&=\overline x^n.
\end{array}
$$
Thus $f_1(p)=q_1$ is a 2-point which has regular parameters 
$$
w_1=w, u_1=\frac{u}{w},v_1=v
$$
so that $u_1w_1=0$ is a local equation of $D_{Y_1}$ and 
\begin{equation}\label{eq413}
\begin{array}{ll}
w_1&=\overline x^n\\
u_1&=\overline x^{a-n}(\overline\lambda(\overline x,y)+\overline x^{c-n}z)\\
v_1&=y
\end{array}
\end{equation}
which has the form 2 (b) of Definition \ref{Def1}.

In case (\ref{eq148}), if $u\mid w$, we have that $f_1(p)$ is a 1-point, and there exists $\alpha\in \bold k$ such that if
$$
u_1=u, v_1=v, w_1=\frac{w}{u}-\alpha,
$$
then $u_1,v_1,w_1$ are regular parameters at $f_1(p)$ and there exists a
series $\overline\lambda$ such that 
\begin{equation}\label{eq301}
\begin{array}{ll}
u_1&=(x^ay^b)^k\\
v_1&=z\\
w_1&=\overline\lambda(x^ay^b,z)+x^{c-ak}y^{d-bk}
\end{array}
\end{equation}
of the form (\ref{eqTF02}) of Definition \ref{torf} at $p$. If $w\mid u$ and $u\not\,\mid w$ 
in case (\ref{eq148}), we either have an expression
\begin{equation}\label{eq194}
w=(x^ay^b)^l(\gamma(x^ay^b,z)+x^{\tilde c}y^{\tilde d})
\end{equation}
at $p$ with $l<k$ and $\gamma$ is a unit series, or we have an expression 
(after possibly making a change of variables in $x$ and $y$, multiplying $x$ and $y$ by appropriate unit series)
\begin{equation}\label{eq195}
w=x^cy^d
\end{equation}
with $ad-bc\ne 0$. Suppose that (\ref{eq194}) holds.  Then  there exist regular parameters $\overline x,\overline y,z$ at $p$
such that
$$
\begin{array}{ll}
u&=(\overline x^a\overline y^b)^k(\overline \gamma(\overline x^a\overline y^b,z)+\overline x^{\overline c}\overline y^{\overline d})\\
v&=z\\
w&=(\overline x^{a}\overline y^{b})^l
\end{array}
$$
for some unit series $\overline\gamma$, and $\overline a\overline d-\overline b\overline c\ne 0$.
Thus $f_1(p)=q_1$ is a 2-point which has regular parameters 
$$
w_1=w, u_1=\frac{u}{w},v_1=v
$$
such that $u_1w_1=0$ is a local equation of $D_{Y_1}$ and 
\begin{equation}\label{eq302}
\begin{array}{ll}
w_1&=(\overline x^a\overline y^b)^l\\
u_1&=(\overline x^a\overline y^b)^{k-l}(\overline \gamma(\overline x^a\overline y^b,z)
+\overline x^{\overline c}\overline y^d)\\
v_1&=z
\end{array}
\end{equation}
which has the form 2 (c) of Definition \ref{Def1}.

If (\ref{eq195}) holds, $f_1(p)=q_1$ is a 2-point which has regular parameters
$$
w_1=w, u_1=\frac{u}{w}, v_1=v
$$
 such that $u_1w_1=0$ is a local equation of $D_{Y_1}$ and 
\begin{equation}\label{eq303}
\begin{array}{ll}
w_1&=x^cy^d\\
u_1&=x^{ak-c}y^{bk-d}\\
v_1&=z
\end{array}
\end{equation}
which has the form (\ref{eqTF21}) of Definition \ref{torf}.

Comparing equations (\ref{eq300}), (\ref{eq413}), (\ref{eq301}), (\ref{eq302}) and (\ref{eq303}), we see that  after possibly replacing $\overline Y$ with a smaller affine neighborhood of $q$, $f_1:\overline X\rightarrow \overline Y_1$ is prepared.
\end{pf}

\begin{Lemma}\label{Lemma1} Suppose that $f:X\rightarrow Y$ is prepared and
$C\subset Y$ is a 2-curve. Then there exists a commutative diagram
$$
\begin{array}{rll}
X_1&\stackrel{f_1}{\rightarrow}&Y_1\\
\Phi_1\downarrow&&\downarrow\Psi_1\\
X&\stackrel{f}{\rightarrow}&Y
\end{array}
$$
where $\Psi_1:Y_1\rightarrow Y$ is the blow up of $C$, $\Phi_1:X_1\rightarrow X$ is a product of blow ups of 2-curves, 
$\Phi_1$ is an isomorphism above $f^{-1}(Y-C)$ and $f_1$ is prepared. If $p_1\in X_1$ is a 3-point, so that $p=\Phi_1(p_1)$ is
necessarily also a 3-point, we have that
$$
\tau_{f_1}(p_1)= \tau_f(p).
$$
Further, if $D_X$ is cuspidal for $f$, then $D_{X_1}$ is cuspidal for $f_1$.
\end{Lemma}

\begin{pf} Let $\Psi_1:Y_1\rightarrow Y$ be the blow up of $C$. By Lemma \ref{Lemma423}, there exists a commutative diagram
$$
\begin{array}{rll}
X_1&\stackrel{f_1}{\rightarrow}&Y_1\\
\Phi_1\downarrow&&\downarrow\Psi_1\\
X&\stackrel{f}{\rightarrow}&Y
\end{array}
$$
such that $\Phi_1:X_1\rightarrow X$ is a sequence of blow ups of 2-curves, $\Phi_1$ is an isomorphism above $f^{-1}(Y-C)$
and $f_1$ is a morphism.

We will show that $f_1:X_1\rightarrow Y_1$ is prepared.

Suppose that $q\in C\subset Y$ and $p\in (f\circ\Phi_1)^{-1}(q)$. Let $q_1=f_1(p)$,
$p'=\Phi_1(p)$.

\noindent{\bf Suppose that $q$ is a 2-point} 

Let $u,v,w$ be permissible parameters at $q$. Suppose that $p'$ satisfies 2 (a) of Definition \ref{Def1}. If $q_1$ is a 2-point,
which has  permissible parameters $u_1,v_1,w_1$ defined by
$$
u=u_1, v=u_1v_1, w=w_1
$$
or
$$
u=u_1v_1, v=v_1, w=w_1,
$$
then $u_1,v_1$ are toroidal forms at $p$. If $q_1$ is a 1-point, which has  
permissible parameters $u_1,v_1,w_1$ defined by
$$
u=u_1, v=u_1(v_1+\alpha), w=w_1,
$$
with $0\ne\alpha$, then $u_1,v_1$ are toroidal forms at $p$.

Further, if $p$ is a 3-point, then $p'$ is a 3-point and $q_1$ is a 2-point.   Thus $\tau_f(p')\ge 0$. After possibly interchanging $u$ and $v$, we have that $u,v,w$ have a form (\ref{eq16}) at
$p'$ in terms of permissible parameters $x,y,z$ for $u,v,w$ at $p'$,  there are $u_1,v_1\in{\cal O}_{Y_1,q_1}$ such that $u_1,v_1,w$ are permissible
parameters at $q_1$, with $u=u_1, v=u_1v_1$ and there are permissible parameters
$x_1,y_1,z_1$ in $\hat{\cal O}_{X_1,p}$ for $u_1,v_1,w$ such that
$$
\begin{array}{ll}
x&=x_1^{a_{11}}y_1^{a_{12}}z_1^{a_{13}}\\
y&=x_1^{a_{21}}y_1^{a_{22}}z_1^{a_{23}}\\
z&=x_1^{a_{31}}y_1^{a_{22}}z_1^{a_{23}}
\end{array}
$$
where $a_{ij}$ are natural numbers with $\text{Det}(a_{ij})=\pm1$.
Thus ${\bf Z}u+{\bf Z}v={\bf Z}u_1+{\bf Z}v_1$, $H_{f_1,p}=H_{f\circ\Phi_1,p}=H_{f,p'}$,
$A_{f_1,p}=A_{f\circ\Phi_1,p}=A_{f,p'}$, $L_{f_1,p}=L_{f\circ\phi_1,p}=L_{f,p'}$ and $\tau_{f_1}(p)=\tau_{f\circ\Phi_1}(p)=\tau_f(p')$.

Now suppose that $u,v,w$ satisfy 2 (b) of Definition \ref{Def1} at $p'$. Then $\Phi_1$ is an isomorphism near $p=p'$ since $p'$ is a 1-point. If permissible parameters at $q_1$ are
$u_1,v_1,w$ with  $u=u_1, v=u_1v_1, w=w_1$, then $u_1,v_1,w_1$ have the form 2 (b) at $p$ also. If
$u=u_1v_1, v=v_1$ then we can interchange $u$ and $v$ 
in 2 (b) and make an appropriate change of permissible parameters at $p$ to get a form  2 (b) for $u_1,v_1,w_1$ at $p$.
If $u=u_1, v=u_1(v_1+\alpha)$, with $\alpha\ne 0$, then $q_1$ is a 1-point,
and $u_1,w_1$ are toroidal forms  at $p$. 

Suppose that $u,v,w$ satisfy 2 (c) of Definition \ref{Def1} at $p'$. Then there are regular parameters $x_1,y_1,z_1$
in $\hat{\cal O}_{X_1,p}$ defined by 
\begin{equation}\label{eq14}
x=x_1^{\overline a}y_1^{\overline b},
y=x_1^{\overline c}y_1^{\overline d},
z=z_1
\end{equation}
with $\overline a\overline d-\overline b\overline c=\pm1$ or 
\begin{equation}\label{eq15}
x=x_1^{\overline a}(y_1+\alpha)^{\overline b},
y=x_1^{\overline c}(y_1+\alpha)^{\overline d},
z=z_1
\end{equation}
with $\overline a\overline d-\overline b\overline c=\pm 1$ and $0\ne\alpha\in \bold k$.

If (\ref{eq14}) holds, then $u,v,w$ satisfy 2 (c) of Definition \ref{Def1} in $\hat{\cal O}_{X_1,p}$. If $q_1$ is a 2-point, which has regular parameters
$u_1,v_1,w_1$ defined by 
\begin{equation}\label{eq362}
u=u_1v_1, v=v_1, w=w_1
\end{equation}
or 
\begin{equation}\label{eq363}
u=u_1,v=u_1v_1, w=w_1
\end{equation}
then at $p$ there is an expression of $u_1,v_1,w_1$ of the form 2 (c) of Definition \ref{Def1}. 
 If $q_1$ is a 1-point, which has regular parameters $u_1,v_1,w_1$ defined by 
\begin{equation}\label{eq364}
u=u_1, v=u_1(v_1+\alpha), w=w_1,
\end{equation}
with $0\ne\alpha\in \bold k$, then $u_1,w_1$ are toroidal forms at $p$.

If (\ref{eq15}) holds, then $p$ is a 1-point and there exist $\overline x_1,\overline y_1\in\hat{\cal O}_{X_1,p}$ such
that $\overline x_1,\overline y_1,z$ are regular parameters in $\hat{\cal O}_{X_1,p}$ and 
\begin{equation}\label{eq304}
x^ay^b=\overline x_1^s,
x^cy^d=\overline x_1^t(\overline y_1+\beta)
\end{equation}
for some $0\ne \beta\in \bold k$, and positive integers $s,t$.  $u,v,w$ thus satisfy  2 (b) of Definition \ref{Def1} at $p$.
 If $q_1$ is a 2-point, which has regular parameters $u_1,v_1,w_1$ defined by (\ref{eq362}) or (\ref{eq363}),
then at $p$ there is an expression of $u_1,v_1,w_1$ of the form 2 (b) of
Definition \ref{Def1}. If $q_1$ is a 1-point, which has permissible parameters
$u_1,v_1,w_1$ defined by (\ref{eq364}),
then $q_1$ is a 1-point and $u_1,w_1$ are toroidal forms at $p$.

\noindent{\bf Suppose that $q$ is a 1-point.} Let $u,v,w$ be permissible parameters at $q$ (satisfying 3 of Definition \ref{Def1}). In this case, 
$\Psi_1$ is an isomorphism near $q_1$, so we can identify $q$ with $q_1$. $u,v,w$ thus have a form (\ref{eqTF01}) or (\ref{eqTF02}) of Definition \ref{torf}
at $p'$. If (\ref{eqTF01}) holds at $p'$, then $p=p'$, so $u,v,w$ have a form (\ref{eqTF01}) at $p$. Suppose that
(\ref{eqTF02}) holds at $p'$. Then  $p$ has regular parameters defined by (\ref{eq14}) or (\ref{eq15}), and we see that $u,v,w$ have a form (\ref{eqTF01}) or (\ref{eqTF02}) at $p$.

\noindent{\bf Suppose that $q$ is a 3-point.} 
Let $u,v,w$ be permissible parameters at $q$. Then 
after possibly permuting $u,v$ and $w$, we have that  $u,v,w$ have one of the forms (\ref{eqTF1}) - (\ref{eqTF3}) of Definition \ref{torf} at $p'$, and thus $u,v,w$ also have one of the forms (\ref{eqTF1}) - (\ref{eqTF3}) of Definition \ref{torf} at $p$.

First assume $q_1$ is a 3-point.  Further assume that
$p$ is a 3-point. Then  $u,v,w$ have the form of (\ref{eq16}) at $p$, and have a form (\ref{eqTF3}) of Definition \ref{torf} at $p$.
Furthermore, if $\tau_f(p')=-\infty$, then $\tau_{f\circ\Phi_1}(p)=-\infty$, and if $\tau_f(p)\ge 0$, then $\tau_{f\circ\Phi_1}(p)\ge 0$ and $H_{f,p'}=H_{f\circ\Phi_1,p}$, $A_{f,p'}=A_{f\circ\Phi_1,p}$ and $\tau_{f}(p')=\tau_{f\circ\Phi_1}(p)$.
After possibly interchanging $u$ and $v$, $q_1$ has permissible parameters $u_1,v_1,w_1$
such that one of the following equations (\ref{eq17}), (\ref{eq365}) or (\ref{eq18}) hold: 
\begin{equation}\label{eq17}
u=u_1,v=u_1v_1, w=w_1
\end{equation}
or 
\begin{equation}\label{eq365}
u=u_1, v=v_1, w=u_1w_1
\end{equation}
 or 
\begin{equation}\label{eq18}
u=u_1w_1, v=v_1,w=w_1.
\end{equation}
Assume that (\ref{eq17}) or (\ref{eq365}) holds. Then $u_1,v_1,w_1$ have a form (\ref{eqTF3}) of Definition \ref{torf} at $p$.
If $\tau_f(p')=-\infty$, then $\tau_{f_1}(p)=-\infty$, and if $\tau_f(p')\ge 0$, then
$$
H_{f\circ\Phi_1,p}=H_{f_1,p},
$$
 $A_{f\circ\Phi_1,p}=A_{f_1,p}$ and
$\tau_{f_1}(p)=\tau_{f\circ\Phi_1}(p)=\tau_f(p')$.

Assume that (\ref{eq18}) holds.  If $\tau_{f}(p')=-\infty$, then we can interchange $u,v$ and $w$ to obtain the case (\ref{eq17}),
which we have already analyzed, 
so we may assume that
$\tau_f(p')\ge 0$.
 In (\ref{eq16}), we have at $p$ an expression
$w=M_0(\gamma+\overline N_0)$ where $M_0$ is a monomial and 
$$
\gamma=\alpha_0+\sum_{i\ge 1}\alpha_i\overline M_i
$$
is such that
$$
\text{rank}(u,v,\overline M_i)=\text{rank}(u,v,M_0)=2
$$
for $1\le i$ and
$$
\text{rank}(u,v,\overline N_0)=3.
$$
Thus
$$
H_{f\circ\Phi_1,p}={\bf Z}u+{\bf Z}v+{\bf Z}M_0+\sum_{i\ge 1}{\bf Z}\overline M_i=H_{f,p'}
$$
and
$$
A_{f\circ\Phi_1,p}={\bf Z}u+{\bf Z}v+{\bf Z}M_0=A_{f,p'}.
$$

If $\text{rank}(v,M_0)=2$, then $w,v$ is a toroidal form  at $p$, so we have,
after a change of variables in $x,y,z$ at $p$, an expression of $w,v,u$ of the form of
(\ref{eq16}), and (by Lemma \ref{Lemma0}) we obtain the same calculation of $H_{f\circ\Phi_1,p}$
and $A_{f\circ\Phi_1,p}$ for these new parameters. Thus (\ref{eq18}) has been transformed into
(\ref{eq365}), from which it follows that $w_1,v_1,u_1$ have a form (\ref{eqTF3}) of Definition \ref{torf}  at $p$ and $\tau_{f_1}(p)=\tau_f(p')$.

If $\text{rank}(v,M_0)=1$, then $\text{rank}(w,u)=2$ and $w,u$ is a toroidal form at $p$. As in the above paragraph, we change variables to obtain a form (\ref{eq17}), from which it follows that $w_1,u_1,v_1$ have a form (\ref{eqTF3}) of Definition \ref{torf} at $p$ and
$\tau_{f_1}(p)=\tau_f(p')$.

Suppose that $q_1$ is a 3-point and   $p$ is a 2-point.
 Then, after possibly permuting $u,v,w$, we have that $u,v,w$ are
such that $u,v$ are toroidal forms of type (\ref{eqTF21})
or type (\ref{eqTF22}) of Definition \ref{torf} at $p$, and $w=M\gamma$ where $M$
is a monomial in $x,y$, and $\gamma$ is a unit series.

First assume  that $u,v$ are of type (\ref{eqTF21}) of Definition \ref{torf} at $p$. After possibly interchanging
$u$ and $v$,
$q_1$ has permissible parameters $u_1,v_1,w_1$ of one of the forms (\ref{eq17}) (\ref{eq365}) or (\ref{eq18}).
In any of these cases, we  have an expression
$$
u_1=M_1\gamma_1, v_1=M_2\gamma_2, w_1=M_3\gamma_3
$$
where $\gamma_1,\gamma_2,\gamma_3$ are unit series and $M_1,M_2,M_3$ are monomials in $x,y$ with
$$
\text{rank}(M_1,M_2,M_3)=2.
$$
 Thus (after possibly interchanging $u_1,v_1,w_1$) $u_1,v_1,w_1$ have a form (\ref{eqTF21}) of Definition \ref{torf}  at $p$.

Now assume that   $u,v$ is of type (\ref{eqTF22}) 
of Definition \ref{torf} at $p$ (and there does not exist a permutation of $u,v,w$ such that $u,v$ are of type (\ref{eqTF21}) at $p$.
We continue to assume that $q_1$ is a 3-point.
Then there are  permissible
parameters $x,y,z$ at $p$ such that
$$
\begin{array}{ll}
u&=(x^ay^b)^k\\
v&=(x^ay^b)^t(\alpha+z)\\
w&=(x^ay^b)^l[\gamma(x^ay^b,z)+x^cy^d]
\end{array}
$$
where $\gamma$ is a unit series and $ad-bc\ne 0$.

If $u=v=0$ are local equations for $C$ at $q$ then 
after possibly permuting $u$ and $v$, we may assume that $q_1$ has permissible parameters $u_1,v_1,w_1$ defined by
$$
u=u_1, v=u_1v_1, w=w_1.
$$
Thus  $u_1,v_1,w_1$ have a form (\ref{eqTF22}) of Definition \ref{torf}
 at $p$.
Otherwise, we can assume, after possibly interchanging $u$ and $v$ that $u=w=0$
are local equations for $C$. If permissible parameters are defined at $q_1$ by
$$
u=u_1, v=v_1, w=u_1w_1,
$$
 then  $u_1,v_1,w_1$ have a form (\ref{eqTF22}) of Definition \ref{torf}
at $p$.

The remaining possibility is that
$u_1,v_1,w_1$ are permissible parameters at $q_1$, where
$$
u=u_1w_1,v=v_1, w=w_1.
$$
Then
$$
\begin{array}{ll}
u_1&=(x^ay^b)^{k-l}[\gamma+x^cy^d]^{-1}\\
v_1&=(x^ay^b)^t(\alpha+z)\\
w_1&=(x^ay^b)^l[\gamma+x^cy^d].
\end{array}
$$
Define new  regular parameters $\overline x,\overline y, z$ at $p$
by $x=\overline x\lambda_x,y=\overline y\lambda_y$, where $\lambda_x,\lambda_y$ are unit series such  that
$$
x^ay^b=(\gamma+x^cy^d)^{\frac{-1}{l}}\overline x^a\overline y^b.
$$
Then
$$
\begin{array}{ll}
u_1&=(\overline x^a\overline y^b)^{k-l}(\gamma+x^cy^d)^{-\frac{k}{l}}\\
v_1&=(\overline x^a\overline y^b)^t[\gamma+x^cy^d]^{\frac{-t}{l}}(\alpha+z)\\
w_1&=(\overline x^a\overline y^b)^l.
\end{array}
$$
If
$$
\frac{\partial \gamma}{\partial z}(0,0,0)\ne 0,
$$
 then $w_1, u_1$ have a form (\ref{eqTF22}) of Definition \ref{torf} at $p$. If
$$
\frac{\partial \gamma}{\partial z}(0,0,0)= 0,
$$
then $w_1,v_1$ have a form (\ref{eqTF22}) of Definition \ref{torf} at $p$. 

There is a similar analysis (to the case when (\ref{eqTF22}) holds at $p$) if $p$ is a 1-point 
(and $q_1$ is a 3-point). After possibly permuting $u,v,w$, we find permissible parameters $u_1,v_1,w_1$ at $q$ such that
one of the forms (\ref{eq17}) -- (\ref{eq18}) hold. As in the above paragraph, we see that (after possibly interchanging $u_1,v_1,w_1$)
 $u_1,v_1,w_1$ have a form (\ref{eqTF1}) of Definition \ref{torf} at $p$.

Still assuming that $q$ is a 3-point, assume that  $q_1$ is a 2-point.  If $p$ is a 3-point then $p'$ is also a 3-point and
there are  permissible parameters
$x,y,z$  for $u,v,w$ at $p$ such that
$$
\begin{array}{ll}
u&=x^ay^bz^c\gamma_1\\
v&=x^dy^ez^f\gamma_2\\
w&=x^ly^mz^n\gamma_3
\end{array}
$$
where $\gamma_1,\gamma_2,\gamma_3$ are  unit series, and
$$
\text{rank}\left(\begin{array}{lll}
a&b&c\\
d&e&f\\
l&m&n
\end{array}
\right)\ge 2.
$$
  After possibly permuting $u, v$ and $w$, we may assume that $q_1$ has
permissible parameters $u_1, v_1, w_1$ defined by 
\begin{equation}\label{eq305}
u=u_1, v=v_1,w=u_1(w_1+\alpha)
\end{equation}
with $\alpha\ne 0$. Thus $(a,b,c)=(l,m,n)$, so that 
$$
\text{rank}\left(\begin{array}{lll}
a&b&c\\
d&e&f
\end{array}
\right)=2
$$
 and $\tau_{f_1}(p)=\tau_{f\circ\Phi_1}(p)=\tau_f(p')\ge 0$. We thus have that $u_1,v_1$ are toroidal forms at $p$, 
 of type (\ref{eqTF3}) of Definition \ref{torf} and
  $\tau_{f_1}(p)=\tau_{f}(p')$.

We have a similar analysis if $q_1$ is a 2-point, $p$ is a 2-point ($q$ is a 3-point), and $u,v$ satisfy (\ref{eqTF21})
of Definition \ref{torf} at $p$.  Then $q_1$ has permissible
parameters $u_1,v_1,w_1$ satisfying (\ref{eq305}), and $u_1,v_1$ are toroidal forms of type (\ref{eqTF21}) of Definition \ref{torf} at $p$.

Now assume that $p$ is a 2-point, $q_1$ is a 2-point, ($q$ is a 3-point) and $u,v$ satisfy (\ref{eqTF22}) of Definition \ref{torf} at $p$ and (\ref{eqTF21}) of
Definition \ref{torf} does not hold at $p$ for any permutation of $u,v,w$. Then we have permissible parameters $x,y,z$ at $p$ such that
$$
\begin{array}{ll}
u&=(x^ay^b)^l\\
v&=(x^ay^b)^t(\beta+z)\\
w&=(x^ay^b)^m(\gamma(x^ay^b,z)+x^cy^d)
\end{array}
$$
where $\gamma$ is a unit series and $\beta\ne 0$. If $u=w=0$ are local equations of $C$ at $q$, then $q_1$ has
regular parameters $u_1,v_1,w_1$ defined by 
\begin{equation}\label{eq306}
u=u_1, v=v_1, w=u_1(w_1+\alpha)
\end{equation}
with $\alpha\ne 0$. Thus $u_1,v_1$ are toroidal forms of type (\ref{eqTF22}) of Definition \ref{torf} at $p$. 

If $u=v=0$ are local equations of $C$ at $q$, then 
$q_1$ has permissible parameters $u_1,v_1,w_1$ defined by 
\begin{equation}\label{eq307}
u=u_1, v=u_1(v_1+\beta), w=w_1
\end{equation}
with $\beta\ne 0$ and we thus have $t=l$. We have
$$
\begin{array}{ll}
u_1&=(x^ay^b)^l\\
w_1&=(x^ay^b)^m(\gamma+x^cy^d)\\
v_1&=z
\end{array}
$$
$u_1,w_1,v_1$ have the  form 2 (c) of Definition \ref{Def1} at $p$.

There is a similar analysis in the case when ($q$ is a 3-point) $q_1$ is a 2-point and $p$ is a 1-point. Then (after possibly permuting
$u,v,w$) $u,v$ satisfy (\ref{eqTF1}) of Definition \ref{torf} at $p$. If $u=w=0$ are local equations of $C$ at $q$, then
$q_1$ has regular parameters $u_1,v_1,w_1$ defined by (\ref{eq306}), and
$u_1,v_1$ are toroidal forms of type (\ref{eqTF1}) of Definition \ref{torf}
at $p$.  If $u=v=0$ are local equations of $C$ at $q$, then $q_1$ has regular parameters $u_1,v_1,w_1$ defined by (\ref{eq307}), and $u_1,w_1,v_1$ have the form 2 (b) of Definition \ref{Def1} at $p$.

Comparing the above expressions, we see that $f_1$ is prepared.

\end{pf}

\begin{Remark}\label{Remark422} Suppose that $f:X\rightarrow Y$ is prepared and $\Phi_1:X_1\rightarrow X$ is either the blow up a 3-point or the blow up
of a 2-curve. Then $f_1=f\circ\Phi_1:X_1\rightarrow Y$ is prepared. If $p_1\in X_1$ is a 3-point then $p=\Phi_1(p_1)$ is a 3-point and 
$\tau_{f_1}(p_1)=\tau_f(p)$. If $D_X$ is cuspidal for $f$, then $D_{X_1}$ is cuspidal for $f_1$.
\end{Remark}

The remark follows by substitution of local forms of $\Phi_1$ into local forms (Definition \ref{Def1}) of the prepared morphism $f$.

\begin{Lemma}\label{Lemma7}  Suppose that $f:X\rightarrow Y$ is prepared, 
$\Omega$ is a set of 2-points of $Y$, and we have assigned to   each $q_0\in \Omega$ permissible parameters $u=u_{q_0},v=v_{q_0},w=w_{q_0}$
at $q_0$. Then there exists a commutative diagram
\begin{equation}\label{eq26}
\begin{array}{rll}
X_1&\stackrel{f_1}{\rightarrow}&Y_1\\
\Phi \downarrow&&\downarrow \Psi\\
X&\stackrel{f}{\rightarrow}&Y
\end{array}
\end{equation}
such that
\begin{enumerate}
\item[1.] $f_1$ is prepared.
\item[2.] $\Phi,\Psi$ are products of blow ups of 2-curves.
\item[3.] Let
$$
\Omega_1=\left\{\begin{array}{l}
q_1\in Y_1 \text{ such that $q_1$ is a 2-point and $q_1=f_1(p_1)$}\\
\text{for some 3-point $p_1\in (f\circ\Phi)^{-1}(\Omega)$}.
\end{array}
\right.
$$
Suppose that $q_1\in\Omega_1$ with $q_0=\Psi(q_1)\in\Omega$. Then there exist
permissible parameters $u_1,v_1,w_1$ at $q_1$ such that
$$
\begin{array}{ll}
u_{q_0}&=u_1^{\overline a}v_1^{\overline b}\\
v_{q_0}&=u_1^{\overline c}v_1^{\overline d}\\
w_{q_0}&=w_1
\end{array}
$$
for some $\overline a,\overline b,\overline c,\overline d\in {\bf N}$ with
$\overline a\overline d-\overline b\overline c=\pm 1$, and if $p_1\in f_1^{-1}(q_1)$ is
a 3-point, then there exist permissible parameters $x,y,z$ at $p_1$ for $u_1,v_1,w_1$
such that 
\begin{equation}\label{eq361}
\begin{array}{ll}
u_1&=x_1^{a_1}y_1^{b_1}z_1^{c_1}\\
v_1&=x_1^{d_1}y_1^{e_1}z_1^{f_1}\\
w_1&=\gamma_1+N_1
\end{array}
\end{equation}
where $N_1=x_1^{g_1}y_1^{h_1}z_1^{i_1}$,
with $\text{rank}(u_1,v_1,N_1)=3$, $\gamma_1=\sum_i\alpha_iM_i$
where $\alpha_i\in \bold k$ and each $M_i$ is a monomial in $x_1,y_1,z_1$ such that there are expressions  
\begin{equation}\label{eq201}
M_i^{e_i}=u_1^{a_i}v_1^{b_i}
\end{equation}
with $a_i,b_i,e_i\in{\bf N}$ and  $\text{gcd}(a_i,b_i,e_i)=1$ for all $i$.
Further, there is a bound $r\in {\bf N}$ such that $e_i\le r$ for all $M_i$ in expressions (\ref{eq201}).
\item[4.] If $D_X$ is cuspidal for $f$, then $D_{X_1}$ is cuspidal for $f_1$.
\item[5.] $\Phi$ is an isomorphism above $f^{-1}(Y-\Sigma(Y))$.
\end{enumerate}
\end{Lemma}

\begin{pf} Suppose that $q_0\in\Omega$. Let the 3-points in $f^{-1}(q_0)$ be $\{p_1,\ldots,p_t\}$. Each $p_i$ has
permissible parameters $x,y,z$ such that there is an expression
of the form (\ref{eq16}), 
\begin{equation}\label{eq22}
\begin{array}{ll}
u_{q_0}&=x^ay^bz^c\\
v_{q_0}&=x^dy^ez^f\\
w_{q_0}&=\gamma+N
\end{array}
\end{equation}
where $N=x^gy^hz^i$, $\gamma=\sum\alpha_iM_i$ with relations 
\begin{equation}\label{eq23}
M_i^{e_i}=u^{a_i}v^{b_i}
\end{equation}
with $a_i,b_i,e_i\in{\bf Z}$, $e_i>0$.

We construct an infinite commutative diagram of morphisms 
\begin{equation}\label{eq416}
\begin{array}{rll}
\vdots&&\vdots\\
\downarrow&&\downarrow\\
X_n&\stackrel{f_n}{\rightarrow}&Y_n\\
\Phi_n\downarrow&&\downarrow\Psi_n\\
\vdots&&\vdots\\
\Phi_2\downarrow&&\downarrow\Psi_2\\
X_1&\stackrel{f_1}{\rightarrow}&Y_1\\
\Phi_1\downarrow&&\downarrow\Psi_1\\
X&\stackrel{f}{\rightarrow}&Y
\end{array}
\end{equation}
as follows. Order the 2-curves of $Y$, and let $\Psi_1:Y_1\rightarrow Y$ be the blow up of the 2-curve $C$ of smallest order. Then construct
(by Lemma \ref{Lemma1}) a commutative diagram 
\begin{equation}\label{eq415}
\begin{array}{rll}
X_1&\stackrel{f_1}{\rightarrow}&Y_1\\
\Phi_1\downarrow&&\downarrow\Psi_1\\
X&\stackrel{f}{\rightarrow}&Y
\end{array}
\end{equation}
where $f_1$ is prepared, $\Phi_1$ is a product of blow up of 2-curves and $\Phi_1$ is an isomorphism above $f^{-1}(Y-C)$. Order the 2-curves of $Y_1$
so that the 2-curves contained in the exceptional divisor of $\Psi_1$ have larger order than the order of the (strict transforms of the) 2-curves of $Y$.

Let $\Psi_2:Y_2\rightarrow Y_1$ be the blow up of the 2-curve $C_1$ on $Y_1$ of smallest order, and construct a commutative diagram
$$
\begin{array}{rll}
X_2&\stackrel{f_2}{\rightarrow}&Y_2\\
\Phi_2\downarrow&&\downarrow\Psi_2\\
X_1&\stackrel{f_1}{\rightarrow}&Y_1
\end{array}
$$
as in (\ref{eq415}). We now iterate to construct (\ref{eq416}). Let $\overline\Psi_n=\Psi_1\circ\cdot\circ\Psi_n:Y_n\rightarrow Y$,
$\overline\Phi_n=\Phi_1\circ\cdots\circ\Phi_n:X_n\rightarrow X$.

Let $\nu$ be a 0-dimensional valuation of ${\bold k}(X)$. Let $p_n$ be the center of
$\nu$ on $X_n$, $q_n=f_n(p_n)$. We will say that $\nu$ is resolved on $X_n$ if one of the
following holds:
\begin{enumerate}
\item[1.] $\overline \Psi_n(q_n)\not\in \Omega$ or
\item[2.] $\overline\Psi_n(q_n)=q_0\in \Omega$ and
\begin{enumerate}
\item $p_n$ is not a 3-point
or
\item $p_n$ is a 3-point such that a form (\ref{eq361}) holds for $p_n$ and $q_n=f_n(p_n)$ so that (\ref{eq201}) holds.
\end{enumerate}
\end{enumerate}
Observe that if $\nu$ is resolved on $X_n$, then there exists a neighborhood $U$ of
the center of $\nu$
 in $X_n$ such that if $\omega$ is a 0-dimensional valuation of ${\bold k}(X)$ whose center is in $U$, then $\omega$
is resolved on $X_n$, and if $n'>n$,
 then $\nu$ is resolved on $X_{n'}$.

We will now show that for any 0-dimensional valuation $\nu$ of ${\bold k}(X)$, there exists $n\in{\bf N}$  such that $\nu$ is resolved on $X_{n}$.

If the center of $\nu$ on $Y$ is not in $\Omega$ or if the center of $\nu$ on $X$ is not a
3-point, then $\nu$ is resolved on $X$, so we may assume that the center
of $\nu$ on $Y$ is $q_0\in \Omega$ and the center  of $\nu$ on $X$ is a 3-point $p$.

Suppose that  $\nu(u_{q_0})$ and $\nu(v_{q_0})$ are rationally dependent. Then there exists $n$ such that the center of $\nu$ on $Y_{n}$
is a 1-point. 
 Thus $\nu$ is resolved on $X_{n}$.

Suppose that $\nu(u_{q_0})$ and $\nu(v_{q_0})$ are rationally independent. At the center $p$ of $\nu$ on
$X$, 
$$
u=u_{q_0}, v=v_{q_0}, w=w_{q_0}
$$
have an expression (\ref{eq16}). We may identify $\nu$ with an extension of $\nu$ to the
quotient field of $\hat{\cal O}_{X,p}$ which dominates $\hat{\cal O}_{X,p}$.
We have 
$$
u=x^ay^bz^c, v=x^dy^ez^f,
$$
and
$$
M_i^{e_i}=u^{k_i}v^{l_i}=x^{a_i}y^{b_i}z^{c_i}
$$
with $k_i,l_i\in{\bf Z}$, $e_i>0$, $a_i,b_i,c_i\in{\bf N}$. Thus, for all $i$,
$(k_i,l_i)\in\sigma$, where
$$
\sigma=\{(k,l)\in{\bf Q}^2\mid ka+ld\ge 0, kb+le\ge 0, kc+lf\ge 0\}.
$$
Since
$$
\text{rank}\left(\begin{array}{lll}
a&b&c\\
d&e&f
\end{array}\right)=2,
$$
$\sigma$ is a rational polyhedral cone which contains no nonzero linear subspaces, and is contained in the (irrational)
half space
$$
\{(k,l)\mid k\nu(u)+l\nu(v)\ge 0\}.
$$
Let $\lambda_1=(m_1,m_2)$, $\lambda_2=(n_1,n_2)$ be integral vectors such that $\sigma={\bf Q}_+\lambda_1+{\bf Q}_+\lambda_2$.
Since $\lambda_1,\lambda_2$ are rational points in $\sigma$, we have $\nu(u^{m_1}v^{m_2})>0$ and $\nu(u^{n_1}v^{n_2})>0$.

Since $\nu(u)$ and $\nu(v)$ are rationally independent, there exists (by Theorem 2.7 \cite{C1}) a sequence of quadratic
transforms ${\bold k}[u,v]\rightarrow {\bold k}[u_1,v_1]$  such that the center of $\nu$ on ${\bold k}[u_1,v_1]$ is $(u_1,v_1)$, there is an expression
$$
u=u_1^{\overline a}v_1^{\overline b}, v=u_1^{\overline c}v_1^{\overline d}
$$
for some $\overline a,\overline b,\overline c,\overline d\in{\bf N}$ with $\overline a\overline d-\overline b\overline c=\pm1$, and
$u^{m_1}v^{m_2},u^{n_1}v^{n_2}\in {\bold k}[u_1,v_1]$. Thus there exists a rational polyhedral cone $\sigma_1\subset{\bf Q}^2$ containing
$\lambda_1$ and $\lambda_2$ such that ${\bold k}[u_1,v_1]={\bold k}[\sigma_1\cap{\bf Z}^2]$. We thus have $M_i^{e_i}\in {\bold k}[u_1,v_1]$ for all $i$, so that
for all $i$, $M_i^{e_i}$ is a monomial in $u_1$ and $v_1$.

 There exists $n$
such that the center of $\nu$ on $Y_{n}$ has permissible parameters $u_1,v_1,w_1$
where 
\begin{equation}\label{eq21}
u=u_1^{\overline a}v_1^{\overline b}, v=u_1^{\overline c}v_1^{\overline d}, w=w_1.
\end{equation}
 We thus have $\tilde a_i,\tilde b_i\in{\bf N}$ such that

\begin{equation}\label{eq117}
M_i^{e_i}=u_1^{\tilde a_i}
v_1^{\tilde b_i}
\end{equation}
 for all $i$. 
 Let $p_n$ be the center of $\nu$ on $X_{n}$. If $p_n$ is not a 3-point,
 then $\nu$ is resolved on $X_{n}$.

If $p_n$ is a 3-point, then $u_1,v_1,w_1$ (defined by (\ref{eq21})) are permissible parameters at $f_{n}(p_n)$. There exist permissible parameters $x_1,y_1,z_1$ at $p_n$
for $u_1,v_1,w_1$ defined by 
$$
\begin{array}{ll}
x&=x_1^{a_{11}}y_1^{a_{12}}z_1^{a_{13}}\\
y&=x_1^{a_{21}}y_1^{a_{22}}z_1^{a_{23}}\\
z&=x_1^{a_{31}}y_1^{a_{132}}z_1^{a_{33}}
\end{array}
$$
where $a_{ij}\in{\bf N}$, and $\text{Det}(a_{ij})=\pm 1$. Substituting into the expression (\ref{eq16}) of $u,v,w$ at $p$, we have
 expressions

$$
u_1=x_1^{\tilde a}y_1^{\tilde b}z_1^{\tilde c},
v_1=x_1^{\tilde d}y_1^{\tilde e}z_1^{\tilde f},
N=x_1^{\tilde g}y_1^{\tilde h}z_1^{\tilde i},
$$
where
$$
\text{Det}\left(\begin{array}{lll}\tilde a&\tilde b&\tilde c\\
\tilde d&\tilde e&\tilde f\\
\tilde g&\tilde h&\tilde i\end{array}\right)\ne 0.
$$

Let $\vec v_1=(\tilde a,\tilde b,\tilde c)$, $\vec v_2=(\tilde d,\tilde e,\tilde f)$, 
$$
\sigma_2={\bf Q}_+\vec v_1+{\bf Q}_+\vec v_2\subset {\bf Q}^2.
$$
 By Gordon's
Lemma, (Proposition 1 \cite{F}) $\sigma_2\cap {\bf Z}^3$ is a finitely generated semi group. Let $\vec w_1,\ldots,\vec w_n\in\sigma_2\cap {\bf Z}^3$ be generators. There
exists $0\ne r\in{\bf N}$ and $\delta_j,\epsilon_j\in{\bf N}$ such that
$$
\vec w_j=\frac{\delta_j}{r}\vec v_1+\frac{\epsilon_j}{r}\vec v_2
$$
for $1\le j\le n$. Since the exponents of 
$$
M_i=u_1^{\frac{\tilde a_i}{e_i}}v_1^{\frac{\tilde b_i}{e_i}}
$$
are in $\sigma_2\cap {\bf Z}^3$ for all $i$, we have an expression
$$
M_i=u_1^{\frac{\overline a_i}{r}}v_1^{\frac{\overline b_i}{r}}
$$
with $\overline a_i,\overline b_i\in{\bf N}$ for all $M_i$  appearing in the expansion (\ref{eq16}) of $w$. Thus $\nu$ is resolved on $X_{n}$.

By compactness of the Zariski-Riemann manifold of $X$ (\cite{Z1}), there exist finitely many $X_i$,
$1\le i\le t$, such that  the center of any 0-dimensional valuation $\nu$ of ${\bold k}(X)$ is resolved on some $X_i$.
Thus $X_t\rightarrow Y_t$ satisfies the conclusions of the Lemma.

\end{pf}

\begin{Definition}\label{Def357} Suppose that $f:X\rightarrow Y$ is prepared, and $q\in Y$ is a 
2-point. Permissible parameters $u,v,w$ at $q$ are {\bf super parameters} for $f$ at $q$ if 
at all
$p\in f^{-1}(q)$, there exist permissible parameters $x,y,z$ for $u,v,w$ at $p$ such that we have one of the forms:
\begin{enumerate}
\item[1.] $p$ is a 1-point
$$
\begin{array}{ll}
u&=x^a\\
v&=x^b(\alpha+y)\\
w&=x^c\gamma(x,y)+x^d(z+\beta)
\end{array}
$$
where $\gamma$ is a unit series (or zero), $0\ne\alpha\in \bold k$ and $\beta\in \bold k$,
\item[2.] $p$ is  a 2-point of the form of  (\ref{eqTF21}) of Definition \ref{torf}
$$
\begin{array}{ll}
u&=x^ay^b\\
v&=x^cy^d\\
w&=x^ey^f\gamma(x,y)+x^gy^h(z+\beta)
\end{array}
$$
where $ad-bc\ne 0$, $\gamma$ is a unit series (or zero), and $\beta\in \bold k$.
\item[3.] $p$ is a 2-point of the form of  (\ref{eqTF22}) of Definition \ref{torf}
$$
\begin{array}{ll}
u&=(x^ay^b)^k\\
v&=(x^ay^b)^t(\alpha+z)\\
w&=(x^ay^b)^l\gamma(x^ay^b,z)+x^cy^d
\end{array}
$$
where $0\ne\alpha\in \bold k$, $ad-bc\ne 0$ and $\gamma$ is a unit series (or zero).
\item[4.] $p$ is a 3-point
$$
\begin{array}{ll}
u&=x^ay^bz^c\\
v&=x^dy^ez^f\\
w&=x^gy^hz^i\gamma+x^jy^kz^l
\end{array}
$$
where $\text{rank}(u,v,x^jy^kz^l)=3$, $\text{rank}(u,v,x^gy^hz^i)=2$
and
$\gamma$ is a unit series in monomials $M$ such that $\text{rank}(u,v,M)=2$ (or $\gamma$ is zero).
\end{enumerate}
\end{Definition}

\begin{Lemma}\label{Lemma358}
Suppose that $f:X\rightarrow Y$ is prepared, $q\in Y$ is a 2-point and $u,v,w$ are permissible
parameters at $q$. Then there exists a commutative diagram 
\begin{equation}\label{eq355}
\begin{array}{rll}
X_1&\stackrel{f_1}{\rightarrow}&Y_1\\
\Phi_1\downarrow&&\downarrow\Psi_1\\
X&\stackrel{f}{\rightarrow}&Y
\end{array}
\end{equation}
such that
\begin{enumerate}
\item[1.] $f_1$ is prepared.
\item[2.] $\Phi_1$ is a product of blow ups of 2-curves and 3-points such that $\Phi_1$ is an isomorphism above $f^{-1}(Y-\Sigma(Y))$. $\Psi_1$ is a product
 of blow ups of 2-curves.
\item[3.] Suppose that $q_1\in\Psi_1^{-1}(q)$ is a 2-point, so that $q_1$ has permissible parameters
$u_1,v_1,w_1$ defined by 
\begin{equation}\label{eq356}
\begin{array}{ll}
u&=u_1^av_1^b\\
v&=u_1^cv_1^d\\
w&=w_1
\end{array}
\end{equation}
for some $a,b,c,d\in{\bf N}$ with $ad-bc\ne 0$. Then $u_1,v_1,w_1$ are super parameters at $q_1$.
\item[4.] If $D_X$ is cuspidal for $f$, then $D_{X_1}$ is cuspidal for $f_1$.
\item[5.] Suppose that $p_1\in X_1$ is a 3-point. Then $p=\Phi_1(p)$ is a 3-point and $\tau_{f_1}(p_1)=\tau_f(p)$.
\end{enumerate}
\end{Lemma}

\begin{pf} By Lemma \ref{Lemma1} and Remark \ref{Remark422}, any diagram (\ref{eq355}) satisfying 2 satisfies 1, 4 and 5.
Further, if $p\in f^{-1}(q)$,  $u,v,w$ are super parameters  at $p$, $p_1\in \Phi_1^{-1}(p)$ is such that $f_1(p_1)=q_1$
is a 2-point,  then the permissible parameters $u_1,v_1,w_1$ of
(\ref{eq356}) at $q_1$ are super parameters at $p_1$.

\vskip .2truein
\noindent{\bf Step 1.} We will show that there exists a sequence of blow ups of 2-curves
and 3-points $\Phi_1:X_1\rightarrow X$ such that $\Phi_1$ is an isomorphism over $f^{-1}(Y-\Sigma(Y))$ and $u,v,w$ are super parameters at all 3-points 
$p\in(f\circ\Phi_1)^{-1}(q)$.

 Suppose that $p\in f^{-1}(q)$ is a 3-point, so that there exist permissible parameters
 $x,y,z\in\hat{\cal O}_{X,p}$ at $p$ for $u,v,w$ such that  
 \begin{equation}\label{eq420}
\begin{array}{ll}
u&=x^ay^bz^c\\
v&=x^dy^ez^f\\
w&=g(x,y,z)+N
\end{array}
\end{equation}
of the form of (\ref{eqTF3}) of Definition \ref{torf} and (\ref{eqTF3w}) of Lemma \ref{Lemmatorf}. There exist
regular parameters $\overline x,\overline y,\overline z$ in ${\cal O}_{X,p}$, and unit
series $\lambda_1,\lambda_2,\lambda_3\in\hat{\cal O}_{X,p}$ such that
$$
x=\overline x\lambda_1, y=\overline y\lambda_2,z=\overline z\lambda_3.
$$
$\overline x\overline y\overline z=0$ is a local equation of $D_X$ at $p$.

There is an expression $g=\sum\alpha_iM_i$ where $\alpha_i\in \bold k$ and $M_i=x^{a_i}y^{b_i}z^{c_i}$ are monomials in $x,y,z$ such that $\text{rank}(u,v,M_i)=2$. Let $I^p$ be the ideal in ${\cal O}_{X,p}$ generated by the $\overline x^{a_i}\overline y^{b_i}\overline z^{c_i}$ for $a_i,b_i,c_i$ appearing in some $M_i$. 
There exists an $r$ such that 
$$
I^p=(\overline x^{a_0}\overline y^{b_0}\overline z^{c_0},\overline x^{a_1}\overline y^{b_1}\overline z^{c_1},\ldots,
\overline x^{a_r}\overline y^{b_r}\overline z^{c_r}).
$$
We have relations $M_i^{e_i}=u^{\alpha_i}v^{\beta_i}$ with $e_i,\alpha_i,\beta_i\in{\bf Z}$ and $e_i>0$ for all $i$.
Thus for $a\in\text{spec}({\cal O}_{X,p})$, $(I^p)_a$ is principal if $u$ or $v$ is not in $a$.

By Lemma \ref{Lemma419}, there exists a sequence of blow ups of 2-curves and 3-points $\Phi_1:X_1\rightarrow X$ such that $\Phi_1$ is an isomorphism over
$f^{-1}(X-\Sigma(Y))$ and $I^p{\cal O}_{X_1,p_1}$ is invertible for all 3-points $p\in X$ such that $f(p)=q$ and $p_1\in\Phi_1^{-1}(p)$.  
 $f\circ\Phi_1:X_1\rightarrow X$ is prepared by Remark \ref{Remark422}.

Suppose that  $p_1\in X_1$ is a 3-point. Then $p=\Phi_1(p_1)$ is also a 3-point, and $\tau_{f\circ\Phi_1}(p_1)=\tau_f(p)$ (by Remark \ref{Remark422}).
Let notation be as in (\ref{eq420}). There exist permissible parameters $x_1,y_1,z_1$ for the permissible parameters $u,v,w$ at $p_1$ such that $x_1,y_1,z_1$ are defined by
\begin{equation}
\begin{array}{ll}
x&= x_1^{a_{11}} y_1^{a_{12}}z_1^{a_{13}}\\
 y&= x_1^{a_{21}} y_1^{a_{22}} z_1^{a_{23}}\\
 z&= x_1^{a_{31}} y_1^{a_{32}} z_1^{a_{33}},
\end{array}
\end{equation}
where $a_{ij}\in\bf N$, $\text{Det}(a_{ij})=\pm 1$.
Thus all of the $M_i$ and $N$ are distinct monomials when expanded in the variables $x_1,y_1,z_1$, so that
$g(x,y,z)$ is a monomial in $x_1,y_1,z_1$ times a unit series (in $x_1,y_1,z_1$).
Thus $u,v,w$ are super parameters at $p_1$.

\vskip .2truein
\noindent{\bf Step 2.} We will show that there exists a sequence of blow ups of 2-curves
$\Phi_2:X_2\rightarrow X_1$, where $\Phi_1:X_1\rightarrow X$ is the map constructed in Step 1,
 such that $\Phi_2$ is an isomorphism over $(f\circ\Phi_1)^{-1}(Y-\Sigma(y))$ and $u,v,w$ are super parameters at all 
$p\in (f\circ\Phi_1\circ\Phi_2)^{-1}(q)$ for which $p$ is a 3-point or $p$ is a 2-point of type
(\ref{eqTF21}) of Definition \ref{torf} for $u,v,w$.

 Let $f_1=f\circ\Phi_1$. Let $\nu$ be a 0-dimensional valuation
of ${\bold k}(X)$. Let $p$ be the center of $\nu$ on $X_1$. Say that $\nu$ is resolved on $X_1$ if
\begin{enumerate}
\item[1.] $f_1(p)\ne q$, or 
\item[2.] $p$ is not a 2-point of type (\ref{eqTF21}) of Definition \ref{torf}, or
\item[3.] $p$ is a 2-point of type (\ref{eqTF21}) and $u,v,w$ are super parameters at $p$.
\end{enumerate}

We construct an infinite sequence of morphisms 
\begin{equation}\label{eq418}
\cdots\rightarrow X_n\stackrel{\Phi_n}{\rightarrow}\cdots\stackrel{\Phi_3}{\rightarrow} X_2\stackrel{\Phi_2}{\rightarrow}X_1
\end{equation}
as follows. Order the 2-curves $C$ of $X_1$ such that $q\in (f\circ\Phi_1)(C)\subset\Sigma(Y)$. Let $\Phi_2:X_2\rightarrow X_1$ be the blow up of the
2-curve $C_1$ on $X_1$ of smallest order. Order the 2-curves $C'$ of $X_2$ such that $q\in (f\circ\Phi_1\circ\Phi_2)(C')\subset \Sigma(Y)$ so that
the 2-curves contained in the exceptional divisor of $\Phi_2$ have order larger than the order of the (strict transform of the) 2-curves $C$ of $X$ such that $q\in f(C)\subset \Sigma(Y)$. Let $\Phi_3:Y_3\rightarrow Y_2$ be the blow ups of the 2-curve $C_2$ on $Y_3$ of smallest order, and repeat to inductively construct
the morphisms $\Phi_n:X_n\rightarrow X_{n-1}$. Let 
$\overline\Phi_n=\Phi_2\circ\cdots\Phi_n:X_n\rightarrow X_1$. The morphisms $f\circ\overline\Phi_n$ are prepared by Remark \ref{Remark422}.

Suppose that $p\in f_1^{-1}(q)$ is a 2-point satisfying (\ref{eqTF21}) of Definition \ref{torf}, so there exist
permissible parameters $x,y,z$ at $p$ for $u,v,w$ such that
$$
\begin{array}{ll}
u&=x^ay^b\\
v&=x^cy^d\\
w&=g(x,y)+x^ey^fz
\end{array}
$$
with $ad-bc\ne 0$. There exist regular parameters $\overline x,\overline y,\overline z$ in 
${\cal O}_{X_1,p}$ and unit series $\lambda_1,\lambda_2\in\hat{\cal O}_{X_1,p}$ such that
$$
x=\overline x\lambda_1, y=\overline y\lambda_2.
$$
$\overline x\overline y=0$ is a local equation of $D_{X_1}$ at $p$.

If $\nu(\overline x),\nu(\overline y)$ are rationally independent, then the center of $\nu$ on $X_n$ is a 2-point for all $n$, there exists an $n$ such that the center of $\nu$ on $X_{n}$ is a 2-point satisfying (\ref{eqTF21}) of Definition \ref{torf}
and $u,v,w$ are super parameters at $p_1$, by embedded resolution of plane curve singularities (cf. Section 3.4 and Exercise 3.3 \cite{C3}) applied to $g(x,y)=0$.
If $\nu(\overline x),\nu(\overline y)$ are rationally dependent, then there exists an $n$ 
 such that the
center $p_1$ of $\nu$ on $X_{n}$ is a 1-point. 

By compactness of the Zariski-Riemann manifold of ${\bold k}(X)$ \cite{Z1}, there exists an $n$ such that all points of $X_n$ are resolved.
Thus  there exists a sequence
of blow ups of 2-curves $\Phi_2:X_2\rightarrow X_1$ such that the conclusions of Step 2 hold.

\vskip .2truein
\noindent{\bf Step 3.} We will show that there exists a diagram (\ref{eq355}) satisfying the
conclusions of the lemma.

After replacing $f$ with $f\circ\Phi_1\circ\Phi_2$, we can assume that $f$ satisfies the
conclusions of Step 2.
We construct a sequence of diagrams (\ref{eq355}) satisfying 1, 2, 4 and 5 of the conclusions of
the lemma as follows.  Let $\Psi_1:Y_1\rightarrow Y$ be the blow up of the 2-curve $C$
containing $q$. We order the two 2-curves in $Y_1$ which dominate $C$. Let $\Psi_2:Y_2\rightarrow Y_1$ be the blow up of the 2-curve of smallest order. Now extend the
ordering to the 2-curves of $Y_2$ which dominate $C$, by requiring that the two 2-curves on
the exceptional divisor of $\Psi_2$ which dominate $C$ have larger order than the order of the (strict transform of the) 2-curve on $Y_1$ dominating $C$ (which was not blown up by $\Psi_2$). Now let 
$\Psi_3:Y_3\rightarrow Y_2$ be the blow up of the 2-curve of smallest order. We continue this
process to construct a sequence of blow ups of 2-curves
$$
\cdots\rightarrow Y_n\stackrel{\Psi_n}{\rightarrow} Y_{n-1}\rightarrow\cdots\rightarrow
Y_1\stackrel{\Psi_1}{\rightarrow}Y.
$$
Let $\overline\Psi_n=\Psi_1\circ\cdots\Psi_{n-1}\circ\Psi_n$. Let
$$
U=\left\{\begin{array}{l}
p\in f^{-1}(q)\text{ such that $u,v,w$ have a form (\ref{eqTF1}) or (\ref{eqTF22}) of}\\
\text{Definition \ref{torf} or of 2 (b) or 2 (c) of Definition \ref{Def1} at $p$}
\end{array}\right\}.
$$
$U$ is an open subset of $f^{-1}(q)$. For each $p\in U$, there exists 
$n(p)$ such that $n\ge n(p)$ implies the rational map $\overline\Psi_n^{-1}\circ f$
is defined at $p$, and $(\overline\Psi_n^{-1}\circ f)(p_1)$ is a 1-point, for $p_1$ in some neighborhood $U_p$ of $p$ in $U$. $\{U_p\}$ is an open cover of  $U$,
so there exists a finite subcover $\{U_{p_1},\ldots, U_{p_m}\}$ of $U$.
Let $n=\text{max}\{n(p_1),\ldots, n(p_m)\}$. We have that $\overline\Psi_n^{-1}\circ f(p)$
is a 1-point if $p\in U$.

By Lemma \ref{Lemma1}, we can now construct a commutative diagram
$$
\begin{array}{rll}
X_n&\stackrel{f_n}{\rightarrow}&Y_n\\
\overline\Phi_n\downarrow&&\downarrow\overline\Psi_n\\
X&\stackrel{f}{\rightarrow}&Y
\end{array}
$$
such that 1, 2, 4  and 5 of the conclusions of the lemma hold.

If $p_1\in X_n$ is such that $f_n(p_1)=q_1\in\overline\Psi_n^{-1}(q)$ is a 2-point, then $\overline\Phi_n(p_1)\in Y-U$,
and thus $u,v,w$ have a form  2 or 4 of Definition \ref{Def357} at $p$. As observed at the beginning of the proof, $p_1$ must have
one of the forms 1 -- 4 of Definition \ref{Def357} with respect to the permissible parameters $u_1,v_1,w_1$ at $q_1$
defined by 3 of the statement of Lemma \ref{Lemma358}. Thus $f_n$ satisfies 3 of the statement of Lemma \ref{Lemma358}.
\end{pf}

\begin{Theorem}\label{Theorem5} Suppose that $f:X\rightarrow Y$ is prepared, $\tau=\tau_f(X)\ge 0$ and
 all 3-points $p$ of $X$ such that $\tau_f(p)=\tau$ map to 2-points of $Y$.
 
 Let
$$
\Omega=\left\{\begin{array}{l}
q\in Y\text{ such that $q$ is a 2-point and}\\
\text{there exists a 3-point $p\in f^{-1}(q)$}
\end{array}\right\}
$$
 
 Then there exists a commutative
diagram
$$
\begin{array}{rll}
X_1&\stackrel{f_1}{\rightarrow}&Y_1\\
\Phi\downarrow&&\Psi\downarrow\\
X&\stackrel{f}{\rightarrow}&Y
\end{array}
$$
such that the following conditions hold:
\begin{enumerate}
\item[1.] $\Phi$ is a product of blow ups of 2-curves and 3-points such that $\Phi$ is an isomorphism above $f^{-1}(Y-\Sigma(Y))$.
$\Psi$ is a product of blow ups of 2-curves.
\item[2.]  $f_1$ is prepared.
\item[3.] All 3-points $p_1$ of $X_1$ such that $\tau_{f_1}(p_1)=\tau$ map to 2-points of $Y_1$, and if $p_1\in X_1$ is a 3-point,
 then $\tau_{f_1}(p_1)=\tau_{f}(\Phi(p_1))$.
\item[4.] Let
$$
\Omega_1=\left\{\begin{array}{l}q\in Y_1 \text{ such that $q$ is a 2-point, $q\in\Psi^{-1}(\Omega)$}\\
 \text{and there exists a
3-point $p\in f_1^{-1}(q)$}.
\end{array}
\right\}.
$$
If $q\in \Omega_1$ is a 2-point and $p_1,\ldots, p_r\in f_1^{-1}(q)$ are the 3 points
in $f_1^{-1}(q)$, then there exist  $u,v\in{\cal O}_{Y_1,q}$ 
and $w_i\in \hat{\cal O}_{Y_1,q}$ for $1\le i\le r$ such that $u,v,w_i$ are (formal) permissible
parameters at $q$ for $1\le i\le r$ and at the point $p_i$ we have permissible parameters
$x,y,z$ for $u,v,w_i$ such that  we have an expression
\begin{enumerate}
\item
$$
\begin{array}{ll}
u&=x^ay^bz^c\\
v&=x^dy^ez^f\\
w_i&=M\gamma
\end{array}
$$
where $\gamma$ is a unit series, $M$ is a monomial in $x,y,z$ and there is a relation
$$
M^{e_i}=u^{a_i}v^{b_i}
$$
with $a_i,b_i,e_i\in{\bf Z}$, $e_i>1$ and $\text{gcd}(a_i,b_i,e_i)=1$
or
\item
$$
\begin{array}{ll}
u&=x^ay^bz^c\\
v&=x^dy^ez^f\\
w_i&=x^gy^hz^i
\end{array}
$$
where
$$
\left(\begin{array}{lll}
a&b&c\\d&e&f\\g&h&i\end{array}\right)
$$
has rank 3.
\end{enumerate}
Further,  there exists $w\in {\cal O}_{Y_1,q}$ and
$\lambda_i(u,v)\in {\bold k}[[u,v]]$ for $1\le i\le r$ such that 
$$
w_i=w-\lambda_i(u,v)
$$
 for
$1\le i\le r$.
\item[5.] For $q\in \Omega_1$  and for $u,v,w_i$  with $1\le i\le r$ in 4 above,
$u,v,w_i$ are super parameters at $q$. 
\item[6.] Suppose that $D_X$ is cuspidal for $f$. Then $D_{X_1}$ is cuspidal for $f_1$.
\item[7.] Suppose that $p_1\in X_1$ is a 3-point. Then $p=\Phi_1(p)$ is a 3-point and $\tau_{f_1}(p_1)=\tau_f(p)$.
\end{enumerate}
\end{Theorem}

\begin{Remark}\label{Remark225} The proof actually produces expressions 4 (a) with
$$
w_i=M\gamma=g_i(u^{\frac{1}{l}},v^{\frac{1}{l}})+N
$$
where $l\in {\bf N}$ and $g_i$ is a series, $N$ is a monomial in $x,y,z$ and $\text{rank}(u,v,N)=3$.
\end{Remark}

\begin{pf} 
 
By Lemma \ref{Lemma7}, there exists a diagram 
\begin{equation}\label{eq19}
\begin{array}{rll}
X_1&\stackrel{f_1}{\rightarrow}&Y_1\\
\Phi\downarrow&&\downarrow\Psi\\
X&\stackrel{f}{\rightarrow}&Y
\end{array}
\end{equation}
where $\Phi$, $\Psi$ are products of blow ups of 2-curves
such that for all 2-points $q\in\Psi^{-1}(\Omega)$, there exist (algebraic) permissible parameters
$u_q,v_q,w_q$ at $q$ such that if $p\in X_1$ is a 3-point such that $f(p)=q$,
then there are permissible
parameters $x,y,z$ at $p$ such that there is an expression 
\begin{equation}\label{eq20}
\begin{array}{ll}
u_q&=x^ay^bz^c\\
v_q&=x^dy^ez^f\\
w_q&=\gamma + N
\end{array}
\end{equation}
where $\text{rank}(u_q,v_q)=2$, $\gamma$ is a (possibly trivial if $\tau_{f_1}(p)=0$) series
$\gamma=\sum\alpha_iM_i$ in monomials $M_i$ in $x,y,z$ such that
$\text{rank}(u_q,v_q,M_i)=2$ for all $i$, and $N$ is a monomial in $x,y,z$ such that $\text{rank}(u_q,v_q,N)=3$.
Further, there are relations 
\begin{equation}\label{eq4}
M_i^{e_i}=u_q^{a_i}v_q^{b_i}
\end{equation}
for all $M_i$, with $a_i,b_i,e_i\in{\bf N}$ and $e_i>0$.
 Let $h_p(u_q,v_q)$ be the series  in the monomials $M_i$ of $\gamma$ such that $e_i=1$ in 
(\ref{eq4}).  Let $w_p=w_q-h_p$. Let
$$
\Omega'=\left\{
\begin{array}{l}
q\in Y_1\text{ such that $q$ is a 2-point, $q\in\Psi^{-1}(\Omega)$}\\
\text{and there exists a 3-point $p\in f^{-1}(q)$}
\end{array}\right\}.
$$
After performing a further sequence of blowups of 2-curves above $X_1$, we then achieve (by Lemma \ref{Lemma353}) that 4 of the conclusions of Theorem \ref{Theorem5} hold.
Now the proof follows from Lemma \ref{Lemma358}, applied successively to all $q\in\Omega'$ and permissible
parameters $u,v,w_p$, for $p\in f_1^{-1}(q)$ a 3-point.

\end{pf}

\section{Relations}
In this section, we suppose that $Y$ is a nonsingular projective 3-fold with toroidal structure $D_Y$, and
$f:X\rightarrow Y$ is a birational morphism of nonsingular projective 3-folds, with toroidal structures $D_Y$ and
$D_X=f^{-1}(D_Y)$, such that $D_Y$ contains the fundamental locus of $f$.

\begin{Definition}\label{Def153}
Suppose that $Y$ is a nonsingular 3-fold with SNC divisor $D_Y$. A 2-point pre-relation
$R$ on $Y$ is an association from a finite set $U(R)$ of 2-points of $Y$. If $q\in U(R)$,
then 
\begin{equation}\label{eq131}
R(q)=\left(\begin{array}{l}
S=S_R(q),E_1=E_{R,1}(q),E_2=E_{R,2}(q),w=w_R(q),u=u_R(q),\\
v=v_R(q),e=e_R(q),a=a_R(q),b=b_R(q),\lambda=\lambda_R(q)
\end{array}\right)
\end{equation}
where $E_1,E_2$ are the components of $D_Y$ containing $q$, $a, b, e\in{\bf Z}$,
$\text{gcd}(a,b,e)=1$ and $e>1$.  $u,v\in{\cal O}_{Y,q}$, $w\in \hat{\cal O}_{Y,q}$
are such that $u,v,w$ are (formal) permissible parameters at $q$, $u=0$ is a local equation
of $E_1$, $v=0$ is a local equation of $E_2$. $S=\text{spec}(\hat{\cal O}_{Y,q}/(w))$,
$0\ne \lambda\in \bold k$.

We will also allow 2-point pre-relations with $a=b=-\infty$, $e=1$ and $\lambda=1$. 
\end{Definition}

Observe that if $R$ is a 2-point pre-relation (and $a,b\ne-\infty$) then $R(q)$ is  determined by the expression 
\begin{equation}\label{eq151}
w^e-\lambda u^av^b.
\end{equation}
Depending on the signs of $a$ and $b$, this expression determines a (formal) germ of an (irreducible)  surface
singularity 
\begin{equation}\label{eq152}
F=F_R(q)=0
\end{equation}
 of one of the following forms:
$$
F=w^e-\lambda u^av^b=0
$$
if $a,b\ge 0$ and $a+b>0$,
$$
F=w^eu^{-a}-\lambda v^b=0
$$
if $a<0$, $b>0$,
$$
F=w^ev^{-b}-\lambda u^a=0
$$
if $b<0$, $a>0$.

In the remaining case, $a,b\le 0$,
$$
F=w^eu^{-a}v^{-b}-\lambda
$$
is a unit in $\hat{\cal O}_{Y,q}$ and $F(q)\ne 0$.

If $a,b=-\infty$, then $R(q)$ is determined by 
\begin{equation}\label{eq359}
F_R(q)=w_R(q)=0.
\end{equation}

\begin{Definition}\label{Def199}
A 2-point pre-relation $R$ on $Y$ is algebraic if
there exists a nonsingular irreducible locally closed surface $\Omega(R)\subset Y$ such that $\Omega(R)$ makes SNCs with $D_Y$, $U(R)\subset\Omega(R)$ and
$S_R(q)$ is the (formal) germ of $\Omega(R)$ at $q$ for all $q\in U(R)$.
\end{Definition}

\begin{Definition}\label{Def155}
A 3-point pre-relation
$R$ on $Y$ is an association from a finite set $U(R)$ of 3-points of $Y$. If $q\in U(R)$
then 
\begin{equation}\label{eq133}
R(q)=\left(
\begin{array}{l}
E_1=E_{R,1}(q),E_2=E_{R,2}(q),E_3=E_{R,3}(q),u=u_R(q),v=v_R(q),\\
w=w_R(q),
a=a_R(q),b=b_R(q),
c=c_R(q),\lambda=\lambda_R(q)
\end{array}
\right)
\end{equation}
where $E_1,E_2,E_3$ are the components of $D_Y$ containing $q$, $a,b,c\in{\bf Z}$,
$\text{gcd}(a,b,c)=1$, $\text{min}\{a,b,c\}<0<\text{max}\{a,b,c\}$.
$u,v,w\in{\cal O}_{Y,q}$
are  permissible parameters at $q$ such that  $u=0$ is a local equation
of $E_1$, $v=0$ is a local equation of $E_2$, $w=0$ is a local equation of $E_3$, and
$0\ne \lambda\in \bold k$.
\end{Definition}

Observe that if $R$ is a 3-point pre-relation then $R(q)$ is uniquely determined by the expression 
\begin{equation}\label{eq156}
u^av^bw^c=\lambda.
\end{equation}
Depending on the signs of $a$, $b$ and $c$, this expression determines a germ of an (irreducible) surface
singularity 
\begin{equation}\label{eq157}
F=F_R(q)=0
\end{equation}
 of one of the following forms: 
\begin{equation}\label{eq136}
\begin{array}{l}
F=w^c-\lambda u^{-a}v^{-b}=0\text{ if }a,b\le0,c> 0\\
F=v^{b}-\lambda u^{-a}w^{-c}=0\text{ if }a,c\le0,b> 0\\
F=u^{a}-\lambda v^{-b}w^{-c}=0\text{ if }b,c\le 0, a> 0\\
F=w^{-c}-\frac{1}{\lambda}u^av^b\text{ if }a,b>0,c<0\\
F=v^{-b}-\frac{1}{\lambda}u^aw^c\text{ if }a,c>0,b<0\\
F=u^{-a}-\frac{1}{\lambda}v^bw^c\text{if }b,c>0,a<0
\end{array}
\end{equation}

A pre-relation $R$ is resolved if $F_R(q)$ is a unit in $\hat{\cal O}_{Y,q}$ for all
$q\in U(R)$ (This includes the case $U(R)=\emptyset$).

\begin{Definition}\label{Def154}
A subvariety $G$ of $Y$ is an admissible center for a 2-point pre-relation $R$ on $Y$ if one of the
following holds:
\begin{enumerate}
\item[1.] $G$ is a 2-point.
\item[2.] $G$ is a 2-curve of $Y$.
\item[3.] $G\subset D_Y$ is a nonsingular curve which contains a 1-point and makes SNCs with $D_Y$. If
$q\in U(R)\cap G$ then $S_R(q)$ contains the germ of $G$ at $q$.
If $R$ is algebraic, then $G$ makes SNCs with $\Omega(R)$.
\end{enumerate}
\end{Definition}

\begin{Definition}\label{Def157}
A curve $C\subset Y$ is an  admissible center for a 3-point pre-relation $R$ on $Y$ if $C$ is a 2-curve.
\end{Definition}

Observe that admissible centers are possible centers.

Suppose that $R$ is a 2-point (or a 3-point) pre-relation on $Y$, $G$ is an admissible center for $R$,
and $\Psi:Y_1\rightarrow Y$ is the blow up of $G$.

If $R$ is a 2-point pre-relation then the transform $R^1$ of
$R$ on $Y_1$ is the 2-point  pre-relation on $Y_1$ defined by the condition that
$U(R^1)$ is the union over $q\in U(R)$ of 2-points $q_1$ in $\Psi^{-1}(q)$ such that
$q_1$ is on the strict transform of $w_R(q)=0$.
For such $q_1$, $R^1(q_1)$ is determined by the strict transform of the form $F_R(q)=0$
(\ref{eq152}) (or (\ref{eq359}))
 on $Y_1$
at $q_1$.

If $q\in U(R)\cap G$, and 
$$
u=u_R(q), v=v_R(q), w=w_R(q),
$$
 then $G$ has local equations of one of the following forms at $q$:
\begin{enumerate}
\item[1.] $u=v=w=0$,
\item[2.] $u=v=0$,
\item[3.] $u=w=0$ or $v=w=0$.
\end{enumerate}

If $q_1\in U(R_1)\cap \Psi^{-1}(q)$, then after possibly interchanging $u$ and $v$,
$$
u_1=u_{R^1}(q_1), v_1=v_{R^1}(q_1),  w_1=w_{R^1}(q_1)
$$
are defined, respectively, by
\begin{enumerate}
\item[1.] $u=u_1, v=u_1v_1, w=u_1w_1$,
\item[2.] $u=u_1, v=u_1v_1, w=w_1$,
\item[3.] $u=u_1,v=v_1, w=u_1w_1$.
\end{enumerate}

If $R$ is algebraic, then the transform $R^1$ of $R$ is algebraic, where $\Omega(R^1)$ is the
strict transform of $\Omega(R)$ by $\Psi$.

If $R$ is a 3-point pre-relation then the transform $R^1$ of
$R$ on $Y_1$ is the 3-point  pre-relation on $Y_1$ defined by the condition that
$U(R^1)$ is the union over $q\in U(R)$ of 3-points $q_1$ in $\Psi^{-1}(q)$ such that
$q_1$ is on the strict transform of  the form $F_R(q)=0$ of (\ref{eq157})
 on $Y_1$.
For such $q_1$, $R^1(q_1)$ is determined by the strict transform of the form $F_R(q)=0$ at $q_1$.

After possibly interchanging $u=u_R(q),v=v_R(q)$ and $w=w_R(q)$, 
$$
u_1=u_{R^1}(q_1), v_1=v_{R^1}(q_1),  w_1=w_{R^1}(q_1)
$$
are defined by 
$$
u=u_1, v=u_1v_1, w=w_1.
$$ 

\begin{Definition}\label{Def156}
Suppose that $f:X\rightarrow Y$ is a prepared morphism. A primitive 2-point
relation $R$ for $f$ is
\begin{enumerate}
\item[1.] A 2-point pre-relation $\overline R$ on $Y$,
\item[2.] A set of 3-points $T=T(R)\subset \cup_{q\in U(\overline R)}f^{-1}(q)$
 such that if $p\in T(R)$
and
$$
\overline R(f(p))=(S,E_1,E_2,w,u,v,e,a,b,\lambda_p)
$$
then   there exist permissible parameters $x,y,z$ at $p$ for $u,v,w$ such that
$$
w^e=u^av^b\overline\Lambda(x,y,z)
$$
where $\overline\Lambda$ is a unit series such that $\overline\Lambda(0,0,0)=\lambda_p$ if $a,b\ne-\infty$, and
$u,v,w$ have a monomial form at $p$ if $a=b=-\infty$.
\end{enumerate}
We define $R(p)=\overline R(f(p))$ if $p\in T(R)$, and denote
$$
R(p)=
\left(\begin{array}{l}
S=S_{ R}(p),E_1(p), E_2(p),w=w_{ R}(p),u=u_{ R}(p),\\
v=v_{ R}(p),
e=e_{ R}(p),
a=a_{ R}(p),b=b_{R}(p),\lambda_p=\lambda_{ R}(p)\end{array}\right).
$$
A 2-point relation $R$ for $f$ is
a finite set of 2-point pre-relations $\{\overline R_i\}$ on $Y$
with associated primitive 2-points relations $R_{i}$ for $f$
such that the sets $T(R_{i})$ are pairwise disjoint.
We denote $U(R)=\cup_iU(\overline R_i)$ and $T(R)=\cup_iT(R_{i})$, and
define
$$
R(p)=R_{i}(p)
$$
if $p\in T(R_{i})$.

We further require that 
$$
u_{\overline R_i}(q)=u_{\overline R_j}(q), v_{\overline R_i}(q)=v_{\overline R_j}(q)
$$
if $q\in U(\overline R_i)\cap U(\overline R_j)$.
We will call $\{\overline R_i\}$  the 2-point pre-relations associated to $R$.
We will say that $R$ is algebraic if each $\overline R_i$ is algebraic and 
\begin{equation}\label{eq253}
\Omega(\overline R_i)\cap U(R)=U(\overline R_i)
\end{equation}
for all $i$. For $p\in T(R_i)\subset T(R)$, we denote
$$
R(p)=\left(\begin{array}{l}
S=S_{R}(p),E_1(p), E_2(p),w=w_{ R}(p),u=u_{R}(p),\\
v=v_{R}(p),
e=e_{R}(p),
a=a_{R}(p),b=b_{R}(p),\lambda_p=\lambda_{R}(p)\end{array}\right).
$$
\end{Definition}

\begin{Definition}\label{Def160}
Suppose that $f:X\rightarrow Y$ is a prepared morphism. A primitive 3-point
relation $R$ for $f$ is
\begin{enumerate}
\item[1.] A 3-point pre-relation $\overline R$ on $Y$,
\item[2.] A set of 3-points $T(R)\subset \cup_{q\in U(\overline R)}f^{-1}(q)$ such that if $p\in T(R)$
and
$$
\overline R(f(p))=(E_1,E_2,E_3,u,v,w,a,b,c,\lambda_p),
$$
then there exist permissible parameters $x,y,z$ at $p$ for $u,v,w$ such that

$$
u^av^bw^c=\Lambda(x,y,z)
$$
where $\Lambda$ is a unit series such that $\Lambda(0,0,0)=\lambda_p$.

We define $R(p)=\overline R(f(p))$ if $p\in T(R)$, and denote
$$
R(p)=
\left(
\begin{array}{l}
E_1=E_1(p),E_2= E_2(p), E_3=E_3(p),u=u_{R}(p),v=v_{R}(p),\\
w=w_{}(p),
a=a_{R}(p),b=b_{R}(p), c=c_{R}(p),\lambda_p=\lambda_{R}(p)
\end{array}\right).
$$
\end{enumerate}

A 3-point relation $R$ for $f$ is
a finite set of 3-point pre-relations $\{\overline R_i\}$ on $Y$
with associated primitive 3-points relations $R_{i}$ for $f$
such that the sets $T(R_{i})$ are pairwise disjoint.
We denote $U(R)=\cup_iU(\overline R_i)$ and $T(R)=\cup_iT(R_{i})$, and
define
$$
R(p)=R_{i}(p)
$$
if $p\in T(R_{i})$.

We further require that we have  equalities
$$
u_{\overline R_i}(q)=u_{\overline R_j}(q),
v_{\overline R_i}(q)=v_{\overline R_j}(q),
w_{\overline R_i}(q)=w_{\overline R_j}(q)
$$
if $q\in U(\overline R_i)\cap U(\overline R_j)$. 
We will say that $\{\overline R_i\}$ are the 3-point pre-relations associated to $R$.
For $p\in T(R_i)\subset T(R)$, we denote
$$
R(p)=
\left(
\begin{array}{l}
E_1=E_1(p),E_2= E_2(p), E_3=E_3(p),u=u_{R}(p),v=v_{R}(p),\\
w=w_{R}(p),
a=a_{R}(p),b=b_{R}(p), c=c_{R}(p),\lambda_p=\lambda_{R}(p)
\end{array}\right).
$$
\end{Definition}

A (2-point or 3-point) relation $R$ is resolved if $T(R)=\emptyset$.

\begin{Definition}\label{Def161}
Suppose that $f:X\rightarrow Y$ is prepared, $R$ is a 2-point (respectively 3-point)
relation for $f$ and
$$
\begin{array}{rll}
X_1&\stackrel{f_1}{\rightarrow}&Y_1\\
\Phi\downarrow&&\downarrow\Psi\\
X&\stackrel{f}{\rightarrow}&Y
\end{array}
$$
is a commutative diagram such that
\begin{enumerate}
\item[1.] $\Psi$ is a product of blow ups which are admissible for all of the pre-relations
$\overline R_i$ associated to $R$ (and their transforms) and $\Phi$ is a
product of blow ups of possible centers
\item[2.] $f_1$ is prepared.
\item[3.] Let $\overline R_i^1$ be the transforms of the $\overline R_i$ on $Y_1$ and let
$$
T_i=\{p\in X_1\mid p\text{ is a 3-point and }
p\in\Phi^{-1}(T(R_{i}))\cap f_1^{-1}(U(\overline R_i^1))\}.
$$
Suppose that the condition 2 of Definition \ref{Def156}
(respectively 2 of Definition \ref{Def160}) are satisfied for
$f_1:X_1\rightarrow Y_1$, and all $\overline{R_i}^1$ and $T_i$.
\end{enumerate}
 Then the transform $R^1$ of $R$ for $f_1$ is the 2-point (respectively 3-point)
relation for $f_1$
defined by Definition \ref{Def156} (respectively Definition \ref{Def160}) as
$$
T(R^1)=\cup T_i,
$$
$$
R^1(p)=\overline R_i^{1}(f_1(p))
$$
for $p\in T_i$.
\end{Definition}

\begin{Theorem}\label{Theorem3}
Suppose that $R$ is a 3-point pre-relation on $Y$. Then there exists a sequence of blow ups of 2-curves
$Y_1\rightarrow Y$, such that if $R^1$ is the transform of $R$ on $Y_1$, then
$R^1$ is resolved.
\end{Theorem}

\begin{pf}
Suppose that $\Phi_1:Y_1\rightarrow Y$ is a sequence of blow ups of 2-curves. Let $R^1$ be the transform of $R$ on $Y_1$.
We will
say that a 0-dimensional valuation $w$ of ${\bold k}(Y)$ is resolved on $Y_1$ if the center $q_1$ of $w$ on $Y_1$  is not 
in $U(R^1)$.

Observe that if $\omega$ is resolved on $Y_1$, then there exists an open neighborhood $\Sigma$ of
$q_1$ in $Y_1$ such that all 0-dimensional valuations $\nu$ of ${\bold k}(Y)$ whose center on $Y_1$ is in $\Sigma$ are resolved on $Y_1$. Further, if $\Phi_2:Y_2\rightarrow Y$ is a sequence of blow ups of 2-curves,
which factors through $Y_1$, and $\omega$ is resolved on $Y_1$, then $\omega$ is resolved on $Y_2$.

 We will show that
for each 0-dimensional valuation $\nu$ of ${\bold k}(Y)$, there exists a sequence of blow ups of
 2-curves $Y_{\nu}\rightarrow Y$ such that the center of $\nu$ is resolved on $Y_{\nu}$.

Let $\nu$ be a 0-dimensional valuation of ${\bold k}(Y)$, and suppose that the center of $\nu$
on $Y$ is $q\in U(R)$. Let $F=F_R(q)$ (with the notation of (\ref{eq136})).    The sequence of
blow ups of 2-curves $Y_{\nu}\rightarrow Y$ such that $\nu$ is resolved on $Y_{\nu}$ is constructed as
follows.

After possibly permuting $u,v,w$ in (\ref{eq136}), possibly replacing $\lambda$
with $\frac{1}{\lambda}$, and observing that $\text{gcd}(a,b,c)=1$, we have an expression 
 \begin{equation}\label{eq313}
 F=w^{\overline c}-\lambda u^{\overline a}v^{\overline b}
 \end{equation}
(with $\overline a=\pm a,\overline b=\pm b\ge 0$, $\overline c=\pm c>0$) and such that if $\overline a=0$ then
$\overline c\le \overline b$ and if
$\overline b=0$ then $\overline c\le \overline a$.

Suppose that  $\overline a+\overline b<\overline c$ in (\ref{eq313}), so that $\overline a>0$ (and $\overline b>0$). Let $\Phi_1:Y_1\rightarrow Y$ be the blow up of the 2-curve with
local equations $u=w=0$ at $q$. Let $R^1$ be the transform of $R$ on $Y_1$.
 There are two 3-points $q_1,q_2\in \Phi_1^{-1}(q)$.
$q_1$ has regular parameters $u_1,v_1,w_1$ defined by
$u=u_1, v=v_1, w=u_1w_1$. The strict transform of $F=0$ has the local equation
$$
F^1=v_1^{\overline b}-\frac{1}{\lambda}u_1^{\overline c-\overline a}w_1^{\overline c}=0
$$
at $q_1$.
We  have a form (\ref{eq313}) for $F^1=F_{R^1}(q_1)$ with a reduction in $\overline c$ to $\overline b$.
$q_2$ has regular parameters $u_1,v_1,w_1$ defined by
$u=u_1w_1, v=v_1, w=w_1$. The strict transform of $F=0$ has the local equation
$$
F^1=w_1^{\overline c-\overline a}-\lambda u_1^{\overline a}v_1^{\overline b}=0
$$
at $q_2$. Thus $q_2\in U(R^1)$ and we have a form (\ref{eq313}) for $F^1=F_{R^1}(q_2)$ with a reduction in $\overline c$
to $\overline c-\overline a$.

Suppose that $\overline a\ge \overline c$ in (\ref{eq313}). Let $\Phi_1:Y_1\rightarrow Y$ be the blow up of the 2-curve with
local equation $u=w=0$ at $q$. Let $R^1$ be the transform of $R$ on $Y_1$.
 There are two 3-points $q_1,q_2\in \Phi_1^{-1}(q)$.
$q_1$ has regular parameters $u_1,v_1,w_1$ defined by
$u=u_1,v=v_1, w=u_1w_1$. The strict transform of $F=0$ has the local equation
$$
F^1=w_1^{\overline c}-\lambda u_1^{\overline a-\overline c}v_1^{\overline b}=0
$$
at $q_1$.
Suppose that  $q_1\in U(R^1)$ (which holds if $\overline b+\overline a-\overline c>0$). If
$\overline a=\overline c$ and $\overline b<\overline c$ we have a reduction in $\overline c$ in the expression of the form (\ref{eq313}).
If $\overline b=0$ and $\overline a-\overline c<\overline c$ we have a reduction in $\overline c$.
Otherwise,
$\overline c$ stays the same, but we have a reduction in $\overline a+\overline b$
in $F^1=F_{R^1}(q_1)$.
$q_2$ has regular parameters $u_1,v_1,w_1$ defined by
$u=u_1w_1,v=v_1, w=w_1$. The strict transform of $F=0$ has the local equation
$$
F^1=u_1^{\overline a}v_1^{\overline b}w_1^{\overline a-\overline c}-\frac{1}{\lambda}=0
$$
at $q_2$. Thus $q_2\not\in U(R^1)$.

We have a similar analysis to the above paragraph if $\overline b\ge \overline c$.  In this case we blow up
the 2-curve with local equations $v=w=0$ at $q$. There is at most a single point
$q_1\in \Phi_1^{-1}(q)\cap U(R^1)$. We have a reduction in a local equation $F^1=F_{\overline R^1}(q_1)=0$ of the strict transform of $F=0$ at $q_1$ of the form (\ref{eq313}) of $\overline c$, or  $\overline c$ stays the same but we have a reduction in
$\overline a+\overline b$.

The final case which we must consider is when $\overline a+\overline b\ge \overline c$
 and $\overline a,\overline b<\overline c$.  Let
$\Phi_1:Y_1\rightarrow Y$ be the blow up of the 2-curve with
local equation $u=w=0$ at $q$. There are two 3-points $q_1,q_2\in \Phi_1^{-1}(q)$.
$q_1$ has regular parameters $u_1,v_1,w_1$ defined by
$u=u_1,v=v_1, w=u_1w_1$. The strict transform of $F=0$ has the local equation
$$
F^1=v_1^{\overline b}-\frac{1}{\lambda}w_1^{\overline c}u_1^{\overline c-\overline a}=0
$$
at $q_1$, so that $q_1\in U(R^1)$ and we have a reduction in $\overline c$ in the expression of $F^1=F_{\overline R^1}(q_1)$ 
of the form of (\ref{eq313}).
$q_2$ has regular parameters $u_1,v_1,w_1$ defined by
$u=u_1w_1,v=v_1, w=w_1$. The strict transform of $F=0$ has the local equation
$$
F^1=w_1^{\overline c-\overline a}-\lambda u_1^{\overline a}v_1^{\overline b}=0
$$
at $q_2$. Thus $q_2\in U(R^1)$ and we have a drop in $\overline c$ in $F^1=F_{\overline R^1}(q_2)$ in (\ref{eq313}).

By descending induction on the above invariants, always performing one of the above 
blow ups if $\nu$ is not resolved, we can construct the desired
morphism $Y_{\nu}\rightarrow Y$ such that $\nu$ is resolved on $Y_{\nu}$.

Now by  compactness of the Zariski-Riemann manifold \cite{Z1}, there exist finitely many
$$
Y_1,\ldots, Y_n\in \{Y_{\nu}\mid \nu\text{ is a 0-dimensional valuation of ${\bold k}(Y)$}\}
$$
 such that if $R^1,\ldots, R^n$ are the respective transforms
of $R$ on $Y_1,\ldots,Y_n$ and $\nu$ is a 0-dimensional valuation of ${\bold k}(Y)$, then the center of $\nu$ on some $Y_i$
is not in $U(R^i)$.
  By Lemma \ref{Lemma353}, there exists a sequence of blow ups of 2-curves
$\overline \Phi:\overline Y\rightarrow Y$ such that  there exist morphisms $\overline \Phi_i:\overline Y\rightarrow Y_i$ for all $i$ which factor $\overline \Phi$. In particular,
the transform of $R$ on $\overline Y$
 is resolved.
\end{pf}

\begin{Theorem}\label{Theorem137} Suppose that $f:X\rightarrow Y$ is prepared
and that $R$ is a 3-point relation for $f$. Then
there exists a commutative diagram, 
$$
\begin{array}{rll}
X_1&\stackrel{f_1}{\rightarrow}&Y_1\\
\Phi\downarrow&&\downarrow \Psi\\
X&\stackrel{f}{\rightarrow}&Y
\end{array}
$$
such that $\Phi$ and $\Psi$ are products of blow ups of 2-curves, $f_1$ is prepared and the transform $R^1$ of $R$ for $f_1$ is defined and is resolved.
In particular, all 3-points of $X_1$ in $\Phi^{-1}(T(R))$ must map to 2-points of $Y_1$. Furthermore,
$\tau_{f_1}(p_1)= \tau_f(\Phi(p_1))$ if $p_1\in X_1$ is a 3-point.
If $D_X$ is cuspidal for $f$ then $D_{X_1}$ is cuspidal for $f_1$.
\end{Theorem}

\begin{pf} Let $\{\overline R_i\}$ be the pre-relations on $Y$ associated to $R$.
By Theorem \ref{Theorem3} there exists a sequence of blow ups of 2-curves
$\Psi:Y_1\rightarrow Y$ such that the pre-relations $\{\overline R_{i}^1\}$
which are the transforms of $\{\overline R_i\}$ on $Y_1$ are  resolved on $Y_1$.
By Lemma \ref{Lemma1} there exists a sequence of blow ups of 2-curves $\Phi:X_1\rightarrow X$
such that $f_1=\Psi^{-1}\circ f\circ\Phi:X_1\rightarrow Y_1$ is a prepared morphism, and if $D_X$ is cuspidal for $f$, then $D_{X_1}$ is cuspidal for $f_1$. Furthermore, $\tau_{f_1}(p_1)=\tau_f(\Phi(p_1))$ if $p_1\in X_1$ is a 3-point.

Suppose that $p_1\in X_1$ is a 3-point. Then $p=\Phi(p_1)$ is a 3-point. Suppose that
$u,v,w$ are permissible parameters  at $q=f(p)$ and $x,y,z$ are permissible parameters for $u,v,w$
at $p$. Then there exist permissible parameters $x_1,y_1,z_1$  for 
$u,v,w$ at $p_1$ such that 
\begin{equation}\label{eq162}
\begin{array}{ll}
x&=x_1^{b_{11}}y_1^{b_{12}}z_1^{b_{13}}\\
y&=x_1^{b_{21}}y_1^{b_{22}}z_1^{b_{23}}\\
z&=x_1^{b_{31}}y_1^{b_{32}}z_1^{b_{33}}
\end{array}
\end{equation}
with $\text{Det}(b_{ij})=\pm1$.
We have (after possibly exchanging $u,v,w$) an expression of the form (\ref{eq16}) at $p$.
 On substitution of (\ref{eq162}) into (\ref{eq16}) we see that
we have an expression
\begin{equation}\label{eq163}
\begin{array}{ll}
u&=x_1^{\overline a}y_1^{\overline b}z_1^{\overline c}\\
v&=x_1^{\overline d}y_1^{\overline e}z_1^{\overline f}\\
w&=\sum \alpha_{i}M_i+N
\end{array}
\end{equation}
at $p_1$ of the form of (\ref{eq16}).

We have a factorization
$$
Y_1=Y_m'\stackrel{\Psi_m'}{\rightarrow}Y_{m-1}'\rightarrow\cdots\rightarrow Y_1'\stackrel{\Psi_1'}{\rightarrow}Y_0'=Y
$$
where each $\Psi_i'$ is the blow up of a 2-curve $C_i$.

 Let $\overline f_i:
X_1\rightarrow Y_i'$ be the resulting maps.

Let $R_i$ be the primitive relations associated to $R$, and let $(\overline R')_i^j$
be the transform of the pre-relation $\overline R_i$ on $Y_j'$. By induction on $j$, we will
show that the transform $(R')^j$ of $R$ for $\overline f_j$ is defined. It suffices to
verify this for $(R')^1$.
Suppose that (with the notation of Definition \ref{Def160}), $p_1\in  \Phi^{-1}(T(R_i))$ is a 3-point. Let $p=\Phi(p_1)$,
$q=f(p)$, $\overline q=\overline f_1(p_1)$. $p=\Phi(p_1)\in T(R_{i})$, so
that $q\in U(\overline R_i)$ is a 3-point. Then
$\tau_f(p)\ge 0$, and
\begin{equation}\label{eq165}
w=x_1^{\overline g}y_1^{\overline h}z_1^{\overline i}\Lambda
\end{equation}
where $\Lambda$ is a unit series in (\ref{eq163}). 

Let $F=F_{\overline R_i}(q)=0$ be the expression of (\ref{eq136})
which determines $\overline R_i(q)$. Then we have (with the notation of (\ref{eq156}),
(\ref{eq163}),  (\ref{eq165}) and Definition \ref{Def160}) that
$$
a(\overline a,\overline b,\overline c)+b(\overline d,\overline e,\overline f)
+c(\overline g,\overline h,\overline i)=(0,0,0)
$$
and $\Lambda(0,0,0)^{c}=\lambda=\lambda_p$. 
After possibly interchanging $u$ and  $v$, permissible parameters at $\overline q$ are $u_1,v_1,w_1$
with 
\begin{equation}\label{eq164}
u=u_1,v=u_1(v_1+\alpha),w=w_1
\end{equation}
for some $\alpha\in \bold k$, or 
\begin{equation}\label{eq408}
u=u_1, v=v_1, w=u_1(w_1+\alpha)
\end{equation}
for some $\alpha\in \bold k$ or 
\begin{equation}\label{eq409}
u=u_1w_1, v=v_1, w=w_1.
\end{equation}

Substituting the expressions (\ref{eq164}), (\ref{eq408}) or (\ref{eq409}) 
 into the expression $F=0$ of (\ref{eq136}) and computing the
strict transform of $F=0$ at $\overline q$, we see that if $\overline q$ is a 3-point,
then $\overline\alpha=0$ in (\ref{eq164}), $\overline q$ is on the strict transform of $F=0$,
$\overline q\in U((\overline R')_i^1)$, and the transform $(R')^1$ of $R$ for
$\overline f_1$ is defined. By induction on the number of 2-curves blown up by $\Psi$,
we see that the transform $R^1$ of $R$ for $f_1$ is defined. Since $U(R^1)=\emptyset$, $R^1$ is resolved.

\end{pf}

\section{well prepared morphisms}

Suppose that $f:X\rightarrow Y$ is a birational, prepared morphism of nonsingular projective 3-folds with
toroidal structures $D_Y$ and $D_X=f^{-1}(D_Y)$. Further suppose that the fundamental locus of $f$ is contained in $D_Y$. If $R$ is a 2-point relation for $f$
with associated 2-point pre-relations $\{\overline R_i\}$,
 then for $p\in T(R_i)$ we have
that 
\begin{equation}\label{eq311}
R(p)=\left(\begin{array}{l} 
S_i=S_R(p),E_1=E_{R,1}(p),E_2=E_{R,2}(p),w_i=w_R(p),
u=u_R(p),\\
v=v_R(p),e_i=e_R(p),a_i=a_R(p),b_i=b_R(p),
\overline\lambda_i=\lambda_R(p)
\end{array}\right)
\end{equation}
which we will abbreviate (as in (\ref{eq151})) as 
\begin{equation}\label{eq168}
R(p)=w_i^{e_i}-\overline \lambda_iu^{a_i}v^{b_i},
\end{equation}
with $e_i>1$, if $a_i,b_i\ne-\infty$, or (as in (\ref{eq359}) 
\begin{equation}\label{eq254}
R(p)=w_i
\end{equation}
if $a_i,b_i=-\infty$. In this case, $u,v,w_i$ have a monomial form at $p$.
Recall that if $p'\in T(R)$ is such that $f(p')=f(p)$, then $u_R(p')=u_R(p)=u$ and $v_R(p')=v_R(p)=v$.
Let $I$ be an index set for the $\{\overline R_i\}$ associated to $R$.

\begin{Definition}\label{Def128}
Suppose that $\tau\ge 0$. A prepared morphism $f:X\rightarrow Y$ is $\tau$-quasi-well prepared with 2-point relation $R$ if: 
\begin{enumerate}
\item[1.] $p\in X$ a 3-point implies $\tau_f(p)\le\tau$.
\item[2.] $T(R)$ is the set of 3-points $p$ on $X$ such that $\tau_f(p)=\tau$.
\item[3.] Suppose that $p\in T(R)$. Then $\tau>0$ implies $R(p)$ has a form (\ref{eq168}), $\tau=0$ implies $R(p)$ has a form (\ref{eq254}).
\item[4.] If $q\in U(\overline R_i)\cap U(\overline R_j)$, then there exists  $\lambda_{ij}(u,v)\in {\bold k}[[u,v]]$, 
with 
$$
u=u_{\overline R_i}(q)=u_{\overline R_j}(q), v=v_{\overline R_i}(q)=v_{\overline R_j}(q),
$$
 such that 
$$
w_j=w_{i}+\lambda_{ij}(u,v).
$$
where $w_i=w_{\overline R_i}(q)$, $w_j=w_{\overline R_j}(q)$.
\item[5.] Suppose that $q\in U(\overline R_i)$,
where $\overline R_i$ is a 2-point relation associated to $R$. Then
 $u=u_{\overline R_i}(q),v=v_{\overline R_i}(q),w_i=w_{\overline R_i}(q)$ are super parameters at $q$ (Definition \ref{Def357}).
\end{enumerate}
\end{Definition}

\begin{Definition}
$f:X\rightarrow Y$ is $\tau$-quasi-well prepared with
2-point relation $R$ and pre-algebraic structure  if $f$ is
 $\tau$-quasi-well prepared with 2-point relation $R$ and
 $w_{\overline R_i}(q)\in{\cal O}_{Y,q}$ for all $\overline R_i$ associated to $R$, and $q\in U(\overline R_i)$.
\end{Definition}

\begin{Definition}\label{Def65} $f:X\rightarrow Y$ is $\tau$-well prepared with
2-point relation $R$   if
\begin{enumerate}
\item[1.] $f$ is $\tau$-quasi-well prepared with  2-point relation
$R$ and pre-algebraic structure.
\item[2.] The primitive pre-relations $\{\overline R_i\}$ associated to $R$ are algebraic, and $R$ is algebraic (Definition \ref{Def156}).
\item[3.] Suppose that  $q\in U(\overline R_i)\cap U(\overline R_j)$.
 Let $w_i=w_{\overline R_i}(q)$ and
 $w_j=w_{\overline R_j}(q)$,  $u=u_{\overline R_i}(q)=u_{\overline R_j}(q), v=v_{\overline R_i}(q)=v_{\overline R_j}(q)$. Then there exists a  unit series
$\phi_{ij}\in {\bold k}[[u,v]]$ and $a_{ij}, b_{ij}\in{\bf N}$ (or $\phi_{ij}=0$ with $a_{ij}=b_{ij}=-\infty$)  with 
\begin{equation}\label{eq64}
w_j=w_i+u^{a_{ij}}v^{b_{ij}}\phi_{ij}.
\end{equation}
\item[4.] For $q\in U(R)$, set  $I_q=\{i\in I\mid q\in U(\overline R_i)\}$. Then the set  
\begin{equation}\label{eq255}
\left\{(a_{ij},b_{ij})\mid i,j\in I_q\right\}
\end{equation}
 from equation (\ref{eq64}) is totally ordered.
\end{enumerate}
\end{Definition}

\begin{Definition}\label{Def66}
Suppose that $f:X\rightarrow Y$ is $\tau$-quasi-well prepared with  2-point relation $R$.
\begin{enumerate}
\item[1.] A 2-point $q\in U(R)$ is prepared for $R$.
\item[2.] A 2-point $q\in Y$ such that $q\not\in U(R)$ is prepared for $R$ if there exist
super parameters $u,v,w$ at $q$ (where $w$ could be formal) 
\item[3.] A 2-curve $C\subset Y$ is prepared for $R$.
\end{enumerate}
\end{Definition}

If $E$ is a component of $D_Y$, $\overline R_i$ is pre-algebraic, and $q\in U(\overline R_i)$, we will denote
$\overline{E\cdot S_{\overline R_i}(q)}$ for the Zariski closure in $Y$ of the curve germ $u=w=0$ at $q$, where $u=0$ is
a local equation of $E$, $w=0$ is an (algebraic) local equation of $S_{\overline R_i}(q)$ at $q$.

\begin{Definition}\label{Def200}
Suppose that $f:X\rightarrow Y$ is $\tau$-well prepared with 2-point relation $R$ for $f$.
 A nonsingular curve $C\subset D_Y$ which makes SNCs with $D_Y$ is prepared for $R$
of type 4 if
\begin{enumerate}
\item[1.] $C=\overline{E_{\alpha}\cdot S_{\overline R_i}(q_{\beta})}$ for some
component $E_{\alpha}$ of $D_Y$,  pre-relation $\overline R_i$ associated to
$R$ and $q_{\beta}\in U(\overline R_i)$.
\item[2.] $\Omega(\overline R_i)$ contains $C$.
\item[3.] If $C'=\overline{E_{\gamma}\cdot S_{\overline R_j}(q_{\delta})}$ is such that
$C'\subset\Omega(\overline R_j)$, $C\ne C'$,
and $q\in C\cap C'$, then $q\in U(\overline R_i)\cap U(\overline R_j)$ and
$C'=\overline{E_{\gamma}\cdot S_{\overline R_j}(q)}$.
\item[4.] If $j\ne i$ and $C=\overline{E_{\gamma}\cdot S_{\overline R_j}(q_{\delta})}$
then $C$ satisfies 1 and 2 and 3 of this definition (for $\overline R_j$). (In this case we have by
(\ref{eq253}) that
$U(\overline R_j)\cap C=U(\overline R_i)\cap C$).
\item[5.]
Let
$$
I_C=\{j\in I\mid C=\overline{E_{\gamma}\cdot S_{\overline R_j}(q_{\delta})}
\text{ for some }\overline R_j, E_{\gamma}, q_{\delta}\in U(\overline R_j)\}.
$$
 Suppose that $q\in C$ is a 1-point or a 2-point such that $q\not \in U(R)$.
Then there exist $u,v\in{\cal O}_{Y,q}$ such that for $j\in I_C$ there exists  $\tilde w_j\in
{\cal O}_{Y,q}$ such that
\begin{enumerate}
\item $\tilde w_j=0$ is a local equation of $\Omega(\overline R_j)$ and
 $u,v,\tilde w_j$ are permissible parameters at $q$ such that $u=\tilde w_j=0$ are local
equations of $C$ at $q$.
\item If $q$ is a 1-point and $p\in f^{-1}(q)$, then there exists a relation of one of the following forms
 for  $u,v,\tilde w_j$ at $p$.

\begin{enumerate}
\item
$p$ a 1-point 
\begin{equation}\label{eq59}
\begin{array}{ll}
u&=x^a\\
v&=y\\
\tilde w_j&=x^c\gamma(x,y)+x^dz
\end{array}
\end{equation}
where $\gamma$ is a unit series  (or zero),
\item
$p$ a 2-point 
\begin{equation}\label{eq61}
\begin{array}{ll}
u&=(x^ay^b)^k\\
v&=z\\
\tilde w_j&=(x^ay^b)^l\gamma(x^ay^b,z)+x^cy^d
\end{array}
\end{equation}
where $\gamma$ is a unit series (or zero) and $ad-bc\ne 0$.
\end{enumerate}
\item If $q$ is a 2-point, then  $u,v,\tilde w_j$ are super parameters at $q$.
\item If $i,j\in I_C$ and $q$ is a 1-point, there exist relations
$$
\tilde w_i-\tilde w_j=u^{c_{ij}}\phi_{ij}(u,v)
$$
where $\phi_{ij}$ is a unit series (or $\phi_{ij}=0$ and $c_{ij}=\infty$).
\item If $i,j\in I_C$ and $q$ is a 2-point (with $q\not\in U(R)$) then
there exist relations
$$
\tilde w_i-\tilde w_j=u^{c_{ij}}v^{d_{ij}}\phi_{ij}(u,v)
$$
where $\phi_{ij}$ is a unit series (or $\phi_{ij}=0$ and $c_{ij}=d_{ij}=\infty$), and the set $\{(c_{ij},d_{ij})\}$ is totally ordered.

\end{enumerate}
\end{enumerate}
\end{Definition}

If $f:X\rightarrow Y$ is $\tau$-well prepared with 2-point relation $R$, and $\overline R_i$ is a pre-relation
associated to $R$, we will feel free to replace $\Omega(\overline R_i)$ with an open subset of $\Omega(\overline R_i)$ containing $U(\overline R_i)$, and all curves $C=\overline{E\cdot S_{\overline R_i}(q)}$ such that $E$ is a component of
$D_Y$, $q\in U(\overline R_i)$ and $C$ is prepared for $R$ of type 4. This convention will allow some simplification of the statements
of the theorems and proofs.

\begin{Definition}\label{Def130}  $f:X\rightarrow Y$ is
$\tau$-very-well prepared  with 2-point relation $R$ if
\begin{enumerate}
\item[1.] $f$ is $\tau$-well prepared with 2-point relation $R$.
\item[2.] If $E$ is a component of $D_Y$ and $q\in U(\overline R_i)\cap E$,
 then $C=\overline{E\cdot S_{\overline R_i}(q)}$
is prepared for $R$
of type 4 (Definition \ref{Def200}).
\item[3.] For all $\overline R_i$ associated to $R$, let
$$
V_i(Y)=\left\{\gamma=\overline{E_{\alpha}\cdot S_{\overline R_{i}}(q_{\gamma})}\mid
q_{\gamma}\in U(\overline R_i), E_{\alpha}\text{ is a component of }D_Y\right\}.
$$
Then
$$
F_i=\sum_{\gamma\in V_i(Y)}\gamma
$$
is a SNC divisor on $\Omega(\overline R_i)$ whose intersection graph
is a tree.
\end{enumerate}
\end{Definition}

If $f:X\rightarrow Y$ is $\tau$-very-well prepared, we will feel free to replace $\Omega(\overline R_i)$ with an open neighborhood of $F_i$ in $\Omega(\overline R_i)$.
This will allow some simplification of the proofs.

\begin{Remark}\label{Remark281} Suppose that $f:X\rightarrow Y$ is $\tau$-very well
prepared. Then it follows from Definition \ref{Def130} and (\ref{eq253}) that $F_i\cap U(R)=U(\overline R_i)$
for all $\overline R_i$ associated to $R$.
\end{Remark}

\begin{Definition}\label{Def289}
Suppose that $f:X\rightarrow Y$ is $\tau$-quasi-well prepared (or $\tau$-well prepared
or $\tau$-very-well prepared)
with  2-point relation $R$.
Let $\{\overline R_i\}$ be the pre-relations  associated to $R$.
Suppose that $G$ is a
 point or nonsingular curve in $Y$ which is an admissible center for all of the $\overline R_i$. Then $G$ is called
a permissible center for $R$ if there exists a commutative diagram 
 \begin{equation}\label{eq30}
\begin{array}{rll}
X_1&\stackrel{f_1}{\rightarrow}&Y_1\\
\Phi\downarrow&&\downarrow\Psi\\
X&\stackrel{f}{\rightarrow}&Y
\end{array}
\end{equation}
where  $\Psi$ is the blow up of $G$ and $\Phi$ is a sequence of blow ups
$$
X_1=\overline X_n\rightarrow \cdots\rightarrow \overline X_1\rightarrow X
$$
of
nonsingular curves and 3-points $\gamma_i$
 which are possible centers such that
\begin{enumerate}
\item[1.] $f_1$ is prepared and the assumptions of Definition \ref{Def161} hold so that the transform $R^1$ of $R$
for $f_1$ is defined.
\item[2.] $f_1:X_1\rightarrow Y_1$ is $\tau$-quasi-well prepared, (or $\tau$-well prepared or
$\tau$-very-well prepared)
with  2-point relation $R^1$.
\end{enumerate}
\end{Definition}

(\ref{eq30}) is called a $\tau$-quasi-well prepared (or $\tau$-well prepared or $\tau$-very-well prepared) diagram of $R$ (and $\Psi$).

\begin{Definition}\label{Def219}
Suppose that  $f:X\rightarrow Y$ is  $\tau$-well prepared (or $\tau$-very-well prepared) with 2-point relation
$R$ and $C\subset Y$ is prepared for $R$ of type 4. Then $C$ is a $*$-permissible center
for $R$ if there exists a commutative diagram 
\begin{equation}\label{eq233}
\begin{array}{rll}
X_1&\stackrel{f_1}{\rightarrow}&Y_1\\
\Phi\downarrow&&\downarrow\Psi\\
X&\stackrel{f}{\rightarrow}&Y
\end{array}
\end{equation}
such that
\begin{enumerate}
\item[1.] $f_1$ is prepared and the assumptions of Definition \ref{Def161} hold so that
the transform $R^1$ of $R$ for $f_1$ is defined.
\item[2.] $f_1:X_1\rightarrow Y_1$ is $\tau$-well prepared (or $\tau$-very-well prepared).
\item[3.] (\ref{eq233}) has a factorization 
\begin{equation}\label{eq220}
\begin{array}{rcl}
X_1=\overline X_m&\stackrel{f_1=\overline f_m}{\rightarrow}&\overline Y_m=Y_1\\
\downarrow&&\downarrow\\
\vdots&&\vdots\\
\downarrow&&\downarrow\\
\overline X_2&\stackrel{\overline f_2}{\rightarrow}&\overline Y_2\\
\overline\Phi_2\downarrow&&\downarrow\overline\Psi_2\\
\overline X_1&\stackrel{\overline f_1}{\rightarrow}&\overline Y_1\\
\overline \Phi_1\downarrow&&\downarrow\overline\Psi_1\\
X&\stackrel{f}{\rightarrow}&Y
\end{array}
\end{equation}
where $\overline\Psi_1$ is the blow up of $C$, 
\begin{equation}\label{eq310}
\begin{array}{rll}
\overline X_1&\stackrel{\overline f_1}{\rightarrow}&\overline Y_1\\
\overline\Phi_1\downarrow&&\downarrow\overline\Psi_1\\
X&\stackrel{f}{\rightarrow}&Y
\end{array}
\end{equation}
is a $\tau$-well prepared diagram of $R$ and $\overline\Psi_1$ of the form (\ref{eq30}), each
$\overline \Psi_{i+1}:\overline Y_{i+1}\rightarrow \overline Y_i$ for $i\ge 1$ is
the blow up of a 2-point $q\in\overline Y_i$ which is prepared for the transform $R^i$
of $R$ on $\overline X_i$ of type 2 of Definition \ref{Def66}, and
$$
\begin{array}{rll}
\overline X_{i+1}&\stackrel{\overline f_{i+1}}{\rightarrow} &\overline Y_{i+1}\\
\overline\Phi_{i+1}\downarrow&&\downarrow\overline\Psi_{i+1}\\
\overline X_{i}&\stackrel{\overline f_i}{\rightarrow}&\overline Y_{i}
\end{array}
$$
is a $\tau$-well prepared diagram of $R^i$ and $\overline\Psi_{i+1}$ of the form of (\ref{eq30}).
\end{enumerate}
\end{Definition}

\begin{Definition}\label{Def396}
Suppose that  $f:X\rightarrow Y$ is  $\tau$-quasi-well prepared (or $\tau$-well prepared or $\tau$-very-well prepared) with 2-point relation $R$.
Suppose that 
\begin{equation}\label{eq405}
\begin{array}{rll}
X_1&\stackrel{f_1}{\rightarrow}&Y_1\\
\Phi\downarrow&&\downarrow\Psi\\
X&\stackrel{f}{\rightarrow}&Y
\end{array}
\end{equation}
is a commutative diagram 
such that there is a factorization
\begin{equation}\label{eq412}
\begin{array}{rcl}
X_1=\overline X_m&\stackrel{f_1=\overline f_m}{\rightarrow}&\overline Y_m=Y_1\\
\downarrow&&\downarrow\\
\vdots&&\vdots\\
\downarrow&&\downarrow\\
\overline X_2&\stackrel{\overline f_2}{\rightarrow}&\overline Y_2\\
\overline\Phi_2\downarrow&&\downarrow\overline\Psi_2\\
\overline X_1&\stackrel{\overline f_1}{\rightarrow}&\overline Y_1\\
\overline \Phi_1\downarrow&&\downarrow\overline\Psi_1\\
X&\stackrel{f}{\rightarrow}&Y
\end{array}
\end{equation}

where each commutative diagram
$$
\begin{array}{rll}
\overline X_{i+1}&\rightarrow &\overline Y_{i+1}\\
\overline\Phi_{i+1}\downarrow&&\downarrow\overline\Psi_{i+1}\\
\overline X_{i}&\rightarrow&\overline Y_{i}
\end{array}
$$
is either of the form (\ref{eq30}) or of the form (\ref{eq233}). Then (\ref{eq405}) is called a $\tau$-quasi-well prepared
(or $\tau$-well prepared or $\tau$-very-well prepared ) diagram of $R$ (and $\Psi$).
\end{Definition}

\begin{Lemma}\label{Lemma31}
Suppose that $f:X\rightarrow Y$ is $\tau$-quasi-well prepared
(or $\tau$-well prepared or $\tau$-very-well prepared) and $C\subset Y$ is a  2-curve. Then $C$ is a permissible center for $R$, and there exists  a $\tau$-quasi-well-prepared
(or $\tau$-well prepared or $\tau$-very-well prepared) diagram
(\ref{eq30}) of $R$ and the blow up $\Psi:Y_1\rightarrow Y$ of $C$
such that $\Phi$ is a product of blow ups of 2-curves. Furthermore,
\begin{enumerate}
\item[1.] If $D_X$ is cuspidal for $f$ then $D_{X_1}$ is cuspidal for $f_1$.
\item[2.] Further suppose that $f$ is $\tau$-well prepared. Then
\begin{enumerate}
\item  Let $E$ be the exceptional divisor for $\Psi$. Suppose that
$q\in U(\overline R_i^1)\cap E$ for some $\overline R_i$ associated to $R$. Let $\gamma_i=\overline{S_{\overline R_i^1}(q)\cdot E}$.
Then $\gamma_i=\Psi^{-1}(\Psi(q))$ is a prepared curve for $R^1$ of type 4.
\item If $\gamma$ is a prepared curve for $R$, then the strict transform of $\gamma$ on $Y_1$ is
a prepared curve for $R^1$.
\end{enumerate}
\item[3.] $\Phi$ is an isomorphism over $f^{-1}(Y-\Sigma(Y))$
\end{enumerate}
\end{Lemma}

\begin{pf} By Lemma \ref{Lemma1}, there exists a commutative diagram
$$
\begin{array}{rll}
X_1&\stackrel{f_1}{\rightarrow}&Y_1\\
\Phi\downarrow&&\downarrow\Psi\\
X&\stackrel{f}{\rightarrow}&Y
\end{array}
$$
where $\Phi$ is a product of blow ups of 2-curves and $f_1$ is prepared, with the
property that $\tau_{f_1}(p_1)=\tau_f(\Phi(p_1))$ if $p_1\in X_1$ is a 3-point.
We further have that $D_{X_1}$ is cuspidal for $f_1$ if $D_X$ is cuspidal for $f$, and $\Phi$ is an isomorphism over $f^{-1}(Y-\Sigma(Y))$.

Let $\{\overline R_i\}$ be the 2-point pre-relations on $Y$ associated to $R$. $C$ is an
admissible center for the $\{\overline R_i\}$ (Definition \ref{Def154}). Let $\{\overline R_i^1\}$ be the transforms of the $\{\overline R_i\}$ on $Y_1$.

We will show that the conditions of Definition \ref{Def161} hold so that we can define the
transform $R^1$ of $R$ for $f_1$. Suppose that $q_1\in U(\overline R_i^1)$, and
$p_1\in f_1^{-1}(q_1)\cap \Phi^{-1}(T( R_i))$ is a 3-point. Let $q=\Psi(q_1)$,
$p=\Phi(p_1)$. There exist permissible parameters $u=u_R(p),v=v_R(p),w=w_R(p)$ at $q$ such that
$R(p)=\overline R_i(q)$ is determined by
$$
w^e-\lambda u^av^b
$$
if $a=a_{ R}(p),b=b_{ R}(p)\ne-\infty$, and by
$$
w=0
$$
if $a=a_{ R}(p)= b=b_{ R}(p)=-\infty$.
There exist permissible parameters $x,y,z$ for $u,v,w$ at $p$ such that an expression (\ref{eq16}) of Definition \ref{Def221}
holds for $u,v,w$ and we have a relation
$$
w^e=u^av^b\Lambda(x,y,z)
$$
where $\Lambda(x,y,z)$ is a unit series with $\Lambda(0,0,0)=\lambda$ if $a,b\ne-\infty$ and $u,v,w$ have a monomial
form in $x,y,z$ if $a=b=-\infty$.  After possibly interchanging $u$ and $v$, we may assume (since $q_1$ is a 2-point)
that $q_1$ has permissible parameters $\overline u,\overline v,w$ such that
$$
u=\overline u,
v=\overline u\overline v.
$$
Since $p_1$ is a 3-point, $\hat{\cal O}_{X_1,p_1}$ has regular parameters $x_1,y_1,z_1$ such that
$$
\begin{array}{ll}
x&=x_1^{a_{11}}y_1^{a_{12}}z_1^{a_{13}}\\
y&=x_1^{a_{21}}y_1^{a_{22}}z_1^{a_{23}}\\
z&=x_1^{a_{31}}y_1^{a_{32}}z_1^{a_{33}}
\end{array}
$$
with $\text{Det}(a_{ij})=\pm 1$. Thus $\overline u,\overline v,w$ has an expression of the form of (\ref{eq16})
in $x_1,y_1,z_1$. If $a,b\ne-\infty$  we have the relation
$$
w^e=\overline u^{a+b}\overline v^b\Lambda,
$$
and $\overline R_i^1(q_1)$ is determined by 
\begin{equation}\label{eq236}
w^e=\overline u^{a+b}\overline v^b\lambda.
\end{equation}
If $a=b=-\infty$, $\overline u,\overline v,w$ have a monomial expression in $x_1,y_1,z_1$.
The transform $R^1$ of $R$ for $f_1$ is thus defined.

Now we will verify that $f_1$ is $\tau$-quasi-well prepared. From the above calculations,
we see that $f_1$ satisfies 1, 2, 3  and 4 of Definition \ref{Def128}.
It remains to verify that 5 of Definition \ref{Def128} holds.

Suppose that $q_1\in U(\overline R_i^1)$ and $p_1\in f_1^{-1}(q_1)$.
Then  
$$
u=u_{\overline R_i}(q), v=v_{\overline R_i}(q), w=w_{\overline R_i}(q)
$$
are super parameters at $q=\Psi(q_1)$, and $p=\Phi(p_1)\in f^{-1}(q)$ has permissible parameters
 $x,y,z$ for $u,v,w$ such that one of the forms  of
Definition \ref{Def357} hold for $u,v,w$ and $x,y,z$.

After possibly interchanging $u$ and $v$, we have 
\begin{equation}\label{eq206}
\begin{array}{ll}
u_{\overline R_i^1}(q_1)&=\overline u=u\\
v_{\overline R_i^1}(q_1)&=\overline v=\frac{v}{u}\\
w_{\overline R_i^1}(q_1)&=w.
\end{array}
\end{equation}
We can verify that there exist permissible parameters $x_1,y_1,z_1$ at $p_1$ such that
$\overline u,\overline v,w$ have one of the forms  of Definition \ref{Def357} in $x_1,y_1,z_1$.
The most difficult case to verify is when $p_1$ is a 1-point and $p$ is a 3-point.
Then $\hat{\cal O}_{X_1,p_1}$ has regular parameters $\overline x_1,\overline y_1,\overline z_1$ defined
by
$$
\begin{array}{ll}
x&=\overline x_1^{\overline a}(\overline y_1+\overline\alpha)^{\overline b}(\overline z_1
+\overline\beta)^{\overline c}\\
y&=\overline x_1^{\overline d}(\overline y_1+\overline\alpha)^{\overline e}(\overline z_1+
\overline\beta)^{\overline f}\\
w&=\overline x_1^{\overline g}(\overline y_1+\overline\alpha)^{\overline h}
(\overline z_1+\overline\beta)^{\overline i}
\end{array}
$$
where $\overline\alpha,\overline\beta\in \bold k$ are nonzero and
$$
\text{Det}\left(\begin{array}{lll}
\overline a&\overline b&\overline c\\
\overline d&\overline e&\overline f\\
\overline g&\overline h&\overline i
\end{array}\right)=\pm 1.
$$
We substitute into
$$
\begin{array}{ll}
u&=x^ay^bz^c\\
v&=x^dy^ez^f\\
w&=x^gy^hz^i\gamma+x^jy^kz^l
\end{array}
$$ of 4 of Definition \ref{Def357}.
Using the fact that
$$
\text{rank}\left(\begin{array}{lll}
a&b&c\\
d&e&f\\
j&k&l
\end{array}\right)=3,
$$
we can make a change of variables in $\overline x_1,\overline y_1,\overline z_1$
 to get permissible parameters $x_1,y_1,z_1$
at $p_1$ satisfying 
\begin{equation}\label{eq205}
\begin{array}{ll}
u&=x_1^{\overline a}\\
v&=x_1^{\overline d}(\overline\gamma+y_1)\\
x^jy^kz^l&=x_1^{\overline g}(\overline\epsilon+z_1)
\end{array}
\end{equation}
with $\overline d>\overline a$ and $0\ne\overline \epsilon,\overline\gamma\in \bold k$.

For each of the monomials $M$ in the series $x^gy^hz^i\gamma$ we have a relation 
\begin{equation}\label{eq234}
M^{\tilde e}=u^{\tilde a}v^{\tilde b}
\end{equation}
with $\tilde e,\tilde a,\tilde b\in{\bf Z}$. On substitution of (\ref{eq206})
and (\ref{eq205}) into (\ref{eq234}) we see that
$$
M=x_1^{e'}\phi_1(y_1)
$$
where $\phi_1$ is a unit series. Thus $\overline u,\overline v,w$ have an expansion of the form
1 of Definition \ref{Def357} in terms of $x_1,y_1,z_1$.

The other cases can be verified by a similar but simpler argument to show that 
$\overline u,\overline v, w$ are super parameters at $q_1$. Thus
5 of
Definition \ref{Def128} holds for $f_1$, so that $f_1$ is $\tau$-quasi-well prepared.

Suppose that $f$ is $\tau$-well prepared. We will verify that $f_1$ is $\tau$-well prepared. 1 and 2
 of Definition \ref{Def65} are immediate.
We must verify that 3 and 4 of Definition \ref{Def65} hold for $f_1$.

Suppose that $q_1\in U(\overline R_i^1)\cap U(\overline R_j^1)$. Let  $q=\Phi(q_1)\in
U(\overline R_i)\cap U(\overline R_j)$, and
$$
\begin{array}{l}
u=u_{\overline R_i}(q)=u_{\overline R_j}(q),\\
v=v_{\overline R_i}(q)=v_{\overline R_j}(q),\\
w_i=w_{\overline R_i}(q), w_j=w_{\overline R_j}(q).
\end{array}
$$
We have a relation 
\begin{equation}\label{eq237}
w_j=w_i+u^{a_{ij}}v^{b_{ij}}\phi_{ij}(u,v),
\end{equation}
from (\ref{eq64}) for $f$.
After possibly interchanging $u$ and $v$ we have permissible parameters
$$
\begin{array}{l}
\overline u=u_{\overline R_i^1}(q_1)=u_{\overline R_j^1}(q_1)\\
\overline v=v_{\overline R_i^1}(q_1)=v_{\overline R_j^1}(q_1)\\
w_i=w_{\overline R_i^1}(q_1), w_j=w_{\overline R_j^1}(q_1)
\end{array}
$$
at $q_1$, where $u=\overline u$, $v=\overline u\overline v$. We have 
\begin{equation}\label{eq208}
w_j=w_i+\overline u^{a_{ij}+b_{ij}}\overline v^{b_{ij}}\phi_{ij}(\overline u,\overline u
\overline v)
\end{equation}
so 3 of Definition \ref{Def65} holds for $f_1$.

Since the set  (\ref{eq255})  of Definition \ref{Def65} is totally ordered for $q\in U(R)$, it
follows from  (\ref{eq208}) that the corresponding set  (\ref{eq255})
for $q_1\in U(R^1)$ is totally ordered. Thus 4 of Definition \ref{Def65} holds for $f_1$ and $R^1$
and $f_1$ is $\tau$-well prepared.

We now verify 2 of Lemma \ref{Lemma31}.
Suppose that $f$ is $\tau$-well prepared and (\ref{eq30}) is a $\tau$-well prepared
diagram of $R$ and $\Psi$. Let $E=\Psi^{-1}(C)$ be the exceptional divisor of $\Psi$.
 Suppose that
$q_1\in U(\overline R_i^1)\cap E$. Let  $q=\Psi(q_1)$. There exist permissible parameters
$$
u=u_{\overline R_i}(q), v=v_{\overline R_i}(q), w_i=w_{\overline R_i}(q)
$$
at $q$ such that $u=v=0$ are local equations of $C$, and
after possibly interchanging $u$ and $v$,
$$
\overline u=u_{\overline R_i^1}(q_1)=u,
\overline v=v_{\overline R_i^1}(q_1)=\frac{v}{u},
\overline w_i=w_{\overline R_i^1}(q_1)=w_i
$$
are permissible parameters at $q_1$.
Since $\overline u$ is a local equation of $E$ at $q_1$,
$\gamma_i=\overline{ S_{\overline R_i^1}(q)\cdot E}=\Psi^{-1}(q)$.
Let $\gamma=\gamma_i$.

We will verify that $\gamma_i$ is prepared for $R^1$ of type 4.
Since $q\in\Omega(\overline R_i)$, $C$ intersects $\Omega(\overline R_i)$ transversally at $q$ (and possibly a finite number of other points),
and $\Omega(\overline R_i^1)$ is the strict transform of $\Omega(\overline R_i)$ by $\Psi$, we have that
$\gamma_i\subset\Omega(\overline R_i^1)$. Thus 2 of Definition \ref{Def200} holds.
Suppose that for some $j$, component $E_{\alpha}$ of $D_{Y_1}$ and $q_{\beta}\in U(\overline R_j^1)$,
$$
\gamma'=\overline{E_{\alpha}\cdot S_{\overline R_j^1}(q_{\beta}))}\subset
\Omega(\overline R_j^1),
$$
$\gamma\ne\gamma'$ and there exists $\overline q\in \gamma\cap\gamma'$. Let $E_1=\Psi(E_{\alpha})$,
a component of $D_Y$. Then $\overline\gamma=\Psi(\gamma')$ is a curve on $D_Y$
through $q$.    Since $\gamma'\subset \Omega(\overline R_j^1)$, we must have
$\overline\gamma\subset
\Omega(\overline R_j)\cap E_1$, so that $q\in U(R)\cap\overline\gamma=U(\overline R_j)\cap \overline \gamma$,  and thus $\overline\gamma=\overline{E_1\cdot S_{\overline R_j}
(q)}$. We thus have that $E_{\alpha}$ is the strict transform of $E_1$,  
$\overline q=E_{\alpha}\cdot \gamma\in U(\overline R_i^1)\cap U(\overline R_j^1)$ and $\gamma'=
\overline{E_{\alpha}\cdot S_{\overline R_j^1}(\overline q)}$. Thus 3 of Definition
\ref{Def200} holds.

Suppose that $\gamma'=\overline{E_{\alpha}\cdot S_{\overline R_j^1}(q_{\delta})}$ and
$\gamma=\gamma'$. Then we must have $\gamma'=\overline{E\cdot S_{\overline R_j^1}(q_1)}$
and 4 of Definition \ref{Def200} holds.

Since for $q\in U(\overline R_j)$, $U(\overline R_j^1)$ contains both 2-points
above $q$ in $Y_1$, we need only verify 5 of Definition \ref{Def200} at 1-points $\overline q\in\gamma$.
Let
$$
I_{\gamma}=\{j\mid \gamma=\overline{E\cdot S_{\overline R_j^1}(q_1)}\}.
$$
At $\overline q$ there exist regular parameters $\tilde u,\tilde v,\tilde w_j$
(for all $j\in I_{\gamma}$) such that 
\begin{equation}\label{eq235}
u=\tilde u,
v=\tilde u(\tilde v+\alpha),
w_j=\tilde w_j
\end{equation}
where $u=u_{\overline R_j}(q)=u_{\overline R_i}(q)$,
$v=v_{\overline R_j}(q)=v_{\overline R_i(q)}$, $w_j=w_{\overline R_j}(q)$
and $0\ne\alpha\in \bold k$.

Thus 5 (a) of Definition \ref{Def200} holds for $\tilde u,\tilde v,\tilde w_j$ at $\overline q$.

As in our verification that $f_1$ is $\tau$-quasi-well prepared, we see that if
$p\in f_1^{-1}(\overline q)$, then there exist permissible parameters $x,y,z$ for $u,v,w_j$ at $p$
such that one of the forms of Definition \ref{Def357} hold for
$u,v,w_j$ and $x,y,z$. Substituting in (\ref{eq235}), we see that (since $\alpha\ne 0$)
$u,v,w_j$ must satisfy a form 1 or 3 of Definition \ref{Def357} at $p$, and $\tilde u,\tilde v,\tilde w_j$
must satisfy one of the forms (i) or (ii) of 5 (b) of Definition \ref{Def200}.
Thus 5 (b) of Definition \ref{Def200} holds at $p$. Substituting (\ref{eq235}) into
(\ref{eq237}), we see that 5 (d) of Definition \ref{Def200} holds at $\overline q$.
Thus $\gamma=\gamma_i$ is prepared for $R^1$ of type 4.

We now verify that if $\gamma$ is a prepared curve for $R$ on $Y$ then the strict
transform $\gamma'$ of $\gamma$ on $Y_1$ is a prepared curve  for $R^1$. We may assume that $\gamma$ is prepared of type 4 and $\gamma\cap C\ne\emptyset$.
The verification that $\gamma'$ is prepared for $f_1$ now follows from a local calculation
at points $q\in C\cap\gamma$.

Suppose that $f$ is $\tau$-very-well prepared. We have seen that $f_1$ is $\tau$-well prepared,
so that 1 of Definition \ref{Def130} holds for $f_1$, and 2 of Definition \ref{Def130}
holds for $f_1$ by our verification of  2 of this lemma. Let
$\overline\Psi_i:\Omega(\overline R_i^1)\rightarrow \Omega(\overline R_i)$ be the restriction
of $\Psi$ to $\Omega(\overline R_i^1)$. Then $\overline\Psi_i$ is the blow up of the
union of nonsingular points $C\cdot \Omega(\overline R_i)$ on the nonsingular surface
$\Omega(\overline R_i)$. Thus since
$$
F_i=\sum_{\gamma\in V_i(Y)}\gamma
$$
is a SNC divisor on $\Omega(\overline R_i)$ whose intersection graph is a tree,
$$
\overline\Psi^{-1}(F_i)=\sum_{\gamma'\in V_i(Y_1)}\gamma'
$$
is a SNC divisor on $\Omega(\overline R_i^1)$ whose intersection graph is a tree
and 3 of Definition \ref{Def130} holds for $f_1$. Thus $f_1$ is $\tau$-very-well prepared.
\end{pf}

\begin{Remark}\label{Remark424} The proof of Lemma \ref{Lemma31} shows that if $f:X\rightarrow Y$ is $\tau$-quasi-well prepared (or $\tau$-well prepared or
$\tau$-very-well prepared), $C\subset Y$ is a 2-curve, $\Psi:Y_1\rightarrow Y$ is the blow up of $C$ and $\Phi:X_1\rightarrow X$ is a sequence of blow ups of
2-curves and 3-points such that the rational map $f_1:X_1\rightarrow Y_1$ is a morphism, then
$$
\begin{array}{rll}
X_1&\stackrel{f_1}{\rightarrow}&Y_1\\
\Phi\downarrow&&\downarrow\Psi\\
X&\stackrel{f}{\rightarrow}&Y
\end{array}
$$
is $\tau$-quasi-well prepared (or $\tau$-well prepared or $\tau$-very-well prepared) for $R$ and $\Psi$. If $D_{X}$ is cuspidal for $f$, then $D_{X_1}$
is cuspidal for $f_1$. In fact, with the above notation, if $f$ satisfies 1 -- 4 of Definition \ref{Def128}, then $f_1$ satisfies 1 -- 4 of Definition 
\ref{Def128}.
\end{Remark}

\begin{Lemma}\label{Lemma32}
Suppose that $f:X\rightarrow Y$ is $\tau$-quasi-well prepared
(or $\tau$-well prepared or $\tau$-very-well prepared) and $q\in U(R)$ is a  2-point
(prepared of type 1 in Definition \ref{Def66}).
 Then $q$ is a permissible center for $R$, and  there exists
 a $\tau$-quasi-well prepared (or $\tau$-well prepared or $\tau$-very-well prepared)
diagram (\ref{eq30}) of $R$ and the blow up $\Psi:Y_1\rightarrow Y$ of $q$ such that:
\begin{enumerate}
\item[1.] Suppose that $D_X$ is cuspidal for $f$. Then $D_{X_1}$ is cuspidal for $f_1$.
\item[2.] Suppose that $f$ is $\tau$-well prepared. Then
\begin{enumerate}
\item
Let $E$ be the exceptional divisor of $\Psi$.
Suppose that $q_1\in U(\overline R_i^1)\cap E$. Let $\gamma_i=\overline{S_{\overline R_i^1}(q_1)
\cdot E}$.
 Then $\gamma_i$ is  a prepared
curve for $R^1$ of type 4. Suppose that $q'\in U(\overline R_j^1)\cap E$.
 Let $\gamma_j=\overline{S_{\overline R_j^1}(q')\cdot E}$. Then either
\begin{enumerate}
\item $\gamma_i=\gamma_j$ or
\item $\gamma_i,\gamma_j$ intersect transversally at a 2 point on $E$ (their tangent spaces have distinct directions at this point and are otherwise disjoint). 
\end{enumerate}
\item If $\gamma$ is a prepared curve on $Y$ then the strict transform of $\gamma$ is
a prepared curve on $Y_1$.
\end{enumerate}
\item[3.] $\Phi$ is an isomorphism over $f^{-1}(Y-\Sigma(Y))$.
\end{enumerate}
\end{Lemma}

\begin{pf}
There exists a pre-relation $\overline R_i$ associated to $R$ such that $q\in U(\overline R_i)\subset U(R)$.
Fix such an $i$.
Let
$$
u=u_{\overline R_i}(q),
v=v_{\overline R_i}(q),
w_i=w_{\overline R_i}(q).
$$
$u,v,w_i$ are super parameters at $q$, and $w_i=0$ is a local equation
of $S_{\overline R_i}(q)$. Let $m_q\subset {\cal O}_{Y,q}$ be the maximal ideal.  

By Lemma \ref{Lemma353}, there exists a morphism $\Phi_0:X_0\rightarrow X$ which is a sequence of blow ups of 2-curves such that $(u,v){\cal O}_{X,p}$
is invertible for all $p\in (f\circ\Phi_0)^{-1}(q)$,  $\Phi_0$ is an isomorphism over $f^{-1}(Y-\Sigma(Y))$, $f\circ\Phi_0$ is prepared and
$u,v,w_i$ are super parameters for $f\circ\Phi_0$.

We will next show that there exists a sequence of blow ups of 2-curves and 3-points $\Phi_1:X_1\rightarrow X_0$
such that $(u,v){\cal O}_{X_1,p}$ is invertible at all $p\in (f\circ\Phi_0\circ\Phi_1)^{-1}(q)$ and
if $m_q{\cal O}_{X_1,p}$ is not invertible, then we have
permissible parameters $x,y,z$ for $u,v,w$ at $p$ of
 one of the following  forms:

$p$ is a 1-point 
\begin{equation}\label{eq43}
u=x^a, v=x^b(\alpha+y), w_i=x^dz
\end{equation}
with $\alpha\ne 0$ and $d<\text{min}\{a,b\}$ or,

$p$ is a  2-point of the type of (\ref{eqTF21}) of Definition \ref{torf} 
\begin{equation}\label{eq44}
u=x^ay^b, v=x^cy^d, w_i=x^gy^h(z+\overline\gamma)
\end{equation}
with $\overline\gamma\in \bold k$, $ad-bc\ne 0$, and $(g,h)\le\text{min}\{(a,b),(c,d)\}$.

In fact, we will construct $\Phi_1:X_1\rightarrow X_0$ such that only (\ref{eq43}) or
(\ref{eq44}) occur at points $p_1$ above $q$ such that $m_q{\cal O}_{X_1,p_1}$ is not
invertible.

Suppose that $p\in (f\circ\Phi_0)^{-1}(q)$ is a 3-point, so that $p$ has a form 4 of Definition \ref{Def357} at $p$. 
There exist $\overline x,\overline y,\overline z\in{\cal O}_{X_0,p}$ and series $\lambda_1,\lambda_2,\lambda_3\in \hat{\cal O}_{X_0,p}$ such that
$x=\overline x\lambda_1$, $y=\overline y\lambda_2$, $z=\overline z\lambda_3$.
 Let $I^p\subset{\cal O}_{X_0,p}$ be the ideal
$$
I^p=(u,v,\overline x^g\overline y^h\overline z^i,\overline x^j\overline y^k\overline z^l).
$$
By Lemma \ref{Lemma419}, there exists a sequence of blow ups of 2-curves and 3-points $\Phi_{1}:X_{1}\rightarrow
X_0$ such that $\Phi_0\circ\Phi_1$ is an isomorphism above $f^{-1}(Y-\Sigma(Y))$ and $I^p{\cal O}_{X_{1},p_1}$ is invertible for all 
3-points $p\in (f\circ\Phi_0)^{-1}(q)$ and $p_1\in\Phi_{1}^{-1}(p)$.
Thus if $p\in (f\circ\Phi_0)^{-1}(q)$ is a 3-point and $p_1\in (f\circ\Phi_0\circ\Phi_1)^{-1}(p)$, then $m_q{\cal O}_{X_1,p_1}$ is invertible or $p_1$ has a form
(\ref{eq43}) or (\ref{eq44}). $f\circ\Phi_0\circ\Phi_1$ is prepared, and 
$u,v,w_i$ are super parameters for $f\circ\Phi_{0}\circ\Phi_1$.

We construct an infinite sequence of morphisms 
\begin{equation}\label{eq421}
\cdots\rightarrow X_n\stackrel{\Phi_n}{\rightarrow}\cdots\stackrel{\Phi_3}{\rightarrow} X_2\stackrel{\Phi_2}{\rightarrow}X_1
\end{equation}
as follows. Order the 2-curves $C$ of $X_1$ such that $q\in (f\circ\Phi_0\circ\Phi_1)(C)\subset\Sigma(Y)$. Let $\Phi_2:X_2\rightarrow X_1$ be the blow up of the
2-curve $C_1$ on $X_1$ of smallest order. Order the 2-curves $C'$ of $X_2$ such that $q\in (f\circ\Phi_0\circ\Phi_1\circ\Phi_2)(C')\subset \Sigma(Y)$ so that
the 2-curves contained in the exceptional divisor of $\Phi_2$ have order larger than the order of the (strict transform of the) 2-curves $C$ of $X$ such that $q\in f(C)\subset \Sigma(Y)$. Let $\Phi_3:Y_3\rightarrow Y_2$ be the blow ups of the 2-curve $C_2$ on $Y_3$ of smallest order. Let 
$\overline\Phi_n=\Phi_2\circ\cdots\Phi_n:X_n\rightarrow X_1$. The morphisms $f\circ\Phi_0\circ\Phi_1\circ\overline\Phi_n$ are prepared, and 
$u,v,w_i$ are super parameters for $f\circ\Phi_0\circ\Phi_1\circ\overline\Phi_n$.

Suppose that $\nu$ is a 0-dimensional valuation of ${\bold k}(X)$.  Let $p_n$ be the center of $\nu$ on $X_n$.
Say that $\nu$ is resolved on $X_n$ if  (at least)
one of the following holds:
\begin{enumerate}
\item[1.] $m_q{\cal O}_{X_n,p_n}$ is invertible.
\item[2.] $p_n$ is a 1-point of the form (\ref{eq43}).
\item[3.] $p_n$ is a 2-point of the form (\ref{eq44}).
\end{enumerate}

If $\nu$ is resolved on $X_n$, with center $p_n$, then there exists an open neighborhood $U$ of
$p_n$ in $X_n$  such that a 0-dimensional valuation $\overline \nu$ of ${\bold k}(X)$
is resolved on $X_n$  if the center  of $\overline\nu$ is in $U$,  since 1, 2 or 3 is an open condition
on $X_n$.
Further, if $n'>n$, then $X_{n'}$ is also
resolved at all 0-dimensional valuations $\overline\nu$ which are resolved on $X_n$.

If the center of $\nu$ on $X_n$ is a 3-point, then the center of $\nu$ on $X_1$ is also a 3-point, so $m_q{\cal O}_{X_n,p_n}$ is invertible.

Suppose that the center of $\nu$ on $X_n$ is a 1-point. Then $u,v,w_i$ have a form 1 of Definition \ref{Def357} at $p_n$, and thus have a form (\ref{eq43}) at $p_n$ if $m_q{\cal O}_{X_n,p_n}$ is not invertible.

Suppose that the center $p$ of $\nu$ on $X_1$ is a 2-point such that $u,v,w_i$ have a form 2 of Definition \ref{Def357} at $p$.
There exist $\overline x,\overline y\in{\cal O}_{X_1,p}$ and series $\lambda_1,\lambda_2\in\hat{\cal O}_{X_1,p}$ such that
$x=\overline x\lambda_1$, $y=\overline y\lambda_2$.
Let $I^p\subset {\cal O}_{X_1,p}$ be the ideal $(u,v,\overline x^e\overline y^f,\overline x^g\overline y^h)$.
There exists an $n$
such that $I^p{\cal O}_{X_{n},p_1}$ is invertible for all $p_1\in\overline\Phi_{n}^{-1}(p)$.
Let $p_1$ be the center of $\nu$ on $X_{n}$. $u,v,w_i$ are super parameters for $f\circ\Phi_0\circ\Phi_1\circ\overline\Phi_{n}$ at $q$. 
 Then $p_1$ is either a 1-point (so that $\nu$ is resolved on $X_n$) or a 2-point of the form 2 of Definition
\ref{Def357}. 
If $p_1$ is a 2-point and $m_q{\cal O}_{X_{n},p_1}$ is not invertible, then $u,v,w_i$ have a form (\ref{eq44})
(with $\overline\gamma=0$).
Thus $\nu$ is resolved on $X_{n}$.

Suppose that the center $p$ of $\nu$ on $X_1$ is a 2-point such that $u,v,w_i$ have a  form 3  of Definition \ref{Def357} at $p$. Then there exist $\overline x,\overline y\in{\cal O}_{X,p}$ and series $\lambda_1,\lambda_2\in\hat{\cal O}_{X,p}$ such that $x=\overline x\lambda_1$, $y=\overline y\lambda_2$. 
 Let $I^p\subset{\cal O}_{X,p}$ be the ideal
$$
I^p=(u,v, (\overline x^a\overline y^b)^l,\overline x^c\overline y^d).
$$
There exists an $n$ such that $I^p{\cal O}_{X_{n},p_1}$ is invertible for all $p_1\in\overline \Phi_{n}^{-1}(p)$.
Let $p_1$ be the center of $\nu$ on $X_{n}$.   
 Then  $m_q{\cal O}_{X_{n},p_1}$ is invertible
or $p_1$ is a 1-point of the form (\ref{eq43}).
Thus $\nu$ is resolved on $X_{n}$.

By compactness of the Zariski-Riemann manifold \cite{Z1}, there exists an $n$ 
such that every valuation $\nu$ of ${\bold k}(X)$ is resolved on  $X_{n}$. Let $X_1=X_n$ and $f_1=f\circ\Phi_0\circ\Phi_1\circ\overline\Phi_n:X_1\rightarrow Y$.
If $m_q{\cal O}_{X_1,p_1}$ is not invertible at some $p_1\in f_1^{-1}(q)$, then one of
the forms (\ref{eq43}) or (\ref{eq44}) must hold at $p_1$.

The locus of points  $p$ on $X_1$ where $m_q{\cal O}_{X_1,p}$ is not invertible is a  (possibly not irreducible) curve $\overline E$ which makes SNCs with the toroidal structure of
$X$. $\overline E$ is supported at points of the form (\ref{eq43}) (with $d<\text{min}(a,b)$) and (\ref{eq44}) (with $\overline\gamma=0$
and $(g,h)<\text{min}\{(a,b),(c,d)\}$).
$x=z=0$ is a local equation of $\overline E$ in (\ref{eq43}). $x=z=0$, $y=z=0$ or $xy=z=0$ are the possible
local equations of $\overline E$ in (\ref{eq44}).

For an irreducible component $C$ of $\overline E$, define an invariant
$$
A(C)=\text{min}\{a,b\}-d>0
$$
computed at a 1-point $p\in C$ (which has an expression (\ref{eq43})).  Let $C$ be a component of $\overline E$ such that
$A(C)=\text{max}_{\overline C\subset E}A(\overline C)$. $C$ is nonsingular and makes SNCs with
$D_{X_1}$.
Let $\Phi_2:X_2\rightarrow X_1$ be the blow up of $C$.

Suppose that $p\in C$, so that $u,v,w_i$ have  the form
(\ref{eq43}) or (\ref{eq44}) at $p$. We may assume that $x=z=0$ are local equations of $C$ at $p$. Suppose that $p_1\in\Phi_2^{-1}(p)$. $p_1$ has (formal) regular parameters
$x_1,y_1,z_1$ defined by 
\begin{equation}\label{eq47}
x=x_1, y=y_1, z=x_1(z_1+\beta)
\end{equation}
with $\beta\in \bold k$ or 
\begin{equation}\label{eq48}
x=x_1z_1,z=z_1.
\end{equation}

Suppose that $p\in C$ is a 1-point, so that (\ref{eq43}) holds for $u,v,w_i$ at $p$.
Under (\ref{eq47}) we have $m_q\hat{\cal O}_{X_2,p_1}$ and thus $m_q{\cal O}_{X_2,p_1}$ is invertible except possibly if $\beta=0$.
If $m_q{\cal O}_{X_2,p_1}$ is not invertible we have
$$
u=x_1^a, v=x_1^b(\alpha+y), w_i=x_1^{d+1}z_1
$$
with
$$
d+1<\text{min}\{a,b\}.
$$
Then the curve $C_1$ with local equations $x_1=z_1$ is a component of the locus where
$m_q{\cal O}_{X_1}$ is not invertible. We have
$$
A(C_1)=\text{min}\{a,b\}-(d+1)<A(C).
$$
Under (\ref{eq48}) we have a 2-point
$$
u=(x_1z_1)^a,
v=(x_1z_1)^b(\alpha+y),
w_i=x_1^dz_1^{d+1}
$$
and a local equation of the toric structure $D_{X_2}$ is $x_1z_1=0$.
Since
$$
d+1\le\text{min}\{a,b\},
$$
$m_q{\cal O}_{X_1,p_1}$ is invertible.

Suppose that $p\in C$ is a 2-point, so that 
(\ref{eq44}) holds at $p$ (with $\overline\gamma=0$).
We may assume that $x=z=0$ is a local equation of $C$ at $p$. Then
$g<\text{min}\{a,c\}$. Let $p_1\in\Phi_2^{-1}(p)$. $p_1$ has regular parameters
$x_1,y_1,z_1$ defined by (\ref{eq47}) or (\ref{eq48}).

Under the substitution (\ref{eq47}) $u,v,w_i$ have  a form (\ref{eq44}) at $p_1$. If  $m_q{\cal O}_{X_2,p_1}$
 is not invertible, we have $\beta=0$ and
$$
u=x_1^ay_1^b, v=x_1^cy_1^d, w_i=x_1^{g+1}y_1^hz_1
$$
which is back in the form (\ref{eq44}) with $\overline\gamma=0$.

Under the substitution (\ref{eq48}), 
\begin{equation}\label{eq212}
u=x_1^ay_1^bz_1^a,
v=x_1^cy_1^dz_1^c,
w_i=x_1^gy_1^hz_1^{g+1}
\end{equation}
so that $p_1$ is a 3-point, and since 
$$
(g+1,h)\le\text{min}\{(a,b),(c,d)\},
$$
$m_q{\cal O}_{X_2,p_1}$ is invertible.

Observe that  $u,v,w_i$ are super parameters for $f_1\circ\Phi_2$. 

By descending induction on $\text{max}_{\overline C\subset E}\{A(\overline C)\}$, we construct a sequence of blow ups
$\Phi_4:X_4\rightarrow X_2$ such that for $\overline f_4=f_1\circ\Phi_2\circ\Phi_4:X_4\rightarrow Y$,
 $u,v,w_i$ are super parameters at $q$ for $\overline f_4$, and
 $m_q{\cal O}_{X_4}$ is  invertible.  Thus
$\overline f_4:X_4\rightarrow Y$  factors through the blow up $\Psi:Y_1\rightarrow Y$
of $q$. Let $f_4:X_4\rightarrow Y_1$, $\overline\Phi:X_4\rightarrow X$ be the resulting maps.

 Let
$q_1\in \Psi^{-1}(q)$.
We obtain permissible
parameters $\overline u, \overline v,\overline w$ at $q_1$ of one of the following forms:
\begin{enumerate}
\item[1.] $q_1$ a  1-point
$$
u=\overline u, v=\overline u(\overline v+\alpha),
w_i=\overline u(\overline w+\beta)
$$
with $\alpha,\beta\in \bold k$, $\alpha\ne 0$. In this case there are no 3-points in $f_4^{-1}(q_1)$ and $\overline u,\overline v$ are toroidal forms at 
all points $p\in f_4^{-1}(q_1)$. 
\item[2.] $q_1$ a 2-point 
\begin{equation}\label{eq118}
u=\overline u,
v=\overline u\overline v,
w_i=\overline u(\overline w_i+\alpha)
\end{equation}
with $\alpha\in \bold k$, or 
\begin{equation}\label{eq172}
u=\overline u\overline v,
v=\overline v,
w_i=\overline v(\overline w_i+\alpha)
\end{equation}
with $\alpha\in \bold k$, or
\item[3.] $q_1$ a 3-point 
\begin{equation}\label{eq238}
u=\overline u\overline w_i,v=\overline v\overline w_i,w_i=\overline w_i.
\end{equation}
\end{enumerate}

If $q_1$ has the form (\ref{eq118}) or (\ref{eq172}) and $p\in f_4^{-1}(q_1)$ then $\overline u,\overline v$ are toroidal forms at $p$.

Suppose that $q_1$ has the form (\ref{eq238}). Let
$p\in f_4^{-1}(q_1)$. $u,v,w_i$ have one of the forms 1 - 4 of Definition \ref{Def357} at $p$.
Since $w_i$ must divide $u$ and $v$, we certainly have that $\overline u,\overline v,\overline w_i$ are monomials in the local equations of the toroidal structure at $p$, times unit series.

Suppose that $u,v,w_i$ have a form 1 of Definition \ref{Def357}  at $p$. We have an expression
$$
\begin{array}{ll}
u&=x^a\\
v&=x^b(\alpha+y)\\
w_i&=x^c(\tilde\gamma(x,y)+x^{\tilde d}z)
\end{array}
$$
where $0\ne\alpha\in \bold k$, $c<a$, $c<b$, $\tilde\gamma$ is a unit series and $\tilde d\ge0$.
Set $\overline x=x(\tilde\gamma+x^{\tilde d}z)^{\frac{1}{c}}$. We have expansions
$$
\begin{array}{ll}
w_i&=\overline x^c\\
u&=\overline x^a(\tilde\gamma+x^{\tilde d}z)^{-\frac{a}{c}}\\
v&=\overline x^b(\tilde\gamma+x^{\tilde d}z)^{-\frac{b}{c}}(\alpha+y)
\end{array}
$$
at $p$.

If $\tilde d>0$ and $\frac{\partial\tilde\gamma}{\partial y}(0,0)=0$, then
there exist $\overline y,\overline z\in \hat{\cal O}_{X_4,p}$ such that $\overline x,
\overline y,\overline z$ are regular parameters in $\hat{\cal O}_{X_4,p}$ and
$$
\begin{array}{ll}
w_i&=\overline x^c\\
v&=\overline x^b(\overline\alpha+\overline y)\\
u&=\overline x^a\hat\gamma(\overline x,\overline y,\overline z)
\end{array}
$$
where $0\ne\overline\alpha\in \bold k$ and $\hat \gamma$ is a unit series. Then

$$
\begin{array}{ll}
\overline w_i&=\overline x^c\\
\overline v&=\overline x^{b-c}(\overline\alpha+\overline y)\\
\overline u&=\overline x^{a-c}\hat\gamma(\overline x,\overline y,\overline z)
\end{array}
$$
and $\overline w_i$, $\overline v$ are toroidal forms at $p$.

If $\tilde d=0$ or $\frac{\partial \tilde\gamma}{\partial y}(0,0)\ne 0$, then
there exist $\overline y,\overline z\in\hat{\cal O}_{X_4,p}$ such that $\overline x,\overline y,
\overline z$ are regular parameters in $\hat{\cal O}_{X_4,p}$ and
$$
\begin{array}{ll}
w_i&=\overline x^c\\
u&=\overline x^a(\overline\alpha+\overline y)\\
v&=\overline x^b\hat\gamma(\overline x,\overline y,\overline z)
\end{array}
$$
where $0\ne\overline\alpha\in \bold k$ and $\hat\gamma$ is a unit series. Then
$$
\begin{array}{ll}
\overline w_i&=\overline x^c\\
\overline u&=\overline x^{a-c}(\overline\alpha+\overline y)\\
\overline v&=\overline x^{b-c}\hat\gamma
\end{array}
$$
and $\overline w_i, \overline u$ are toroidal forms  at $p$.

Suppose that $u,v,w_i$ have a form 3 of Definition \ref{Def357} at $p$. There are
two cases. Either 
\begin{equation}\label{eq239}
\begin{array}{ll}
u&=(x^ay^b)^k\\
v&=(x^ay^b)^t(\alpha+z)\\
w_i&=x^cy^d
\end{array}
\end{equation}
with $0\ne\alpha\in \bold k$ and $ad-bc\ne 0$,  or 
\begin{equation}\label{eq240}
\begin{array}{ll}
u&=(x^ay^b)^k\\
v&=(x^ay^b)^t(\alpha+z)\\
w_i&=(x^ay^b)^l\tilde \gamma
\end{array}
\end{equation}
with $0\ne\alpha\in \bold k$, $l\le\text{min}\{k,t\}$ and $\tilde\gamma$ is a unit series.

If (\ref{eq239}) holds then $\overline u,\overline w_i$ 
are toroidal forms at $p$
of the form of (\ref{eqTF21}) of Definition \ref{torf}. 

Assume that (\ref{eq240}) holds. We then have
$$
\begin{array}{ll}
\overline u&=(x^ay^b)^{k-l}\tilde\gamma^{-1}\\
\overline v&=(x^ay^b)^{t-l}\tilde\gamma^{-1}(\alpha+z)\\
\overline w_i&=(x^ay^b)^l\tilde\gamma.
\end{array}
$$

If $\frac{\partial \tilde\gamma}{\partial z}(0,0,0)\ne 0$, then there exist regular
parameters $\overline x,\overline y,\overline z$ at $p$ such that
$$
\begin{array}{ll}
\overline u&=(\overline x^a\overline y^b)^{k-l}\\
\overline w_i&=(\overline x^a\overline y^b)^l(\overline\beta+\overline z)\\
\overline v&=(\overline x^a\overline y^b)^{t-l}\hat\gamma(\overline x,\overline y,\overline z)
\end{array}
$$
where $0\ne\overline\beta\in \bold k$ and $\hat\gamma$ is a unit series. 
 Thus $\overline u,\overline w_i$ are toroidal forms at $p$

If $\frac{\partial \tilde\gamma}{\partial z}(0,0,0)=0$ then there exist regular parameters
$\overline x,\overline y,\overline z$ at $p$ such that
$$
\begin{array}{ll}
\overline u&=(\overline x^a\overline y^b)^{k-l}\\
\overline v&=(\overline x^a\overline y^b)^{t-l}(\overline\beta+\overline z)\\
\overline w_i&=(\overline x^a\overline y^b)^l\hat\gamma(\overline x,\overline y,\overline z)
\end{array}
$$
where $0\ne\overline\beta\in \bold k$ and $\hat\gamma$ is a unit series. Thus $\overline u,\overline v$ are toroidal forms at $p$.

If $u,v,w_i$ have a form 2 or 4 of Definition \ref{Def357} at $p$, then a simpler analysis shows that two of $\overline u,\overline v,\overline w_i$ are toroidal forms at $p$, and $\overline u,\overline v,\overline w_i$ are monomials
in local equations of the toroidal structure at $p$ times unit series.

We conclude that the morphism $f_4$ is prepared.

We will now verify that $f_4$ is $\tau$-quasi-well prepared.

For $j$ such that $q\in U(\overline R_j)$, let $w_j=w_{\overline R_j}(q)$.
We will analyze our construction of $X_4\rightarrow X$ to show that $u,v,w_j$
are super parameters for $\overline f_4$.

Since $f$ is $\tau$-quasi-well prepared,  for $j$ such that
$q\in U(\overline R_j)$, there exists a series $\lambda_{ij}(u,v)$
such that $w_j=w_i+\lambda_{ij}(u,v)$. 

The morphism $\Phi_0\circ\Phi_1:X_1\rightarrow X$
which we constructed is a product of blow ups of 2-curves and 3-points. Since $u,v,w_j$ are super parameters for $f$, $u,v,w_j$
are thus also super parameters for $f\circ\Phi_0\circ\Phi_1$. 

 $X_4\rightarrow X_1$ is a sequence of
blow ups  of curves $C$ such that $A(C)>0$.
Suppose that $C$ is such a curve, and $p\in C$. 
It suffices to analyze the blow up $\Phi_2:X_2\rightarrow X_1$ of 
$C$ in our construction.
We saw that $u,v,w_i$ have a form (\ref{eq43}) at $p$ with
$d<\text{min}\{a,b\}$ or a form (\ref{eq44}) with $\overline\gamma=0$, $(g,h)\le\text{min}\{(a,b),(c,d)\}$  and
$g<\text{min}\{a,c\}$. In either case $x=z=0$ are local equations of $C$. We will
show that for all $j$ such that $q\in U(\overline R_j)$, $u,v,w_j$ are super parameters for $f\circ\Phi_0\circ\Phi_1\circ\Phi_2$.

First assume that $p\in C$ is a 1-point. Then,   we have
$$
u=x^a, v=x^b(\alpha+y), w_i=x^dz
$$
with $0\ne\alpha\in \bold k$ and $d<\text{min}\{a,b\}$. Thus, since
$$
w_j=w_i+\lambda_{ij}(u,v),
$$
 $w_j=x^d\overline z$ where $\overline z=z+x\Omega(x,y)$
for some series $\Omega$. It follows that $x=\overline z=0$ are local equations of $C$
at $p$, and $u,v,w_j$ are super parameters for $f\circ\Phi_1\circ\Phi_2$.

Now assume that $p\in C$ is a 2-point. Then  we have
$$
u=x^ay^b, v=x^cy^d, w_i=x^gy^hz
$$
where after possibly interchanging $u$ and $v$, we have
$(g,h)<(a,b)\le(c,d)$ (recall that $(u,v){\cal O}_{X_1,p}$ is invertible). 
Since $w_j-w_i\in {\bold k}[[u,v]]$, we have an expression 
\begin{equation}\label{eq243}
w_j=x^gy^h\overline z
\end{equation}
at $p$ where $x=\overline z=0$ are local equations of $C$ at $p$.
it follows that $u,v,w_j$ are super parameters for $f\circ\Phi_0\circ\Phi_1\circ\Phi_2$.

 By induction, $u,v,w_j$ are super parameters for $\overline f_4=f\circ\Phi_0\circ\Phi_1\circ\Phi_2\circ\Phi_4$.

$q$ is an admissible center for all 2-point relations $\overline R_j$ associated to $R$
(Definition \ref{Def154}). For all $j$, let $\overline R_j^1$ be the transform of $\overline R_j$ on $Y_1$.
Suppose that $p_1\in\overline f_4^{-1}(q)\cap\overline\Phi^{-1}(T(\overline R_j))$
is a 3-point.
$p=\overline\Phi(p_1)\in T(\overline R_j)$ is a 3-point with permissible parameters $x,y,z$
such that 
\begin{equation}\label{eq209}
\begin{array}{ll}
u&=x^ay^bz^c\\
v&=x^dy^ez^f\\
w_j&=M_0\gamma
\end{array}
\end{equation}
where $\text{rank}(u,v)=2$, $\gamma(x,y,z)$ is a unit series, $M_0$ is a monomial in $x,y,z$ and $w_j=w_{\overline R_j}(q)$. Let $w_j^{e_j}=\overline\lambda_ju^{a_j}v^{b_j}$
define $\overline R_j(q)$ if $\tau>0$,
$w_j=0$ define $\overline R_j(q)$ if $\tau=0$. If $\tau>0$, then 
\begin{equation}\label{eq210}
M_0^{e_j}=u^{a_j}v^{b_j}
\text{and }\gamma^{e_j}(0,0,0)=\overline\lambda_j,
\end{equation}
and $\text{rank}(u,v,M_0)=3$ if $\tau=0$.
By its construction, $\overline\Phi$ is a sequence of blow ups
of 2-curves and 3-points above $p$.

Thus since $p_1$ is a 3-point, we have permissible parameters $x_1,y_1,z_1$ at $p_1$ such that
$$
\begin{array}{ll}
x&=x_1^{a_{11}}y_1^{a_{12}}z_1^{a_{13}}\\
y&=x_1^{a_{21}}y_1^{a_{22}}z_1^{a_{23}}\\
z&=x_1^{a_{31}}y_1^{a_{32}}z_1^{a_{33}}
\end{array}
$$
and $\text{Det}(a_{ij})=\pm1$.
On substitution into (\ref{eq209}) we see that an expression 
\begin{equation}\label{eq370}
\begin{array}{ll}
u&=x_1^{\overline a}y_1^{\overline b}z_1^{\overline c}\\
v&=x_1^{\overline d}y_1^{\overline e}z_1^{\overline f}\\
w_j&=M_0\gamma
\end{array}
\end{equation}
where $\text{rank}_{(x_1,y_1,z_1)}(u,v)=2$,
  holds at $p_1$
for $u,v,w_j$, and the relation (\ref{eq210}) holds at $p_1$ if $\tau>0$, and $\text{rank}_{(x_1,y_1,z_1)}(u,v,M_0)=3$
if $\tau=0$.

Suppose that $q_1\in U(\overline R_j^1)\cap \Psi^{-1}(q)$. After possibly interchanging $u$ and $v$,
$q_1$ has permissible parameters $\overline u,\overline v,\overline w_j$ with
\begin{equation}\label{eq245}
u=\overline u,
v=\overline u\overline v,
w_j=\overline u\overline w_j.
\end{equation}
The pre-relation $\overline R_j^1(q_1)$ is then defined if $\tau>0$ by 
\begin{equation}\label{eq211}
\overline w_j^{e_j}=\overline\lambda_j \overline u^{a_j+b_j-e_j}\overline v^{b_j},
\end{equation}
and $\overline R_j^1(q_1)$ is defined by $\overline w_j=0$ if $\tau=0$.
We have seen that if $p_1\in \overline f_4^{-1}(q)\cap \overline\Phi^{-1}(T(\overline R_j))$ is a 3-point,
then there are permissible parameters $x_1,y_1,z_1$ at $p_1$ such that an expansion of the form (\ref{eq370})  holds for $u,v,w_j$, and (\ref{eq210}) holds if $\tau>0$. If $\tau=0$, then $u,v,w_j$ is a monomial form at $p$. If we also have that $p_1\in f_4^{-1}(q_1)$
then we have the expression
$$
\begin{array}{ll}
\overline u&=x_1^{\overline a}y_1^{\overline b}z_1^{\overline c}\\
\overline v&=x_1^{\overline d-\overline a}y_1^{\overline e-\overline b}z_1^{\overline f-\overline c}\\
\overline w_j&=\frac{M_0}{\overline u}\gamma
\end{array}
$$
at $p_1$, and (if $\tau>0$) (\ref{eq210}) becomes
$$
(\frac{M_0}{\overline u})^{e_j}=\overline u^{a_j+b_j-e_j}\overline v^{b_j}.
$$
 Thus the transform $R^1$ of $R$ for
$f_4$ (Definition \ref{Def161}) is defined.
We further have $\tau_{f_4}(p_1)=\tau_f(\overline\Phi(p_1))=\tau$.

If $p\in\overline f_4^{-1}(q)\cap \overline\Phi^{-1}(T(\overline R_j))$
is a 3-point and
$p\not\in f_4^{-1}(U(\overline R_j^1))$ then $f_4(p)=q_1$ where (after possibly interchanging
$u$ and $v$) $q_1$ has permissible parameters 
\begin{equation}\label{eq246}
\begin{array}{ll}
u&=\overline u\\
v&=\overline u(\overline v+\alpha)\\
w_j&=\overline u(\overline w_j+\beta)
\end{array}
\end{equation}
$\alpha,\beta\in \bold k$ and and at least one of $\alpha,\beta$ is zero, or 
\begin{equation}\label{eq213}
\begin{array}{ll}
u&=\overline u\overline w_j\\
v&=\overline v\overline w_j\\
w_j&=\overline w_j.
\end{array}
\end{equation}

Suppose that (\ref{eq246}) holds at $q_1$. Since $\text{rank}(u,v)=2$ in (\ref{eq370}), we must have
$\alpha=0$ and $0\ne\beta$. But then $M_0=u$, a contradiction to the assumption that $e_j>1$ (and $\text{gcd}(a_j,b_j,e_j)=1$)  in (\ref{eq210}) if $\tau>0$, or to the assumption that
$\text{rank}(u,v,M_0)=3$ if $\tau=0$.

Suppose that (\ref{eq213}) holds at $q_1$. Then $q_1=f_4(p)$ is a 3-point. From equations (\ref{eq213})
and (\ref{eq370}) we have (in the notation of Definition \ref{Def221}) that
$\tau_{f_4}(p)=-\infty$ if $\tau=0$ and if $\tau>0$ then
$$
H_{f_4,p}=H_{\overline f_4,p}=H_{f,\overline\Phi(p)},
$$
$$
A_{f_4,p}=A_{f,\overline \Phi(p)}+{\bf Z}M_0
$$
since $q=\overline f_4(p)=f(\overline\Phi(p))$ is a 2-point. Thus, since $e_j>1$ in (\ref{eq210}), we have
$$
\tau_{f_4}(p)=\mid H_{f_4,p}/A_{f_4,p}\mid
<\mid H_{f,\overline \Phi(p)}/A_{f,\overline\Phi(p)}\mid=\tau.
$$

Finally, suppose that $p\in \overline f_4^{-1}(q)-\cup_j\overline\Phi^{-1}
(T(\overline R_j))$ is a 3-point.  Suppose that $\tilde p=\Phi_2\circ\Phi_4(p)$ is a 2-point. Then $u,v,w_i$
have a form (\ref{eq44}) at $\tilde p$, and (with $\overline\gamma=0$) $u,v,w_i$ have a form (\ref{eq212}) at $p$, with $(g+1,h)\le\text{min}\{(a,b),(c,d)\}$.
We see (since $(u,v){\cal O}_{X_4, p}$ is invertible) that $w_i\mid u$, $w_i\mid v$ at $p$, and thus $f_4(p)$
is a 3-point with permissible parameters $\overline u,\overline v,\overline w_i$ defined by
$$
u=\overline u\overline w_i, v=\overline v\overline w_i, w_i=\overline w_i.
$$
Thus, $\overline u,\overline v,\overline w_i$ have a toroidal form at $p$, so that $\tau_{f_4}(p)=-\infty<\tau$.

Suppose that $\tilde p =\Phi_2\circ\Phi_4(p)$ is a 3-point. Then $\Phi_2\circ\Phi_4$ is the identity near $p$, and thus $\overline\Phi(p)$ is a 3-point and $\overline\Phi$ factors as a sequence of blow ups of 2-curves and 3-points near $p$. Thus
$\tau_{f_4}(p)\le \tau_f(\overline\Phi(p))<\tau$.
Thus 1, 2 and 3 of Definition \ref{Def128} hold for $f_4$.

Now we verify 4 of Definition \ref{Def128}. Suppose that $q_1\in U( R^1)$. Let
$q=\Psi(q_1)\in U(R)$. Since $f$ is $\tau$-quasi-well prepared, there exists
$w_q\in{\cal O}_{Y,q}$ satisfying 4 of Definition \ref{Def128} for $f$. If
$q_1\in U(\overline R_i^1)\cap U(\overline R_j^1)$ then (after possibly interchanging $u$ and $v$) we have
$$
\begin{array}{ll}
u&=\overline u\\
v&=\overline u\overline v\\
w_i&=\overline u\overline w_i\\
w_j&=\overline u\overline w_j
\end{array}
$$
where
$$
u=u_{\overline R_i}(q)=u_{\overline R_j}(q), v=v_{\overline R_i}(q)=v_{\overline R_j}(q), w_i=w_{\overline R_i}(q), w_j=w_{\overline R_j}(q)
$$
and
$$
\overline u=u_{\overline R_i^1}(q_1)=u_{\overline R_j^1}(q_1), \overline v=v_{\overline R_i^1}(q_1)=u_{\overline R_j^1}(q_1),
 \overline w_i=w_{\overline R_i^1}(q_1), \overline w_j=w_{\overline R_j^1}(q_1).
$$

Since $f$ is $\tau$-quasi-well prepared, there exists a series $\lambda_{ij}(u,v)$ such that
$$
w_j=w_i+\lambda_{ij}(u,v).
$$
Since $\lambda_{ij}(0,0)= 0$, $\overline u\mid\lambda_{ij}(\overline u,\overline u\overline v)$,
and 
$$
\overline w_j=\overline w_i+\frac{\lambda_{ij}(\overline u,\overline u\overline v)}{\overline u}.
$$
Thus 4 of Definition \ref{Def128} holds for $f_4$.

Earlier in the proof we verified that if $q\in U(\overline R_j)$, for some $\overline R_j$ associated to $R$, then 
$$
u=u_{\overline R_j}(q), v=v_{\overline R_j}(q), w_j=w_{\overline R_j}(q)
$$
are super parameters for $\overline f_4$.
If $q_1\in U(\overline R_j^1)\cap \Psi^{-1}(q)$, then (after possibly interchanging
$u$ and $v$) we have permissible parameters
$$
\overline u=u_{\overline R_j^1}(q_1), \overline v=v_{\overline R_j^1}(q_1),
\overline  w_j=w_{\overline R_j^1}(q_1)
$$
such that
$$
u=\overline u,
v=\overline u\overline v,
w_j=\overline u\overline w_j.
$$
Substituting into the forms of Definition \ref{Def357}, we see that $\overline u,\overline v,\overline w_j$ are
super parameters for $f_4$ at $q_1$. Thus 5 of Definition
\ref{Def128} holds for $f_4$ and $f_4$ is $\tau$-quasi-well prepared.

We now verify that $D_{X_4}$ is cuspidal for $f_4$ if $D_X$ is cuspidal for $f$  (this is 1 of the conclusions of the lemma). Since the property of being cuspidal is
stable under blow ups of 2-curves and 3-points, it suffices to show that if $C$ is a component of
$\overline E$ on $D_{X_1}$, such that $A(C)>0$, $C$ contains no 2-points, and $p\in C$ then there
exists a Zariski open neighborhood $U$ of $p$ in $X_1$ such that $f_4$ is toroidal on
$(\Phi_2\circ\Phi_4)^{-1}(U)$.

There exist permissible
parameters $x,y,z$ at $p$ such that (after possibly interchanging $u$ and $v$) we have expressions 
\begin{equation}\label{eq312}
\begin{array}{ll}
u&=u_{\overline R_i}(q)=x^a\\ 
v&=v_{\overline R_i}(q)=x^b(\alpha+y)\\
w_i&=w_{\overline R_i}(q)=x^dz
\end{array}
\end{equation}
with $0\ne\alpha$, $d<a\le b$  and $x=z=0$ are local equations of $C$ at $p$.

We consider the effect of the blow up of $C$, $\Phi_2:X_2\rightarrow
X_1$.  If $p_1\in \Phi_2^{-1}(p)$ then $p_1$ has regular parameters $x_1,y_1,z_1$
of one of the forms 
\begin{equation}\label{eq249}
x=x_1, y=y_1, z=x_1(z_1+\gamma)
\end{equation}
with $\gamma\in \bold k$, or 
\begin{equation}\label{eq250}
x=x_1z_1, y=y_1, z=z_1.
\end{equation}
Under (\ref{eq250}) we have a local factorization of the rational map $X_1\rightarrow Y_1$, at $p$,
by
$$
\begin{array}{ll}
\frac{u}{w_i}&=x_1^{a-d}z_1^{a-d-1}\\
w_i&=x_1^dz_1^{d+1}\\
\frac{v}{w_i}&=x_1^{b-d}z_1^{b-d-1}(\alpha+y_1)
\end{array}
$$
which is toroidal (of the form 2 of Definition \ref{Def274}).

Under (\ref{eq249}) we obtain a form 
$$
u=x_1^a,
v=x_1^b(\alpha+y),
w_i=x_1^{d+1}(z_1+\gamma)
$$
 which gives a toroidal factorization of
$X_1\rightarrow Y_1$, at $p$ (of the form 3 or 5 following Definition \ref{Def274}), except if $\gamma=0$ and $x_1=z_1=0$ 
are local equations of a curve $C_1$ with
$0<A(C_1)<A(C)$, and such that $C_1$ contains no 2-points (or 3-points). By successive blowing up of curves in the construction of
$\Phi_4:X_4\rightarrow X_2$, we see that  $X_4\rightarrow Y_1$
is a toroidal morphism  on $(\Phi_2\circ\Phi_4)^{-1}(U)$
for some Zariski open neighborhood $U$ of $p$.

Now suppose that $f$ is $\tau$-well prepared. We will verify that $f_4$ is $\tau$-well
prepared. 1 and 2 of Definition \ref{Def65} are immediate. 

We will verify that 3 of Definition \ref{Def65} holds for $f_4$.
Suppose that $q_1\in U(\overline R_i^1)\cap U(\overline R_j^1)$ and
$$
\overline u=u_{\overline R_i^1}(q_1)=u_{\overline R_j^1}(q_1),
\overline v=v_{\overline R_i^1}(q_1)=v_{\overline R_j^1}(q_1),
$$
$$
\overline w_i=w_{\overline R_i^1}(q_1), \overline w_j=w_{\overline R_j^1}(q_1).
$$
Let $q=\Psi(q_1)$,
$$
u=u_{\overline R_i}(q)=u_{\overline R_j}(q),
v=v_{\overline R_i}(q)=v_{\overline R_j}(q),
$$
$$
w_i=w_{\overline R_i}(q), w_j=w_{\overline R_j}(q).
$$
Since $f$ is $\tau$-well prepared, there exist unit series $\phi_{ij}(u,v)$ such that
$$
w_j=w_i+u^{a_{ij}}v^{b_{ij}}\phi_{ij}
$$
(or $\phi_{ij}=0$ and $a_{ij}=b_{ij}=\infty$).
After possibly interchanging $u$ and $v$, we may assume that
$$
u=\overline u, v=\overline u\overline v,
w_i=\overline u\overline w_i,
w_j=\overline u\overline w_j.
$$
Let
$$
\overline\phi_{ij}(\overline u,\overline v)=\phi_{ij}(\overline u,\overline u\overline v).
$$
Then we have 
\begin{equation}\label{eq214}
\overline w_j=\overline w_i+\overline u^{\overline a_{ij}}\overline v^{\overline b_{ij}}
\overline\phi_{ij}
\end{equation}
where 
\begin{equation}\label{eq244}
\overline a_{ij}=a_{ij}+b_{ij}-1,
\overline b_{ij}=b_{ij}.
\end{equation}
Thus 3 of Definition \ref{Def65} holds for $f_4$. 

Since the set  (\ref{eq255})   associated to $q$ and $R$ is totally ordered, the set (\ref{eq255}) 
associated to $q_1$ and $R^1$ is also totally ordered. Thus 4 of Definition \ref{Def65} holds
for $f_4$, and we see that $f_4$ is $\tau$-well prepared.

 We now verify 2  of Lemma \ref{Lemma32}.
Suppose that $q_1\in U(\overline R_i^1)\cap E$ for some $\overline R_i^1$ associated to $R^1$.
Continuing with the notation we used in the verification that $f_4$ is $\tau$-well prepared,
let $\gamma_i=\overline{S_{\overline R_i^1}(q_1)\cdot E}$. $\gamma_i$ is covered by two affine charts,
with uniformizing parameters $\overline u,\overline v,\overline w$ defined by 
\begin{equation}\label{eq76}
u=\overline u,
v=\overline u\overline v,
w_i=\overline u\overline w_i
\end{equation}
and 
\begin{equation}\label{eq77}
u=\tilde u\tilde v,
v=\tilde v,
w_i=\tilde v\tilde w_i.
\end{equation}

In the chart (\ref{eq76}), $\overline u=0$ is a local equation for $E$,
$\overline w_i=w_{\overline R_i^1}(q_1)=0$ is a local
equation for $S_{\overline R_i^1}(q_1)$ and $\overline v=0$ is a local equation for the strict transform of the
component $E_2$ of $D_Y$ with local equation $v=0$ at $q$. In (\ref{eq77}), $\tilde v=0$ is a local equation of $E$ and $\tilde u=0$ is a local equation of the
strict transform of the component $E_1$ of $D_Y$
with local equation $u=0$ at $q$. In the chart defined by (\ref{eq76}) $\overline u=\overline w_i=0$ are local equations of $\gamma_i$ and in the chart defined by (\ref{eq77}) $\tilde v=\tilde w_i=0$ are local equations of $\gamma_i$.
Thus $\gamma_i$ makes SNCs with $D_{Y_1}$, and $\gamma_i$ is a line on $E\cong {\bf P}^2$.

If $q\in U(\overline R_j)$ for some  $j\ne i$, then we see from (\ref{eq214}) and (\ref{eq244})
 that
$\gamma_j=\overline{S_{\overline R_j^1}(q')\cdot E}$ (where $q'\in \Psi^{-1}(q)\cap U(\overline R_j^1)$) has  local equations $\overline u=\overline w_j=0$ in the chart  (\ref{eq76})
where
$$
\overline w_j=\overline w_i+\overline u^{a_{ij}+b_{ij}-1}\overline v^{b_{ij}}\overline \phi_{ij}.
$$
In the chart (\ref{eq77}), $\gamma_j$ has local equations $\tilde v=\tilde w_j=0$
where $\tilde w_j=\tilde w_i+\tilde u^{a_{ij}}\tilde v^{a_{ij}+b_{ij}-1}\tilde \phi_{ij}$.
If $a_{ij}+b_{ij}>1$, then $\gamma_i=\gamma_j$ so that 2 (a) (i) holds in the statement of Lemma \ref{Lemma32}.
If $a_{ij}+b_{ij}=1$ then 2 (a) (ii) of the statement of Lemma \ref{Lemma32} holds.

We will now verify that $\gamma_i$ is a prepared curve for $R^1$ of type 4. 1 -- 4 of Definition \ref{Def200} are immediate from the above calculation.

We now verify that 5 of Definition \ref{Def200} holds for $\gamma_i$. Let
$$
I_{\gamma_i}=\{j\mid q\in U(\overline R_j)\text{ and }\gamma_j=\gamma_i\}.
$$
We have permissible parameters $u, v,w_j=w_{\overline R_j}(q)$ at $q$  for $j\in I_{\gamma_i}$.
By construction, $\gamma_i\cap U(\overline R_j^1)$ is the set of 2-points in $\gamma_i$ for all $j\in I_{\gamma_i}$.
Suppose that $q_1\in\gamma_i$ is a 1-point. Then we have (after possibly interchanging $u$ and $v$) permissible parameters $\overline u,\overline v,\overline w_j$
at $q_1$ for $j\in I_{\gamma_i}$ defined by 
\begin{equation}\label{eq290}
u=\overline u,
v=\overline u(\overline v+\alpha),
w_j=\overline u\overline w_j
\end{equation}
for some $0\ne\alpha\in \bold k$, where $\overline w_j=0$ is a local equation of
$\Omega(\overline R_j^1)$
at $q_1$.  Thus 5 (a) of Definition \ref{Def200} holds. From the relation
$$
w_j-w_i=u^{a_{ij}}v^{b_{ij}}\phi_{ij}(u,v),
$$
we have
$$
\overline w_j-\overline w_i=\overline u^{a_{ij}+b_{ij}-1}(\overline v+\alpha)^{b_{ij}}
\phi_{ij}(\overline u,\overline u(\overline v+\alpha)).
$$
Thus 5 (d) of Definition \ref{Def200} holds.

1 or 3 of Definition \ref{Def357}  hold at all
$$
p\in f_4^{-1}(q_1)\subset \overline f_4^{-1}(q)
$$
for $u,v,w_j$, and after substitution of (\ref{eq290}) into this form, we see that 5 (b)
of Definition \ref{Def200} holds. Thus $\gamma_i$ is prepared for $R^1$ of type 4.
We have completed the verification of 2 (a) of Lemma \ref{Lemma32}.

We now verify 2 (b) of Lemma \ref{Lemma32}. Suppose that $\gamma$ is prepared for $R$,
 $q\in\gamma$, and $\gamma$ is prepared for $R$ of
type 4. Then $\gamma=\overline{E_2\cdot S_{\overline R_i}(q)}$
for some $\overline R_i$ and component $E_2$ of $D_Y$ containing $q$ and $\gamma\subset\Omega(\overline R_i)$. Let $u=u_{\overline R_i}(q),v=v_{\overline R_i}(q),w_i=w_{\overline R_i}(q)$. We may assume that
$v=0$ is a local equation of $E_2$. Then $v=w_i=0$ are local equations at $q$ of $\gamma$.
Let $\gamma'$ be the strict transform of $\gamma$ on $Y_1$. $q_1=\gamma'\cdot E$ has
permissible parameters
$$
\overline u=u_{\overline R_i^1}(q_1),
\overline v=v_{\overline R_i^1}(q_1),
\overline w_i=w_{\overline R_i^1}(q_1),
$$
where
$$
u=\overline u,
v=\overline u\overline v,
w_i=\overline u\overline w_i,
$$
and $\overline v=\overline w_i=0$ are local equations of $\gamma'$.
Let $\tilde E_2$ be the strict transform of $E_2$. Since $q_1\in U(\overline R_i^1)$, we have $\gamma'=\overline{\tilde E_2\cdot S_{\overline R_i^1}(q_1)}$ and $\gamma'\subset\Omega(\overline R_i^1)$.  Thus the
conditions of Definition \ref{Def200} hold for $\gamma'$, and $\gamma'$ is prepared for $R^1$.
If $\gamma$ is a 2-curve, then the strict transform $\gamma'$ of $\gamma$ on $Y_1$ is a 2-curve so $\gamma'$ is prepared for $R^1$.

Finally, suppose that $f$ is $\tau$-very-well prepared. We have shown that 1 and 2 of
Definition \ref{Def130} hold for $f_4$.
Since whenever $\overline R_i$ is a pre-relation associated to $R$ containing $q$, $\Omega(\overline R_i^1)\rightarrow \Omega(\overline R_i)$ is the blow up of a point on a nonsingular
surface, and $V_i(Y)$ satisfies 3 of Definition \ref{Def130}, 3 of Definition \ref{Def130}
holds for $V_i(Y_1)$. Thus $f_4$ is $\tau$-very-well prepared.
\end{pf}

\begin{Lemma}\label{Lemma171} Suppose that $f:X\rightarrow Y$ is $\tau$-quasi-well
prepared (or $\tau$-well prepared or $\tau$-very-well prepared)
with  2-point relation $R$. Suppose that $q\in Y$ is
a 2-point such that $q\not\in U(R)$ and $q$ is prepared (of type 2 of Definition \ref{Def66}) for $R$. Then $q$ is a permissible center for $R$ and there exists  a $\tau$-quasi-well prepared (or $\tau$-well prepared
or $\tau$-very-well prepared)
diagram (\ref{eq30})of $R$ and the blow up $\Psi:Y_1\rightarrow Y$ of $q$ such that:
\begin{enumerate}
\item[1.] Suppose that $D_X$ is cuspidal for $f$. Then $D_{X_1}$ is cuspidal for $f_1$.
\item[2.] Suppose that $f$ is $\tau$-well prepared.
If $\gamma$ is a prepared curve on $Y$ then the strict transform of $\gamma$ is a prepared
curve on $Y_1$.
\item[3.] $\Phi$ is an isomorphism over $f^{-1}(Y-\Sigma(Y))$
\end{enumerate}
\end{Lemma}

\begin{pf} The proof is a simplification of the proof of Lemma \ref{Lemma32}.
\end{pf}

\begin{Lemma}\label{Lemma67}
Suppose that $f:X\rightarrow Y$ is $\tau$-very-well prepared
  and $C\subset Y$ is a prepared curve of type 4
(of Definition \ref{Def200}). Further suppose that
$q_{\delta}\in C\cap U(\overline R_j)$ for some $\overline R_j$ associated to $R$ implies $C=\overline{E\cdot S_{\overline R_j}(q_{\delta})}$ for some component $E$ of $D_X$. 
Then $C$ is a *-permissible center for $R$, and   there exists a $\tau$-very-well prepared diagram
$$
\begin{array}{rll}
X_1&\stackrel{f_1}{\rightarrow}&Y_1\\
\Phi_1\downarrow&&\downarrow\Psi_1\\
X&\stackrel{f}{\rightarrow}&Y
\end{array}
$$
of $R$  of the form of (\ref{eq233}). If $D_X$ is cuspidal for $f$ then  $D_{X_1}$ is cuspidal for $f_1$.
$\Phi_1$ is an isomorphism over $f^{-1}(Y-(\Sigma(Y)\cup C))$. If $C$ is contained in the fundamental locus of $f$ then $\Phi_1$ is an isomorphism over
$f^{-1}(Y-\Sigma(Y))$.
\end{Lemma}

\begin{pf} If $C$ is contained in the fundamental locus of $f$, then ${\cal I}_C{\cal O}_{X,p}$ is invertible if $p\in f^{-1}(Y-\Sigma(Y))$ by Lemma
\ref{Lemma143}.

Let $C=\overline{E_{\alpha}\cdot S_{\overline R_i}(q_{\beta})}$.
For $\overline q\in C$, we have permissible parameters 
\begin{equation}\label{eq215}
u,v,\tilde w_i
\end{equation}
such that $u=\tilde w_i=0$  are local equations at $\overline q$ for $C$,
with the notation of 5  of Definition \ref{Def200},
if $\overline q\not\in U(\overline R_i)$, and  $u=u_{\overline R_i}(\overline q)$,
$v=v_{\overline R_i}(\overline q)$, $\tilde w_i=w_{\overline R_i}(\overline q)$ if
$\overline q\in U(\overline R_i)$.
As in the proof of Lemma \ref{Lemma32}, after blowing up 2-curves and 3-points above $X$, by a
morphism $\Phi_0\circ\Phi_1:X_1\rightarrow X$, with associated morphism $f_1=f\circ\Phi_0\circ\Phi_1:X_1
\rightarrow Y$, we have that the following holds:
Suppose that $\overline q\in C$ is a 2-point, and  $\overline p\in f_1^{-1}(\overline q)$ and ${\cal I}_C{\cal O}_{X_1,\overline p}$ is not invertible,  then one of the forms (\ref{eq43}) or (\ref{eq44}) hold at $\overline p\in f_1^{-1}(\overline q)$ (with $d<a$ in (\ref{eq43}), $(g,h)< (a,b)$ and $\overline\gamma=0$ in (\ref{eq44})). We have that $\Phi_0\circ\Phi_1$ is an isomorphism over $f^{-1}(Y-(\Sigma(Y)\cup C))$ and if $C$ is contained in the fundamental locus of $f$, then
$\Phi_0\circ\Phi_1$ is an isomorphism over $f^{-1}(Y-\Sigma(Y))$.

If $\overline q\in C$ is a 1 point then  a form (\ref{eq60}) 
below holds at $\overline p\in f_1^{-1}(\overline q)$ if ${\cal I}_C{\cal O}_{X_1,\overline p}$
is not invertible.
\begin{equation}\label{eq60}
\begin{array}{ll}
u&=x^a\\
v&=y\\
\tilde w_i&=x^dz
\end{array}
\end{equation}
where $d<a$ and $\overline p$ is a 1-point.

As in the proof of Lemma \ref{Lemma32}, the locus of points in $X_1$ where ${\cal I}_C
{\cal O}_{X_1}$ is not invertible is a (possibly
reducible) curve $\overline E$ which makes SNCs with the toroidal structure of $X_1$.
As in the proof of Lemma \ref{Lemma32}, we can construct a sequence of blow ups of
sections over components of $\overline E$, $X_4\rightarrow X_1$,  such that the resulting map $\overline f_4:X_4\rightarrow Y$ factors
through the blow up $\Psi_1:Y_1\rightarrow Y$ of $C$. Let $\overline\Phi:X_4\rightarrow X$
be the composite map. By construction, if $u,v,\tilde w_i$ are our permissible parameters
at $\overline q\in C$ of (\ref{eq215}), and $\overline p\in \overline f_4^{-1}(\overline q)$,
we have permissible parameters at $\overline p$ such that a form 5 (b)  of
Definition \ref{Def200} holds for $u,v,\tilde w_i$ at $\overline p$ if $\overline q$ is a 1-point and $u,v,\tilde w_i$ are super parameters for $\overline f_4$ if $\overline q$ is  a 2-point.

Let $f_4:X_4\rightarrow Y_1$ be the resulting morphism. As in the proof of Lemma
\ref{Lemma32}, we see that $f_4$ is prepared.

Suppose that $\overline q\in C\cap U(\overline R_j)$ for some $j$.
Then $\overline q\in U(\overline R_i)$, since $C$ is prepared of type 4. By the hypothesis
of this lemma, the germ of $C$ at $\overline q$ is contained in $S_{\overline R_j}(\overline q)$.
Since $C\subset E_{\alpha}$, the germ of $C$ at $\overline q$ is $E_{\alpha}\cdot
S_{\overline R_j}(\overline q)$. Thus in the forms of (\ref{eq64}) of Definition \ref{Def65} for
$R$, we have $a_{jk}>0$ for all $j,k\in I_{\overline q}$.

 $\Psi_1^{-1}(\overline q)$ is covered by 2 affine charts. The first chart has
uniformizing parameters $\overline u,\overline v,\overline w_i$ defined by 
\begin{equation}\label{eq216}
u=\overline u,
v=\overline v,
\tilde w_i=\overline u\overline w_i.
\end{equation}
The second chart has uniformizing parameters $u', v', w_i'$ defined by
\begin{equation}\label{eq217}
u=u'w_i',
v=v',
\tilde w_i=w_i'.
\end{equation}

For $j\in I_{\overline q}$ we have a relation
$$
\tilde w_j=\tilde w_i+u^{a_{ij}}v^{b_{ij}}\phi_{ij}(u,v)
$$
with $a_{ij}>0$ (where $\tilde w_j=w_{\overline R_j}(\overline q)$).
As in the proof of Lemma \ref{Lemma32}, it follows that the transform $R^1$ of $R$ for $f_4$ is defined,  $f_4$ is $\tau$-quasi-well prepared, and $D_{X_4}$ is cuspidal for $f_4$ if $D_X$ is cuspidal for $f$.

We now verify that $f_4$ is $\tau$-well prepared. 1 of Definition \ref{Def65} is 
immediate.

Suppose that $C\cap U(\overline R_j)\ne\emptyset$ for some $j$.
Then $C=\overline{E_{\alpha}\cdot \overline R_j(q_{\beta})}\subset \Omega(\overline R_j)$. Since $\Omega(\overline R_j)$ is nonsingular and
makes SNCs with $D_{Y}$, $\Omega(\overline R_j^1)$ is nonsingular, makes SNCs with $D_{Y_1}$,
 contains $U(\overline R_j^1)$, and $\Omega(\overline R_j^1)\cap U(R^1)=U(\overline R_j^1)$. 
 If $C\cap U(\overline R_j)=\emptyset$, then after possibly replacing $\Omega(\overline R_j)$ with a neighborhood
 of $F_j$ in $\Omega(\overline R_j)$ (with the notation of Definition \ref{Def130}, and following our convention on 
 $\Omega(\overline R_j)$ stated after Definition \ref{Def130}), we have that $\Omega(\overline R_j)\cap C
 =\emptyset$.
 Thus 2 of Definition \ref{Def65} holds.

If $\overline q\in C\cap U(\overline R_j)$, then $q_j=U(\overline R_j^1)\cap f_4^{-1}(\overline q)$ is a 2-point
in the chart (\ref{eq216}), and 
\begin{equation}\label{eq218}
w_{\overline R_j^1}(q_j)=\overline w_j=\frac{\tilde w_j}{u}
=\overline w_i+\overline u^{a_{ij}-1}\overline v^{b_{ij}}\phi_{ij}(\overline u,\overline v).
\end{equation}
Thus $q_j\in U(\overline R_i^1)$ if and only if $a_{ij}-1+b_{ij}>0$.
Let $q_i=U(\overline R_i^1)\cap f_4^{-1}(\overline q)$,
$$
I_{q_i}=\{j\mid q_i\in U(\overline R_j^1)\}.
$$
$j\in I_{q_i}$ if and only if $j\in I_{\overline q}$ and $a_{ij}+b_{ij}>1$. Suppose that $\tau>0$, and
$j\in I_{q_i}$. If $\overline R_j(\overline q)$ is defined by an expression
$$
\tilde w_j^{e_j}=\lambda_ju^{a_j}v^{b_j},
$$
then $\overline R_j^1(q_i)$ is defined by the expression 
\begin{equation}\label{eq219}
\overline w_j^{e_j}=\lambda_j\overline u^{a_j-e_j}\overline v^{b_j}.
\end{equation}
From equation (\ref{eq218})  we see that the set  (\ref{eq255}) of Definition
\ref{Def65} corresponding to $q_1$ and $R^1$ is totally ordered, since the set  (\ref{eq255}) 
of Definition \ref{Def65} corresponding to $\overline q$ and $R$ is totally ordered. In particular, we see
that 3 and 4 of Definition \ref{Def65} hold for $f_4$.
We have completed
the verification that $f_4$ is $\tau$-well prepared.

Let $E=\Psi_1^{-1}(C)$ and if $C\cap U(\overline R_j)\ne\emptyset$, let
$\gamma_j=\Omega(\overline R_j^1)\cdot E$. $\gamma_j$ is also
nonsingular and makes SNCs with $D_{Y_1}$. We have that
$\gamma_j=\overline{E\cdot S_{\overline R_j^1}(q_1)}$ for all $q_1\in E\cap
U(\overline R_j^1)$, and $\gamma_j$ is a section of $E$ over $C$. In particular, 1 and 2 of
Definition \ref{Def200} hold for $\gamma_j$. 

Suppose that for some $\overline R_k$ associated to $R$, there exists $\gamma'\subset\Omega(\overline R_k^1)$, with $\gamma'=\overline{E_{\beta}\cdot
S_{\overline R_k^1}(q_{\delta})}$, and $\gamma'\cap \gamma_j\ne\emptyset$.

First suppose that $\gamma'\subset E$. Then $U(\overline R_k^1)\cap E\ne\emptyset$ and
$\gamma'=\gamma_k=\overline{E\cdot S_{\overline R_k^1}(q_2)}$ for all $q_2\in U(\overline
R_k^1)\cap E$. Let $F=E_{\alpha}$ be the component of $D_Y$ containing $C$. Since $C$ is prepared
and $\Psi_1(\gamma_j)=\Psi_1(\gamma_k)=C$, we have that
$C\cap U(\overline R_j)=C\cap U(\overline R_k)$ and $C=\overline{F\cdot S_{\overline R_j}(\overline q)}
=\overline{F\cdot S_{\overline R_k}(\overline q)}$ for $\overline q\in C\cap U(\overline R_j)$.

Suppose that $q_2\in \gamma_j\cap \gamma_k$. Let $\overline q=\Psi_1(q_2)\in C$.
Let $u,v,\tilde w_i$ and $u,v,\tilde w_k$ be the permissible parameters at $\overline q$ of (\ref{eq215}).

Suppose that $\overline q$ is a 1-point.
We have
$$
\tilde w_j=\tilde w_k+u^{c_{jk}}\phi_{jk}(u,v),
$$
where $\phi_{jk}$ is a unit series (or $\phi_{jk}=0$)
by 5 (d) of Definition \ref{Def200}, since $C$ is prepared for $R$ of type 4. We have
that $q_2$ has permissible parameters
$\overline u,\overline v,\overline w_j$ where 
\begin{equation}\label{eq251}
u=\overline u,
v=\overline v,
\tilde w_j=\overline w_j\overline u.
\end{equation}
Define $\overline w_k$ by $\tilde w_k=\overline w_k\overline u$. We have that
$$
\overline w_j=\overline w_k+\overline u^{c_{jk}-1}\phi_{jk}(\overline u,\overline v).
$$
$\overline u=\overline w_j=0$ are local equations of $\gamma_j$ at $q_2$ and
$\overline u=\overline w_k=0$ are local equations of $\gamma_k$ at $q_2$. Thus
$q_2$ is a 1-point with $c_{jk}-1>0$ and $\gamma_j=\gamma_k$.

Now suppose that $q_2\in \gamma_j\cap\gamma_k$ and $\overline q=\Psi_1(q_2)\in C$
is a 2-point. First suppose that $\overline q\not\in C\cap U(\overline R_j)
=C\cap U(\overline R_k)$. We have a relation
$$
\tilde w_j=\tilde w_k+u^{c_{jk}}v^{d_{jk}}\phi_{jk}(u,v)
$$
where $\phi_{ij}$ is a unit series (or $\phi_{jk}=0$) by 5 (e) of Definition \ref{Def200} for $C$. We have $c_{jk}\ge 1$.
$q_2$ has permissible parameters $\overline u,\overline v,\overline w_j$ where
\begin{equation}\label{eq252}
u=\overline u,
v=\overline v,
\tilde w_j=\overline w_j\overline u.
\end{equation}
Define $\overline w_k$ by $\tilde w_k=\overline w_k\overline u$. We then have 
\begin{equation}\label{eq380}
\overline w_j=\overline w_k+\overline u^{c_{jk}-1}\overline v^{d_{jk}}\phi_{jk}(\overline u,\overline v).
\end{equation}
Thus $q_2$ is a 2-point (and $q_2\not\in U( R^1)$).

Finally, suppose that $q_2\in \gamma_j\cap \gamma_k$ and
$\overline q=\Psi_1(q_2)\in C\cap U(\overline R_j)= C\cap U(\overline R_k)$. Then  $q_2\in U(\overline R_j^1)
\cap U(\overline R_k^1)$, and we have equations (\ref{eq252}) and (\ref{eq380}).

Now suppose that $\gamma'=\overline{E_{\beta}\cdot S_{\overline R_k^1}(q_{\delta})}\not\subset E$,
and $\gamma'\cap \gamma_j\ne\emptyset$. Then there exists a component $G$ of $D_Y$ such that
$\overline\gamma=\Psi_1(\gamma')\subset G$, and $E_{\beta}$ is the strict transform of $G$.
Suppose that $q_2\in \gamma'\cap \gamma_j$.
$\overline q=\Psi_1(q_2)\in \overline\gamma\cap C$ implies $\overline q\in U(\overline R_j)\cap U(\overline R_k)$ and
$\overline\gamma=\overline{G\cdot S_{\overline R_k}(\overline q)}$ by 3 of Definition \ref{Def200}
for $R$. Thus $q_2\in U(\overline R_j^1)\cap U(\overline R_k^1)$.

Suppose that $q_2\in \gamma_j$, $\overline q=\Psi_1(q_2)\in C$
and $p\in f_4^{-1}(q_2)$. Let $u,v,\tilde w_j$ be the permissible parameters at $\overline q$ of (\ref{eq215}).
Further suppose that $q_2$
is a 1-point. Let $\overline u,\overline v,\overline w_j$ be the permissible parameters at
$q_2$ of (\ref{eq251}). $u,v,\tilde w_j$ satisfy a form 5 (b) of Definition \ref{Def200}
at $p$. Substituting (\ref{eq251}) into these forms, we see that
$\overline u,\overline v,\overline w_j$ satisfy a form of 5 (b) of Definition
\ref{Def200} at $p$. Now suppose that $q_2\in \gamma_j$ is a 2-point, but $q_2\not\in
U(\overline R_j^1)$. Let $\overline u, \overline v, \overline w_j$ be the permissible parameters
at $q_2$ of (\ref{eq252}). $u,v,\tilde w_j$ satisfy a form 5 (c) of Definition \ref{Def200}
at $p$. Substituting (\ref{eq252}) into these forms, we see that $\overline u,\overline v, \overline w_j$ satisfy a form 5 (c) of Definition \ref{Def200} at $p$.
Thus 5 of Definition \ref{Def200} holds for $\gamma_j$.

We have seen that the curves $\gamma_j$ only fail to be prepared of type 4 for $f_4$ at a
finite set of 2-points $T_1\subset E$, where condition 3 of Definition \ref{Def200} fails.
If $q\in T_1$, then $q\not\in U(R^1)$ and   there exist $\gamma_j$ and $\gamma_k$ such that 
$\gamma_j\ne\gamma_k$ and $q\in\gamma_j\cap\gamma_k$. 
 For $q\in T_1$, let
$$
J_q=\{j\mid q\in \gamma_j=\overline{E\cdot S_{\overline R_j^1}(q_{\epsilon})}\text{ for some }q_{\epsilon}
\in U(\overline R_j^1)\}.
$$
Observe that:
\begin{enumerate}
\item[1.] For $q\in T_1$, $j\in J_q$ and $p\in f_4^{-1}(q)$, 
the permissible parameters $\overline u,\overline v,\overline w_j$ at $q$ defined by  (\ref{eq216}) have a form 5 (c) of Definition \ref{Def200}.
\item[2.] For $i,j\in J_q$ there exists a relation of the form 5 (e) of Definition
\ref{Def200} for $\overline w_i$ and $\overline w_j$, and the set $\{(a_{ij},b_{ij})\}$ is totally ordered.
\end{enumerate}
In particular, the points in $T_1$ are prepared for $R^1$ of type 2 of Definition \ref{Def66}.
Let $\Psi_2:Y_2\rightarrow Y_1$
be the blow up of $T_1$, and let
$$
\begin{array}{rll}
X_5&\stackrel{f_5}{\rightarrow}&Y_2\\
\downarrow&&\downarrow\Psi_2\\
X_4&\stackrel{f_4}{\rightarrow}&Y_1
\end{array}
$$
be the $\tau$-well prepared diagram of Lemma \ref{Lemma171}.
Let
$$
T_2=\left\{
\begin{array}{l}  q\in \Psi_2^{-1}(T_1) \text{ such that }q\in \gamma_i^2\cap \gamma_k^2\\
\text{ where }\gamma_i^2,\gamma_j^2\text{ are the strict transforms of some}\\
\gamma_i,\gamma_k\text{ such that }\gamma_i\ne\gamma_k\text{ and }\gamma_i\cap \gamma_k\ne \emptyset
\end{array}\right\}.
$$
The points of $T_2$ must again be prepared for the transform $R^2$ of $R$ on $X_5$ of type 2 of Definition \ref{Def66}, and the points of $T_2$ satisfy the corresponding statements 1 and 2 above that the points of $T_1$ and $J_q$ satisfy.

We can iterate this process a finite number of times to produce a $\tau$-well prepared
diagram
$$
\begin{array}{rll}
\overline X_m&\stackrel{\overline f_m}{\rightarrow}&\overline Y_m\\
\downarrow&&\downarrow\\
\overline X_1=X_4&\stackrel{f_4}{\rightarrow}&\overline Y_1=Y_1\\
\downarrow&&\downarrow\\
X&\stackrel{f}{\rightarrow}&Y
\end{array}
$$
of the form of (\ref{eq220}) such that the strict transform of the $\gamma_i$ are
disjoint on $\overline Y_m$ above $T_1$.
It follows that the strict transforms of the $\gamma_i$ are prepared of type 4 for the transform $R^m$ of $R$ on $\overline X_m$.

To show that $\overline f_m$ is $\tau$-very-well prepared, it only remains to verify that the
strict transform $\gamma'$ of a curve $\gamma=\overline{E_{\beta}\cdot S_{\overline R_i}(q)}$
on $Y$ (with $\gamma\ne C$ and $\gamma\cap C\ne\emptyset$) is prepared on $\overline Y_m$. 
Since $\gamma$ is prepared of type 4, we have that $\gamma\cap C\subset U(R_i)$.  By our previous
analysis, we then know that $\gamma'\cap E\subset U(R^1)$, so that $\overline Y_m\rightarrow Y_1$ is an isomorphism in a neighborhood of $\gamma'$. It thus suffices to check
that $\gamma'$ is prepared of type 4 on $Y_1$. This follows by a local
analysis.

\end{pf}

\begin{Remark}\label{Remark293}
\begin{enumerate}
\item[1.] Suppose that  $f:X\rightarrow Y$ is $\tau$-quasi-well prepared and $C\subset D_Y$ is a nonsingular (integral) curve which makes SNCs with $D_Y$ and contains a 1-point such that 
\begin{enumerate}
\item $q\in C\cap U(\overline R_i)$ for some 2-point pre-relation $\overline R_i$ associated to $R$ implies the (formal) germ of $C$ at $q$ is contained in $S_{\overline R_i}(q)$, and 
\item Suppose that $q\in C-U(R)$. Then there exist permissible parameters
$u,v,w$ at $q$ such that $u=w=0$ are local equations of $C$
at $q$, and
\begin{enumerate}
\item If $q$ is a 1-point, $C$ is not a component of the fundamental locus of $f$, and $p\in f^{-1}(q)$, then there exists a relation of one of the following forms
 for  $u,v,w$ at $p$.

\begin{enumerate}
\item
$p$ a 1-point 
\begin{equation}\label{eq59*}
\begin{array}{ll}
u&=x^a\\
v&=y\\
w&=x^c\gamma(x,y)+x^dz
\end{array}
\end{equation}
where $\gamma$ is a unit series  (or zero),
\item
$p$ a 2-point 
\begin{equation}\label{eq61*}
\begin{array}{ll}
u&=(x^ay^b)^k\\
v&=z\\
w&=(x^ay^b)^l\gamma(x^ay^b,z)+x^cy^d
\end{array}
\end{equation}
where $\gamma$ is a unit (or zero) and $ad-bc\ne 0$.
\end{enumerate}
\item If $q$ is a 2-point,  $u,v,w$ are super parameters for $f$ at $q$.
\end{enumerate}
\end{enumerate}
Then there exists a $\tau$-quasi-well prepared diagram
$$
\begin{array}{rll}
\overline X_1&\stackrel{\overline f_1}{\rightarrow}&\overline Y_1\\
\overline\Phi_1\downarrow&&\downarrow\overline\Psi_1\\
X&\stackrel{f}{\rightarrow} &Y
\end{array}
$$
where $\overline\Psi_1$ is the blow up of $C$. $\overline\Phi_1$ is an isomorphism over $f^{-1}(Y-(\Sigma(Y)\cup C))$.
If $C$ is contained in the fundamental locus of $f$, then $\overline\Phi_1$ is an isomorphism over $f^{-1}(Y-\Sigma(Y))$.
\item[2.] Further suppose that 
$f:X\rightarrow Y$ is $\tau$-well prepared, and if $\gamma=\overline{E\cdot R_k(q_{\alpha})}$ is prepared for
$R$ of type 4, then either $C=\gamma$ or $q\in C\cap \gamma$ implies $q\in U(\overline R_k)$ and the germ of $C$ at $q$ is contained in $S_{R_k}(q)$. 
Then
there exists a $\tau$-well prepared diagram
$$
\begin{array}{rll}
X_1&\stackrel{f_1}{\rightarrow}&Y_1\\
\Phi\downarrow&&\downarrow\Psi\\
X&\stackrel{f}{\rightarrow}&Y
\end{array}
$$
 where $\Psi$ is the blow up $\overline\Psi_1$ of
$C$, possibly followed by  blow ups of 2-points which are prepared for
the transform of $R$ (of type 2 of Definition \ref{Def66}) if $C$ is prepared of type 4 for $R$, such that
\begin{enumerate}
\item If $\gamma\subset Y$ is prepared for $f$, (and $\gamma\ne C$) then the strict transform of $\gamma$
is prepared for $f_1$.
\item If $C\subset Y$ is prepared for $f$ (of type 4) and $q\in U(\overline R_i)\cap C$ for some $i$, then $\overline{E\cdot S_{\overline R_i^1}(q')}$ is prepared for $f_1$ (of type 4) for all $q'\in E\cap U(\overline R_i^1)$, where $E$ is the component of $D_{Y_1}$ dominating $C$.
\item Suppose that $F$ is a component of $D_{Y_1}$ such that $\Psi(F)$
is a point. If $q\in F\cap U(\overline R_i^1)$ for some $i$, then
$\overline{F\cdot S_{\overline R_i^1}(q)}$ is a prepared curve of
type 4 for $R^1$.
\item If $D_X$ is cuspidal for $f$, then $D_{X_1}$ is cuspidal for $f_1$.
\item $\overline\Phi_1$ is an isomorphism over $f^{-1}(Y-(\Sigma(Y)\cup C))$.
If $C$ is contained in the fundamental locus of $f$, then $\overline\Phi_1$ is an isomorphism over $f^{-1}(Y-\Sigma(Y))$.
\end{enumerate}
\end{enumerate}
\end{Remark}

The proof of Remark \ref{Remark293} is a variation of the proof of Lemma
\ref{Lemma67}, using Lemma \ref{Lemma143} if $C$ is a component of the fundamental locus of $f$.

\section{Existence of a $\tau$-very-well prepared morphism}
Suppose that $f:X\rightarrow Y$ is a birational morphism of nonsingular projective 3-folds, with toroidal structures
$D_Y$ and $D_X=f^{-1}(D_Y)$.

\begin{Theorem}\label{Theorem169} Suppose that $f:X\rightarrow Y$ is prepared.
Let $\tau=\tau_f(X)$. Suppose that $\tau\ge 0$ and if $p\in X$ is a 3-point and
$\tau_f(p)=\tau$ then $f(p)$ is a 2-point on $Y$. Further suppose that $D_X$ is cuspidal for $f$.
 Then there exists a commutative diagram
$$
\begin{array}{rll}
X_1&\stackrel{f_1}{\rightarrow}&Y_1\\
\Phi\downarrow&&\downarrow\Psi\\
X&\stackrel{f}{\rightarrow}&Y
\end{array}
$$
where $\Phi$, $\Psi$ are products of blow ups of 2-curves,
 and there exists a 2-point relation $R^1$ for $f_1$
such that $f_1$ is $\tau$-quasi-well
prepared with 2-point relation $R^1$. Further, $D_{X_1}$ is cuspidal for $f_1$.
\end{Theorem}

\begin{pf} Let
$$
\begin{array}{rll}
X_1&\stackrel{f_1}{\rightarrow}&Y_1\\
\downarrow&&\downarrow\\
X&\stackrel{f}{\rightarrow}&Y
\end{array}
$$
be a diagram
 satisfying the conclusions of Theorem \ref{Theorem5}.

Let $T$ be the 3-points $p\in X_1$ with $\tau_{f_1}(p)=\tau$. Let $I$
be an index set of $T$, and let $U=\{f_1(p)\mid p\in T\}$.
We necessarily have that $U$ consists of 2-points. Suppose that $q\in U$ and $p_i\in T\cap f_1^{-1}(q)$. We will define a 2-point pre-relation $\overline R^1_{p_i}$ on $Y_1$ which has the property that $U(\overline R^1_{p_i})=\{q\}$.  
Let
$u,v,w_i$  (with $w_i\in\hat{\cal O}_{Y_1,q}$)
be the permissible parameters at $q$ of 4 of Theorem \ref{Theorem5}, which have the property  that 4 (a) or (b) of Theorem \ref{Theorem5} holds for permissible parameters $x,y,z$ at $p_i$.
If $\tau>0$, define $\overline R^1_{p_i}(p_i)$ from the expression
$$
w_i^{e_i}=\overline\lambda_i u^{a_i}v^{b_i}
$$
of 4 (a) of Theorem \ref{Theorem5} where
$\overline\lambda_i=\gamma(0,0,0)^{e_i}$. We have $\text{gcd}(a_i,b_i,e_i)=1$ and $e_i>1$. If $\tau=0$, then an expression as a monomial form 4 (b) of Theorem \ref{Theorem5} holds at $p_i$, and we define $\overline R^1_{p_i}(p_i)$ by  $a_i=b_i=-\infty$ and the expression $w_i=0$. 
We now define a primitive 2-point relation $R^1_{p_i}$ on $X_1$ by $T(R^1_{p_i})=\{p_i\}$ and $R_{p_i}^1(p_i)=\overline R_{p_i}(p_i)$. We can thus define a 2-point relation $R^1$ on $X_1$ with $T(R^1)=T$, $U(R^1)=U$,
and where the $ R^1_{p_i}$ defined above are the primitive 2-point relations
associated to $R^1$.
It follows from Theorem \ref{Theorem5} that $f_1$ is $\tau$-quasi-well prepared for $R^1$, and $D_{X_1}$ is cuspidal for $f_1$.
\end{pf}

\begin{Lemma}\label{Lemma139} Suppose that $Y$ is a nonsingular projective 3-fold with
toroidal structure $D_Y$ and 2-point pre-relations $R_1,\ldots ,R_n$  and $\gamma\subset D_Y$ is a reduced (but possibly not
irreducible) curve such that $\gamma$ has no components which are 2-curves and $\gamma$ is nonsingular at 1-points.
If the $R_i$ are algebraic, let $\Omega(R_i)$ be the locally closed subset of Definition \ref{Def199}, and assume
that $\overline\gamma\cap\Omega(R_i)$ is a union of 2-points, if $\overline\gamma$ is an irreducible component of $\gamma$
such that $\overline\gamma\cap\Omega(R_i)$ is a finite set.
Consider the following algorithm.
\begin{enumerate}
\item[1.] Perform a sequence of blow ups of 2-curves $\overline Y_1\rightarrow Y$ so that
the strict transform $\overline \gamma_1$ of $\gamma$ on $\overline Y_1$ contains no
3-points.
\item[2.] Perform an arbitrary sequence of blow ups $\overline Y_2\rightarrow \overline Y_1$
of 2-curves.  Let $\overline \gamma_2$ be the strict transform of $\gamma$ on $\overline Y_2$.
Let $\overline R_1^2,\ldots,\overline R_n^2$ be
the transforms of $R_1,\ldots, R_n$ on $\overline Y_2$.
\item[3.]  Blow up all 2-points $q\in\overline\gamma_2$ such that $\overline\gamma_2$ does not
make SNCs with $D_Y$ at $q$ or $q\in U(\overline R_i^2)$ for some $i$ and
 the germ of $\overline\gamma_2$ at $q$ is not
contained in $S_{\overline R_i^2}(q)$, or $q\in \Omega(\overline R_i^2)$ for some $i$, and the germ of 
$\overline\gamma_2$ at $q$ is not contained in the germ of $\Omega(\overline R_i^2)$ at $q$.
\item[4.] Iterate steps 1. - 3.
\end{enumerate}

Then, after finitely many iterations, we produce a sequence of admissible  blow ups
$\Psi_1:Y_1\rightarrow Y$
 such that
the strict transform $\gamma_1$ of $\gamma$ on $Y_1$ is nonsingular and makes SNCs with $D_{Y_1}$.
Further, if $R_{i}^1$ is the transform of $R_i$ on $Y_1$ for $1\le i\le n$ and $q\in U(R_{i}^1)
\cap \gamma_1$ for some $i$, then the germ of $\gamma_1$ at $q$ is contained in $S_{R_i^1}(q)$.
If $q\in \Omega(R_i^1)\cap \gamma_1$ for some $i$, then the germ of $\gamma_1$ at $q$ is contained in the
germ of $\Omega(R_i^1)$ at $q$.
\end{Lemma}

\begin{pf} This follows from embedded resolution of plane curve singularities (cf. Section 3.4 and Exercise 3.13 \cite{C3}).
\end{pf}

\begin{Lemma}\label{Lemma323} Suppose that $f:X\rightarrow Y$ is prepared, and
$C$ is a reduced (but possibly not irreducible) curve, consisting of components of the fundamental locus of $f$ which contain a 1-point. Then 
\begin{enumerate}
\item[1.] $C$ is nonsingular at 1-points of $Y$.
\item[2.] If $f$ is $\tau$-well prepared, and $\gamma\subset Y$ is prepared for $R$
of type 4  then either
$\gamma$ is a component of $C$ or $C\cap \gamma$ contains no 1-points of $Y$.
\end{enumerate}
\end{Lemma}

\begin{pf}
1 follows from Lemma \ref{Lemma141}. 2 follows from Lemma \ref{Lemma141} and 5 (b) of Definition \ref{Def200}.
\end{pf}

\begin{Theorem}\label{Theorem4}  Suppose that $\tau\ge0$, $f:X\rightarrow Y$
is $\tau$-quasi-well  prepared (or $\tau$-well prepared) with
 2-point relation $R$, and $C\subset Y$ is a reduced (but possibly not irreducible) curve
consisting of components of the fundamental locus of $f$ which contain a 1-point of $Y$. Assume that $D_X$ is cuspidal for $f$.
 Then there exists a sequence of blow ups of 2-curves and 2-points $\Psi_1:Y_1\rightarrow Y$ and a
$\tau$-quasi-well prepared (or $\tau$-well prepared) diagram of $R$ and $\Psi_1$ 
\begin{equation}\label{eq140}
\begin{array}{rll}
X_1&\stackrel{f_1}{\rightarrow}&Y_1\\
\Phi_1\downarrow&&\Psi_1\downarrow\\
X&\stackrel{f}{\rightarrow}&Y
\end{array}
\end{equation}
such that
\begin{enumerate}
\item[1.]  The strict transform $\overline C$ of $C$
on $Y_1$ is nonsingular and
makes SNCs with $D_{Y_1}$.
\item[2.] If $C_j$ is an irreducible component of
 $\overline C$ and $q\in U(\overline R_i^1)$ for some $\overline R_i^1$ associated to $R^1$ is such that $q\in C_j$ then the germ of $C_j$ at $q$ is contained in $S_{\overline R_i^1}(q)$.
\item[3.] If $f$ is $\tau$-well prepared, $C_j$ is an irreducible  component of $\overline C$
and $\gamma=\overline{E\cdot R_k^1(q_{\alpha})}$ is prepared for $R^1$ of type 4,  then either $C_j=\gamma$ or
$q\in C_j\cap \gamma$ implies $q\in U(\overline R_k^1)$ (and thus the germ of $C_j$ at $q$ is contained in $S_{R_k}(q)$).
\item[4.]  $D_{X_1}$ is cuspidal for $f_1$.
\item[5.] $\Phi_1$ is an isomorphism over $f^{-1}(Y-\Sigma(Y))$.
\end{enumerate}
\end{Theorem}

\begin{pf}  $C\subset D_Y$ and
$C$ is nonsingular  at 1-points   by Lemma \ref{Lemma323}.
We now  follow the algorithm of Lemma \ref{Lemma139}.  If $f$ is $\tau$-well prepared, then
by our convention on the $\Omega(\overline R_i)$ (following Definition \ref{Def200}), we may replace $\Omega(\overline R_i)$ with an open subset,
so that if $C_j$ is a component of $C$, then either $C_j$ is disjoint from $\Omega(\overline R_i)$, or intersects $\Omega(\overline R_i)$ in a
finite set of 2-points which are in  $U(\overline R_i)$,
or on a curve $\gamma\subset\Omega(\overline R_i)$ which is
prepared for $R$ of type 4. By Lemma \ref{Lemma323}, if $C_j$ is a component of $C$ and $\gamma\subset\Omega(\overline R_i)$ is
prepared of type 4 for $R$, such that $C_j\ne\gamma$, then $\gamma\cap C_j$ is a set of 2-points.

We first construct a sequence of blow ups of 2-curves, $\Psi_1:Y_1\rightarrow Y$ so that the strict transform $C_1$ of $C$ on $Y_1$ contains no 3-points. Let
$$
\begin{array}{rll}
X_1&\stackrel{f_1}{\rightarrow}&Y_1\\
\Phi_1\downarrow&&\downarrow\Psi_1\\
X&\stackrel{f}{\rightarrow}&Y
\end{array}
$$
be a $\tau$-quasi-well prepared (or $\tau$-well prepared) diagram of $R$ and $\Psi_1$ such that $\Phi_1$ is an isomorphism above $f^{-1}(Y-\Sigma(Y))$ (this exists by Lemma \ref{Lemma31}).

 If $f$ (and thus also $f_1$) is $\tau$-well prepared, and $q$ is a 2-point on a curve $\gamma\subset Y_1$
 which is prepared for $R^1$ of  type 4, then $q$ is prepared for $f_1$ (of type 1 or 2 of Definition \ref{Def66}).

Consider the 2-points 
$$
\Omega=\left\{\begin{array}{l}
\text{2-points }q\in C_1\text{ which are not in $U(R^1)$
and are not on a curve }\\
\text{which is prepared of type 4 for $f_1$ (if $f$ is $\tau$-well prepared)}\end{array}\right\}.
$$
  By Lemma \ref{Lemma358} and Remark \ref{Remark424}, there exists a $\tau$-quasi-well
prepared (or $\tau$-well prepared) diagram
$$
\begin{array}{rll}
 X_2&\stackrel{ f_2}{\rightarrow}& Y_2\\
\Phi_2\downarrow&&\downarrow\Psi_2\\
X_1&\stackrel{f_1}{\rightarrow}& Y_1
\end{array}
$$
where $\Phi_2$ is a product of blow ups of 2-curves and 3-points and  $\Psi_2$ is a  product of blow ups of 2-curves such that 
$\Phi_2$ is an isomorphism above $f^{-1}(Y_1-\Sigma(Y_1))$ and at all 2-points $\overline q_1\in\Psi_2^{-1}(\Omega)$,
there exist
permissible parameters $\overline u,\overline v,w$, at $\overline q_1$,
such that $\overline u,\overline v,w$ are super parameters for $f_2$ at $\overline q_1$, and thus $\overline q_1$ is prepared for $R^2$ with respect to $\overline u,
\overline v,w$.
In particular, all 2-points on the strict transform $ C_2$ of $C$ on $Y_2$
are prepared for $R^2$ (of type 1 or 2 of Definition \ref{Def66}). The map $Y_2\rightarrow Y_1$
is Step 2 of the algorithm of Lemma \ref{Lemma139}.

Now perform Step 3 of the algorithm, blowing up all 2-points $q$ on the strict transform
$ C_2$ of $C$ on $Y_2$ where 
$ C_2$ is singular, or
$ C_2$ does not make SNCs with $D_{ Y_2}$ at $q$,  or $q\in\Omega(\overline R_i^2)$ for some 
$\overline R_i^2$ associated to $R^2$
 and an irreducible component
$C'$ of $C_2$ contains $q$ but
 the  germ of $C'$ at $q$ is not contained in the germ of $\Omega(\overline R_i^2)$ at $q$.
  Let $\Psi_3:Y_3\rightarrow Y_2$ be the resulting map. By Lemmas \ref{Lemma32} and \ref{Lemma171},
there exists a $\tau$-quasi-well prepared (or $\tau$-well prepared) diagram
$$
\begin{array}{rll}
X_3&\stackrel{f_3}{\rightarrow}&Y_3\\
\Phi_3\downarrow&&\downarrow\Psi_3\\
X_2&\stackrel{f_2}{\rightarrow}&Y_2
\end{array}
$$
such that $\Phi_3$ is an isomorphism above $f_2^{-1}(Y_2-\Sigma(Y_2))$.

Now by Step 4 of the algorithm of Lemma \ref{Lemma139}, we can iterate this process to
construct a $\tau$-quasi-well prepared (or $\tau$-well prepared) diagram (\ref{eq140})
such that the conclusions of Theorem \ref{Theorem4} hold. By our construction, and Lemmas \ref{Lemma31}, \ref{Lemma32} 
and \ref{Lemma171} and Remark \ref{Remark424}, $D_{X_1}$ is cuspidal for $f_1$. 
\end{pf}

\begin{Theorem}\label{Theorem6}  Suppose that
$f:X\rightarrow Y$, $C$ and $f_1:X_1\rightarrow Y_1$ are as in the assumptions and
conclusions of Theorem \ref{Theorem4}. Then there exists a sequence of blow ups of 2-curves
$\Psi_2:Y_2\rightarrow Y_1$, and a $\tau$-quasi-well prepared (or $\tau$-well prepared)
diagram 
$$
\begin{array}{rll}
X_2&\stackrel{f_2}{\rightarrow}&Y_2\\
\Phi_2\downarrow&&\Psi_2\downarrow\\
X_1&\stackrel{f_1}{\rightarrow}& Y_1
\end{array}
$$
 such that
\begin{enumerate}
\item[1.] The conclusions of 1, 2 and 3  of Theorem \ref{Theorem4} hold for the strict transform
$C_2$ of $C$ on $Y_2$.
\item[2.] $D_{X_2}$ is cuspidal for $f_2$.
\item[3.] The components $\gamma$ of $C_2$  are
permissible centers (or *-permissible centers if $\gamma$ is prepared of type 4) for $R^2$.
\item[4.] $\Phi_2$ is an isomorphism over $f_1^{-1}(Y_1-\Sigma(Y_1))$.
\end{enumerate}

In the resulting $\tau$-quasi-well prepared (or $\tau$-well prepared) diagram 
$$
\begin{array}{rll}
X_3&\stackrel{f_3}{\rightarrow}&Y_3\\
\Phi_3\downarrow&&\downarrow\Psi_3\\
X_2&\stackrel{f_2}{\rightarrow}&Y_2,
\end{array}
$$
$\Psi_3$ is the blow up of $C_2$, possibly followed by blow ups of 2-points which are
prepared for the transform of $R$ if $f$ is $\tau$-well prepared, and $C$ contains a component which is prepared of type 4 for $R$. Further,  $D_{X_3}$ is cuspidal for $f_3$ and $\Phi_3$ is an isomorphism above $f_2^{-1}(Y_2-\Sigma(Y_2))$.
\end{Theorem}

\begin{pf}
For all $q\in \overline C$ a 2-point such that $q\not\in U(R^1)$,
 let $u_q,v_q,w_q$ be permissible
parameters at $q$ such that $u_q=w_q=0$ are local equations for $\overline C$ at $q$.
By Lemma \ref{Lemma358} and Remark \ref{Remark424}, there exists a $\tau$-quasi-well prepared (or $\tau$-well prepared)
diagram
$$
\begin{array}{rll}
X_2&\stackrel{f_2}{\rightarrow}&Y_2\\
\Phi_2\downarrow&&\Psi_2\downarrow\\
X_1&\stackrel{f_1}{\rightarrow}& Y_1
\end{array}
$$
where $\Phi_2$ is a product of blow ups of 2-curves and 3-points, $\Psi_2$ is a  product of blow ups of 2-curves, such that 
$\Phi_2$ is an isomorphism over $f_1^{-1}(Y_1-\Sigma(Y_1))$ and for all 2-points $q\in\overline C$, 
 which are not in $U(R^1)$,
 at 2-points $q_1\in \Psi_2^{-1}(q)$ we have permissible
parameters $\overline u,\overline v,w$ at $q_1$ such that $q_1$ is prepared (of type 2 of Definition \ref{Def66}) for $R^2$ with respect
to the parameters $\overline u,\overline v,w$. In particular, this is true for the point
$q_1\in\Psi_2^{-1}(q)$ on the strict transform $C_2$ of $C$ on $Y_2$. At this $q_1$, $\overline u,
\overline v,w$  satisfy
$$
u_q=\overline u\overline v^n, v_q=\overline v, w_q=w
$$
for some $n$,  $\overline u=w=0$ are local equations of $C_2$ at $q_1$, and $\overline u,\overline v,w$ are
super parameters at $q_1$.

Thus, the hypothesis of Remark \ref{Remark293}
are satisfied, and the conclusions of Theorem \ref{Theorem6} hold.
\end{pf}

\begin{Lemma}\label{Lemma348}  Suppose that $\tau\ge0$, $f:X\rightarrow Y$
is $\tau$-quasi-well prepared with  2-point relation $R$ and $D_X$ is cuspidal for $f$.
Further suppose there exists a $\tau$-quasi-well prepared diagram
\begin{equation}\label{eq336}
\begin{array}{rll}
\tilde X&\stackrel{\tilde f}{\rightarrow}&\tilde Y\\
\tilde\Phi\downarrow&&\downarrow\tilde\Psi\\
X&\stackrel{f}{\rightarrow}&Y,
\end{array}
\end{equation}
where $\tilde R$ is the transform of $R$ on $\tilde X$, such that if $q_1\in U(\tilde R)$ is on a component $E$ of $D_{\tilde Y}$ such that $\tilde\Psi(E)$ is not a point, then
$T(\tilde R)\cap \tilde f^{-1}(q_1)=\emptyset$. Further suppose that    $D_{\tilde X}$ is cuspidal for $\tilde f$.
Then there exists a commutative diagram
$$
\begin{array}{rll}
X_1&\stackrel{f_1}{\rightarrow}&Y_1\\
\Phi\downarrow&&\downarrow\Psi\\
\tilde X&\stackrel{\tilde f}{\rightarrow}&\tilde Y
\end{array}
$$
such that $\Phi$, $\Psi$ are products of blow ups of possible centers, and
$f_1$ is $\tau$-quasi-well prepared with  2-point relation
$R^1$ and pre-algebraic structure. ($R^1$ will in general not be the transform of $\tilde R$.)
Further,  $D_{X_1}$ is cuspidal for $f_1$.
\end{Lemma}

\begin{pf}
Given a diagram (\ref{eq336}), we will define a new 2-point relation $\tilde R'$
on $\tilde X$ for $\tilde f$. This is accomplished as follows. Suppose that
$q_1\in U(\tilde R)$ is such that $\tilde f^{-1}(q_1)\cap T(\tilde R)\ne\emptyset$. Let
$$
J_{q_1}=\{i\mid T(\tilde R_i)\cap\tilde f^{-1}(q_1)\ne \emptyset\}.
$$
For $j\in J_{q_1}$, let $q=\tilde\Psi(q_1)$,
$$
u=u_{\overline R_j}(q), v=v_{\overline R_j}(q), w_j=w_{\overline R_j}(q).
$$
 Let
$$
u_1=u_{\overline{\tilde R_j}}(q_1),
v_1=v_{\overline{\tilde R_j}}(q_1),
w_{j,1}=w_{\overline{\tilde R_j}}(q_1).
$$
Since $\tilde\Psi$ is a composition of admissible blow ups for the transforms of the pre-relations $\overline R_i$ on $Y$,
we have relations 
\begin{equation}\label{eq381}
u=u_1^{\tilde a}v_1^{\tilde b},
v=u_1^{\tilde c}v_1^{\tilde d},
w_j=u_1^{\tilde e}v_1^{\tilde f}w_{j,1}
\end{equation}
with $\tilde a\tilde d-\tilde b\tilde c=\pm 1$.
Since if $q_1$ is on a component $E$ of $D_{\tilde Y}$ we must have $\tilde\Psi(E)$ is a point,  we have $\tilde a,\tilde b,\tilde c,\tilde d$ all nonzero.

For $p_1\in T(\tilde R_j)\cap \tilde f^{-1}(q_1)$, there exist permissible
parameters $x_1,y_1,z_1$ for $u_1,v_1,w_{j,1}$ (and $u,v,w_j$) at $p_1$ and an expression
\begin{equation}\label{eq338}
\begin{array}{ll}
u&=x_1^ay_1^bz_1^c\\
v&=x_1^dy_1^ez_1^f\\
w_{j}&=x_1^gy_1^hz_1^i\gamma_j
\end{array}
\end{equation}
where $\gamma_j$ is a unit series in $\hat{\cal O}_{\tilde X,p_1}$.
Since $\tilde a\tilde d-\tilde b\tilde c=\pm 1$, there exist (after possibly interchanging $u$ and $v$)
$m>0$ and a factorization of the expression of $u$ and $v$ in
(\ref{eq381}) by the three successive substitutions: 
\begin{equation}\label{eq337}
\begin{array}{l}
u=\overline u, v=\overline u^m\overline v\\
\overline u=\tilde u\tilde v, \overline v=\tilde v\\
\tilde u=u_1^{\overline a}v_1^{\overline b}, \tilde v=u_1^{\overline c}v_1^{\overline d}
\end{array}
\end{equation}
for some $\overline a,\overline b,\overline c,\overline d\in {\bf N}$, with $\overline a\overline d-\overline b\overline c=\pm1$.
substituting (\ref{eq337}) into (\ref{eq338}), we see that
$$
(a,b,c)>(d,e,f)-m(a,b,c)>0.
$$
Thus 
\begin{equation}\label{eq339}
m(a,b,c)<(d,e,f)<(m+1)(a,b,c).
\end{equation}
 If there does not exist a natural number $r$ such that $w_j|u^r$ and $w_j|v^r$ in $\hat{\cal O}_{\tilde X,p_1}$, then  there exists
an expression by (\ref{eq339}) (after possibly interchanging $u$ and $v$, and $x_1,y_1,z_1$) 
\begin{equation}\label{eq349}
\begin{array}{ll}
u&=x_1^ay_1^b\\
v&=x_1^dy_1^e\\
w_j&=x_1^gy_1^hz_1^i\gamma_j
\end{array}
\end{equation}
with $i\ne 0$. But we now obtain a contradiction, since $uv=0$ is a local equation of $D_{\tilde X}$ at $p_1$.

In conclusion, there exists a natural number $r$ such that $w_j$ divides $u^r$ and $v^r$
in $\hat{\cal O}_{\tilde X,p_1}$.

Set $\eta(p_1)=\text{max}\{2r,\tilde e,\tilde f\}$, and $\eta=\text{max}\{\eta(p_1)\mid p_1\in T(\tilde R)\cap \tilde f^{-1}(q_1)\}$.
Fix $j\in J_{q_1}$ and $p_1\in T(\tilde R_j)\cap \tilde f^{-1}(q_1)$. There exists $\sigma(u,v,w_j)\in {\bold k}[[u,v,w_j]]=\hat{\cal O}_{Y,q}$ such that
the order of the series $\sigma$ is greater than $\eta$ and $w_j+\sigma\in{\cal O}_{Y,q}$. Let 
$$
w_{p_1}^*=w_j+\sigma(u,v,w_j).
$$
In $\hat{\cal O}_{\tilde X,p_1}$, we have
$$
w_{p_1}^*=w_j\gamma_{p_1}''=x_1^gy_1^hz_1^i\gamma_{p_1}'
$$
where $\gamma_{p_1}', \gamma_{p_1}''$ are unit series. Set
$$
w_{p_1}=\frac{w_{p_1}^*}{u_1^{\tilde e}v_1^{\tilde f}}=w_{j,1}+\frac{\sigma(u_1^{\tilde a}v_1^{\tilde b},u_1^{\tilde c}v_1^{\tilde d},
u_1^{\tilde e}v_1^{\tilde f}w_{j,1})}{u_1^{\tilde e}v_1^{\tilde f}}
\in\hat{\cal O}_{\tilde Y,q_1}\cap {\bold k}(Y)={\cal O}_{\tilde Y,q_1}.
$$
We further have 
$$
w_{p_1}=w_{j,1}\gamma_{p_1}''.
$$
For $k\in J_{q_1}$, there exists $\lambda_{jk}(u,v)\in {\bold k}[[u,v]]$ such that $w_k=w_j+\lambda_{jk}(u,v)$.
Suppose that $p_2\in T(\tilde R_k)\cap \tilde f^{-1}(q_1)$. Write
$$
\lambda_{jk}(u,v)=\alpha_{p_2}(u,v)+h_{p_2}(u,v)
$$
 where $\alpha_{p_2}(u,v)$ is a polynomial, and $h_{p_2}(u,v)$  is a series of order greater than $\eta$. Set
$$
w_{p_2}^*=w_j+\sigma(u,v,w_j)+\lambda_{jk}(u,v)-h_{p_2}(u,v)\in{\cal O}_{Y,q}
$$
$$
\begin{array}{ll}
w_{p_2}^*&=w_k+\sigma(u,v,w_k-\lambda_{jk}(u,v))-h_{p_2}(u,v)\\
&=w_k+\overline\sigma_{p_2}(u,v,w_k)
\end{array}
$$
where $\overline\sigma_{p_2}$ is a series of order greater than $\eta$. Set
$$
w_{p_2}=\frac{w_{p_2}^*}{u_1^{\tilde e}v_1^{\tilde f}}.
$$
From (\ref{eq381}) we see that $u_1,v_1,w_{p_2}$ are permissible parameters at $q_1$ and $w_{p_2}\in{\cal O}_{\tilde Y,q_1}$.
We further have that  $w_{p_2}=w_{k,1}\gamma_{p_2}''$ for some unit series $\gamma_{p_2}''\in \hat{\cal O}_{\tilde X,p_2}$.

We further have
$$
w_{p_2}-w_{p_1}=\frac{w^*_{p_2}-w^*_{p_1}}{u_1^{\tilde e}v_1^{\tilde f}}
=\frac{\lambda_{jk}(u,v)-h_{p_2}(u,v)}{u_1^{\tilde e}v_1^{\tilde f}}
\in {\bold k}((u_1,v_1))\cap {\bold k}[[u_1,v_1,w_{j,1}]]={\bold k}[[u_1,v_1]].
$$
We now define the new 2-point relation $\tilde R'$ on $\tilde X$ for $\tilde f$.
Set $T=T(\tilde R)$, $U=\tilde f(T(\tilde R))$. For $p_1\in T$ we define
a primitive 2-point relation $R_{p_1}$ by $U(R_{p_1})=\{q_1=\tilde f(p_1)\}$,
$T(R_{p_1})=\{p_1\}$, 
$$
u_{R_{p_1}}(p_1)=u_{\tilde R}(p_1),
v_{R_{p_1}}(p_1)=v_{\tilde R}(p_1),
w_{R_{p_1}}(p_1)=w_{p_1}.
$$
If $\tau>0$, we define 
$$
a_{R_{p_1}}(p_1)=a_{\tilde R}(p_1), b_{R_{p_1}}(p_1)=b_{\tilde R}(p_1), e_{R_{p_1}}(p_1)=e_{\tilde R}(p_1)
$$
and
$$
\lambda_{R_{p_1}}(p_1)=\lambda_{\tilde R}(p_1)\gamma_{p_1}''(0,0,0)^{e_{\tilde R}(p_1)}.
$$
If $\tau=0$, we define
$$
a_{R_{p_1}}(p_1)=b_{R_{p_1}}(p_1)=-\infty.
$$

Let $\tilde R'$ be the 2-point relation associated to the $R_{p_1}$ for all
$p_1\in T$.
We have $T(\tilde R')=T(\tilde R)$ and $U(\tilde R')=\tilde f(T(\tilde R))$.
From the above calculations, we see that $f:\tilde X\rightarrow \tilde Y$ with the relation $\tilde R'$ satisfies 1 - 4 of the
conditions of Definition \ref{Def128} of a $\tau$-quasi-well prepared
morphism. Finally, by Lemma \ref{Lemma358} and Remark \ref{Remark424},  there exists a commutative diagram
$$
\begin{array}{rll}
\tilde X_1&\stackrel{\tilde f_1}{\rightarrow}&\tilde Y_1\\
\tilde\Phi_1\downarrow&&\downarrow\tilde\Psi_1\\
\tilde X&\stackrel{\tilde f}{\rightarrow}&\tilde Y
\end{array}
$$
where $\tilde\Phi_1$ is a product of blow ups of 2-curves and 3-points and $\tilde\Psi_1$ is a product of blow ups of 2-curves such that
the transform $\tilde R^1$ of $\tilde R'$ for $\tilde f_1$ is defined
and $\tilde R^1$ satisfies 5 of the conditions of Definition \ref{Def128}
(as well as 1 -- 4). Then $\tilde f_1$ is $\tau$-quasi-well prepared with
pre-algebraic structure. Thus the conclusions of Lemma \ref{Lemma348} hold.
\end{pf}

\begin{Theorem}\label{Theorem80} Suppose that $\tau\ge0$, $f:X\rightarrow Y$
is $\tau$-quasi-well prepared with  2-point relation $R$ and $D_X$ is cuspidal for $f$. 
 Then there exists a commutative diagram
$$
\begin{array}{rll}
X_1&\stackrel{f_1}{\rightarrow}&Y_1\\
\Phi\downarrow&&\downarrow\Psi\\
X&\stackrel{f}{\rightarrow}&Y
\end{array}
$$
such that $\Phi$ and $\Psi$ are products of blow ups of possible centers and 
$f_1$ is $\tau$-quasi-well prepared with  2-point relation
$R^1$ and pre-algebraic structure. ($R^1$ will in general not be the transform of $R$.)
Further, $D_{X_1}$ is cuspidal for $f_1$.
\end{Theorem}

\begin{pf}

We will show that there exists a $\tau$-quasi-well prepared diagram (\ref{eq336})
as in the hypothesis of Lemma \ref{Lemma348}. Then Lemma \ref{Lemma348} implies that the conclusions of Theorem \ref{Theorem80} hold.

\noindent{\bf Step 1.}
Let $A_0$ be the set of 2-points  $q\in Y$ such that
$q\in U(\overline R_i)$ for some $\overline R_i$ associated to $R$ and  $f^{-1}(q)\cap T(R_i)\ne\emptyset$.

For $q\in A_0\cap U(\overline R_i)$,
set 
\begin{equation}\label{eq345}
u=u_{\overline R_i}(q), v=v_{\overline R_i}(q),
w_i = w_{\overline R_i}(q).
\end{equation}
Let $\Psi_1:Y_1\rightarrow Y$ be the blowup of all $q\in A_0$, and let
$$
\begin{array}{rll}
X_1&\stackrel{f_1}{\rightarrow}&Y_1\\
\Phi_1\downarrow&&\downarrow\Psi_1\\
X&\stackrel{f}{\rightarrow}&Y
\end{array}
$$
be a $\tau$-quasi-well prepared diagram of $R$ and $\Psi_1$ (such a diagram exits by
Lemma \ref{Lemma32}).
Let $A_1$ be the set of all 2-points
$q_1\in Y_1$ such that for some $i$, $q_1\in U(\overline R_i^1)$, $f_1^{-1}(q_1)\cap T(R_i^1)\ne\emptyset$ and $q_1$
 is on a component $E$ of $D_{Y_1}$ such that $\Psi_1(E)$ is not a 2-point. We have $A_1\subset \Psi_1^{-1}(A_0)$.
Let $\Psi_2:Y_2\rightarrow Y_1$ be the blowup of all
$q_1\in A_1$, and let
$$
\begin{array}{rll}
X_2&\stackrel{f_2}{\rightarrow}&Y_2\\
\Phi_2\downarrow&&\downarrow\Psi_2\\
X_1&\stackrel{f_1}{\rightarrow}&Y_1
\end{array}
$$
be a $\tau$-quasi-well prepared diagram of $R^1$ and $\Psi_2$.  Continue in this way to construct a sequence of $n$ blow ups of 2-points
$\Phi_{k+1}:X_{k+1}\rightarrow Y_{k+1}$ for $0\le k\le n-1$ with $\tau$-quasi-well prepared diagrams
$$
\begin{array}{rll}
X_{k+1}&\stackrel{f_{k+1}}{\rightarrow}&Y_{k+1}\\
\Phi_{k+1}\downarrow&&\downarrow\Psi_{k+1}\\
X_{k}&\stackrel{f_k}{\rightarrow}&Y_{k}
\end{array}
$$
of $R^k$ and $\Psi_{k+1}$.
We have a resulting $\tau$-quasi-well prepared diagram of $R$ 

\begin{equation}\label{eq98}
\begin{array}{rll}
X_{n}&\stackrel{f_{n}}{\rightarrow}&Y_{n}\\
\Phi\downarrow&&\downarrow\Psi\\
X&\stackrel{f}{\rightarrow}&Y.
\end{array}
\end{equation}
Suppose that $q_{n}\in Y_{n}$ is a 2-point such that  $q_n$ is on
a component $E$ of $D_{Y_n}$ such that $\Psi(E)$ is  not a point, and
for some $i$, $q_{n}\in U(\overline R_i^{n})$ and $f_{n}^{-1}(q_{n})\cap T(R_i^{n})\ne\emptyset$. We have permissible parameters 
\begin{equation}\label{eq414}
u_1=u_{\overline R_i^{n}}(q_{n}),v_1=v_{\overline R_i^{n}}(q_{n}),
w_{i,1}=w_{\overline R_i^{n}}(q_{n})
\end{equation}
 at $q_{n}$ such that for $\Psi(q_n)=q$ and with notation of (\ref{eq345}),

\begin{equation}\label{eq96}
\begin{array}{ll}
u&=u_1\\
v&=u_1^nv_1\\
w_i&=u_1^nw_{i,1}
\end{array}
\end{equation}
or

$$
\begin{array}{ll}
u&=u_1v_1^n\\
v&=v_1\\
w_i&=v_1^nw_{i,1}.
\end{array}
$$

Suppose that $p\in X$ is a 3-point such that
$p\in T(R_i)\cap f^{-1}(q)$. Then
 there are
 permissible parameters $x,y,z$ for $u,v,w_i$ at $p$ such that
\begin{equation}\label{eq97}
\begin{array}{ll}
u&=x^ay^bz^c\\
v&=x^dy^ez^f\\
w_i&=x^gy^hz^i\gamma
\end{array}
\end{equation}
where $\gamma$ is a unit series.

We will show that we can choose $n$ sufficiently large in the diagram (\ref{eq98}), so that if $p_{n}\in X_{n}$ is a 3-point such that
$p_{n}\in\Phi^{-1}(p)\cap T(R^{n})$ and $q_n=f_{n}(p_{n})$ is on a component $E$ of $D_{Y_n}$ such that $\Psi(E)$ is  not a point,
  then (\ref{eq97}) must have one of the following
forms (after possibly interchanging $u,v$ and $x,y,z$): 
\begin{equation}\label{eq99}
\begin{array}{ll}
u&=x^ay^b\\
v&=x^dy^ez^f\\
w_i&=x^gy^hz^i\gamma
\end{array}
\end{equation}
where $b\ne 0$, $f\ne 0$ and $i\ne 0$ or 
\begin{equation}\label{eq100}
\begin{array}{ll}
u&=x^a\\
v&=x^dy^ez^f\\
w_i&=x^gy^hz^i\gamma
\end{array}
\end{equation}
with $e$ and $f\ne 0$, and $h$ or $i\ne 0$. 

We will now prove this statement. 
Since $p$ is a 3-point, there exist regular parameters $\overline x,\overline y,\overline z$ in ${\cal O}_{X,p}$
and unit series $\lambda_1,\lambda_2,\lambda_3\in\hat{\cal O}_{X,p}$ such that $\overline x=x\lambda_1$, $\overline y=y\lambda_2$, $\overline z=z\lambda_3$.
Let $\nu$ be any valuation of ${\bold k}(X)$ which has center
 $p_{n}$ on $X_n$. $q_n$
has permissible parameters (\ref{eq414}).

After possibly interchanging $u$ and $v$,  we have a relation (\ref{eq96}), so that
$\nu(v)>n\nu(u)$. We can reindex $x,y,z$ so that $0<\nu(\overline x)\le\nu(\overline y)\le\nu(\overline z)$. Then
$$
(f+e+d-nc)\nu(\overline z)\ge (f-nc)\nu(\overline z)+(e-nb)\nu(\overline y)+(d-na)\nu(\overline x)>0.
$$
Thus if $c\ne 0$, and $n>f+e+d$, we have a contradiction.
Thus taking  $n>f+e+d$ for all $d,e,f$ in local forms (\ref{eq97}) for 3-points 
$p\in T(R)$,
we achieve that $c=0$ in all local forms (\ref{eq97})   which are the images of 3-points  $p_n\in T(R^{n})$ which
map to a point $q_{n}$ of $Y_{n}$ which is on a component $E$ of $D_{Y_n}$ such that $\Psi(E)$ is not a point. 

If $i=0$ (and $c=0$) in (\ref{eq97}) we have
$$
(h+g-nb)\nu(\overline y)\ge (h-nb)\nu(\overline y)+(g-na)\nu(\overline x)>0
$$

so that if $b\ne 0$ and $n>h+g$ we have a contradiction.  Thus, by taking $n\gg 0$ in (\ref{eq98}), we see that if $b\ne 0$, then a form (\ref{eq99}) must hold at $p$ (since $uv=0$ is a local equation of $D_X$ at $p$).
If $b=c=0$  in (\ref{eq97}), then a similar calculation shows that a form (\ref{eq100}) must hold at $p$ (for $n \gg0$).

We observe that in (\ref{eq99})   we have 
\begin{equation}\label{eq101}
(z)\cap\hat{\cal O}_{Y,q}=(v,w_i).
\end{equation}

Suppose that (\ref{eq100}) holds. If $i\ne 0$ then 
\begin{equation}\label{eq399}
(z)\cap\hat{\cal O}_{Y,q}=(v,w_i).
\end{equation}
If $h\ne 0$, then 
\begin{equation}\label{eq400}
(y)\cap \hat{\cal O}_{Y,q}=(v,w_i).
\end{equation}

We will show that in (\ref{eq99}), $v=w_i=0$ is a formal branch of an algebraic curve $C$ in the fundamental locus of
$f:X\rightarrow Y$.
Let $R={\cal O}_{Y,q}$, $S={\cal O}_{X,p}$. $\overline z=0$ is a  local equation for
a component of $D_X$.  We have that $v\in (\overline z)\cap R$ and $u\not\in(\overline z)\cap R$
so that $(\overline z)\cap R=(v)$ or $(\overline z)\cap R=a$ where $a$ is a height two prime
containing $v$.
We have $(z\hat S)\cap\hat R=(v,w_i)$. Suppose that $(\overline z)\cap R=(v)$. We then have an induced morphism
$$
\hat R/(v)\rightarrow \hat S/(z)
$$
which is an inclusion by the Zariski Subspace Theorem (Theorem 10.14 \cite{Ab}). This is impossible, so that $a$ is a height 2 prime in $R$, and defines a curve $C$, which is necessarily in the fundamental locus of $f$
since $\overline z=0$ is a local equation at $p$ of a component of $D_{X}$ which dominates $C$.
A similar argument shows that in (\ref{eq100}), $v=w_i=0$ is a formal branch of an algebraic curve $C$ in the fundamental
locus of $f$.

\vskip .2truein
\noindent{\bf Step 2.} Let $C$ be the reduced curve in $Y$ whose components are the curves
in the fundamental locus of $f$ which are not 2-curves.
Let $\overline C$ be the reduced curve in $Y_{n}$ which is the strict transform of $C$.
The components of $\overline C$ are then in the fundamental locus of $f_n$.
 By Theorems \ref{Theorem4}
and \ref{Theorem6}, we can perform a sequence of blow ups
of 2-points and 2-curves $\Psi':Y'\rightarrow Y_{n}$ so that we can construct a
$\tau$-quasi-well prepared diagram of $\Psi'$ and $R^n$ 
\begin{equation}\label{eq105}
\begin{array}{rll}
X'&\stackrel{f'}{\rightarrow}&Y'\\
\Phi'\downarrow&&\downarrow\Psi'\\
X_{n}&\rightarrow&Y_{n}
\end{array}
\end{equation}
where $R'$ is the transform of $R^n$ on $X'$,
such that  the strict transform $\tilde C$ of $\overline C$ on $Y'$ is  nonsingular,  
if $q'\in U(R_i')\cap \tilde C$ for some $i$ then
 the germ at $q'$ of $\tilde C$ is contained in $S_{R_i'}(q')$, and
  the (disjoint) components of $\tilde C$ are permissible centers for $R'$.

Let $\Psi(1):Y(1)\rightarrow Y'$ be the blow up of  
$\tilde C$.  By Theorem \ref{Theorem6}, 
we have a $\tau$-quasi-well prepared diagram of $\Psi(1)$ and $R'$ 
\begin{equation}\label{eq113}
\begin{array}{rll}
X(1)&\stackrel{f(1)}{\rightarrow}& Y(1)\\
\Phi(1)\downarrow&& \downarrow\Psi(1)\\
X'&\stackrel{f'}{\rightarrow}&Y'.
\end{array}
\end{equation}

Let $R(1)$ be the transform of $R'$ on $X(1)$.
Suppose that $\overline q\in U(R_i(1))\subset Y(1)$ is a 2-point on the exceptional divisor of $\Psi(1)$, $\overline q$ is on 
a component $E$ of $D_{Y(1)}$ such that $(\Psi\circ\Psi'\circ\Psi(1))(E)$ is not a point of $Y$,
 and there exists a 3-point
$\overline p_i\in (f(1))^{-1}(\overline q)\cap T(R_i(1))\subset  X(1)$. Let $q=\Psi\circ\Psi'\circ\Psi(1)(\overline q)$, $\tilde q=\Psi(1)(\overline q)$.
Let
$$
\tilde u=u_{\overline R'_i}(\tilde q), \tilde v=v_{\overline R_i'}(\tilde q), \tilde w_i=w_{\overline R_i'}(\tilde q),
$$
$$
u=u_{\overline R}(q), v=v_{\overline R}(q), w_i=w_{\overline R_i}(q).
$$
We have an expression, after possibly interchanging $u$ and $v$, 
\begin{equation}\label{eq342}
u=\tilde u, v=\tilde u^e\tilde v, w_i=\tilde u^f\tilde w_i
\end{equation}
for some $e,f>0$. $\tilde v=0$ is a local equation of the strict transform of $D_Y$ at $\tilde q$, and
$\tilde v=\tilde w_i=0$ are local equations of  $\tilde C$ at $\tilde q$ (by (\ref{eq101}) (\ref{eq399}) or (\ref{eq400})). $\Psi(1)$ is the blow up of $(\tilde v,\tilde w_i)$
above $\tilde q$. 
Since $\overline q\in U(R_i(1))$, we must have
$$
\tilde u=u_{\overline R_i(1)}(\overline q), \tilde v=v_{\overline R_i(1)}(\overline q),
\tilde w_i=v_{\overline R_i(1)}(\overline q)w_{\overline R_i(1)}(\overline q).
$$

Substituting into (\ref{eq342}), we have 
\begin{equation}\label{a*}
u=u_{\overline R_i(1)}(\overline q),
v=u_{\overline R_i(1)}(\overline q)^ev_{\overline R_i(1)}(\overline q),
w_i=u_{\overline R_i(1)}(\overline q)^fv_{\overline R_i(1)}(\overline q)w_{\overline R_i(1)}(\overline q)
\end{equation}
with $e,f>0$.

We now apply steps 1 and 2 of the proof to $f(1):X(1)\rightarrow Y(1)$ and $R(1)$. We construct a $\tau$-quasi-well prepared diagram
$$
\begin{array}{rll}
X(2)&\stackrel{f(2)}{\rightarrow}&Y(2)\\
\overline \Phi(2)\downarrow&&\downarrow\overline\Psi(2)\\
X(1)&\rightarrow &Y(1),
\end{array}
$$
where $R(2)$ is the transform of $R(1)$ on $X(2)$ such that if $q_2\in U(R_i(2))\subset Y(2)$ is a 2-point such that $q_2$ is on a component $E$ of $D_{Y(2)}$ such that $\overline\Psi(2))(E)$ is not a point of $Y(1)$ and there exists a 3-point $p_2\in f(2)^{-1}(q_2)\cap T(R_i(2))\subset X(2)$ and $q_1=\overline\Psi(2)(q_2)$ then we have an expression:
\begin{equation}\label{b*}
\begin{array}{ll}
u_{\overline R_i(1)}(q_1)&=u_{\overline R_i(2)}(q_2)\\
v_{\overline R_i(1)}(q_1)&=u_{\overline R_i(2)}(q_2)^ev_{\overline R_i(2)}(q_2)\\
w_{\overline R_i(1)}(q_1)&=u_{\overline R_i(2)}(q_2)^fv_{\overline R_i(2)}(q_2)w_{\overline R_i(2)}(q_2)
\end{array}
\end{equation}
or
\begin{equation}\label{c*}
\begin{array}{ll}
u_{\overline R_i(1)}(q_1)&=u_{\overline R_i(2)}(q_2)v_{\overline R_i(2)}(q_2)^e\\
v_{\overline R_i(1)}(q_1)&=v_{\overline R_i(2)}(q_2)\\
w_{\overline R_i(1)}(q_1)&=u_{\overline R_i(2)}(q_2)v_{\overline R_i(2)}(q_2)^fw_{\overline R_i(2)}(q_2)
\end{array}
\end{equation}
with $e,f>0$.

Let $q=(\Psi\circ\Psi'\circ\Psi(1)\circ\overline\Psi(2))(q_2)$.
Substituting (\ref{b*}) or (\ref{c*}) into (\ref{a*}), we see that if $q_2$ is on a component $E$ of $D_{Y(2)}$ which does not contract to $q$, then we have an expression
\begin{equation}\label{d*}
\begin{array}{ll}
u=u_{R_i}(q)&=u_{\overline R_i(2)}(q_2)\\
v=v_{R_i}(q)&=u_{\overline R_i(2)}(q_2)^{e_2}v_{\overline R_i(2)}(q_2)\\
w=w_{R_i}(q)&=u_{\overline R_i(2)}(q_2)^{f_2}v_{\overline R_i(2)}(q_2)^2w_{\overline R_i(2)}(q_2)
\end{array}
\end{equation}
with $e_2,f_2\ge 2$.

Iterating steps 1 and 2, we construct a sequence of $\tau$-quasi-well prepared diagrams

\begin{equation}\label{g*}
\begin{array}{rll}
\vdots&&\vdots\\
\downarrow&&\downarrow\\
X(n)&\stackrel{f(n)}{\rightarrow}&Y(n)\\
\overline \Phi(n)\downarrow&&\downarrow\overline\Psi(n)\\
X(n-1)&
\stackrel{f(n-1)}{\rightarrow} &Y(n-1)\\
\downarrow&&\downarrow\\
\vdots&&\vdots\\
\downarrow&&\downarrow\\
X(1)&\stackrel{f(1)}{\rightarrow}&Y(1)\\
\overline \Phi(1)\downarrow&&\downarrow\overline\Psi(1)\\
X&\stackrel{f}{\rightarrow} &Y.
\end{array}
\end{equation}

This algorithm continues as long as there exists $q_n\in U(R_i(n))$ for some $i$ such that $q_n$
is on a component $E$ of $D_{Y(n)}$ which does not contract to a point of $Y$, and $f(n)^{-1}(q_n)\cap T(R_i(n))\ne\emptyset$.

Suppose that the algorithm never terminates. 

Let $\nu$ be a 0-dimensional valuation of ${\bf k}(X)$. We will say that $\nu$ is resolved on $X(n)$ if the center of $\nu$ on $X(n)$ is at a point $p_n$ of $X(n)$ such that either $p_n\not\in T(R(n))$ or $p_n\in T(R(n))$ and all components $E$ of $D_{Y(n)}$ containing $q_n=f(n)(p_n)$ contract to a point of $Y$.

By our construction, if $\nu$ is resolved on $X(n)$, then $\nu$ is resolved on $X(m)$ for all $m\ge n$. Further, the set of $\nu$ in the Zariski-Riemann manifold $\Omega(X)$ of $X$ which are resolved on $X$ is an open set.

Suppose that  $\nu$ is a 0-dimensional valuation of ${\bf k}(X)$ such that $\nu$ is not resolved on $X(n)$ for all $n$. 
Let $p_n$ be the center of $\nu$ on $X(n)$, $q_n$ be the center of $\nu$ on $Y(n)$.

 There exists an $i$ such that for all $n$,  $q_n\in U(R_i(n))$ and $p_n\in f(n)^{-1}(q_n)\cap T(R_i(n))$. 
We have expressions
\begin{equation}\label{e*}
\begin{array}{ll}
u=u_{R_i}(q)&=u_{\overline R_i(n)}(q_n)\\
v=v_{R_i}(q)&=u_{\overline R_i(n)}(q_n)^{e_n}v_{\overline R_i(n)}(q_n)\\
w=w_{R_i}(q)&=u_{\overline R_i(n)}(q_n)^{f_n}v_{\overline R_i(n)}(q_n)^nw_{\overline R_i(n)}(q_n)
\end{array}
\end{equation}
with $e_n,f_n\ge n$ for all $n$.

 From (\ref{e*}), we see that $\nu(v)>n\nu(u)$ for all $n\in{\bf N}$. Thus $\nu$ is a composite valuation, and there exists a prime ideal $P$ of the valuation ring $V$ of
$\nu$ such that $v\in P$, $u\not\in P$. Let $\nu_1$ be a valuation whose valuation ring is $V_P$.
We have $\nu_1(u)=0$, $\nu_1(v)>0$. From (\ref{e*}) we see that
\begin{equation}\label{f*}
\nu_1(w)-n\nu_1(v)>0
\end{equation}
for all $n\in{\bf N}$.

At $p=p_0\in X$, we have a form (\ref{eq99}) or (\ref{eq100}). In (\ref{eq99}) we have $\nu_1(y)=0$.    $uv=0$ is a local equation of $D_X$ at $p$. Thus either $a>0$ or $d>0$. If $a>0$ then $\nu_1(x)=0$ and $\nu_1(z)>0$, a contradiction to (\ref{f*}) since $f\ne 0$. If $d>0$, we again have a contradiction to (\ref{f*}).
In (\ref{eq100}) we have $\nu_1(x)=0$ and  $\nu_1(y),\nu_1(z)\ge 0$, a contradiction to (\ref{f*}), since $e,f>0$.

 We have shown that for all 0-dimensional valuations $\nu$ of ${\bf k}(X)$, there exists $n$ such that $\nu$ is resolved on $X(n)$.

By compactness of the Zariski-Riemann manifold \cite{Z1} there exists $N$ such that  all $\nu\in \Omega(X)$ are resolved on $X(N)$, a contradiction to our assumption that (\ref{g*}) is of infinite length. The diagram
$$
\begin{array}{lll}
X(N)&\rightarrow &Y(N)\\
\downarrow&&\downarrow\\
X&\rightarrow&Y
\end{array}
$$
thus satisfies the hypothesis of (\ref{eq336}) of Lemma \ref{Lemma348}, so that the conclusions of Theorem \ref{Theorem80} hold.

\end{pf}

\begin{Lemma}\label{Lemma33} Suppose that $\tau\ge 0$, $f:X\rightarrow Y$ is $\tau$-quasi-well prepared 
with 2-point relation $R$ and pre-algebraic structure (or $\tau$-well prepared with 2-point relation $R$),
$q\in U(\overline R_i)\subset U(R)$  and $p\in f^{-1}(q)\cap T(R_i)$ is a 3-point.  Suppose that $E$ is a component of $D_Y$
containing $q$. Let $C=\overline{E\cdot S_{\overline R_i}(q)}$.

Let $\Psi:Y_n\rightarrow Y$ be obtained by blowing up $q$,  then blowing up the point $q_1$ which is the intersection of
the exceptional divisor over $q$ and  the strict transform of $C$ on $Y_1$, and iterating this procedure $n$ times, blowing up the intersection point of the last exceptional divisor with the strict transform of $C$. Let
$$
\begin{array}{rll}
X_n&\stackrel{f_n}{\rightarrow}&Y_n\\
\Phi\downarrow&&\downarrow\Psi\\
X&\stackrel{f}{\rightarrow}&Y
\end{array}
$$
be the $\tau$-quasi-well prepared (or $\tau$-well prepared) diagram of  $R$ and $\Psi$ obtained from Lemma \ref{Lemma32}
(so that $\Phi$ is an isomorphism above $f^{-1}(Y-\Sigma(Y)))$.

Suppose that for all $n>>0$ there exists a 3-point $p_n\in\Phi^{-1}(p)\cap
T(R_{i}^n)$ such that
$f_n(p_n)=q_n\in \Psi^{-1}(q)\cap C_n$, where $C_n$ is the strict transform of $C$ on $Y_n$. Then
$C$ is a component of the fundamental locus of $f$.
 \end{Lemma}

\begin{pf} At $p$ there is an expression
$$
\begin{array}{ll}
u&=x^ay^bz^c\\
v&=x^dy^ez^f\\
w_i&=x^gy^hz^i\gamma
\end{array}
$$
where 
$$
u=u_{\overline R_i}(p),
v=v_{\overline R_i}(p),
w_i=w_{\overline R_i}(p),
$$
$x,y,z$ are permissible parameters at $p$ for $u,v,w_i$, $\gamma$ is a unit series, and $v=0$ is a local equation of
$E$ at $q$. Then $v=w_i=0$ are local equations of $C$ at $q$.

 By our construction of $\Psi$,  we  have that 
$$
u=u_{\overline R_i^n}(p_n),
v_1=v_{\overline R_i^n}(p_n),
w_{i,1}=w_{\overline R_i^n}(p_n)
$$
are defined by
$$
u_1=u,
v_1=u^nv_1,
w_{i1}=u^nw_{i,1}.
$$

Now the conclusions of the lemma follow from the argument from (\ref{eq97}) to the end of Step 1
in the proof of Theorem \ref{Theorem80}.
\end{pf}

\begin{Theorem}\label{Lemma145} Suppose that $\tau\ge 0$,  $f:X\rightarrow Y$ is
$\tau$-quasi-well
prepared with  2-point relation $R$ and pre-algebraic structure, and $D_X$ is cuspidal for $f$.  
 Then there exists a commutative
diagram
$$
\begin{array}{rll}
X_1&\stackrel{f_1}{\rightarrow}&Y_1\\
\Phi\downarrow&&\downarrow\Psi\\
X&\stackrel{f}{\rightarrow}&Y
\end{array}
$$
such that $\Phi$,$\Psi$ are products of blow ups of possible centers and $f_1$ is
$\tau$-very-well prepared with  2-point relation $R^1$.
(In general, $R^1$ is not the transform of $R$). Further,  $D_{X_1}$ is cuspidal for $f_1$.
\end{Theorem}

\begin{pf}
After modifying $R$ by
replacing the primitive 2-point relations $\{\overline R_i\}$ associated to $R$ with
pre-relations such that each $U(\overline R_i)$ is a single point, we may assume that each $\overline R_i$ is
algebraic (Definition \ref{Def199}), and $R$ is algebraic (Definition \ref{Def156}).

There exists a sequence of blow ups of 2-curves $\Psi_1:Y_1\rightarrow Y$ such that
3 and 4 of Definition \ref{Def65} hold for the transforms $\{\overline R_i^1\}$ of the
$\{\overline R_i\}$ on $Y_1$, by embedded resolution of plane curve
singularities (cf. Section 3.4, Exercise 3.13 \cite{C3}) and by Lemma 5.14 \cite{C3}).  By Lemma \ref{Lemma31}, there exists a $\tau$-quasi-well
prepared diagram
$$
\begin{array}{rll}
X_1&\stackrel{f_1}{\rightarrow}&Y_1\\
\Phi_1\downarrow&&\downarrow\Psi_1\\
X&\stackrel{f}{\rightarrow}&Y
\end{array}
$$
of $R$ and $\Psi_1$
where $\Phi_1$, $\Psi_1$ are products of blow ups of 2-curves. Thus
 $f_1$ is $\tau$-well prepared (with respect to the transform $R^1$ of $R$).
We may thus assume that $f$ is $\tau$-well prepared.
Let
$$
\overline V_0=\left\{
\begin{array}{l}
\gamma_i=\overline{E\cdot S_R(p_i)}\text{ such that  }E\text{ is a component of $D_{Y}$, $p_i\in T(R)$,}\\
 \text{and $\gamma_i$ is a component of the fundamental locus of $f$}.
\end{array}\right\}
$$
Suppose that $\gamma_i\in \overline V_0$. Let $\eta_i$ be a general point of
$\gamma_i$. In a neighborhood of $\eta_i$, $f:X\rightarrow Y$ can be factored by blowing up finitely many curves
which dominate $\gamma_i$ (\cite{Ab1} or \cite{D}). Let $r(0)$ be the total number of components of $D_X$ which dominate
the curves $\gamma_i\in \overline V_0$.
\vskip .2truein
\noindent {\bf Step 1.} By Lemma \ref{Lemma33}, there exists a $\tau$-well prepared diagram of $R$ and $\Psi_1$
$$
\begin{array}{rll}
X_1&\stackrel{f_1}{\rightarrow}& Y_1\\
\Phi_1\downarrow&&\downarrow \Psi_1\\
X&\stackrel{f}{\rightarrow}&Y
\end{array}
$$
where $\Psi_1$ is a product of blow ups of prepared 2-points (of type 1 of Definition \ref{Def66}) such that if
$$
\gamma_i(1)\in V_1=\left\{\begin{array}{l}
\gamma_i(1)=\overline{E\cdot S_{R^1}(p_i)}\text{ such that }
E\text{ is a component of }D_{Y_1},\\
 p_i\in T(R^1)\text{ and }\Psi_1(\gamma_i(1))
\text{ is not a point}
\end{array}
\right\}
$$
then $\gamma_i(1)$ dominates a component $\gamma_i$ of $\overline V_0$, and 
$\Phi_1$ is an  isomorphism over $f^{-1}(Y-\Sigma(Y))$ and thus is an isomorphism over the preimage by $f$ of a general point $\eta_j$ of $\gamma_j$ for all $\gamma_j\in \overline V_0$.
Let $r(1)$ be the number of components of $D_{X_1}$ which dominate some $\gamma_i\in\overline V_0$ and are exceptional for $f_1$. We have $r(1)= r(0)$.
\vskip .2truein
\noindent {\bf Step 2.} By Theorems \ref{Theorem4} and \ref{Theorem6} there exists a
$\tau$-well prepared diagram of  $R^1$
$$
\begin{array}{rll}
X_2&\stackrel{f_2}{\rightarrow}&Y_2\\
\Phi_2\downarrow&&\downarrow\Psi_2\\
X_1&\stackrel{f_1}{\rightarrow}&Y_1
\end{array}
$$
where $\Psi_2$ is a product of blow ups of permissible 2-curves and 2-points,
followed by the blow ups of the strict transforms of the $\gamma_i$ in $\overline V_0$, and possibly by blow ups of more 2-points above the $\gamma_i$.
$\Phi_1\circ\Phi_2$ is an isomorphism over $f^{-1}(Y-\Sigma(Y))$, and thus is an isomorphism over the preimage by $f$ of a general point $\eta_j$ of all $\gamma_j\in \overline V_0$. $\Psi_2$ is
the blow up of the strict transform $\gamma_j(1)$ of $\gamma_j$ over a general point $\eta_j$ of $\gamma_j$ for all $\gamma_j\in \overline V_0$.
Let $r(2)$ be the number of components of $D_{X_2}$ which dominate some $\gamma_i\in\overline V_0$, and are exceptional for $f_2$. We have $r(2)<r(0)$.

Let
$$
V_2=\left\{\begin{array}{l}
\gamma_i(2)=\overline{E\cdot S_{R^2}(p_i)}\text{ such that }
E\text{ is a component of }D_{Y_2},\\
 p_i\in T(R^2)\text{ and }(\Psi_1\circ\Psi_2)(\gamma_i(1))
\text{ is not a point}
\end{array}
\right\}
$$
By construction, if $\gamma_i(2)\in V_2$, then $\gamma_i(2)$ dominates a component of $\overline V_0$.
\vskip .2truein
\noindent {\bf Step 3.}
We now repeat Step 1. As in the construction of $f_1$, there exists a $\tau$-well-prepared
diagram
$$
\begin{array}{rll}
X_3&\stackrel{f_3}{\rightarrow}&Y_3\\
\Phi_3\downarrow&&\downarrow\Psi_3\\
X_2&\stackrel{f_2}{\rightarrow}&Y_2
\end{array}
$$
where $\Psi_3$ is a product of blow ups of 2-points such that if
$$
\gamma_i(3)\in V_3=\left\{\begin{array}{l}
\gamma_i(3)=\overline{E\cdot S_{R^3}(p_i)}\mid
E\text{ is a component of }D_{Y_3},\\
p_i\in T(R^3)\text{ and }
\Psi_1\circ\Psi_2\circ\Psi_3(\gamma_i(3))\text{ is not a point}
\end{array}\right\}
$$
then $\gamma_i(3)$ is in the fundamental locus of $f_3$, $\gamma_i(3)$ dominates a
component $\gamma_i$ of $\overline V_0$ and $\Phi_3$ is an isomorphism over $f_2^{-1}(Y-\Sigma(Y_2))$. 

The total number $r(3)$ of components of $D_{X_3}$ which dominate some $\gamma_i\in \overline V_0$
is  equal to $r(2)$. We now repeat Step 2  to get a reduction $r(4)<r(3)$, where $r(4)$ is the number of components of $D_{X_4}$ which dominate some $\gamma_i\in\overline V_0$
and are exceptional for $f_4$. We see that by induction on $r(2n)$, iterating Step 3, we can construct a $\tau$-well-prepared diagram for  $R$
$$
\begin{array}{rll}
X_4&\stackrel{f_4}{\rightarrow}&Y_4\\
\Phi\downarrow&&\downarrow\Psi\\
X&\stackrel{f}{\rightarrow}&Y
\end{array}
$$
such that if $E$ is a component of $D_{Y_4}$ and $S=S_{R^4}(p)$ for some $p\in T(R^4)$ intersects $E$,
then
$\gamma=\overline{E\cdot S}$ is either contracted to a point
by $\Psi$ or $f_4$ is an isomorphism over the generic
point of $\gamma$.

By Lemma \ref{Lemma33}, there exists a $\tau$-well prepared diagram for  $R^4$
$$
\begin{array}{rll}
X_5&\stackrel{f_5}{\rightarrow}& Y_5\\
\Phi_5\downarrow&&\downarrow\Psi_5\\
X_4&\stackrel{f_4}{\rightarrow}&Y_4
\end{array}
$$
such that if $E$ is a component of $D_{Y_5}$ and $S=S_{R^5}(p)$ for some
$p\in T(R^5)$ intersects $E$, then $\gamma=\overline{E\cdot S}$ is exceptional for
$\Psi\circ\Psi_5$.

Let
$$
W_5=\left\{\begin{array}{l} \gamma=\overline{S_{\overline R_i^5}(q)\cdot E}\mid 
E\text{ is a component of $D_{Y_5}$},\\
\text{$R_i^5$ is associated to $R^5$ and $q\in f_5(T(R_i^5))$}\end{array}\right\}
$$
and let
$$
Z_5=\left\{
\begin{array}{l}
q\in U(R^5)-f_5(T(R^5))\text{ such that there exist }\gamma_i,\gamma_j\in W_5\\
\text{ such that }\gamma_i\ne\gamma_j\text{ and }q\in \gamma_i\cap\gamma_j.
\end{array}\right\}
$$
The points of $Z_5$ are prepared 2-points for $R^5$ (of type 1 of Definition \ref{Def66}). Let $\Psi_6:Y_6\rightarrow Y_5$ be the
blow up of $Z_5$. By Lemma \ref{Lemma32}, there exists a $\tau$-well prepared diagram
$$
\begin{array}{rll}
X_6&\stackrel{f_6}{\rightarrow}&Y_6\\
\Phi_6\downarrow&&\downarrow\Psi_6\\
X_5&\stackrel{f_5}{\rightarrow}&Y_5
\end{array}
$$
of  $R^5$ and $\Psi_6$.

Define
$$
W_6=\left\{\begin{array}{l} \gamma=\overline{S_{\overline R_i^6}(q)\cdot E}\mid 
E\text{ is a component of }D_{Y_6},\\
\text{$R_i^5$ is associated to $R^5$ and $q\in \Psi_6^{-1}(f_5(T(R_i^5)))\cap U(\overline R_i^6)$}\end{array}\right\},
$$
$$
Z_6=\left\{
\begin{array}{l}
q\in U(R^6)-\Psi_6^{-1}(f_5(T(R^5)))\cap U(R^6)\text{ such that there exist }\gamma_i,\gamma_j\in W_6\\
\text{ such that }\gamma_i\ne\gamma_j\text{ and }q\in \gamma_i\cap\gamma_j.
\end{array}\right\}
$$

  We necessarily
have that the curves in $W_6$ are strict transforms of curves in $W_5$. We can
iterate, blowing up $Z_6$, and constructing a $\tau$-well prepared diagram, and repeating
until we eventually construct a $\tau$-well prepared diagram
$$
\begin{array}{rll}
X_7&\stackrel{f_7}{\rightarrow}&Y_7\\
\Phi_7\downarrow&&\downarrow\Psi_7\\
X_6&\stackrel{f_6}{\rightarrow}&Y_6
\end{array}
$$
such that $\Psi_7$ is a sequence of blow ups of prepared 2-points (of type 1 of Definition \ref{Def66}) and if $\gamma_1=\overline{S_{R^7}(p_1)\cdot E_1}$, $\gamma_2=\overline{S_{R^7}(p_2)\cdot E_2}$ for $p_1,p_2\in T(R^7)$ and $E_1, E_2$
components of $D_{Y_7}$, are such that $\gamma_1\ne\gamma_2$, then $\gamma_1\cap\gamma_2
\subset U(R^7)\cap (\Psi_6\circ\Psi_7)^{-1}(f_5(T(R^5))$.

Let $\beta=\Psi\circ\Psi_5\circ\Psi_6\circ\Psi_7$.
We now construct pre-relations $\overline R_{i}'$ on $Y_7$ with associated
primitive relations $R_{i}'$ for $f_7$.

If $\overline R_i$ is a 2-point pre-relation associated to $R$, 
let $T(R_{i}')= T(R_i^7)$
and let 
\begin{equation}\label{eq382}
U(\overline R_{i}')=\left\{
\begin{array}{l} q\in U(R^7)\cap (\Psi_6\circ\Psi_7)^{-1}(f_5(T(R^5))\text{ such that }q\in\overline{E\cdot S_{R^7}(p)}\\
\text{ for some component $E$ of }D_{Y_7}, p\in T(R_i^7)\end{array}\right\}.
\end{equation}

For $q'\in U(\overline R_{i}')$, define $\overline R_{i}'(q')=\overline R_i^7(q')$.
For $p\in T(R_{i}')$ define $R_{i}'(p)=\overline R_{i}^7(f_7(p))$.
Let $R'$ be the 2-point relation for $f_7$ defined by the $\{R_{i}'\}$.
Let $\Omega(\overline R_i')$ be an open subset of $\Omega(\overline R_i^7)$ which contains all $\overline{E\cdot S_{\overline R_i'}(q)}$ for $q\in U(\overline R_i')$. Recall that these curves are all exceptional for $\beta$ and we can take $\Omega(\overline R_i')$ so that $\Omega(\overline R_i')\cap U(\overline R')=U(\overline R_i')$.

$f_7$ is $\tau$-well prepared with 2-point relation $R'$. For all $\overline R_i'$, let 
$$
V_i(Y_7)=\left\{
\gamma=\overline{E_{\alpha}\cdot S_{\overline R_i'}(q)}\text{ such that $q\in U(\overline R_i')$, $E_{\alpha}$
is a component of $D_{Y_7}$}\right\}.
$$
By our construction, Lemmas \ref{Lemma31}, \ref{Lemma32} and \ref{Lemma171},
and 2 of Remark \ref{Remark293}, every curve $\gamma\in V_i(Y_7)$ is prepared for $R^7$ of type 4. By (\ref{eq382}), we
now conclude that every curve $\gamma\in V_i(Y_7)$ is prepared for $R'$ of type 4. 1 and 2 of Definition \ref{Def130}
thus  hold for $f_7$ and $R'$.
3 of Definition \ref{Def130} holds for $f_7$ and $R'$ since for all $\overline R_i'$, $V_i(Y_7)$ consists of exceptional
curves of $\Omega(\overline R_{i}')$ contracting to a nonsingular point $q_i\in\Omega(\overline R_i)$.
Thus $f_7$ is $\tau$-very-well prepared with 2-point relation
 $R'$.

\end{pf}

\begin{Theorem}\label{Theorem268} Suppose that $f:X\rightarrow Y$ is prepared,
$\tau_f(X)=\tau\ge 0$ and $D_X$ is cuspidal for $f$.
 Then there exists a commutative diagram
$$
\begin{array}{rll}
X_1&\stackrel{f_1}{\rightarrow}&Y_1\\
\Phi\downarrow&&\downarrow\Psi\\
X&\stackrel{f}{\rightarrow}&Y
\end{array}
$$
such that $f_1$ is prepared, $\Phi$, $\Psi$ are products of blowups of 2-curves,
$\tau_{f_1}(X_1)=\tau$, $D_{X_1}$ is cuspidal for $f_1$ and all 3-points $p\in X_1$ such that $\tau_{f_1}(p)=\tau$ map to 2-points of $Y_1$.
\end{Theorem}

\begin{pf} We first define a 3-point relation $R$ on $X$. Let
$$
T=\left\{
p\in X\mid p\text{ is a 3-point, $q=f(p)$ is a 3-point and $\tau_f(p)=\tau$}\right\}.
$$
For $p\in T$, since $\tau\ge 0$, there exist permissible parameters $u,v,w$ at $q=f(p)$ and permissible
parameters $x,y,z$ at $p$ such that
$$
\begin{array}{ll}
u&=x^{\overline a}y^{\overline b}z^{\overline c}\\
v&=x^{\overline d}y^{\overline e}z^{\overline f}\\
w&=x^{\overline g}y^{\overline h}z^{\overline i}\gamma
\end{array}
$$
where $\gamma$ is a unit series, $\text{rank}(u,v)=2$ and $\text{rank}(u,v,w)=2$.
Thus there exist $a,b,c\in {\bf Z}$ such that 
\begin{equation}\label{eq262}
(x^{\overline a}y^{\overline b}z^{\overline c})^a
(x^{\overline d}y^{\overline e}z^{\overline f})^b
(x^{\overline g}y^{\overline h}z^{\overline i})^c=1,
\end{equation}
with $\text{gcd}(a,b,c)=1$ and
$$
\text{min}\{a,b,c\}<0<\text{max}\{a,b,c\}.
$$

For $p\in T$, define 3-point pre-relations $\overline R_p$ on $Y$ by $U(\overline R_p)=\{q=f(p)\}$,
and $\overline R(q)$ is defined (with the notation of (\ref{eq133})) so that
$u=0$ is a local equation of $E_1$, $v=0$ is a local equation of $E_2$, $w=0$ is a local
equation of $E_3$, $a,b,c$ are defined by (\ref{eq262}) and $\lambda=\lambda_p=\gamma(0,0,0)^c$.

For $p\in T$, we now define  primitive 3-point relations $R_p$ for $f$ by $T(R_p)=\{p\}$, with associated
3-point pre-relation $\overline R_p$. We define $R$ to be the associated 3-point relation
for $f$ with $T(R)=\cup_{p\in T}T(R_p)=T$.

 By Theorem
\ref{Theorem137}, there exists a commutative diagram
$$
\begin{array}{rll}
X_1&\stackrel{f_1}{\rightarrow}& Y_1\\
\Phi\downarrow&&\downarrow\Psi\\
X&\stackrel{f}{\rightarrow}&Y
\end{array}
$$
where $\Phi$, $\Psi$ are products of blow ups of 2-curves
such that $f_1$ is prepared, $\tau_{f_1}(X_1)=\tau$,  the transform $R^1$ of $R$ for
$f_1$ is resolved, and $D_{X_1}$ is cuspidal for $f_1$. Thus all 3-points $p\in X_1$ with $\tau_{f_1}(p)=\tau$  map to 2-points of $Y_1$. 

\end{pf}

\begin{Theorem}\label{Theorem269} Suppose that $f:X\rightarrow Y$ is prepared,
$\tau_{f}(X)=\tau\ge0$ and $D_{X}$ is cuspidal for $f$.   Then there exists a commutative diagram
$$
\begin{array}{rll}
X_1&\stackrel{f_1}{\rightarrow}&Y_1\\
\Phi_1\downarrow&&\downarrow\Psi_1\\
X&\stackrel{f_1}{\rightarrow}&Y
\end{array}
$$
such that $\Phi_1$ and $\Psi_1$ are products of blow ups of possible centers, $f_1$ is $\tau$-very-well prepared with 2-point relation $R^1$,
 and $D_{X_1}$ is cuspidal for $f_1$.
\end{Theorem}

\begin{pf} 
By Theorem \ref{Theorem268}, there exists a commutative diagram  
$$
\begin{array}{rll}
X_1&\stackrel{f_1}{\rightarrow}&Y_1\\
\Phi\downarrow&&\downarrow \Psi\\
X&\stackrel{f}{\rightarrow}&Y
\end{array}
$$
such that $\Phi$ and $\Psi$ are products of 2-curves, $f_1$ is prepared, $\tau_{f_1}(X_1)\le\tau_f(X)$, $D_{X_1}$
is cuspidal for $f_1$ and all 3-points $p$ of $X_1$
such that $\tau_{f_1}(p)=\tau$ map to 2-points of $Y_1$. Now by Theorems \ref{Theorem169},
 \ref{Theorem80}  and \ref{Lemma145} there exists a commutative diagram
$$
\begin{array}{rll}
X_2&\stackrel{f_2}{\rightarrow}&Y_2\\
\downarrow&&\downarrow\\
X_1&\stackrel{f_1}{\rightarrow}&Y_1
\end{array}
$$
where the vertical arrows are products of blow ups of possible centers  such that $f_2$ is $\tau$-very-well prepared and $D_{X_2}$ is cuspidal for $f_2$.

\end{pf}

\section{Toroidalization}
Suppose that $f:X\rightarrow Y$ is a birational morphism of nonsingular projective 3-folds with toroidal structures $D_Y$ and $D_X=f^{-1}(D_Y)$, such that $D_X$ contains the singular locus of $f$.

\begin{Theorem}\label{Theorem270} Suppose that $\tau\ge 0$ and $f:X\rightarrow Y$  is $\tau$-very-well prepared with 2-point relation $R$.
Further suppose that  $D_X$ is cuspidal for $f$.
Then there exists a $\tau$-very-well prepared diagram
$$
\begin{array}{rll}
X_1&\stackrel{f_1}{\rightarrow}&Y_1\\
\downarrow&&\downarrow\\
X&\stackrel{f}{\rightarrow}&Y
\end{array}
$$
such that the transform $R^1$ of $R$ is resolved and $D_{X_1}$ is cuspidal for $f_1$.
In particular, $f_1$ is prepared, $\tau_{f_1}(X_1)<\tau$, and $D_{X_1}$ is cuspidal for $f_1$. 
\end{Theorem}

\begin{pf} 
Fix a  pre-relation $\overline R_t$ associated to $R$  on $Y$, with associated primitive relation $R_t$.
By induction on the number of pre-relations associated to $R$, it suffices to resolve $R_t$
by a $\tau$-very-well prepared diagram (for $R$).

Recall (Definition \ref{Def130})
$$
V_t(Y)=\left\{\begin{array}{l}\overline{E\cdot S}\text{ such that } E\text{ is a component of }D_Y,\\
S=S_{\overline R_t}(q)\text{ for some }q\in U(\overline R_t)
\end{array}\right\}.
$$

$F_t=\sum_{\gamma\in V_t(Y)}\gamma$ is a SNC divisor on $\Omega(\overline R_t)$ whose intersection graph is a tree.

If $\gamma_1=\overline{E_1\cdot S_{\overline R_t}(q_1)}\in V_0$ and $q\in \gamma_1$, say that $\gamma_1$ is good at $q$ if
whenever $q\in U(\overline R_i)$ for some $i$,
 then  $S_{\overline R_i}(q)$ contains the germ of  $\gamma_1$ at $q$
(so that $\gamma_1=\overline{E_1\cdot S_{\overline R_i}(q)}\subset \Omega(\overline R_i)$).
 Otherwise, say that $\gamma_1$ is bad at $q$. Say that $\gamma_1$ is good
if $\gamma_1$ is good at $q$ for all $q\in\gamma_1$.

Let $Y_0= Y$, $X_0=X$, $f_0=f$. We will show that there exists a sequence of
$\tau$-very-well prepared diagrams 
\begin{equation}\label{eq287}
\begin{array}{rll}
X_{i+1}&\stackrel{f_{i+1}}{\rightarrow}&Y_{i+1}\\
\Phi_{i+1}\downarrow&&\Psi_{i+1}\downarrow\\
X_i&\stackrel{f_i}{\rightarrow}&Y_i
\end{array}
\end{equation}
for $0\le i\le m-1$ such that   the transform $R_t^m$ of $ R_t$ on
$X_m$ is resolved. 

Suppose that  $\gamma_1\in V_0$ and $q\in \gamma_1$
is a bad point. By Remark \ref{Remark281}, we have that $q\in U(\overline R_t)$.
Suppose that  $E_1,E_2$ are the two components of $D_{Y}$
containing $q$, and $\gamma_1=\overline{E_1\cdot S_{\overline R_t}(q)}$. Let
$\gamma_2=\overline{E_2\cdot S_{\overline R_t}(q)}$.

We will show that $q$ is a good point of $\gamma_2$.

$\gamma_1$ not good at $q$  implies there exists
$j\ne t$ such that $q\in U(\overline R_j)$ and the germ of $\gamma_1$ at $q$ is
not contained in $S_{\overline R_j}(q)$.
Let
$$
u=u_{\overline R_t}(q),
v=v_{\overline R_t}(q),
w_t=w_{\overline R_t}(q).
$$
After possibly interchanging $u$ and $v$ we have that $u=w_t=0$ are local equations of
$\gamma_1$, $v=w_t=0$ are local equations of $\gamma_2$ at $q$. Let $w_j= w_{\overline R_j}(q)$.
In the equation
$$
w_j=w_t+u^{a_{tj}}v^{b_{tj}}\phi_{tj}
$$
of 3 of Definition \ref{Def65} we have $a_{tj}=0$.

If $q$ is not a good point for $\gamma_2$ then there exists $k\ne t$ such that
$q\in U(\overline R_k)$ and the germ of $\gamma_2$ at $q$ is not contained in $S_{\overline R_k}(q)$. Let $w_k=w_{\overline R_k}(q)$. In the equation
$$
w_k=w_t+u^{a_{tk}}v^{b_{tk}}\phi_{tk}
$$
we thus have $b_{tk}=0$. But we must have
$$
(0,b_{tj})\le (a_{tk},0)\text{ or }(a_{tk},0)\le (0,b_{tk})
$$
by 4 of Definition \ref{Def65}, which is impossible. Thus $q$ is a good point for $\gamma_2$.

Suppose that  all $\gamma\in V_t(Y)$ are bad.
Pick $\gamma_1\in V_t(Y)$. Since $\gamma_1$ is bad there exists
 $\gamma_2\in V_t(Y)
-\{\gamma_1\}$ such that $\gamma_2$ is good at $q_1=\gamma_1\cap\gamma_2$
(as shown above).
$\gamma_1\cap\gamma_2$ is a single point since $V_t(Y)$ is a tree.
Since $\gamma_2$ is bad and $V_t(Y)$ is a tree, there exists $\gamma_3\in V_t(Y)$ which intersects $\gamma_2$
at a single point $q_2$ and is disjoint from $\gamma_1$ such that $\gamma_3$ is good at
$q_2$. Since $V_t(Y)$ is a finite set, we must eventually find a curve which is good,
a contradiction.

Let $\gamma\in V_t(Y)$ be a good curve, so that it is prepared for $R$
of type 4, and is a *-permissible center (Lemma \ref{Lemma67}) and let $\Psi_1':
Y_1'\rightarrow Y$ be the blow up of $\gamma$.

By Lemma \ref{Lemma67} we can construct a $\tau$-very-well prepared diagram
of the form of (\ref{eq233}) of Definition \ref{Def219} 
\begin{equation}\label{eq321}
\begin{array}{rll}
X_1&\stackrel{f_2}{\rightarrow}&Y_1\\
\downarrow&&\downarrow\\
\downarrow&&Y'_1\\
\downarrow&&\downarrow\Psi_1'\\
X&\rightarrow&Y.
\end{array}
\end{equation}
where $Y_1\rightarrow Y'_1$ is a sequence of blow ups of 2-points which are prepared
for the transform of $R$ of type 2 of Definition \ref{Def66}.
Observe that if $\gamma_1\in V_t(Y)$ is a good curve, with $\gamma_1\ne\gamma$, then the strict transform of $\gamma_1$ is a good curve in $V_t(Y_1)$.

 We now
iterate this process. 
We order the curves in $V_t(Y)$,  and choose $\gamma=\overline{E\cdot S_{R_t}(q)}\in V_t(Y)$ in the construction of the diagram (\ref{eq321}) so that it is the minimum
good curve in $V_t(Y)$.

We inductively define a sequence of $\tau$-very well prepared diagrams (\ref{eq287}) by blowing up the good curve in $V_t(Y_i)$
with smallest order, and then constructing a very well prepared diagram (\ref{eq287}) of the form of (\ref{eq321}).
Then we define the total ordering on $V_t(Y_{i+1})$ so that the ordering of strict transforms
of elements of $V_t(Y_{i})$ is preserved, and these strict transforms have smaller
order than the element of $V_t(Y_{i+1})$ which is not a strict transform of an element of $V_t(Y_{i})$.
We repeat, as long as $R_t^i$ is not resolved ($T(R_t^i)\ne\emptyset$).

Suppose that the algorithm does not converge in the construction of $f_m:X_m\rightarrow Y_m$ such that the transform  $R_t^m$ of $R_t$ is resolved. Then there exists a diagram 
\begin{equation}\label{eq232}
\begin{array}{rll}
\vdots&&\vdots\\
\downarrow&&\downarrow\\
X_n&\stackrel{f_n}{\rightarrow}&Y_n\\
\Phi_n\downarrow&&\downarrow\Psi_n\\
X_{n-1}&\stackrel{f_{n-1}}{\rightarrow}&Y_{n-1}\\
\downarrow&&\downarrow\\
\vdots&&\vdots\\
\downarrow&&\downarrow\\
X_0=X&\stackrel{f_0=f}{\rightarrow}&Y_0=Y
\end{array}
\end{equation}
constructed by infinitely many iterations of the algorithm such that $T(R_t^n)\ne\emptyset$ for all $n$.

Suppose that $q_n\in U(\overline R_t^n)$ is an infinite sequence of points such that $\Psi_n(q_n)=q_{n-1}$ for all $n$ and $\Psi_n$
is not an isomorphism for infinitely many $n$.

By construction, the restriction of
 $\Psi_n$ to $S_{\overline R_t^{n}}(q_{n})$ is an isomorphism onto  $S_{\overline R_t^{n-1}}(q_{n-1})$ for all $n$. Thus
the restriction
$$
\overline\Psi_n=\Psi_1\circ\cdots\circ\Psi_n:S_{\overline R_t^n}(q_n)
\rightarrow S_{\overline R_t}(q)
$$
is an isomorphism, where $q=q_0=\Psi_1\circ\cdots\circ\Psi_n(q_n)$. Without loss of generality, we may assume that no $\Psi_n$ is an
isomorphism (on $Y_n$) at $q_n$. We have permissible parameters
$u_i=u_{\overline R_t^i}(q_i),v_i=v_{\overline R_t^i}(q_i),w_{t,i}=w_{\overline R_t^i}(q_i)$
at $q_i$ for all $i$ such that either 
\begin{equation}\label{eq230}
u_i=u_{i+1},
v_i=v_{i+1},
w_{t,i}=u_{i+1}w_{t,i+1}
\end{equation}
or 
\begin{equation}\label{eq231}
u_i=u_{i+1},
v_i=v_{i+1},
w_{t,i}=v_{i+1}w_{t,i+1}.
\end{equation}

Suppose there exists $k\ne t$  such that $q_n\in U(\overline R_k^n)$ for all $n$.

Let $w_{k,i}=w_{\overline R_k^i}(q_i)$ for $i\ge 0$.

The relations
$$
w_{k,i}-w_{t,i}=u_i^{a_{tk}}v_i^{b_{tk}}\phi_{t,k}
$$
of 3 of Definition \ref{Def65}
transform to
$$
w_{k,i+1}-w_{t,i+1}=u_{i+1}^{a_{tk}-1}v_{i+1}^{b_{tk}}\phi_{t,k}
$$
under (\ref{eq230}), and transform to
$$
w_{k,i+1}-w_{t,i+1}=u_{i+1}^{a_{tk}}v_{i+1}^{b_{tk}-1}\phi_{t,k}
$$
under (\ref{eq231}). But we see that after a finite number of iterations $q_n\not\in U(\overline R_k^n)$, unless $a_{tk}=b_{tk}=\infty$.
Thus there exists $n_0$, such that whenever $n\ge n_0$, $q_n\not\in U(\overline R_k^n)$ if $k\ne t$ and $a_{kt},b_{kt}\ne\infty$ .

By our ordering, we have that there exists an $n_0$ such that if $n\ge n_0$, $\gamma\in V_t(Y_n)$ is good and if $k$ is such that $\gamma\cap U(\overline R_k^n)\ne\emptyset$ then
$a_{kt},b_{kt}=\infty$, so that the Zariski closures of $\Omega(\overline R_k^n)$ and $\Omega(\overline R_t^n)$ are the same.
Thus all elements of $V_t(Y_n)$ are good for $n\ge n_0$, since otherwise, there would be a bad curve $\gamma_1\in V(Y_n)$ which intersects
a good curve $\gamma_2$ at a point $q'$ at which $\gamma_1$ is not good. But then we must have that there exists $k\ne t$ such that 
the Zariski closure of $\Omega(\overline R_k^n)$ is not equal to the Zariski closure of $\Omega(\overline R_t^n)$,
 and 
$q'\in U(\overline R_k^n)$, so that $\gamma_2\cap U(\overline R_k^n)\ne\emptyset$, a contradiction.

Our birational morphism  of $\Omega(\overline R_t^n)$ to $\Omega(\overline R_t)$ is an isomorphism in a neighborhood of $U(\overline R_t^n)$. Thus we have a
natural identification of $V_t(Y_n)$ and $V_t(Y)$, and we see that for $n\ge n_0$, the $\Psi_n$ cyclically blow up the different curves of $V_t(Y)$.

Since $R_t^n$ is (by assumption) not resolved for all $n$, there are 3-points
$p_n\in T(R_t^n)\subset X_n$ such that $\Phi_n(p_n)=p_{n-1}$, and $f_n(p_n)=q_n\in U(\overline R_t^n)$  for all $n$.
Without loss of generality, we may assume that no $\Psi_n$ is an
isomorphism  at $q_n$.

With the above notation at $q_n=f_n(p_n)$, we have that  (\ref{eq230}) and (\ref{eq231}) must alternate in the diagram (\ref{eq232}) for $n\ge n_0$, by
our ordering  of $V_t(Y_n)$. Let $p=p_0\in X=X_0$, $q=q_0=f(p)$.

We have 
\begin{equation}\label{eq322}
u=u_n,
v=v_n,
w_t=u_n^{a_n}v_n^{b_n}w_{t,n}
\end{equation}
where 
$$
u=u_{\overline R_t}(q), v=v_{\overline R_t}(q), w_t=w_{\overline R_t}(q)
$$
and
 $a_n$, $b_n$ are positive integers which both go to infinity as $n$ goes to infinity.

There exists (by Theorem 4 of  Section 4, Chapter VI \cite{ZS}) a  valuation $\nu$ of ${\bold k}(X)$ which dominates the  
(non-Noetherian) local ring $\cup_{n\ge 0}{\cal O}_{X_n,p_n}$, and thus dominates the local rings
${\cal O}_{X_n,p_n}$ for all $n$.  Let $x,y,z$ be permissible parameters for $u,v,w_t$ at $p$. Write (in $\hat{\cal O}_{X,p}$) 
\begin{equation}\label{eq385}
\begin{array}{ll}
u&=x^ay^bz^c\\
v&=x^dy^ez^f\\
w_t&=x^gy^hz^i\gamma
\end{array}
\end{equation}
where $xyz=0$ is a local equation of $D_X$ at $p$ and $\gamma$ is a unit series.

 There exist regular parameters $\overline x,\overline y, \overline z$
in ${\cal O}_{X,p}$ and unit series $\lambda_1,\lambda_2,\lambda_3\in\hat{\cal O}_{X,p}$ such that
$\overline x=x\lambda_1$, $\overline y=y\lambda_2$, $\overline z=z\lambda_3$. We may permute $x,y,z$ so that $0<\nu(\overline x)\le\nu(\overline y)\le \nu(\overline z)$. We have (from (\ref{eq322}))
$$
\nu(w_t)-n\nu(u)-n\nu(v)>0
$$
for all $n\in {\bf N}$. Thus
$$
0<(g-na-nd)\nu(\overline x)+(h-ne-nb)\nu(\overline y)+(i-nf-nc)\nu(\overline z)
\le ((g+h+i)-nf-nc)\nu(\overline z)
$$
for all $n$. Thus $f=c=0$, but this is impossible, since $uv=0$ is a local equation of $D_X$ at $p$.

Thus the algorithm converges in a morphism $f_m:X_m\rightarrow Y_m$ such that $T(R_t^m)= \emptyset$, and after iterating for each primitive relation
associated to $R$, we obtain the construction of $f_1:X_1\rightarrow Y_1$, as in the conclusions of the theorem, such  that
$f_1$ is prepared,  cuspidal for $D_{X_1}$ and $\tau_{f_1}(X_1)<\tau$.

\end{pf}

\begin{Theorem}\label{Theorem391} Suppose that $f:X\rightarrow Y$ is prepared, $\tau_f(X)=-\infty$ and
$D_X$ is cuspidal for $f$. Then $f$ is toroidal.
\end{Theorem}

\begin{pf} Suppose that $E$ is a component of $D_X$ and $E$ contains a 3-point $\overline p$. Let $f(\overline p)=\overline q$. $\overline q$ is a 3-point, and if $\overline u,\overline v,\overline w$ are permissible
parameters at $\overline q$, then there exists an expression 
\begin{equation}\label{eq317}
\begin{array}{ll}
\overline u&=\overline x^{a_{11}}\overline y^{a_{12}}\overline z^{a_{13}}\\
\overline v&=\overline x^{a_{21}}\overline y^{a_{22}}\overline z^{a_{23}}\\
\overline w&=\overline x^{a_{31}}\overline y^{a_{32}}\overline w^{a_{33}}
\end{array}
\end{equation}
at $\overline p$, where $\overline x,\overline y,\overline z$ are permissible parameters at $\overline p$
for $\overline u,\overline v,\overline w$ and $\overline x=0$ is a local equation of $E$. In particular, $f$
has a toroidal form of type 1 following Definition \ref{Def274} at $\overline p$. Thus $f(E)=D$ is a component of
$D_Y$, $f(E)=C$ is a 2-curve or $f(E)=\overline q$. Since $D_X$ is cuspidal for $f$, all components $E$ of $D_X$ map to a 3-point,
a 2-curve or a component of $D_Y$.

Suppose that $E$ is a component of $D_X$, $p\in E$, $q=f(p)$, $u,v,w$ are permissible parameters at $q$,
and $x,y,z$ are permissible parameters for $u,v,w$ at $p$, such that $x=0$ is a local equation of $E$. For $g\in\hat{\cal O}_{X,p}$, let $\text{ord}_Eg$ be the largest power of $x$ which divides $g$ in $\hat{\cal O}_{X,p}$ Let 
$$
\text{Jac}(f)=\text{Det}\left(
\begin{array}{lll}
\frac{\partial u}{\partial x}&\frac{\partial u}{\partial y}&\frac{\partial u}{\partial z}\\
\frac{\partial v}{\partial x}&\frac{\partial v}{\partial y}&\frac{\partial v}{\partial z}\\
\frac{\partial w}{\partial x}&\frac{\partial w}{\partial y}&\frac{\partial w}{\partial z}\\
\end{array}\right).
$$

Then we define
$$
\lambda(E)=\text{ord}_Euvw-\text{ord}_E\text{Jac}(f).
$$
$\lambda(E)$ is an invariant of $E$, since $f(E)$ is a component of $D_Y$, a 2-curve or a 3-point. Notice that
$$
\lambda(E)=\text{ord}_Euv-\text{ord}_E\text{Jac}(f)
$$
if $q\in f(E)$ is a 2-point, and
$$
\lambda(E)=\text{ord}_Eu-\text{ord}_E\text{Jac}(f)
$$
if $q\in f(E)$ is a 1-point.
We compute $\lambda(E)$ in (\ref{eq317}) or in a toroidal form following Definition \ref{Def274}
if $E$ does not contain a 3-point (recall that $D_X$ is cuspidal for $f$), to see that
$$
\lambda(E)=1
$$
 for all components $E$ of $D_X$.
 
 We will now verify that if $p\in D_X$, then $f$ is toroidal at $p$, from which the conclusions of the theorem follow.
 \vskip .2truein
 \noindent{\bf Case 1. Suppose that $f(p)=q$ is a 3-point.}
 Let $u,v,w$ be permissible parameters at $q$, and $x,y,z$ be permissible parameters for $u, v, w$ at $p$.
 
  First suppose that $p$ is a 1-point. Then (since $f$
 is prepared),  we have an expression (after possibly interchanging $u,v,w$) 
 \begin{equation}\label{eq318}
 \begin{array}{ll}
 u&=x^a\\
 v&=x^b(\alpha+y)\\
 w&=x^c(\gamma(x,y)+x^dz)
 \end{array}
 \end{equation}
  with $0\ne\alpha$, $a,b,c>0$, $\gamma$  a unit series, $d\ge 0$.
 
 Let $E$ be the component of $D_X$ containing $p$. Computing $\lambda(E)$ at $p$ from (\ref{eq318}),
 we get
 $$
 1=(a+b+c)-(a+b+c+d-1)=1-d.
 $$
 Thus $d=0$ and $u,v, w$ have the toroidal form 3 following Definition \ref{Def274} at $p$.
 
 Suppose that $p$ is a 2-point. Let $E$ and $E'$ be the  components of $D_X$ containing $p$. Then we have
 (after possibly interchanging $u,v,w$ and $x,y,z$),
 a form 
 \begin{equation}\label{eq325}
 u=x^ay^b, v=x^cy^d, w=x^my^n(g(x,y)+x^ey^fz)
 \end{equation}
 where $g(x,y)$ is a unit series,  $x=0$ is a local equation of $E$, and $y=0$ is a local equation of $E'$,
 or 
 \begin{equation}\label{eq326}
 u=(x^ay^b)^k,  v=(x^ay^b)^t(\alpha+z),
 w=(x^ay^b)^m(g(x^ay^b,z)+x^cy^d)
\end{equation}
 where $g$ is a unit series, $\alpha\ne 0$, $ad-bc\ne 0$, $x=0$ is a local equation of $E$, and $y=0$ is a local equation of $E'$.

 In (\ref{eq325}), we compute, 
 $$
 1=\lambda(E)=a+c+m-[a+c+e+m-1]=1-e
 $$
 implies $e=0$. We also compute
 $$
 1=\lambda(E')=b+d+n-[b+d+n+f-1]=1-f
 $$
 which implies $f=0$. Thus we have a toroidal form 2 following Definition \ref{Def274} at $p$.
 
 Suppose that $p$ satisfies (\ref{eq326}). Then
 $$
 1=\lambda(E)=a(k+t+m)-[a(k+t+m)+c-1]=1-c
 $$
 implies $c=0$.
 $$
 1=\lambda(E')=b(k+t+m)-[b(k+t+m)+d-1]=1-d
 $$
 implies $d=0$. This is impossible, so (\ref{eq326}) cannot occur.
 
 \vskip .2truein
 \noindent{\bf Case 2.
 Suppose that $f(p)=q$ is a 2-point.} Let $u,v,w$ be permissible parameters at $q$, and $x,y,z$ be permissible
 parameters for $u,v,w$ at $p$.
 
 Suppose that $p$ is a 1-point, and let $E$ be the component of $D_X$ containing $p$. We have an expression
 
 \begin{equation}\label{eq319}
 \begin{array}{ll}
 u&=x^a,\\
 v&=x^b(\alpha+y),\\
 w&=g(x,y)+x^cz
 \end{array}
 \end{equation}
 where $x=0$ is a local equation of $E$, and $0\ne\alpha$,
 or 
 \begin{equation}\label{eq329}
 \begin{array}{ll}
 u&=x^a,\\
 v&=x^c(\gamma(x,y)+x^dz),\\
 w&=y
 \end{array}
 \end{equation}
 where $\gamma$ is a unit series.
 
 Computing $\lambda(E)$ in (\ref{eq319}), we have
 $$
 1=\lambda(E)=a+b-(a+b+c-1)=1-c
 $$
  implies $c=0$. Thus after  making an appropriate change of variables, we have a toroidal form 5 following
  Definition \ref{Def274} at $p$.
  
  Computing $\lambda(E)$ in (\ref{eq329}), we have
  $$
  1=\lambda(E)=(a+c)-[a+c+d-1]=1-d
  $$
  implies $d=0$. Thus $f$ has a toroidal form 5 following Definition \ref{Def274} at $p$.
  
  Suppose that $p$ is a 2-point. Let $E$ and $E'$ be the  components of $D_X$
  containing $p$. Then we have a form 
  \begin{equation}\label{eq323}
  u=x^ay^b,
  v=x^cy^d,
  w=g(x,y)+x^ey^fz
  \end{equation}
  where $x=0$ is a local equation of $E$, $y=0$ is a local equation of $E'$, 
  and $ad-bc\ne 0$, or
  \begin{equation}\label{eq324}
  u=(x^ay^b)^k, v=(x^ay^b)^t(\alpha+z),
  w=g(x^ay^b,z)+x^cy^d
  \end{equation}
  where $x=0$ is a local equation of $E$, $y=0$ is a local equation of $E'$, $0\ne\alpha\in \bold k$
  and $ad-bc\ne 0$, or 
  \begin{equation}\label{eq330}
  u=(x^ay^b)^k,
  v=(x^ay^b)^l(\gamma(x^ay^b,z)+x^cy^d), w=z
  \end{equation}
  where $\gamma$ is a unit series, $x=0$ is a local equation of $E$, $y=0$ is a local equation of $E'$, $a,b>0$ and $ad-bc\ne 0$. 
  
  Computing $\lambda(E), \lambda(E')$ in (\ref{eq323}), we get 
  $$
  1=\lambda(E)=a+c-(a+c+e-1)=1-e
  $$
  which implies $e=0$.
  $$
  1=\lambda(E')=b+d-(b+d+f-1)=1-f
  $$
  implies $f=0$. Thus after an appropriate change of variables, we have a toroidal form  4 
  following Definition \ref{Def274} at $p$.
  
  Computing $\lambda(E), \lambda(E')$ in (\ref{eq324}), we get
  $$
  1=\lambda(E)=a(k+t)-[a(k+t)+c-1]=1-c
  $$
  which implies $c=0$, and
  $$
  1=\lambda(E')=b(k+t)-[b(k+t)+d-1]=1-d
  $$
  which implies $d=0$. Since $c=d=0$ is impossible in (\ref{eq324}), we see that this form cannot occur.
  
  Computing $\lambda(E),\lambda(E')$ in (\ref{eq330}), we get
  $$
  1=\lambda(E)=a(k+l)-[a(k+l)+c-1]=1-c
  $$
  implies $c=0$,
  $$
  1=\lambda(E')=b(k+l)-[b(k+l)+d-1]=1-d
  $$
  implies $d=0$. Thus (\ref{eq330}) cannot hold.

  \vskip .2truein
  \noindent{\bf Case 3.  Suppose that $f(p)=q$ is a 1-point.} Let $u,v,w$ be permissible parameters at $q$
  satisfying 3 of Definition \ref{Def1}. Suppose that $p$ is a 1-point, and let $E$ be the component of $D_X$
  containing $p$. We have an expression 
  \begin{equation}\label{eq327}
  u=x^a, v=y, w=g(x,y)+x^cz
  \end{equation}
  at $p$.
  $$
  1=\lambda(E)=a-[a+c-1]=1-c
  $$
  implies $c=0$ and thus $f$ has a toroidal form 6 following Definition \ref{Def274} at $p$.
  
  Suppose that $p$ is a 2-point, and let $E,E'$ be the components of $D_X$ containing $p$. We have an expression 
\begin{equation}\label{eq328}
u=(x^ay^b)^k,
v=z, w=g(x^ay^b,z)+x^cy^d
\end{equation}
at $p$, where $x=0$ is a local equation of $E$, $y=0$ is a local equation of $E'$, $a,b>0$ 
and $ad-bc\ne 0$.
$$
1=\lambda(E)=ak-[ak+c-1]=1-c
$$
implies $c=0$.
$$
1=\lambda(E')=bk-[bk+d-1]=1-d
$$
implies $d=0$. Thus (\ref{eq328}) cannot occur.
\end{pf}

\noindent{\bf Proof of Theorem \ref{TheoremA}} By resolution of singularities and resolution of indeterminacy \cite{H} (cf. Section 6.8 \cite{C3}),
and by \cite{M}, there exists a commutative diagram
$$
\begin{array}{rll}
X_1&\stackrel{f_1}{\rightarrow}&Y_1\\
\Phi_1\downarrow&&\downarrow\Psi_1\\
X&\stackrel{f}{\rightarrow}&Y
\end{array}
$$
where $\Phi_1$, $\Psi_1$ are products of blow ups of points and nonsingular curves, such that $X_1$ and $Y_1$
are nonsingular and projective. By  Lemma \ref{Lemma280}, Remark \ref{Remark390} and
Theorem \ref{Theorem190},  we can construct a commutative diagram
$$
\begin{array}{rll}
X_2&\stackrel{f_2}{\rightarrow}&Y_2\\
\Phi_2\downarrow&&\downarrow\Psi_2\\
X_1&\stackrel{f_1}{\rightarrow}&Y_1
\end{array}
$$
such that $\Phi_2$ and $\Psi_2$ are products of blow ups of points and nonsingular curves, such that $f_2$ is prepared
and $D_{X_2}$ is cuspidal for $f_2$.

Now by descending induction on $\tau(X_2)$ and Theorems \ref{Theorem269} and \ref{Theorem270}, there
exists a commutative diagram
$$
\begin{array}{rll}
X_3&\stackrel{f_3}{\rightarrow}&Y_3\\
\Phi_3\downarrow&&\downarrow\Psi_3\\
X_2&\stackrel{f_2}{\rightarrow}&Y_2
\end{array}
$$
such that $\Phi_2$ and $\Psi_3$ are products of blow ups of possible centers, $f_3$ is prepared, $D_{X_3}$ is cuspidal for
$f_3$ and $\tau_{f_3}(X_3)=-\infty$.

By Theorem \ref{Theorem391}, $f_3$ is toroidal, and the conclusions of the theorem follow.

\end{document}